\documentclass[10pt,a4paper, russian]{article}

\usepackage[OT1]{fontenc}                         
\usepackage[utf8]{inputenc}                       
\usepackage[russian]{babel}        

\usepackage{amssymb, amsmath, amsthm, 
amscd, 
graphicx, 
blackdvi,
cite
}

\usepackage{comment}
\usepackage[left=2cm,right=2cm,top=2cm,bottom=2cm]{geometry}
\usepackage{enumitem,multicol}
\usepackage{float}
\usepackage{amsfonts}
\usepackage{cite}

\usepackage{geometry,comment}
\usepackage{indentfirst}
\usepackage{float}

\usepackage{textcomp}
\usepackage{ifpdf}
\ifpdf
\fi

\graphicspath{{figures/}}

\numberwithin{equation}{section}

\newtheorem{theorem}{Теорема}[section]

\newtheorem{corollary}{Следствие}[section]
\newtheorem{proposition}{Предложение}[section]

\theoremstyle{definition}

\theoremstyle{remark}
\newtheorem*{remark}{Замечание}     

\newcommand{\pder}[2]{\frac{\partial \, #1}{\partial \, #2} }
\newcommand{\map}[3]{#1 \, : \, #2 \to #3}
\newcommand{\mapto}[3]{#1 \, : \, #2 \mapsto #3}
\newcommand{\restr}[2]{\left. #1 \right|_{#2}}
\newcommand{\plan}[1]{{}} 
\newcommand{\ddef}[1]{{\em #1}}
\newcommand{\vect}[1]{\left( \begin{array}{c} #1 \end{array} \right)}
\newcommand{\eq}[1]{$(\protect\ref{#1})$}
\newcommand{\be}[1]{\begin{equation}\label{#1}}
\newcommand{\ee}{\end{equation}}
\newcommand{\sn}{\operatorname{sn}\nolimits}
\newcommand{\cn}{\operatorname{cn}\nolimits}
\newcommand{\dn}{\operatorname{dn}\nolimits}
\newcommand{\E}{\operatorname{E}\nolimits}
\newcommand{\Sym}{\operatorname{Sym}\nolimits}
\newcommand{\tangh}{\operatorname{th}\nolimits}
\newcommand{\am}{\operatorname{am}\nolimits}
\newcommand{\glob}{\operatorname{glob}\nolimits}
\newcommand{\loc}{\operatorname{loc}\nolimits}

\renewcommand{\tanh}{\operatorname{th}}
\renewcommand{\cosh}{\operatorname{ch}}
\renewcommand{\sinh}{\operatorname{sh}}
\renewcommand{\tan}{\operatorname{tg}}

\renewcommand{\leq}{\leqslant}
\renewcommand{\geq}{\geqslant}
\newcommand{\I}{\operatorname{\mathcal{I}}\nolimits}
\newcommand{\CI}{\operatorname{\mathcal{CI}}\nolimits}
\newcommand{\CN}{\operatorname{\mathcal{CN}}\nolimits}
\newcommand{\NN}{\operatorname{\mathcal{N}}\nolimits}
\newcommand{\EE}{\operatorname{\mathcal{E}}\nolimits}
\newcommand{\cl}{\operatorname{cl}\nolimits}
\newcommand{\A}{\operatorname{\mathcal{A}}\nolimits}
\newcommand{\Lie}{\operatorname{Lie}}
\newcommand{\rank}{\operatorname{rank}}
\newcommand{\Abn}{\operatorname{Abn}}
\newcommand{ \Gr}{\operatorname{ Gr}}
\newcommand{ \sing}{\operatorname{sing}}
\newcommand{ \diag}{\operatorname{diag}}
\newcommand{ \Kil}{\operatorname{Kil}}
\newcommand{ \type}{\operatorname{type}}

\newcommand{\sll}{\operatorname{sl}}
\newcommand{\Isom}{\operatorname{Isom}}
\newcommand{ \diam}{\operatorname{diam}}
\newcommand{ \GL}{\operatorname{GL}}
\newcommand{ \Tr}{\operatorname{Tr}}

\newcommand{\ct}{\cos{\tau}}
\newcommand{\st}{\sin{\tau}}
\newcommand{\se}{\sin{(\tau \eta \pp)}}
\newcommand{\ce}{\cos{(\tau \eta \pp)}}
\newcommand{\pp}{\bar{p}_3}
\newcommand{\stl}{\sin{\argl}}
\newcommand{\argl}{\frac{t \eta p_3}{2 I_1}}
\newcommand{\ctl}{\cos{\argl}}
\newcommand{\sht}{\sinh{\tau}}
\newcommand{\cht}{\cosh{\tau}}
\newcommand{\sea}{\mathfrak{se}}

\newcommand{\so}{\mathfrak{so}}
\newcommand{\su}{\mathfrak{su}}
\newcommand{\sha}{\mathfrak{sh}}
\newcommand{\sla}{\mathfrak{sl}}

\renewcommand{\O}{\operatorname{O}}
\renewcommand{\Im}{\operatorname{Im}}
\renewcommand{\Re}{\operatorname{Re}}

\def\vH{\vec H}
\def\lan{\langle}
\def\ran{\rangle}
\def\vh{\vec h}
\def\ad{\operatorname{ad}}
\def\const{\operatorname{const}}
\def\spann{\operatorname{span}}
\def\sgn{\operatorname{sgn}}
\def\ds{\displaystyle}

\def\Id{\operatorname{Id}}
\def\Exp{\operatorname{Exp}}
\def\conj{\operatorname{conj}}
\def\cut{\operatorname{cut}}
\def\SO{\operatorname{SO}}
\def\SE{\operatorname{SE}}
\def\SH{\operatorname{SH}}
\def\SL{\operatorname{SL}}
\def\PSL{\operatorname{PSL}}
\def\SU{\operatorname{SU}}
\def\gl{\mathfrak{gl}}
\def\Sym{\operatorname{Sym}}

\def\VVec{\operatorname{Vec}}
\def\MAX{\operatorname{MAX}}
\def\Max{\operatorname{Max}}
\def\Cut{\operatorname{Cut}}
\def\Conj{\operatorname{Conj}}
\def\Lip{\operatorname{Lip}}
\def\SOLV{\operatorname{SOLV}}
\def\SU{\operatorname {SU}}
\def\tr{\operatorname {tr}}
\def\sr{\sqrt{r}}
\def\sM{\sqrt{M}}
\def\Fb{\F^{\b}}

\def\intt{\operatorname{int}}
\def\tmax{t_{\max}}
\def\tconj{t_{\conj}^1}
\def\tcut{t_{\cut}}

\def\then{\quad\Rightarrow\quad}
\def\tmax{t_{\MAX}^1}
\def\ss{\sn p}
\def\cc{\cn p}
\def\dd{\dn p}
\def\tS{\widetilde{S}}
\def\tG{\widetilde{G}}
\def\tN{\widetilde{N}}
\def\k{\mathbf{k}}
\def\ppp{\mathbf{p}}
\def\ddd{\mathbf{d}}
\def\sss{\mathbf{s}}
\def\bx{\bar{x}}
\def\by{\bar{y}}

\def\x{\mathbf{x}}
\def\y{\mathbf{y}}
\def\ab{\mathbf{a}}

\def\Z{{\mathbb Z}}

\def\R{\mathbb{R}}
\def\Q{\mathbb{Q}}
\def\C{\mathbb{C}}
\def\N{\mathbb{N}}
\def\H{\mathbb{H}}

\def\gg{\mathfrak{g}}
\def\T{\mathbb{T}}

\def\F{\Phi}

\def\p{\psi}
\def\a{\alpha}
\def\b{\beta}
\def\d{\delta}
\def\D{\Delta}
\def\G{\Gamma}
\def\g{\gamma}
\def\eps{\varepsilon}
\def\f{\varphi}
\def\s{\sigma}
\def\t{\theta}
\def\Om{\Omega}

\def\lam{\lambda}

\def\r{\rho}

\def\dc{\dot c}

\def\dx{\dot x}
\def\dy{\dot y}
\def\dz{\dot z}

\def\dq{\dot q}
\def\dh{\dot h}
\def\dg{\dot g}

\newcommand{\onefiglabelsize}[4]
{
\begin{figure}[htbp]
\begin{center}
\includegraphics[width=#4\textwidth]{#1}
\\
\parbox[t]{#4\textwidth}{\caption{#2}\label{#3}}
\end{center}
\end{figure}
}

\newcommand{\twofiglabelsize}[8]
{
\begin{figure}[htbp]
\includegraphics[height=#4\textwidth]{#1}
\hfill
\includegraphics[height=#8\textwidth]{#5}
\\
\hfill
\parbox[t]{0.45\textwidth}{\caption{#2}\label{#3}}
\hfill
\parbox[t]{0.45\textwidth}{\caption{#6}\label{#7}}
\hfill
\end{figure}
}

\newcommand{\onefiglabelsizen}[4]
{
\begin{figure}[htbp]
\begin{center}
\includegraphics[height=#4cm]{#1}
\\
\parbox[t]{0.7\textwidth}{\caption{#2}\label{#3}}
\end{center}
\end{figure}
}

\newcommand{\figout}[1]
{#1}

\author{Ю.Л. Сачков\\
Институт программных систем РАН\\
Переславль-Залесский\\
yusachkov@gmail.com}

\title{Левоинвариантные задачи \\оптимального управления \\на группах Ли\footnote{Исследование выполнено при финансовой поддержке РФФИ в рамках научного проекта No. 20-11-50114}}

\begin{document}


\maketitle

\begin{abstract}
Левоинвариантные задачи оптимального управления на группах Ли образуют важный класс задач с большой группой симметрий. Они интересны в теоретическом плане, так как часто допускают полное исследование, и на этих модельных задачах можно изучить общие закономерности. В частности, задачи на нильпотентных группах Ли доставляют фундаментальную нильпотентную аппроксимацию общих задач. Левоинвариантные задачи также часто возникают в приложениях: в классической и квантовой механике, геометрии, робототехнике, моделях зрения и обработке изображений.

Цель данной работы --- дать обзор основных понятий, методов и результатов, относящихся к левоинвариантным задачам оптимального управления на группах Ли. Основное внимание уделено 
описанию экстремальных траекторий и их оптимальности, времени разреза и множества разреза, 
оптимального синтеза.   

Библиография: 238 названий.
\end{abstract}

{\bf Ключевые слова:} оптимальное управление, геометрическая теория управления, левоинвариантные задачи, субриманова геометрия, группы Ли, оптимальный синтез


\tableofcontents

\newpage

\section{Предисловие}

Исследование инвариантных управляемых систем на группах Ли и однородных пространствах является одной из центральных тем геометрической теории управления. С теоретической точки зрения, это --- естественный и важный класс систем, для которого возможна содержательная глобальная теория (именно такие системы возникают, например, при локальной нильпотентной аппроксимации гладких систем). С другой стороны, такие системы моделируют целый ряд прикладных задач (вращение и качение тел, движение роботов, квантовая механика, компьютерное видение). 

Хорошо известно, что получить точное решение глобальной нелинейной задачи управления (например, задачи управляемости или оптимального управления) представляется очень сложным, если задача не имеет большой группы симметрий. Для инвариантных задач на группах Ли (и их проекций на однородные пространства) точное решение  часто можно найти на основе методов геометрической теории управления с использованием техники дифференциальной геометрии, теории групп и алгебр Ли. Полученное решение инвариантной задачи может дать хорошую аппроксимацию соответствующей нелинейной задачи. Например, инвариантная субриманова геометрия на группе Гейзенберга служит краеугольным камнем всей субримановой геометрии.

Основные задачи, рассматривавшиеся для левоинвариантных систем на группах Ли, --- задача управляемости и задача оптимального управления. По задаче управляемости имеется обширная литература; она описана, например, в обзоре \cite{JMS07}.

В данном обзоре предпринята попытка полного описания имеющихся результатов по левоинвариантным задачам оптимального управления на группах Ли. Некоторые из рассматриваемых задач исследовались классиками без привлечения аппарата групп Ли. Так, Леонард Эйлер изучал вращение твердого тела в пространстве и стационарные конфигурации упругого стержня. Левоинвариантные задачи оптимального управления на группах Ли находятся в сфере пристального внимания геометрической теории управления начиная с 1990-ых годов. Эти работы составляют содержание данного обзора.

Обзор имеет следующую структуру.
В разделе \ref{sec:intro} приводятся базовые сведения геометрической теории управления, относящиеся к группам Ли и левоинвариантным задачам оптимального управления.
Раздел \ref{sec:class}  посвящен классификации трехмерных и четырехмерных левоинвариантных субримановых задач.
В разделах \ref{sec:elementary} и \ref{sec:elliptic}
 рассматриваются задачи, интегрируемые соответственно в элементарных и эллиптических функциях.
Некоторые вопросы, оставшиеся неохваченными по причине объема обзора, перечислены в разделе \ref{sec:uncovered}.
Библиография структурирована по типам рассматриваемых задач.

\medskip

Автор благодарит А.А. Аграчева, А.В. Подобряева, А.П. Маштакова, А.А. Ардентова и И.Ю. Бесчастного за полезные советы по содержанию и изложению в данной работе.

Также автор благодарен Е.Ф. Сачковой за помощь в наборе обзора и постоянную поддержку при работе над ним. 
\section{Введение: левоинвариантные задачи оптимального управления}\label{sec:intro}
\subsection{Определения и постановки задач}
\subsubsection{Группы Ли, левоинвариантные управляемые системы и задачи оптимального \\управления}

\ddef{Группа Ли} $G$ --- это гладкое многообразие, снабженное такой групповой структурой, что отображения
\begin{align*}
 &(g_1, g_2) \mapsto g_1 g_2, \quad G \times G \to G,\\
& g \mapsto g^{-1}, \quad G \to G,
\end{align*}   
являются гладкими. Таким образом, \ddef{левый сдвиг} на  любой элемент $g \in G$,
$$
L_g : h \mapsto g h, \quad G \to G,
$$
есть диффеоморфизм. Обозначим его дифференциал через 
$$
L_{g*} : T_h G \to T_{g h} G.
$$

Векторное поле $X\in \VVec(G)$ называется \ddef{левоинвариантным}, если 
$$
L_{g*} X(h) = X(gh), \quad g, h \in G.
$$
Алгебра Ли левоинвариантных векторных полей на $G$ называется \ddef{алгеброй Ли} 
$ \gg$ группы Ли $G$. Эта алгебра Ли изоморфна  касательному  пространству $T_{\Id}G$ в единичном элементе $\Id$ группы Ли $G$.

Управляемая система 
$$
\dg = f(g, u), \quad g \in G, \quad u \in U
$$
называется \ddef{левоинвариантной},  если
$$
L_{h*}f(g, u) = f(hg, u), \quad g, h \in G, \quad  u \in U.
$$
Например, аффинная  по управлениям система
$$
\dg =  X_0(g) + \sum_{i=1}^k u_i X_i(g), \quad g \in G, \quad u=(u_1, \dots, u_k) \in U \subset \R^k,
$$
является левоинвариантной, если таковыми являются поля $X_0$, $X_1$, $\dots, X_k$.

Задача оптимального управления
\begin{align}
&\dg = f(g, u), \quad g \in G, \quad u \in U, \label{def_sys}\\
&g(0) = g_0, \quad g(t_1) = g_1, \label{def_g01}\\
&J = \int_{0}^{t_1} \f(u) dt \to \min  \label{def_J}
\end{align}
называется левоинвариантной, если таковой является система \eq{def_sys}.

Для левоинвариантных задач оптимального управления можно положить $g_0 = \Id$.

\subsubsection{Субримановы задачи, точки разреза, сопряженные точки, кратчайшие и сферы}

\ddef{Субриманова структура} на  гладком многообразии $M$ --- это   распределение (подрасслоение касательного расслоения $TM$)
$$
\D   = \{\D_q \subset T_qM \mid q \in M\}, \quad \dim \D_q \equiv \const, 
$$
снабженное   скалярным произведением
\begin{align*}
&\langle \cdot, \cdot\rangle = \{\langle \cdot, \cdot\rangle_q - \text{скалярное произведение в} \,\,\D_q \mid q \in M\}.
\end{align*}
В частном случае $\D_q = T_qM$, $q \in M$,  получаем \ddef{риманову структуру} $\langle \cdot, \cdot\rangle$  на многообразии $M$.

Субриманова структура на группе Ли $G$  называется \ddef{левоинвариантной}, если распределение и скалярное произведение сохраняются левыми сдвигами на $G$.
Для такой структуры существует глобальный \ddef{ортонормированный репер} из левоинвариантных полей:
\begin{align*}
&X_1, \dots, X_k \in \gg, \quad k = \dim \D_g,\\
&\D_g = \spann(X_1(g), \dots, X_k(g)),\\
&\langle X_i(g), X_j(g) \rangle = \d_{ij}, \quad g \in G.
\end{align*}
Кривая $g \in \Lip([0, t_1], G)$ называется \ddef{допустимой} для распределения $\D$, если 
$$
\dg(t) = \sum_{i=1}^k u_i(t)X_i(g(t))
$$
для некоторых управлений $u_i \in L^{\infty}([0, t_1])$.
\ddef{Субримановой длиной} допустимой кривой называется число
$$
l(g(\cdot)) = \int_0^{t_1}\left(\sum_{i=1}^k u_i^2(t)\right)^{1/2}dt.
$$
\ddef{Субриманово расстояние} (\ddef{расстояние Карно-Каратеодори}) между точками $g_0, \, g_1 \in G$ есть 
$$
d(g_0, g_1) = \inf\{l(g(\cdot)) \mid g(\cdot) - \text{допустимая кривая, соединяющая}\,\, g_0 \,\,\text{и}\,\, g_1  \} .
$$
\ddef{Субриманова кратчайшая} --- это допустимая кривая $g \in \Lip([0, t_1], \, G)$, для которой
$$
l(g(\cdot)) = d(g(0), g(t_1)).
$$
Такая кривая есть решение задачи оптимального управления
\begin{align}
&\dg = \sum_{i=1}^k u_i X_i(g), \quad g \in G, \quad u \in \R^k, \label{def_SR1}\\
&g(0) = g_0, \quad g(t_1)=g_1, \label{def_SR2}\\
&l(g(\cdot)) = \int_0^{t_1}\left(\sum_{i=1}^k u_i^2(t)\right)^{1/2}dt \to \min. \label{def_SR3}
\end{align}
Здесь терминальное время $t_1$ может быть закрепленным или свободным. Из  неравенства Коши-Буняковского следует, что минимизация длины \eq{def_SR3} эквивалентна минимизации \ddef{энергии}
\begin{align}
&J = \dfrac 12 \int_0^{t_1} \sum_{i=1}^k u_i^2(t) dt \to \min \label{def_SR4}
\end{align}
с фиксированным временем $t_1$.

Допустимая кривая $g \in \Lip([0, t_1], G)$ называется  \ddef{субримановой геодезической}, если 
она натурально параметризована (т.е. $\sum_{i=1}^k u_i^2(t) \equiv 1$)  и 
для любого $\tau \in [0, t_1]$ существует отрезок $I \subset \R$, $\tau \in \intt I$, такой, что сужение 
$g|_{I \cap [0, t_1]}$ есть кратчайшая.

\ddef{Временем разреза} вдоль геодезической $g(\cdot)$ называется величина 
$$
t_{\cut}(g(\cdot)) = \sup\{T>0 \,\mid \,g|_{[0,T]}\,\, \text{есть кратчайшая}\} \in (0, \, +\infty].
$$
Соответствующая точка $g(t_{\cut})$ при $t_{\cut} < + \infty$
называется \ddef{точкой разреза} вдоль геодезической $g(\cdot)$.
\ddef{Множеством разреза} для левоинвариантной субримановой задачи, соответствующей начальной  точке $g_0 = \Id$, называется множество $\Cut \subset G $ точек разреза вдоль всех геодезических, выходящих из точки $g_0$.

Аналогично определяются допустимые  траектории, геодезические, время разреза, точки разреза и множество разреза для общих задач оптимального управления \eq{def_sys}--\eq{def_J}.

Если группа Ли $G$ связна, а распределение $\D$ вполне неголономно, то по теореме Рашевского-Чжоу (см.~далее теорему~\ref{rash}), любые точки $g_0, g_1 \in G$ соединимы допустимой кривой, а субриманово расстояние~$d$ превращает $G$ в метрическое пространство. В этом случае \ddef{субримановы сферы} определяются как в произвольном метрическом пространстве:
$$
S_R(g_0) = \{g \in G \, \mid \, d(g_0, g) = R  \}.
$$
Для левоинвариантных субримановых задач $S_R(g_0)= L_{g_0}(S_R(\Id))$, поэтому достаточно исследовать только центрированные в единице сферы $S_R = S_R(\Id)$.

\ddef{Метрической прямой} называется такая геодезическая $g(t)$, $t \in \R$, что для любых $a, b \in \R$  сужение $\restr{g}{[a,b]}$  есть кратчайшая.

\subsubsection{Группы Карно}
\ddef{Алгебра Карно} $\gg$  есть стратифицированная нильпотентная алгебра Ли, порожденная первым слоем~$\gg^{(1)}$:
\begin{align}
&\gg = \gg^{(1)} \oplus \cdots  \oplus \gg^{(s)} , \label{intro_carnot1}\\
& [\gg^{(1)},  \gg^{(k)}] =  \gg^{(k+1)}, \qquad k = 1, \dots, s -1, \label{intro_carnot2}\\
&  [\gg^{(1)},  \gg^{(s)}] = \{0\}.  \label{intro_carnot3}
\end{align}
Наименьшее $s$, для которого выполнены условия \eq{intro_carnot1}--\eq{intro_carnot3}, называется \ddef{глубиной} (или \ddef{ступенью}) алгебры Карно $\gg$. Размерность первого слоя~$\gg^{(1)}$ называется \ddef{рангом} алгебры Карно.
\ddef{Группа Карно} есть связная односвязная группа Ли, алгебра Ли которой есть алгебра Карно.

Если алгебра Карно свободная нильпотентная, то она называется \ddef{свободной алгеброй Карно}, а соответствующая группа Ли --- \ddef{свободной группой Карно}.

Если на первом слое $\gg^{(1)}$  задано скалярное произведение, то левые сдвиги этого слоя и скалярного произведения задают на соответствующей группе Карно левоинвариантную субриманову структуру. Такие левоинвариантные субримановы структуры на группах Карно возникают как нильпотентные аппроксимации общих субримановых структур в точках общего положения.
В этом обзоре рассматриваются несколько таких субримановых структур:
\begin{itemize}
\item
на группе Гейзенберга, это свободная группа Карно ранга 2, глубины 2 (раздел \ref{subsec:heis}),
\item
задача с вектором роста (3,6), соответствующая группа Ли есть свободная группа Карно ранга 3, глубины 2 (раздел \ref{subsec:36}),
\item
двухступенные свободные нильпотентные группы Ли  (раздел \ref{subsec:twostep}),
\item
двухступенные задачи коранга 1, соответствующая группа Ли есть $(2k+1)$-мерная группа Гейзенберга, это несвободная группа Карно ранга $2k$, глубины 2  (раздел \ref{subsec:corank1}),
\item
двухступенные задачи коранга 2, соответствующая  $(k+2)$-мерная группа Ли есть   несвободная группа Карно ранга $k$, глубины 2  (раздел \ref{subsec:corank2}),
\item
на группе Энгеля, это несвободная группа Карно ранга 2, глубины 3 (раздел \ref{subsec:engel}),
\item
на группе Картана, это свободная группа Карно ранга 2, глубины 3 (раздел \ref{subsec:cartan}).
\end{itemize}

\subsubsection{Библиографические комментарии}  В этом обзоре используются лишь базовые сведения  о гладких многообразиях, группах Ли и алгебрах Ли, см.,~например, \cite{warner}. Первоначальные определения субримановой геометрии содержатся в любом источнике \cite{jurd_book, montgomery_book, notes, ABB_book, intro, agr_UMN, versh_gersh}.

\subsection{Элементы геометрической теории управления}\label{subsec:GCT}

\subsubsection{Теорема Рашевского-Чжоу}
Рассмотрим левоинвариантную субриманову структуру $(\D, \lan \cdot, \cdot \ran)$ на группе Ли $G$ с левоинвариантным ортонормированным репером $X_1, \dots, X_k$. Обозначим через $\Lie(X_1, \dots, X_k)$ подалгебру Ли в $\gg$, 
порожденную полями $X_1, \dots, X_k$.
\begin{theorem}[Рашевский-Чжоу]\label{rash}
Пусть группа Ли $G$ связна, а распределение $\D$ вполне неголономно:
$$
\Lie(X_1, \dots, X_k) = \gg.
$$
Тогда
\begin{itemize}
\item[$(1)$]
любые точки в $G$ соединимы допустимой кривой,
\item[$(2)$]
$(G, d)$ есть метрическое пространство,
\item[$(3)$]
топология на $G$, индуцированная метрикой $d$, эквивалентна  топологии многообразия.
\end{itemize}
\end{theorem}

Распределение $\D$ на многообразии $M$  называется \ddef{интегрируемым}, если через каждую точку $q \in M$  проходит гладкое многообразие $N_q$  такое, что $T_qN_q = \D_q$ (интегральное многообразие распределения~$\D$), в этом случае любая субриманова структура $(\D, \langle \cdot, \cdot \rangle)$  также называется \ddef{интегрируемой}. Это равносильно тому, что для любого локального базиса $X_1, \dots, X_k$  распределения $\D$  выполнено равенство $\Lie(X_1, \dots, X_k)(q) = \D_q$, $q \in M$.

\subsubsection{Теорема Филиппова}
Стандартные условия существования оптимального управления в задачах оптимального управления даются теоремой Филиппова \cite{filippov, notes}.

\subsubsection{Принцип максимума Понтрягина на группах Ли}\label{subsubsec:PMP}
Рассмотрим левоинвариантную задачу оптимального управления \eq{def_sys}--\eq{def_J} на группе Ли $G$ с фиксированным терминальным временем $t_1$. Введем \ddef{гамильтониан  принципа максимума Понтрягина}:
\begin{align}
&h_u^{\nu}(\lam) = \lan \lam, f(g, u)\ran + \nu \f(u),\\
&\lam \in T^*_g G \subset T^*G, \quad u \in U, \quad g \in G, \quad \nu \in \R.
\end{align}
Для фиксированных $u \in U$, $\nu \in \R$ обозначим гамильтоново векторное поле $\Vec{h}_u^{\nu} \in \VVec(T^*G)$, соответствующее гамильтониану $h_u^{\nu} \in C^{\infty}(T^*G)$.
\begin{theorem}[Принцип максимума Понтрягина]
Если $g(t)$ есть оптимальная траектория, соответствующая управлению $u(t)$, то существуют кривая $\lam \in \Lip([0, t_1], T^*G)$, $\lam_t \in T^*_{g(t)}G$, и число $\nu \in \{ -1, 0\}$, для которых выполнены условия:
\begin{itemize}
\item[$(1)$]
$\dot{\lam}_t = \Vec{h}_{u(t)}^{\nu}(\lam_t) $ для п.~в. $t \in [0, t_1]$,
\item[$(2)$]
$\ds {h}_{u(t)}^{\nu}(\lam_t) = \max_{v \in U}h_v^{\nu}(\lam_t)$ для всех $t \in [0, t_1]$,
\item[$(3)$]
$(\lam_t, \nu) \neq (0, 0)$ для всех $t \in [0, t_1]$.

Если терминальное время $t_1$ свободно, то к условиям $(1)$--$(3)$ присоединяется условие
\item[$(4)$] ${h}_{u(t)}^{\nu}(\lam_t) \equiv 0$.
\end{itemize}
\end{theorem}

Траектория $g(t)$ и управление $u(t)$, удовлетворяющие принципу максимума Понтрягина, называются \ddef{экстремальными}, а соответствующая кривая  $\lam_t$  --- \ddef{экстремалью}.

Для левоинвариантной  субримановой задачи \eq{def_SR1}, \eq{def_SR2}, \eq{def_SR4}  принцип максимума Понтрягина детализируется следующим образом. Обозначим гамильтонианы 
 $h_i(\lam) = \lan \lam, X_i(g)\ran$, $\ds H(\lam) = \dfrac 12 \sum_{i=1}^{k} h_i^2(\lam)$, $\lam \in T^*G$.

\begin{corollary}
Пусть $g(t)$ есть субриманова кратчайшая, соответствующая  управлению $u(t)$. Тогда существует кривая $\lam \in \Lip([0, t_1], T^*G)$, $\lam_t \in T^*_{g(t)}G $, для  которой выполнено одно, и только одно, из условий:
\begin{align}
& \dot{\lam}_t = \sum_{i=1}^k u_i(t) \Vec h_i (\lam_t), \quad H(\lam_t) \equiv 0, \quad \lam_t \neq 0, \tag{$A$}\\
&\dot{\lam}_t = \Vec{ H} (\lam_t), \quad u_i(t) = h_i (\lam_t). \tag{$N$}
\end{align}
\end{corollary}
Случай $(A)$  называется \ddef{анормальным}, а случай $(N)$ --- \ddef{нормальным}.

Дополним ортонормированный репер субримановой структуры $X_1, \dots, X_k$ до левоинвариантного  репера $X_1, \dots, X_n$ на группе Ли $G$ и введем гамильтонианы  
$h_i(\lam) = \lan \lam, X_i(g) \ran$. Тогда нормальная гамильтонова система $\dot{\lam} = \Vec{H}(\lam)$ может быть представлена в виде:
\begin{align}
&\dh_i = \{H, h_i \}, \quad i = 1, \dots, n, \label{GCT_vert}\\
& \dg = \sum_{i=1}^k h_i X_i(g),\nonumber
\end{align}
где $\{H, h_i  \}$ есть скобка Пуассона гамильтонианов.

\ddef{Вертикальную подсистему} \eq{GCT_vert} можно  рассматривать как систему дифференциальных уравнений на коалгебре Ли $\gg^*$ после тривиализации  кокасательного расслоения $T^*G \cong G\times \gg^*$ левыми сдвигами.

Вдоль непостоянных нормальных экстремалей $H(\lam_t) \equiv \const >0$, и их можно параметризовать натурально  (длиной дуги), то есть так,  чтобы $H(\lam_t) \equiv  1/2$.  Обозначим цилиндр $C = \gg^*\cap\{H =  1/2\}$, тогда натурально параметризованные нормальные экстремальные траектории задаются с помощью \ddef{экспоненциального отображения}
\begin{align}
& \Exp: (\lam_0, t) \mapsto g(t) = \pi \circ e^{t\Vec H}(\lam_0), \label{elem_exp}\\
&\Exp: C\times\R_+ \to G, \nonumber
\end{align}
где экспонента справа в \eq{elem_exp} обозначает поток гамильтонова поля, а $\pi: T^*G \to G$ есть каноническая проекция.

\ddef{Волновым фронтом} за время $t>0$, соответствующим начальной точке $\Id$, называется множество
$$
W_t = \{\Exp(\lambda, t) \mid \lambda \in C\}.
$$
Очевидно включение $S_t \subset W_t$.

\ddef{Анормальным множеством} для левоинвариантной субримановой структуры на группе Ли $G$, соответствующим начальной точке $g_0 = \Id \in G$, называется множество
$$
\Abn = \{g(t) \mid g(\cdot) \text{  анормальная траектория}, \ t > 0, \ g(0) = \Id\} \subset G.
$$

\subsubsection{Условия оптимальности второго порядка}
\begin{theorem}[Условие Лежандра]
Короткие дуги нормальных  экстремальных траекторий оптимальны.
\end{theorem}
Поэтому нормальные экстремальные траектории являются геодезическими. Момент времени $\widehat t > 0$ называется \ddef{сопряженным  временем} для нормальной геодезической $\Exp(\lam, t)$, $\lam \in C$, если 
$ (\lam, \widehat t)$ есть критическая точка экспоненциального отображения, то есть дифференциал $\Exp_{*(\lam, \widehat t)}: T_{(\lam, \widehat t)}(C \times \R_+) \to T_{\widehat g} G$ вырожден, где 
$\widehat g = \Exp(\lam, \widehat t)$. При этом точка $\widehat g$ называется  \ddef{сопряженной точкой}. \ddef{Первое сопряженное время } вдоль геодезической $g(t)$ есть 
$t^1_{\conj} = \inf \{t > 0 \mid t -\text{сопряженное время вдоль}\,\, g(\cdot) \}$.

\begin{theorem}[Условие Якоби]
Пусть $g: [0,\, t_1] \to G$ есть нормальная геодезическая, не содержащая  анормальных дуг. Тогда:
\begin{itemize}
\item[$(1)$]
$t^1_{\conj} > 0$,
\item[$(2)$]
для любого $\tau \in  (0, t^1_{\conj})$ геодезическая $g|_{[0, \tau]}$ есть локально кратчайшая в топологии $W^{1, 2}$ на пространстве горизонтальных  кривых с теми же граничными точками,
\item[$(3)$]
для любого $\tau > t^1_{\conj}$ геодезическая $g|_{[0, \tau]}$ не является кратчайшей.
\end{itemize}
\end{theorem}

\ddef{Первой каустикой} называется множество
$$
\Conj^1 = \{\Exp(\lam, t) \mid \lam \in C, \,\, t = t^1_{\conj}(\lam)   \}.
$$

\begin{theorem}[Условие Гоха]
Пусть $(\D, \langle \cdot, \cdot\rangle)$  есть вполне неголономная субриманова структура на гладком многообразии~$M$. Если любые субримановы шары компактны и
$$
\D_q + [\D,\D]_q = T_qM, \qquad q \in M,
$$
то любая кратчайшая нормальна.
\end{theorem}

\subsubsection{Библиографические  комментарии}
Этот раздел содержит стандартный материал геометрической  теории управления, см.~\cite{jurd_book, notes, montgomery_book, ABB_book, intro}.

\subsection{Симметрийный метод построения оптимального синтеза}
Отыскание оптимальных траекторий в задачах оптимального управления обычно  состоит из следующих шагов: 
\begin{itemize}
\item[$(1)$]
доказательство существования оптимальных траекторий,
\item[$(2)$]
описание экстремальных траекторий,
\item[$(3)$]
выбор оптимальных траекторий из экстремальных.
\end{itemize}

Шаг $(1)$ обычно выполняется с помощью общих методов  теории управления. Например, для задач субримановой геометрии условия  существования  кратчайших  даются теоремами Рашевского-Чжоу и Филиппова.

Шаг $(2)$, как правило, выполняется  с помощью принципа максимума Понтрягина. Доказывается интегрируемость гамильтоновой  системы принципа максимума Понтрягина, и ее решения параметризуются  в явном виде.

Шаг $(3)$ наиболее сложен. Локальную оптимальность экстремальных траекторий обычно можно исследовать с помощью оценок сопряженного времени. Для изучения глобальной оптимальности в задачах  с большой группой симметрий (в частности,  в левоинвариантных задачах) часто применяется следующий \ddef{симметрийный метод}.
\begin{itemize}
\item[$(3.1)$]
Описываются дискретные  и непрерывные  симметрии экспоненциального отображения.
\item[$(3.2)$]
Отыскиваются  \ddef{точки Максвелла}, соответствующие симметриям (то есть точки, куда несколько симметричных экстремальных траекторий приходят в одно и то же время). Эти точки (и их прообразы  относительно экспоненциального отображения) образуют \ddef{страты Максвелла} в образе (соответственно прообразе) экспоненциального отображения. На каждой экстремальной траектории отыскивается  первое время Максвелла, соответствующее симметриям (то есть первое время, когда экстремальные траектории пересекают страты Максвелла). При достаточно общих  условиях,  экстремальная траектория не может быть оптимальной  после точки Максвелла \cite{max1}.
\item[$(3.3)$]
Доказывается, что на любой геодезической  первое сопряженное время не меньше первого времени Максвелла, соответствующего симметриям. Для этого можно использовать прямые оценки якобиана экспоненциального отображения или гомотопическую инвариантность индекса Маслова (количества сопряженных точек на экстремальной траектории \cite{cime, cartan_conj}).
\item[$(3.4)$]
Рассматривается ограничение экспоненциального отображения на подобласти, вырезаемые  в прообразе и образе экспоненциального отображения стратами Максвелла, соответствующими симметриям. С помощью теоремы Адамара о глобальном диффеоморфизме \cite{hadamard} доказывается, что это ограничение есть диффеоморфизм. 
\item[$(3.5)$]
На основе описанной таким образом глобальной структуры  
экспоненциального отображения часто можно доказать, что время разреза на экстремальных траекториях равно первому времени Максвелла, соответствующему симметриям. Более того, таким образом можно доказать, что для любой конечной точки в указанных  подобластях  в образе экспоненциального отображения существует  единственная оптимальная траектория, которую можно вычислить, обращая экспоненциальное отображение в этих подобластях.
\item[$(3.6)$]
 Наконец, для задач невысокой размерности с большой группой симметрий иногда удается построить полный \ddef{оптимальный синтез}, то есть закон, сопоставляющий каждой конечной точке в пространстве состояний одну или несколько оптимальных  траекторий, приходящих в эту точку. 
\end{itemize}
Многие  задачи оптимального управления, описанные в этом обзоре, исследованы с помощью этого симметрийного метода.

\subsubsection{Библиографические комментарии}
Симметрийный метод есть обобщение классического  метода Адамара в римановой геометрии,  примененного им, в частности,  к  исследованию  оптимального синтеза на поверхностях отрицательной кривизны~\cite{hadamard}. В описанном виде он применялся к левоинвариантным задачам оптимального управления в работах~\cite{martinet, myasn36, s2r2_sym, cut_sre2, el_cut, engel_cut, sh2_3, so3, sl2}.


\section{Классификации левоинвариантных субримановых задач}\label{sec:class}
\subsection{Задачи на трехмерных группах Ли}\label{subsec:class3}
В этом разделе описана классификация, с точностью до локальных изометрий и дилатаций, всех неинтегрируемых левоинвариантных субримановых структур   ранга 2 на трехмерных группах Ли.
\subsubsection{Трехмерные алгебры Ли}\label{subsubsec:Liealg3}
Все трехмерные алгебры Ли, в которых существует двумерное подпространство, не являющееся подалгеброй, суть алгебры Ли следующих групп Ли:
\begin{itemize}
\item
группа Гейзенберга $H_3$,
\item
$\A^+(\R) \oplus \R$, где $\A^+(\R)$ есть группа сохраняющих ориентацию аффинных функций на $\R$,
\item
$\SOLV$, группы Ли, алгебры Ли которых разрешимы и имеют двумерную производную подалгебру,
\item
группы $\SE(2)$ и $\SH(2)$ сохраняющих  ориентацию евклидовых и гиперболических движений плоскости соответственно,
\item
трехмерные простые группы Ли $\SL(2)$ и $\SU(2)$.
\end{itemize}

\subsubsection{Субримановы структуры}
Пусть $G$ --- трехмерная группа Ли и $(\D, \lan\cdot, \cdot  \ran)$ --- неинтегрируемая  левоинвариантная  субриманова структура на $G$ ранга 2.
\begin{proposition}
Пусть группа Ли $G$ односвязна. Существует левоинвариантный репер $(X_0, X_1, X_2)$ на группе Ли $G$ такой, что $(X_1, X_2)$ есть ортонормированный репер для субримановой структуры $(\D, \lan\cdot, \cdot  \ran)$, в котором таблица умножения есть либо
\begin{align}
&[X_1, X_0] = c_{01}^2 X_2, \label{class3_X10}\\
&[X_2, X_0] = c_{02}^1 X_1, \label{class3_X20}\\
&[X_2, X_1] = c_{12}^1 X_1 + c_{12}^2 X_2 + X_0, \label{class3_X21}
\end{align}
либо
\begin{align}
&[X_1, X_0] = \kappa X_2, \label{class3_X102}\\
&[X_2, X_0] = - \kappa X_1, \label{class3_X202}\\
&[X_2, X_1] = X_0.\label{class3_X212}
\end{align}
\end{proposition}
В случае \eq{class3_X10}--\eq{class3_X21} обозначим
$$
\chi = \dfrac{c_{01}^2 +c_{02}^1 }{2}, \quad \kappa = -(c_{12}^1)^2-(c_{12}^2)^2 + \dfrac{c_{01}^2 - c_{02}^1 }{2},
$$
а в  случае \eq{class3_X102}--\eq{class3_X212} $\chi = 0$.

\subsubsection{Классификация трехмерных субримановых структур}
При растяжениях ортонормированного репера $(X_1, X_2)$ инварианты $\chi$ и $\kappa$ умножаются на ненулевую  константу, поэтому их можно нормировать условием
$$
\chi = \kappa = 0 \quad \text{или} \quad \chi^2 + \kappa^2 = 1, \quad \chi \geq 0.
$$

\ddef{Субриманова изометрия}  между двумя субримановыми многообразиями $(M, \D, \langle\cdot, \cdot\rangle)$  и $(M', \D', \langle\cdot, \cdot\rangle')$  есть диффеоморфизм $\map{f}{M}{M'}$, удовлетворяющий условиям:
\begin{itemize}
\item[(1)]
$f_*(\D) = \D'$, 
\item[(2)]
$\langle X_1, X_2\rangle = \langle f_* X_1, f_*X_2\rangle'$  для любых векторных полей $X_1, X_2$, касающихся распределения $\D$.
\end{itemize}

\begin{theorem}\label{th:class3}
Все неинтегрируемые левоинвариантные  субримановы структуры ранга $2$  на $3$-мерных группах Ли классифицируются с точностью до локальных изометрий и дилатаций как на Рис.~$\ref{fig:class3}$,  где каждая структура обозначена точкой 
$(\kappa, \chi)$, и разные точки обозначают локально неизометричные структуры.

Более того,
\begin{itemize}
\item[$(1)$]
если $\chi = \kappa = 0$, то структура локально изометрична субримановой структуре на группе Гейзенберга (см.~раздел $\ref{subsec:heis}$);
\item[$(2)$]
если $\chi^2 + \kappa^2 = 1$, то  существуют не более трех локально неизометричных нормализованных субримановых структур с этими инвариантами; в частности, на каждой унимодулярной группе Ли для любых 
$(\chi, \kappa)$ существует единственная нормализованная структура;
\item[$(3)$]
если $\chi \neq 0$ или $\chi = 0$ и $\kappa \geq 0$, то две структуры с заданными $(\chi, \kappa)$ локально изометричны тогда и только тогда, когда их алгебры Ли изоморфны.
\end{itemize}
\end{theorem}

\figout{
\onefiglabelsize{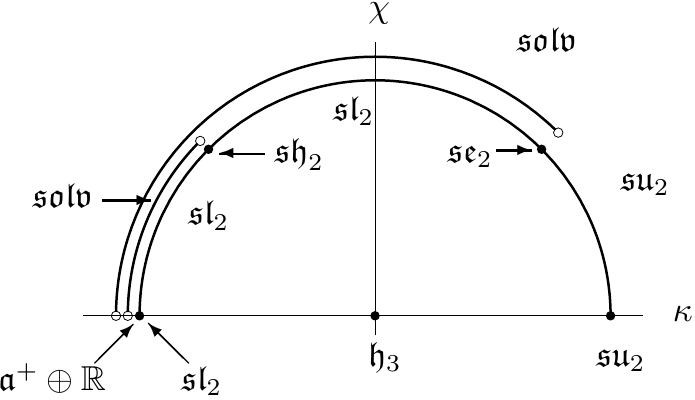}{Классификация 3-мерных контактных левоинвариантных субримановых структур}{fig:class3}{0.5}
}

\subsubsection{Изометрия между $\A^+(\R)\times S^1$ и $\SL(2)$ }
Существуют неизоморфные  группы Ли  с локально изометричными субримановыми структурами: из теоремы~\ref{th:class3} следует, что существует единственная нормализованная левоинвариантная структура на группе 
$\A^+(\R) \oplus \R$ с $\chi = 0$, $\kappa = -1$. Эта структура локально изометрична субримановой структуре на $\SL(2)$, определяемой формой Киллинга.

Группа Ли $\A^+(\R) \oplus \R$ представляется матрицами:
\begin{align*}
&\A^+(\R) \oplus \R = \left\{
\begin{pmatrix}   
a & 0 & b\\
0 & 1 & c \\
0 & 0 & 1
\end{pmatrix} \mid a > 0, \quad b, c \in \R
\right \},
\end{align*}
где действие на вектор $(x, y) \in \R^2$ задается как
\begin{align*}
&\begin{pmatrix}   
a & 0 & b\\
0 & 1 & c \\
0 & 0 & 1
\end{pmatrix}
\begin{pmatrix}
x\\
y\\
1
\end{pmatrix}=
\begin{pmatrix}
a x + b\\
y + c\\
1
\end{pmatrix}.
\end{align*}
Алгебра Ли этой группы Ли  порождена матрицами 
\begin{align*}
&e_1 = 
\begin{pmatrix}  
0 & 0 & 1\\
0 & 0 & 0 \\
0 & 0 & 0 
\end{pmatrix}, \quad
e_2=\begin{pmatrix}  
-1 & 0 & 0\\
0 & 0 & 0 \\
0 & 0 & 0 
\end{pmatrix},\quad
e_3=\begin{pmatrix}  
0 & 0 & 0\\
0 & 0 & 1 \\
0 & 0 & 0 
\end{pmatrix}.
\end{align*}

Рассмотрим субриманову структуру  на  $\A^+(\R)\oplus \R$ с ортонормированным репером $(e_2, e_1 + e_3)$.

Подгруппа $\A^+(\R)$ диффеоморфна полуплоскости $\{(a, b) \in \R^2 \mid a > 0 \}$, которая задается в стандартных полярных координатах как $\{(\r, \t) \mid \r > 0, -\pi/2 < \t < \pi/2  \}$.

Рассмотрим субриманову  структуру на $\A^+(\R)\times S^1$, заданную проецированием структуры на $\A^+(\R)\times \R$.

\begin{theorem}
Диффеоморфизм $F:\A^+(\R)\times S^1 \to  \SL(2)$, заданный как
\begin{align*}
&F(\r, \t, \f) = \dfrac{1}{\sqrt{\r\cos \t}}
\begin{pmatrix}
\cos \f & \sin \f\\
\r \sin(\t - \f) & \r \cos(\t - \f)
\end{pmatrix},
\end{align*}
где $(\r, \t) \in \A^+(\R)$ и $\f \in S^1$, есть глобальная субриманова изометрия.
\end{theorem}

\subsubsection{Библиографические комментарии}
Изложение в этом разделе опирается на работу~\cite{agr_bar}.

Отметим, что в более ранней работе~\cite{FG} получена полная классификация субримановых однородных пространств, то есть субримановых структур, имеющих транзитивную группу изометрий, гладко действующих на многообразии. Использованные в этой работе инварианты $\tau_0$ и $K$ совпадают, с точностью  до нормирующего множителя, с инвариантами $\chi$ и $\kappa$ этого раздела.

Классификация контактных левоинвариантных субримановых метрик на трехмерных группах Ли с точностью до автоморфизмов алгебры Ли получена в работах \cite{versh_gersh88, versh_gersh}.

\subsection{Задачи на четырехмерных группах Ли}\label{subsec:class4}
\subsubsection{Распределения Энгеля}
Пусть $M$ есть
четырехмерное многообразие.   Распределение  $\D \subset TM$ ранга 2 называется \ddef{распределением Энгеля}, если 
$$
\rank([\D, \D])= 3, \quad \rank([\D, [\D, \D]]) = 4,
$$
где $[\D, \D]$ состоит из касательных векторов, которые можно получить с помощью коммутаторов локальных сечений распределения $\D$. Иными словами, распределение  $\D$ имеет вектор роста $(2, 3, 4)$.

\subsubsection{Классификация левоинвариантных энгелевых субримановых структур}

\begin{proposition}
По любой левоинвариантной  энгелевой субримановой структуре можно найти левоинвариантный репер $(X_1, \dots, X_4)$ на соответствующей группе Ли, такой, что $(X_1, X_2)$ есть ортонормированный репер структуры, с таблицей умножения 
\begin{align*}
	&[X_1,X_2]= X_3,\,\, [X_1,X_3] = X_4,\\
  &[X_1,X_4]=\frac12 A X_1+T_5 X_2+T_3 X_3+T_1 X_4, \\
	&[X_2,X_3]=T_6 X_1 + T_4 X_2 + T_2 X_3, \\
	&[X_2,X_4]=T_4 X_3 + T_2 X_4, \\
	&[X_3,X_4]=C X_1+B X_2-\frac12 A X_3+T_4 X_4,
\end{align*}
где $A=T_1 T_4$, $B=T_2 T_5-T_3 T_4$, $C=\frac12 T_1 T_2 T_4-T_3 T_6$.
\end{proposition}

\begin{theorem}
Любая левоинвариантная энгелева субриманова структура однозначно локально определяется структурными константами $T_i$ и принадлежит по крайней мере одному семейству  из таблицы {\em \ref{tab:class4}}. В этой таблице перечислены ограничения на $T_i$, определяющие семейства, а также соответствующие нетривиальные структурные уравнения.
\end{theorem}
\begin{table}[h]
\caption{Классификация левоинвариантных энгелевых субримановых структур}\label{tab:class4}
\begin{tabular}{| l | l|}  \hline        
\itshape Ограничения & \itshape Структурные уравнения, за исключением $[X_1,X_2]=X_3,\,\,[X_1,X_3]=X_4$  \\  \hline
I.  $T_2 = T_4 = T_6 = 0$  & $ [X_1,X_4] = T_5 X_2 + T_3 X_3 + T_1 X_4 $\\  \hline
II. $ T_4 = T_6 = T_5 = 0$ & $ [X_1,X_4] = T_3 X_3 + T_1 X_4$, \\ 
                           & $ [X_2,X_3] = T_2 X_3$,\\
                           & $ [X_2,X_4]= T_2 X_4$ 	\\ \hline
III. $ T_1 = T_2 = T_5 = 0$& $[X_1,X_4] = T_3 X_3$,\\
                           & $[X_2,X_3] =T_6 X_1 + T_4 X_2$, \\
                           & $  [X_2,X_4] = T_4 X_3$,\\
                           & $[X_3,X_4] = -T_6 T_3 X_1 -T_4 T_3 X_2 + T_4 X_4$\\		\hline
IV. $ T_1 = T_3 = 0$,      & $[X_2,X_3] =T_6 X_1 + T_2 X_3$,  \\
    $\quad \,\,\,\, T_4 = T_5 = 0$ & $[X_2,X_4] =  T_2 X_4$ \\		\hline
V.  $T_1\neq 0$,           & $[X_1, X_4] = T_1X_4-\frac{T_1^3+8T_5}{4T_1}X_3 +T_5X_2-\frac{2T_2T_5}{T_1}X_1$,\\
    $\quad \,\,\, T_4 = \frac12 \frac{T_2(T_1^2+4T_3)}{T_1}$,&$[X_2, X_3] = T_2X_3-\frac{4T_2T_5}{T_1^2}X_2+\frac{8T_2^2T_5}{T_1^3}X_1$, \\
$\quad \,\,\, T_5 = -\frac18 T_1^3-\frac12 	 T_1 T_3$,&$[X_2, X_4] = T_2X_4-\frac{4T_2T_5}{T_1^2}X_3$,\\
$\quad \,\,\, T_6 = -\frac{T_2^2(T_1^2+4T_3)}{T_1^2}$  &$ [X_3, X_4] = \frac{2T_2T_5}{T_1}\left(X_3 -\frac{2}{T_1} X_4-\frac{4T_5}{T_1^2} X_2+\frac{8T_2T_5}{T_1^3} X_1\right)$\\\hline												
\end{tabular}
\end{table}
\subsubsection{Интегрируемость и строгая анормальность}
\begin{theorem}
Рассмотрим левоинвариантную  энгелеву субриманову структуру типа {\em III} на группе Ли с алгеброй Ли
\begin{align*}
&[X_1,X_2]=X_3, & &[X_1,X_3]=X_4, 
\\
&[X_1,X_4] = T_3 X_3, & &[X_2,X_3] =T_6 X_1 + T_4 X_2, 
\\
&[X_2,X_4] = T_4 X_3, & &[X_3,X_4] = -T_6 T_3 X_1 -T_4 T_3 X_2 + T_4 X_4,
\end{align*}
где векторные поля $X_1$, $X_2$ образуют ортонормированный репер. Нормальный гамильтонов поток этой структуры суперинтегрируем  в том смысле, что он имеет четыре независимых коммутирующих первых интеграла, включая нормальный гамильтониан $H$, и еще один независимый первый интеграл, коммутирующий с $H$. Если $T_4 \neq 0$, то анормальные геодезические  этой структуры строгие, т.е. не являются нормальными.
\end{theorem}
\subsubsection{Сопряженные точки}
\begin{proposition}
Пусть $g(t)$ есть анормальная геодезическая левоинвариантной  энгелевой субримановой структуры, и пусть $\d =  T_6 - \frac 14(T_2)^2$. 

Если $\d > 0$, то все сопряженные времена имеют вид 
$$
t_{\conj} = \dfrac{\pi k}{\sqrt {\d}}, \quad k \in \N.
$$
Если $g(t)$ строго анормальна, то ограничение $g|_{[0, \tau]}$ есть $C^0$-локально кратчайшая  при  $\tau < \dfrac{\pi}{\sqrt{\d}}$, и не  является таковой при  $\tau > \dfrac{\pi}{\sqrt{\d}}$. 

Если $g(t)$ строго анормальна и  $\d \leq 0$, то ограничение 
$g|_{[0, \tau]}$ есть $C^0$-локально кратчайшая для любого $\tau >0$.
\end{proposition}

\subsubsection{Библиографические комментарии}
Результаты этого раздела получены в работе~\cite{BM}.

\section{Задачи, интегрируемые в элементарных функциях}\label{sec:elementary}
\subsection{Субриманова задача на группе Гейзенберга}\label{subsec:heis}
\subsubsection{Постановка задачи}
\paragraph{Задача Дидоны}
Рассмотрим следующую формализацию древнейшей    задачи оптимизации, восходящей к  IX веку до н.~э.~\cite{tikhomirov_stories}, \ddef{задачи Дидоны}. 

Пусть на евклидовой плоскости заданы точки $a_0, a_1 \in \R^2$, соединенные кривой $\g_0 \subset \R^2$. Пусть также задано число $S \in \R$. Требуется соединить точки $a_0$, $a_1$ кратчайшей кривой $\g \subset \R^2$ так, чтобы кривые $\g_0$ и $\g$ ограничивали на плоскости область алгебраической площади $S$.
\paragraph{Задача оптимального управления}
Эту геометрическую задачу можно переформулировать как задачу оптимального управления
\begin{align}
&\dg = u_1 X_1(g) + u_2 X_2(g), \quad g = (x, y, z) \in \R^3, \label{heisdg}\\
&g(0) = g_0, \quad g(t_1) = g_1, \label{heisg0}\\
&l=\int_0^{t_1}\sqrt{u_1^2 + u_2^2}\,dt \to \min, \label{heisl}\\
&X_1 = \dfrac{\partial}{\partial x} - \dfrac{y}{2}\dfrac{\partial}{\partial z} , \quad 
X_2 = \dfrac{\partial}{\partial y} + \dfrac{x}{2}\dfrac{\partial}{\partial z}. \label{heisX12}
\end{align}
Это субриманова задача    для субримановой структуры на $\R^3$, заданной ортонормированным репером $X_1$, $X_2$.
\paragraph{Алгебра Гейзенберга и группа Гейзенберга}
\ddef{Алгеброй Гейзенберга } называется трехмерная свободная нильпотентная алгебра Ли $\gg$ с двумя образующими, глубины 2.  Существует  базис $\gg = \spann(X_1, X_2, X_3)$, в котором единственная ненулевая скобка Ли есть
$$
[X_1, X_2] = X_3.
$$
Алгебра Гейзенберга имеет градуировку
$\gg = \gg^{(1)} \oplus\gg^{(2)} $, $ \gg^{(1)} = \spann(X_1, X_2)$, $ \gg^{(2)} = \R X_3$,  $[\gg^{(1)}, \gg^{(i)}] =  \gg^{(i+1)} $, 
$\gg^{(3)} = \{0\}$,
поэтому она является алгеброй Карно. Соответствующая связная односвязная группа Ли $G$ называется \ddef{группой Гейзенберга}.

Группа  Гейзенберга имеет линейное представление
\begin{align*}
G=\left\{
\begin{pmatrix}
1 & x & z +\frac{xy}{2} \\
0 & 1 & y\\
0 & 0 & 1 
\end{pmatrix}
\mid x, y, z \in \R^3
\right\},
\end{align*}
дающее закон умножения в этой группе:
\begin{align*}
\begin{pmatrix}
x_1\\
y_1\\
z_1
\end{pmatrix} \cdot
\begin{pmatrix}
x_2\\
y_2\\
z_2
\end{pmatrix} = 
\begin{pmatrix}
x_1 + x_2\\
y_1 + y_2\\
z_1 + z_2 + (x_1 y_2 - x_2 y_1)/2
\end{pmatrix}.
\end{align*}
Векторные поля \eq{heisX12} левоинвариантны на группе Ли $G$. Поэтому задача \eq{heisdg}--\eq{heisX12} есть левоинвариантная субриманова задача на группе Гейзенберга. Это --- простейшая субриманова задача, не являющаяся римановой.

Неголономная левоинвариантная субриманова задача на группе Гейзенберга единственна, с точностью до изоморфизма этой группы~\cite{symmetry}.

В силу левоинвариантности задачи можно считать $g_0 = \Id = (0, 0, 0)$ в \eq{heisg0}.

\subsubsection{Симметрии распределения и субримановой структуры }

Пусть $(\D, \langle\cdot, \cdot\rangle)$  есть субриманова структура на гладком многообразии $M$.  Векторное поле $X \in \mathrm{Vec}(M)$  называется \ddef{инфинитезимальной симметрией}:
\begin{itemize}
\item[(1)]
 распределения $\D$, если его поток $\map{e^{tX}}{M}{M}$  сохраняет $\D$:
$$
\left(e^{tX}\right)_* \D = \D, \qquad t \in \R;
$$
\item[(2)]
субримановой структуры $(\D, \langle\cdot, \cdot\rangle)$, если его поток  сохраняет $\D$ и $\langle\cdot, \cdot\rangle$:
$$
\left(e^{tX}\right)_* \D = \D, 
\quad
\left(e^{tX}\right)^* \langle\cdot, \cdot\rangle = \langle\cdot, \cdot\rangle,
\qquad t \in \R.
$$
\end{itemize}
Пространство инфинитезимальных симметрий распределения (субримановой структуры) есть алгебра Ли.

\begin{theorem}
\begin{itemize}
\item[$(1)$]
Алгебра Ли инфинитезимальных симметрий распределения $\spann(X_1, X_2)$ параметризуется гладкими функциями трех переменных.
\item[$(2)$]
Алгебра Ли инфинитезимальных симметрий субримановой структуры с ортонормированным репером $X_1$, $X_2$ есть четырехмерная алгебра Ли $\spann{(X_0, Y_1, Y_2, Y_3)}$ с таблицей умножения
$$
[X_0, Y_1] = - Y_2, \quad [X_0, Y_2]= Y_1, \quad [Y_1, Y_2] = Y_3.
$$
Векторные поля $Y_1$, $Y_2$, $Y_3$ образуют правоинвариантный репер на группе Гейзенберга, а поле $X_0$ определяет вращение:
$$
X_0 = -y\dfrac{\partial}{\partial x} + x \dfrac{\partial}{\partial y}.
$$
\end{itemize}
\end{theorem}
\subsubsection{Геодезические}
Существование оптимальных управлений в задаче \eq{heisdg}--\eq{heisX12} следует из теорем Рашевского-Чжоу и Филиппова.

Анормальные траектории постоянны.

Для параметризации нормальных экстремалей введем линейные на слоях $T^*G$ гамильтонианы $h_i(\lam) = \lan \lam, X_i \ran$, $i = 1, 2, 3$, $X_3 = [X_1, X_2] = \dfrac{\partial}{\partial z}$. 
Нормальные экстремали суть траектории гамильтонова поля $\Vec H$, где $H = \frac 12 (h_1^2 + h_2^2)$. Натурально параметризованные экстремали принадлежат поверхности уровня $\{H = \frac 12  \}$. Введем на этой поверхности координату $\t$:
$$
h_1 = \cos \t, \quad h_2 = \sin \t.
$$
Натурально параметризованные экстремали
в случае $h_3 = 0$  имеют вид
\begin{align}
&\theta \equiv \theta_0, \qquad h_3 \equiv 0, \label{heisg1}\\
&x = t \cos \theta_0, \qquad y = t \sin \theta_0, \qquad z = 0, \label{heisg2}
\end{align}
а в случае $h_3\neq 0$
\begin{align}
&\theta = \theta_0 + h_3 t, \qquad h_3 \equiv \const, \label{heisg3}\\
&x = (\sin (\theta_0 + h_3 t) - \sin \theta_0)/h_3, \label{heisg4}\\
&y = (\cos \theta_0 - \cos (\theta_0 + h_3 t))/h_3, \label{heisg5}\\
&z = (h_3 t - \sin h_3 t)/(2h_3^2). \label{heisg6}
\end{align}
При $h_3 = 0$ геодезические \eq{heisg2} суть прямые в плоскости $\{ z = 0 \}$, а в случае $h_3\neq 0$ геодезические \eq{heisg4}--\eq{heisg6} суть спирали с переменным наклоном, проецирующиеся на плоскость $(x, y)$ в окружности.

Формулы \eq{heisg1}--\eq{heisg6} дают параметризацию экспоненциального отображения 
\begin{align*}
&\Exp: \, (\t_0, h_3, t) \mapsto (x, y, z),\\
&\Exp: \, C\times \R_+ \to G, \quad C = \gg^* \cap\{ H = 1/2 \}.
\end{align*}
В случае $h_3 = 0$ геодезические $g(t)$, $t \in [0, t_1]$, оптимальны для любого $t_1 > 0$. Эти геодезические $g(t)$, $t \in \R$, и только они, являются  метрическими прямыми.

\subsubsection{Сопряженные времена}
\begin{theorem}
Пусть $(\t_0, h_3) \in C$  и $g(t) = \Exp(\t_0, h_3, t)$.

Если $h_3 = 0$, то на геодезической $g(t)$, $t > 0$,  нет сопряженных точек.

Если $h_3 \neq 0$, то
сопряженные времена вдоль геодезической $g(t)$, $t > 0$, имеют вид:
\begin{align*}
&\dfrac{2 \pi k}{|h_3|} \quad \text{ и } \quad 
\dfrac{2 p_k}{|h_3|}, \qquad k \in \N,
\end{align*}
где $p_k \in (\pi k, \, \frac{\pi}{2} + \pi k) $ есть $k$-й положительный корень уравнения
$$
(2p - \sin 2p)\cos p - (1 - \cos 2p)\sin p = 0.
$$
\end{theorem}

Поэтому первое сопряженное время равно
$$
\tconj = \dfrac{2\pi}{|h_3|},
$$
а
первая каустика есть
$$
\Conj^1 = \{x = y= 0, \quad z \neq 0  \}.
$$
\subsubsection{Время разреза и множество разреза}
\begin{theorem}
Время разреза вдоль геодезической $g(t) = \Exp(\t_0, h_3, t)$ имеет вид:
\begin{align*}
&t_{\cut} = +\infty \quad \text{при} \quad h_3 = 0,\\
&t_{\cut} =\dfrac{2 \pi}{|h_3|} \quad \text{при} \quad h_3 \neq 0.
\end{align*}
\end{theorem}

Поэтому
$$
\tcut(\lam) = \tconj(\lam), \qquad \lam \in C,
$$
и
множество разреза совпадает с первой каустикой:
$$
\Cut = \Conj^1 = \{x = y= 0, \quad z \neq 0  \}.
$$

Геодезические-спирали теряют оптимальность (как локальную, так и глобальную) при первом после начала пересечении  с осью $z$; то есть они оптимальны вплоть до первого витка окружности $(x(t), y(t))$.
\subsubsection{Оптимальный синтез}
Пусть $g_1 = (x_1, y_1, z_1) \neq (0, 0, 0)$ в граничных условиях \eq{heisg0}. Опишем соответствующие решения задачи \eq{heisdg}--\eq{heisX12}.

Если $z_1 = 0$, $x_1^2 + y_1^2 \neq 0$, то кратчайшая есть прямолинейный  отрезок \eq{heisg2}, где
$$
t \in [0, t_1], \quad x_1 = t_1 \cos \theta_0, \quad y_1 = t_1 \sin \theta_0.
$$

Если $z_1 \neq 0$, $x_1^2 + y_1^2 \neq 0$, то кратчайшая есть спираль \eq{heisg4}--\eq{heisg6}, где $t \in [0, t_1]$,
\begin{align*}
&\frac{x_1^2 + y_1^2}{|z_1|} = \frac{8 \sin^2 p_1}{2 p_1 - \sin 2 p_1}, \qquad p_1 \in (0, \pi),\\
&h_3 = \sgn z_1 \sqrt{\frac{2 p_1 - \sin 2 p_1}{2 |z_1|}},\\
&t_1 = \frac{2 p_1}{|h_3|},\\
&x_1 = \frac{2 \sin p_1}{h_3} \cos \tau_1, \qquad y_1 = \frac{2 \sin p_1}{h_3} \sin \tau_1,\\
&\theta_0 = \tau_1 - p_1.
\end{align*}

Наконец, если $z_1 \neq 0$, $x_1^2 + y_1^2 = 0$,  то существует однопараметрическое семейство кратчайших \eq{heisg4}--\eq{heisg6}, где
$t \in \left[0, \frac{2 \pi}{|h_3|}\right]$, $h_3 = \sgn z_1 \sqrt{\frac{\pi}{|z_1|}}$, $\theta_0 \in S^1$.
\subsubsection{Субриманово расстояние и сферы}
Пусть $g = (x, y, z) \in G$. Тогда субриманово расстояние $d_0(g) = d(\Id, g)$ представляется следующим образом.

Если $z = 0$, то $d_0(g) = \sqrt{x^2 + y^2}$.

Если $z \neq 0$, $x^2 + y^2 \neq 0$, то 
\begin{align}
&d_0(g) = \frac{p}{\sin p} \sqrt{x^2 + y^2},\\
&\frac{2 p - \sin 2 p}{8 \sin^2 p} = \frac{|z|}{x^2 + y^2}, \qquad p \in (0, \pi).
\end{align}

Если $z \neq 0$, $x^2 + y^2 = 0$, то $d_0(g) = 2 \sqrt{\pi|z|}$.

Субриманова сфера радиуса $R$ с центром в единице есть поверхность вращения с параметрическими уравнениями 
$$
x = R \frac{\sin p}{p} \cos \tau, \quad y = R \frac{\sin p}{p} \sin \tau, \quad z = R^2 \frac{2 p - \sin 2 p}{8 p^2},
\qquad p \in [- \pi, \pi], \quad \tau \in S^1, 
$$
она похожа на яблоко  и имеет две особые конические точки $(x,y,z) = (0, 0, \pm R^2/(4\pi))$.
Сферы сохраняются вращениями $X_0$ и растягиваются дилатациями $Y = x\dfrac{\partial}{\partial x} + y\dfrac{\partial}{\partial y} + 2z\dfrac{\partial}{\partial z}$:
$$
e^{r X_0} (S_{\Id} (R)) = S_{\Id} (R), \qquad e^{r Y} (S_{\Id} (R)) = S_{\Id} (e^r R).
$$
Единичная сфера $S_{\Id} (1)$ и ее половина изображены на Рис.  \ref{fig:heis2}  и  \ref{fig:heis3}  соответственно.

\figout{
\twofiglabelsize{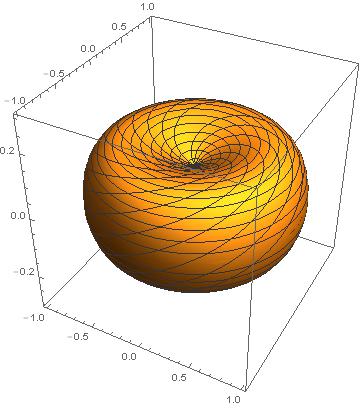}{Субриманова сфера  на группе Гейзенберга}{fig:heis2}{0.5}
{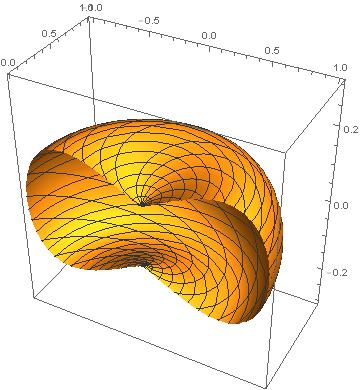}{Субриманова полусфера  на группе Гейзенберга}{fig:heis3}{0.5}
}
\subsubsection{Библиографические комментарии}
Субриманова задача на группе Гейзенберга  описана практически в каждой книге или обзоре по субримановой и неголономной геометрии, см.~\cite{jurd_book, SRG, montgomery_book, capogna, rifford, ABB_book, intro}. По-видимому, первое детальное исследование этой задачи было выполнено в работах~\cite{gaveau, brock, versh_gersh}.

В работе \cite{versh_gersh} исследован геодезический поток для контактных левоинвариантных субримановых структур на трехмерных группах Ли (в том числе на группе Гейзенберга) методами теории динамических систем.

\subsection{Машина Маркова-Дубинса}\label{subsec:markov}
\subsubsection{Постановка задачи}

Рассмотрим модель машины, движущейся по плоскости. Состояние машины задается ее положением и ориентацией на плоскости. 
Машина может ехать вперед  с постоянной линейной скоростью и одновременно поворачиваться с ограниченной угловой скоростью.  Требуется перевести машину из заданного начального состояния в заданное конечное состояние за минимальное время.

После выбора подходящих единиц измерения задача формулируется как задача быстродействия
\begin{align}
&\dot x = \cos \t, \qquad g = (x,y,\t) \in \R^2 \times S^1, \label{markov1}\\
&\dot y = \sin \t, \qquad u \in [-1, 1], \label{markov2}\\
&\dot \t = u,   \label{markov3}\\
&g(0) = g_0, \quad g(t_1) = g_1,  \label{markov4}\\
&t_1 \to \min.  \label{markov5}
 \end{align}
Это левоинвариантная задача на группе Ли $G = \SE(2) \cong \R^2 \ltimes S^1$,  см. раздел \ref{subsec:se2},  поэтому можно положить $g_0 = \Id = (0, 0, 0)$. В терминах левоинвариантных векторных полей на этой группе Ли
\be{markovX12}
X_1 = \cos \t \pder{}{x} + \sin \t \pder{}{y}, \qquad X_2 = \pder{}{\t}
\ee
управляемая система \eq{markov1}--\eq{markov3}  записывается как
\be{markov6}
\dot g = X_1 + u X_2, \qquad g \in G, \quad u \in [-1, 1].
\ee

\subsubsection{Управляемость}
\begin{theorem}
Множество достижимости системы \eq{markov6}  из точки $\Id$  есть вся группа $\SE(2)$.
\end{theorem}

\subsubsection{Оптимальные траектории}

Существование решений в задаче \eq{markov1}--\eq{markov4} следует из теоремы Филиппова.

Введем линейные на слоях $T^*G$  гамильтонианы $h_i(\lam) = \lan \lam, X_i\ran$, $i = 1, 2, 3$, $\lam \in T^*G$, $X_3 = [X_1, X_2]$. 

Из принципа максимума Понтрягина получаем гамильтонову систему
\begin{align*}
&\dot h_1 = - u h_3, \quad \dot h_2 = h_3, \quad \dot h_3 = u h_1, \\
&\dot g = X_1 + u X_2
\end{align*}
 и условие максимума
$$
u(t) h_2(t) = \max_{|v| \leq 1} v h_2(t).
$$

Если на некотором промежутке $h_2(t) \neq 0$, то $u \equiv \sgn h_2(t) = \pm 1$, и экстремальная кривая есть дуга окружности.

Если на некотором промежутке $h_2(t) \equiv 0$, то $u \equiv 0$, и экстремальная кривая есть  прямолинейный отрезок.

\begin{theorem}\label{th:markov}
Оптимальные траектории могут быть одного из следующих двух типов:
\begin{itemize}
\item[$(1)$]
конкатенация дуги окружности единичного радиуса, прямолинейного отрезка, и дуги окружности единичного радиуса,
\item[$(2)$]
конкатенация не более чем трех дуг окружностей единичного радиуса.
\end{itemize}
В случае $(2)$, если $a,b,c$  суть времена движения по дугам, то $\pi < b < 2 \pi$, $\min(a,c) < b - \pi$, $\max(a,c) <b$.
\end{theorem}

\subsubsection{Библиографические комментарии}
Впервые версию задачи \eq{markov1}--\eq{markov5}  рассмотрел в 1887 г.   А.А. Марков \cite{markov}. 

Детально эту задачу изучил в 1957 г. Л. Дубинс \cite{dubins}.  Он показал, что оптимальная траектория принадлежит одному из 6 типов конкатенаций дуг единичных окружностей и прямолинейных отрезков.

Подробный анализ этой задачи методами геометрической теории управления, включая теорему \ref{th:markov},  приведен в \cite{sussmann_tang}.

Анализ и приложение задачи Маркова-Дубинса к управлению движением самолетов см. в работе \cite{pecsvaradi}.

Программная реализация оптимального синтеза в этой задаче описана в работе \cite{ard_gub}.

\subsection{Машина Ридса-Шеппа}
\subsubsection{Постановка задачи}

Рассмотрим вариацию модели машины из раздела \ref{subsec:markov}.
Пусть теперь
машина может ехать вперед  или назад с постоянной линейной скоростью и одновременно поворачиваться с ограниченной угловой скоростью.  Требуется перевести машину из заданного начального состояния в заданное конечное состояние за минимальное время.

Задача формулируется как задача быстродействия
\begin{align}
&\dot x = v \cos \t, \qquad g = (x,y,\t) \in \R^2 \times S^1, \label{reeds1}\\
&\dot y = v \sin \t, \qquad  |v| = 1, \quad u \in [-1, 1], \label{reeds2}\\
&\dot \t = u,   \label{reeds3}\\
&g(0) = g_0, \quad g(t_1) = g_1,  \label{reeds4}\\
&t_1 \to \min.  \label{reeds5}
 \end{align}
Это левоинвариантная задача на группе Ли $G = \SE(2) \cong \R^2 \ltimes S^1$.  Управляемая система \eq{reeds1}--\eq{reeds3} в терминах левоинвариантных векторных полей \eq{markovX12}
  записывается как
\be{reeds6}
\dot g = v X_1 + u X_2, \qquad g \in G, \quad  |v| = 1, \quad u \in [-1, 1].
\ee

\subsubsection{Управляемость}
\begin{theorem}
Множество достижимости системы \eq{reeds6}  из точки $\Id$  есть вся группа $\SE(2)$.
\end{theorem}

\subsubsection{Существование оптимальных траекторий}

Множество значений управляющего параметра для машины Ридса-Шеппа невыпукло (см. \eq{reeds2}), поэтому общие  теоремы существования оптимальных управлений для нее неприменимы. Однако  теорема существования справедлива.

\begin{theorem}
В задаче \eq{reeds1}--\eq{reeds5} оптимальное управление существует.
\end{theorem}

\subsubsection{Оптимальные траектории}
\begin{theorem}
Для любой точки $g_1 \in \SE(2)$ оптимальная траектория
может быть выбрана одного из следующих двух типов:
\begin{itemize}
\item[$(1)$]
конкатенация не более двух дуг окружностей единичного радиуса, прямолинейного отрезка, и не более двух дуг окружностей единичного радиуса,
\item[$(2)$]
конкатенация не более чем четырех дуг окружностей единичного радиуса.
\end{itemize}
\end{theorem}

\subsubsection{Библиографические комментарии}
Задачу \eq{reeds1}--\eq{reeds5} впервые рассмотрели Дж. Ридс и Л. Шепп в работе \cite{reeds}. 
Они показали, что оптимальная траектория принадлежит одному из 48 типов конкатенаций дуг единичных окружностей и прямолинейных отрезков.

Эта задача детально исследована в  работе~\cite{sussmann_tang}:
приведено 
семейство траекторий, содержащее оптимальные траектории в задаче \eq{reeds1}--\eq{reeds5}, удовлетворяющие принципу максимума Понтрягина и условиям оптимальности высших порядков.
В этой работе количество типов оптимальных траекторий уменьшено до 46.

Полный оптимальный синтез построен в работе \cite{soueres_laumond}, см. также \cite{soueres}.

Этой задаче посвящена также работа \cite{boissonat}.

Программная реализация оптимального синтеза в задаче Ридса-Шеппа описана в работе \cite{ard_gub}.

\subsection{Субриманова задача с вектором роста $(3,6)$}\label{subsec:36}
\subsubsection{Постановка задачи и две модели}\label{subsubsec:36_state}
Левоинвариантная субриманова задача с вектором роста $(3,6)$ ставится следующим образом:
\begin{align}
&\dx = u, \quad x, u \in \R^3, \label{36_1}\\
&\dy = x\land u, \quad y \in \R^3  \land \R^3, \label{36_2}\\
&x(0) = 0, \quad y(0) = 0,\quad x(t_1) = x_1, \quad y(t_1) = y_1, \label{36_3}\\
&\int_0^{t_1} (u_1^2 + u_2^2 + u_3^2)^{1/2} dt \to \min.\label{36_4}
\end{align}
Отождествляя  пространство $\R^3  \land \R^3$  с $\R^3$ с помощью оператора $* :\R^3  \land \R^3 \to \R^3$, $y \mapsto z$, где $y \land z$ есть форма объема на $\R^3$, можно заменить внешнее произведение $x \land u$ на векторное  произведение $[x, u]$ в $\R^3$, и получить постановку 
\begin{align*}
&\dx = u, \quad x, y, u \in \R^3, \\
&\dy = [x, u], \\
&x(0) = 0, \quad y(0) = 0,\quad x(t_1) = x_1, \quad y(t_1) = y_1,\\
&\int_0^{t_1} (u_1^2 + u_2^2 + u_3^2)^{1/2} dt \to \min.
\end{align*}
В координатах $x = (x_1, x_2, x_3)$, $y=(y_1, y_2, y_3) \in \R^3$ ортонормированный репер для соответствующей субримановой структуры  имеет вид:
\begin{align}
&X_1 = \dfrac{\partial}{\partial x_1} + x_3\dfrac{\partial}{\partial y_2}  - x_2 \dfrac{\partial}{\partial y_3}, \label{36X1}\\
&X_2 = \dfrac{\partial}{\partial x_2} + x_1 \dfrac{\partial}{\partial y_3} - x_3 \dfrac{\partial}{\partial y_1},\label{36X2}\\
&X_3 = \dfrac{\partial}{\partial x_3  } + x_2 \dfrac{\partial}{\partial y_1} - x_1 \dfrac{\partial}{\partial y_2}. \label{36X3}
\end{align}

Ненулевые скобки Ли в алгебре Ли, порожденной полями $X_i$, имеют вид:
\begin{align}
&[X_1, X_2] = X_{12}, \quad [X_2, X_3]= X_{23}, \quad [X_3, X_1] = X_{31}, \label{36_table}
\end{align}
где $X_{12} = 2 \dfrac{\partial}{\partial y_3}$, $X_{23} = 2 \dfrac{\partial}{\partial y_1}$,  $X_{31} = 2 \dfrac{\partial}{\partial y_2}$.
Будем называть ортонормированный репер \eq{36X1}--\eq{36X3}  первой моделью субримановой $(3, 6)$-структуры.

Вторая модель дается векторными полями
\begin{align}
&\widetilde X_1 = \dfrac{\partial}{\partial x_1} + \dfrac{x_3}{2}\dfrac{\partial}{\partial y_2}  - \dfrac{x_2}{2} \dfrac{\partial}{\partial y_3}, \label{36tX1}\\
&\widetilde X_2 = \dfrac{\partial}{\partial x_2} +\dfrac{x_1}{2} \dfrac{\partial}{\partial y_3} - \dfrac{x_3}{2} \dfrac{\partial}{\partial y_1},\label{36tX2}\\
&\widetilde X_3 = \dfrac{\partial}{\partial x_3  } + \dfrac{x_2}{2} \dfrac{\partial}{\partial y_1} - \dfrac{x_1}{2} \dfrac{\partial}{\partial y_2}. \label{36tX3}
\end{align}

\subsubsection{Симметрии}\label{subsubsec:36_symmetr}
Задача \eq{36_1}--\eq{36_4} инвариантна относительно группы $\SO(3)$, естественно действующей на управления и состояния:
\begin{align*}
&g : (u, x, v\land w)\mapsto (g(u), g(x), g(v)\land g(w)),\\
&g : \R^3\times \R^3\times(\R^3\land \R^3) \to \R^3 \times \R^3 \times (\R^3\land \R^3),\\
&g \in \SO(3).
\end{align*}

\subsubsection{Геодезические в первой модели }\label{subsubsec:36_geod1}
Из таблицы умножения \eq{36_table} следует, что распределение $\D = \spann(X_1, X_2, X_3)$ удовлетворяет условию Гоха $\D^2_{q}= T_qG$, $q \in  G \cong  \R^6$. Поэтому все субримановы кратчайшие  нормальны.
В нормальном случае гамильтониан принципа максимума Понтрягина имеет вид
$$
H = \dfrac{p^2}{2} + (\r, [x, p]) + \dfrac{x^2 \r^2}{2} - \dfrac{(\r, x)^2}{2}, \quad (x, y, p, \r)\in T^* G,
$$
где $(\cdot, \cdot)$ и $[\cdot, \cdot]$ суть скалярное и векторное произведения в $\R^3$. Соответствующая гамильтонова система есть
\begin{align*}
&\dx = p + [\r, x],\\
&\dy = [x, p]+ x^2\cdot \r - (\r,x)\cdot x,\\
&\dot p = -[p, \r] - \r^2\cdot x + (\r, x)\cdot \r,\\
&\dot\r = 0,\\
&x(0) = y(0) = 0, \quad p_0^2 = 1.
\end{align*}

Поэтому $\r \equiv \r_0$, $p = [\r_0, x] + p_0$.

Выберем декартовы координаты $(x_1, x_2, x_3) \in \R^3$ так, чтобы ось $x_3$ была сонаправлена оси $\r_0$, ось $x_1$ принадлежала плоскости $\spann(\r_0, p_0)$, а ось $x_2$ так, чтобы система координат была правосторонней. Тогда горизонтальная подсистема нормальной гамильтоновой системы принимает форму 
\begin{align*}
&\dx_1 = \sin \f - r x_2,\\
&\dx_2 = r x_1,\\
&\dx_3 = \cos \f,\\
&\dy_1 = x_2 \cos \f - r x_1 x_3,\\
&\dy_2 = x_3\sin \f - x_1 \cos\f - r x_2x_3,\\
&\dy_3 = -x_2 \sin \f  + (x_1^2 + x_2^2)\,r,\\
&r = 2|\r_0|, \quad x(0) = y(0) = 0.
\end{align*}

Поэтому в первой модели \eq{36X1}--\eq{36X3}
\begin{align*}
&x_1 = \dfrac{1}{r}\sin \f \sin (r t),\\
&x_2 = \dfrac{1}{r}\sin \f(1 - \cos(r t)),\\
&x_3 = t \cos \f,\\
&y_1 = \dfrac{1}{r^2}\sin\f\cos\f(r t(1 + \cos(rt)) - 2\sin(rt)),\\
&y_2 = \dfrac{1}{r^2}\sin\f\cos\f(r t\sin(rt) + 2\cos(rt)-2),\\
&y_3 =\dfrac{1}{r^2}\sin^2\f\,(r t - \sin(rt)).
\end{align*}

\subsubsection{Геодезические во второй модели }\label{subsubsec:36_geod2}
Во второй модели \eq{36tX1}--\eq{36tX3} экстремальные управления, соответствующие геодезическим постоянной скорости, имеют вид
\begin{align*}
&u(t) = a \cos(ct) + b \sin(ct) +z, \quad |a| = |b| > 0, \quad c \geq 0,
\end{align*}
где $ a, b, z $ взаимно ортогональные векторы в $\R^3$.
Соответствующие  геодезические  суть
\begin{multline*}
q(t) = \left(\dfrac{\sin(ct)}{c} a + \dfrac{1 - \cos(c t)}{c} b + t z \right., \\ 
\left. \dfrac{c t - \sin(c t)}{2 c^2} a \land b + \dfrac{2(1 - \cos(c t)) -  c t \sin (c t)}{2 c^2}a \land z + \dfrac{c t( 1 + \cos(c t)) - 2 \sin(c t)}{2 c^2}b \land z  \right), \,\,c > 0,
\end{multline*}
и
\begin{align}
&q(t) = (t(a + z), 0), \quad c = 0. \label{36abnorm}
\end{align}
При $c = 0$ геодезические суть однопараметрические подгруппы в группе Карно $G = \R^3 \times (\R^3 \land \R^3)$, лежащие в ее первом слое $\R^3 \times \{0\}$.

\subsubsection{Анормальные геодезические и анормальное множество}\label{subsubsec:36_abn}
С точностью до перепараметризации, анормальные  управления постоянны, а анормальные геодезические \eq{36abnorm} суть однопараметрические подгруппы  в $\R^3 \times (\R^3\land\R^3)$, лежащие в первом слое 
$\R^3\times \{ 0\}$.

Анормальное множество, соответствующее единичному элементу $\Id \in \R^3\times (\R^3\land\R^3)$, есть первый слой группы Карно:
$$
\Abn = \R^3\times \{ 0\} = \{(x, 0) \in  \R^3\times (\R^3\land\R^3) \}.
$$
\subsubsection{Время разреза}\label{subsubsec:36_tcut}
Определим следующие функции:
\begin{align*}
&S(\t) = \dfrac{\sin \t}{\t}, \quad U(\t) = \dfrac{\t - \sin \t \cos \t}{4 \t^2}, \quad V(\t) = \dfrac{\sin \t - \t\cos \t}{2 \t^2}, \quad Q(\t) = - \dfrac{U(\t)S(\t)}{V(\t)},
\end{align*}
и пусть $\f_1 \in(\pi, \frac 32 \pi)$ есть первый положительный корень функции $V$. Функция $Q(\t)$ неотрицательна и строго возрастающая при $\t \in [\pi, \f_1)$. Более того, $Q(\pi) = 0$ и $Q(\f_1 - 0) = +\infty$. Поэтому определена обратная функция
$$
Q^{-1}: [0, \, +\infty) \to [\pi, \f_1).
$$

\begin{theorem}
Время разреза для геодезических равно
\begin{align*}
&t_{\cut}(a, b, z, c) = \dfrac{2}{c}\,Q^{-1}\left(\dfrac{|z|^2}{|a|^2}\right), \quad c > 0,\\
&t_{\cut}(a, b, z, 0) = +\infty.
\end{align*}
\end{theorem}
Поэтому метрическими прямыми являются лишь анормальные геодезические \eq{36abnorm} --- однопараметрические подгруппы в первом слое группы Карно.

\subsubsection{Множество разреза}\label{subsubsec:36_cut}
\begin{theorem} \label{th:36_cut}
Множество разреза в данной задаче есть
$$
\Cut = \{(x, y) \in \R^3\times \R^3 \mid y \neq 0, \quad \exists a \in \R: x = a y   \}.
$$
\end{theorem}
Определим следующие функции:
\begin{align*}
&P(\t) = - \dfrac{S(\t)}{V(\t)}\, \sqrt{\dfrac{W(\t)}{U(\t)}}, \quad R(\t) = \dfrac{1 - S^2(\t)}{\sqrt{U(\t)W(\t)}},\\
&W(\t) = U(\t) - S(\t)V(\t).
\end{align*}
Функция $P: [\pi, \f_1) \to [0, + \infty)$  есть возрастающая биекция, поэтому определена обратная  функция $P^{-1}:[0, +\infty) \to [\pi, \f_1)$.
\begin{theorem} \label{th:36_dist}
Пусть $(x,y) \in \Cut$. Тогда 
\begin{align*}
& d^2((0, 0), (x, y))= |x|^2 + R(\t)|y|,\\
&\t = P^{-1}\left(\dfrac{|x|^2}{|y|}  \right) \in [\pi, \f_1).
\end{align*}
\end{theorem}
\subsubsection{Библиографические комментарии}
Изложение в этом разделе опирается на независимые работы~\cite{myasn36} и \cite{montanari_morbidelli}. Пункты~\ref{subsubsec:36_state} (первая модель), \ref{subsubsec:36_symmetr}, \ref{subsubsec:36_geod1}, 
\ref{subsubsec:36_cut} (теорема~\ref{th:36_cut}) опираются на работу~\cite{myasn36}. Пункты \ref{subsubsec:36_state} (вторая модель), \ref{subsubsec:36_geod2}, \ref{subsubsec:36_abn},
\ref{subsubsec:36_tcut}, \ref{subsubsec:36_cut} (теорема~\ref{th:36_dist}) опираются на работу~\cite{montanari_morbidelli}.

\subsection[Субриманова задача  на двухступенных свободных нильпотентных группах Ли]{Субриманова задача  на двухступенных свободных нильпотентных \\группах Ли}\label{subsec:twostep}
\subsubsection{Постановка задачи}\label{subsubsec:twostep_state}
\paragraph{Алгебры Ли и группы Ли }
Двухступенная свободная нильпотентная алгебра Ли $\gg$ удовлетворяет соотношениям 
$$
\gg = \gg^{(1)} \oplus\gg^{(2)}, \quad [\gg^{(1)}, \gg^{(1)}] = \gg^{(2)}, \quad [\gg^{(1)}, \gg^{(2)}] = [\gg^{(2)}, \gg^{(2)}] = \{0\}.
$$
Она имеет базис
$$
\gg = \spann\{X_i, X_{jk}  \mid 1\leq i \leq n, \quad 1 \leq j < k \leq n\},
$$
в котором таблица умножения имеет вид
$$
[X_i, X_j] = X_{ij}, \quad 1 \leq i < j \leq n.
$$
Эта алгебра Ли имеет размерность $n(n+1)/2$, где $n = \dim \gg^{(1)}$.

Обозначим через $G$ связную односвязную  группу Ли с алгеброй Ли $\gg$.  Эта группа Ли моделируется  пространством $\R^n \times(\R^n\land\R^n)\cong\R^n\times \so(n)$ с законом умножения
$$
(x^1, y^1)\cdot (x^2, y^2) = (x^1 + x^2, y^1+y^2 - x^1\land x^2), \quad (x^i, y^i)\in \R^n\times\so(n),
$$
где $x^1\land x^2 = x^1\otimes (x^2)^T - x^2\otimes (x^1)^T$. Поэтому будем далее считать $G = \R^n \times\so(n)$, так что любой элемент $(x, y) \in G$ имеет координатное представление 
$(x_1, \dots, x_n; y_{12}, \dots, y_{(n-1)n})$. В этих координатах следующие векторные поля образуют левоинвариантный репер на $G$:
\begin{align*}
&X_i = \dfrac{\partial}{\partial x_i} - \sum_{j=1, j\neq i}^{n} x_j \dfrac{\partial}{\partial y_{ij}}, \quad i = 1,\dots,n,\\
&X_{ij} = 2\dfrac{\partial}{\partial y_{ij}}, \quad 1\leq i < j \leq n.
\end{align*}

\paragraph{Субриманова задача}
Рассмотрим субриманову структуру на группе Ли $G$ с ортонормированным репером $(X_1, \dots, X_n)$. Соответствующая задача оптимального управления имеет вид
\begin{align}
&\dx = u, \qquad x, u \in \R^n, \label{n(n+1)x}\\
&\dy = x \land u, \qquad y \in \so(n), \label{n(n+1)y}\\
& x(0) =0, \quad y(0) = 0, \quad x(t_1) = x^1, \quad y(t_1) = y^1,\label{n(n+1)in}\\
&\int_0^{t_1}|u|\, dt \to \min.\label{n(n+1)l}
\end{align}
В координатах система \eq{n(n+1)x}, \eq{n(n+1)y} имеет вид
\begin{align*}
&\dx_i = u_i, \qquad i = 1,\dots n,\\
&\dy_{ij}  = x_i u_j - x_ju_i,\qquad 1\leq i < j \leq n.
\end{align*}

Задача \eq{n(n+1)x}--\eq{n(n+1)l} впервые рассматривалась Б.~Гаво~\cite{gaveau}, а затем Р.~Брокеттом~\cite{brock80},  поэтому она называется в литературе  \ddef{задачей Гаво-Брокетта}.

\paragraph{Существование решений}
Субримановы кратчайшие существуют по теоремам Рашевского-Чжоу и Филиппова.

\paragraph{Симметрии}
Задача инвариантна  относительно естественного действия группы $\O(n)$:
$$
g: (u, x, v\land w) \mapsto (g(u), g(x), g(v)\land g(w)), \quad  u, x, v, w \in \R^n, \quad g \in \O(n).
$$
Или, иными словами,
\begin{align}
&g: (u, x, y) \mapsto \left(g(u), g(x), gyg^T\right), \qquad u, x \in \R^n, \quad y \in \so(n).\label{n(n+1)O(n)}
\end{align}

\subsubsection{Экстремали}\label{subsubsec:twostep_extrem}
Введем линейные на слоях $T^*G$ гамильтонианы, соответствующие базисным полям:
$$
h_i(\lam) = \lan \lam, X_i  \ran, \quad h_{ij}(\lam)=\lan\lam, X_{ij}\ran,
$$
и положим
$$
h = (h_1, \dots, h_n)\in \R^n, \quad \H = (h_{ij})\in\so(n).
$$

\begin{theorem}[ПМП]
Пусть $\lam_t = (x(t), y(t), h(t), \H(t))$ есть экстремаль, соответствующая оптимальному управлению $u(t)$.
\begin{itemize}
 \item[$(1)$] 
Если $\lam_t$	 анормальна, то $h\equiv 0$, $\H\equiv \const \neq 0$ и $\H u(t)\equiv 0 $.

Более того, если $n$ четно и все собственные значения $\H$ отличны от нуля, то $(x(t), y(t))\equiv (0, 0)$. Если $n$ нечетно и $\H$ имеет только одно нулевое собственное значение, то $u(t)$ есть постоянный вектор, с точностью до перепараметризации времени.
\item[$(2)$] 
Если $\lam_t$	 нормальна, то $u(t) \equiv h(t)$, $\H\equiv \const \neq 0$, и 
\begin{align*}
&\dg = \sum_{i=1}^n h_i X_i(g), \quad g = (x, y) \in G,\\
&\dh = \H h.
\end{align*}
\end{itemize}	
\end{theorem}
Согласно условию Гоха, все локально оптимальные анормальные траектории нормальны. 

Впрочем, в примере~14~\cite{brock80} показано, что существуют неоптимальные строго анормальные траектории.

Будем далее рассматривать случай общего положения: пуcть при четном $n$ матрица $\H$ невырождена и имеет $n/2$ разных собственных значений, а при нечетном $n$ матрица $\H$ имеет одно нулевое собственное значение и разные все остальные собственные значения.

Запишем ненулевые собственные значения матрицы $\H$ в виде:
$$
\{i\lam_1, -i\lam_1, \dots, i\lam_{[n/2]}, -i\lam_{[n/2]}  \},
$$
где $\lam_k > 0$ для всех $k$. Обозначим соответствующие собственные векторы матрицы $\H$:
$$
\{v_1, v_{-1}, \dots, v_{[n/2]}, v_{-[n/2]} \}.
$$
В случае нечетного $n$ обозначим через $v_0$ вещественный  собственный вектор, соответствующий нулевому собственному значению. Обозначим через $\lan v, w \ran = v^T \bar w$ эрмитово скалярное произведение в $\C^n$, и определим ортогональные векторы
$$
\a_k = 2 \Im(\lan h_0, v_k\ran v_k),\quad  \b_k = 2 \Re(\lan h_0, v_k\ran v_k), \quad \g_0 = \lan h_0, v_0\ran v_0.
$$

\begin{theorem}
Пусть $(x(t), y(t)) \in \R^n\times\so(n)$ есть геодезическая. Если $n$ четно, то
\begin{align}
&x = \sum_{i = 1}^{[n/2]}\dfrac{1}{\lam_i}(\cos(\lam_i t) - 1)\a_i + \dfrac{1}{\lam_i}\sin(\lam_i t)\b_i, \label{n(n+1)x1}\\
&y = \sum_{i, j = 1}^{[n/2]} A_{ij}\,\a_i \land\a_j + B_{ij}\,\a_i\land\b_j +  C_{ij}\,\b_i\land\b_j,\label{n(n+1)y1}
\end{align}
где
\begin{align*}
&A_{ij}=  \dfrac{1 - \cos(\lam_i - \lam_j)t}{2\lam_i(\lam_i - \lam_j)} + \dfrac{\cos(\lam_i + \lam_j)t - 1}{2\lam_i(\lam_i + \lam_j)} -\dfrac{\cos(\lam_j t) - 1}{\lam_i\lam_j}, \quad i \neq j,\\
&B_{ij}=  \dfrac{(\lam_j - \lam_i)\sin(\lam_i + \lam_j)t}{2\lam_i\lam_j(\lam_i + \lam_j)} + \dfrac{(\lam_i + \lam_j)\sin(\lam_i - \lam_j)t}{2\lam_i\lam_j(\lam_i - \lam_j)} -\dfrac{\sin(\lam_j t)}{\lam_i\lam_j}, \quad i \neq j,\\
&B_{ii }= \dfrac{t}{\lam_i}  -\dfrac{\sin (\lam_i t)}{\lam_i^2}, \\
&C_{ij} = \dfrac{1 - \cos(\lam_i - \lam_j)t}{2\lam_i(\lam_i - \lam_j)} - \dfrac{\cos(\lam_i + \lam_j)t - 1}{2\lam_i(\lam_i + \lam_j)}.
\end{align*}

Если $n$ нечетно,  то к суммам  \eq{n(n+1)x1} и \eq{n(n+1)y1} нужно прибавить  соответственно  дополнительные слагаемые $t \g_0 $ и 
\begin{align*}
\sum_{i = 1}^{[n/2]}\left[\dfrac{-t}{\lam_i}(\cos(\lam_i t) + 1) + \dfrac{2\sin (\lam_i t)}{\lam_i^2}\right]\a_i \land \g_0 
+\sum_{i = 1}^{[n/2]}\left[\dfrac{1}{\lam_i^2}(2(1 - \cos(\lam_i t)) - \lam_i t \sin(\lam_i t))\right]\b_i \land \g_0.
\end{align*}
\end{theorem}

\subsubsection{Нижняя оценка множества разреза}\label{subsubsec:twostep_cut}
Рассмотрим следующее подмножество группы $G = \R^n \times \so(n)$:
\begin{align}\label{n(n+1)st}
&C_n = \{(x, y) \in  \R^n \times \so(n) \mid y \neq 0, \quad \exists \,\, 0 \neq M\in \O(n): M x = x, \quad M y M^T = y, \quad M|_{\ker y} = \Id  \}.
\end{align}
Условие \eq{n(n+1)st} означает, что элемент $M \in\O(n)$ стабилизирует $(x, y)$ относительно действия \eq{n(n+1)O(n)}.
\begin{proposition}
Для любого $n \geq 2$ множество $C_n \subset G$ есть полуалгебраическое множество коразмерности $2$.
\end{proposition}
\begin{theorem}
Для всех $n \geq 2$ имеет место включение:
\begin{align}
&C_n \subset \Cut. \label{n(n+1)Cn}
\end{align}
\end{theorem}
\begin{remark}
При $n = 2, 3$ включение \eq{n(n+1)Cn} превращается в равенство, см.~разделы~\ref{subsec:heis}, \ref{subsec:36}. Верно ли 
это при $n \geq 4$, неизвестно.
\end{remark}

\subsubsection{Анормальное множество}\label{subsubsec:twostep_abnorm}
Приведем известные описания и свойства анормального множества $\Abn$.
\begin{theorem}\label{th:n(n+1)abn1}
Имеет место равенство
\begin{align}
&\Abn = \bigcup_{W \in \Gr(n, n-2)}W\times(W \land W),
\end{align}
где $ \Gr(n, n-2)$ есть грассманиан $(n-2)$-мерных подпространств в $\R^n$.
\end{theorem}

Рангом элемента $y \in\so(n)$ назовем размерность образа оператора $y: \R^n \to \R^n$. Для открытого плотного подмножества в $\so(n)$ ранг принимает максимальное значение: $n$ при четном $n$, и $n-1$ при нечетном~$n$.

Обозначим через $\so(n)_{\sing}$ множество элементов $y \in \so(n)$, имеющих ранг меньше максимального,  через $\so(n)_{d}$ множество элементов ранга $d$,  и через $\so(n)_{< d}$ множество элементов ранга меньше $d$.
\begin{theorem}\label{th:n(n+1)abn2}
Если $n$ нечетно, то $\Abn = \R^n\times \so(n)_{\sing}$.
\end{theorem}
\begin{theorem}\label{th:n(n+1)abn3}
Если $n$ четно, то $\Abn$ есть объединение $Y_1 \cup Y_2$ двух квазипроективных подмногообразий 
\begin{align*}
&Y_1 = \{(x, y)\in \R^n\times \so(n) \mid x \in\Im y, \quad y \in \so(n)_{n-2}  \},\\
&Y_2 = \R^n\times \so(n)_{< n-2}.
\end{align*}
В частности, $\Abn$ есть особое алгебраическое многообразие коразмерности $3$.
\end{theorem}
\begin{theorem}\label{th:n(n+1)abn4}
$\Abn \subset G$ есть полуалгебраическое множество коразмерности не меньше $3$.
\end{theorem}
\begin{theorem}\label{th:n(n+1)abn5}
Для всех $k \geq 2$ имеют место включения
$$
\bar C_n\setminus C_n \subset \Abn  \subsetneq \bar C_n.
$$
Для любого $n \geq 4$ первое включение строгое, поэтому $\Cut \cap \Abn \neq \emptyset$.
Более того, существуют анормальные геодезические с конечным временем разреза.
\end{theorem}

\subsubsection{Библиографические комментарии}
Рассматриваемая в данном разделе субриманова задача исследовалась в известных работах Б.~Гаво~\cite{gaveau}, Р.~Брокетта~\cite{brock80} и В.~Лиу и Х.~Суссмана~\cite{liusus}. Изложение в этом разделе опирается на следующие источники: пункты~\ref{subsubsec:twostep_state}, \ref{subsubsec:twostep_extrem} --- статья~\cite{monroy_meneses}, пункт~\ref{subsubsec:twostep_cut} --- статья~\cite{rizzi_serres},
 пункт~\ref{subsubsec:twostep_abnorm}: теоремы~\ref{th:n(n+1)abn1}--\ref{th:n(n+1)abn4} --- статья~\cite{ledonne_sard}, теорема~\ref{th:n(n+1)abn5} --- статья~\cite{rizzi_serres}.

\subsection{Двухступенная субриманова задача коранга 1}\label{subsec:corank1}
\subsubsection{Алгебра Ли и группа Ли}
Рассмотрим алгебру Ли 
$$
\gg = \spann(X_1, \dots, X_k, Y_1, \dots, Y_k, Z), \quad k \in \N,
$$
с таблицей умножения 
\begin{align*}
&[X_i, Y_j] = -\d_{ij}b_i Z, \qquad b_i > 0, \quad i, j = 1, \dots, k,\\
&[X_i, Z] = [Y_i, Z] = 0, \qquad  i = 1, \dots, k.
\end{align*}

Пусть $G$ есть связная односвязная группа Ли с алгеброй Ли $\gg$. Тогда $G \cong \R^{2k + 1}$ и можно выбрать координаты
$$
g = (x, y, z) \in \R^{2k + 1}, \quad x = (x_1, \dots, x_k) \in\R^k, \quad y = (y_1, \dots, y_k)\in \R^k, \quad z \in \R,
$$
в которых умножение в $G$ принимает вид
$$
g\cdot g' = \left(x + x', y + y', z + z' - \frac 12 \sum_{i = 1}^k b_i(x_i\cdot x_i' - y_i\cdot y_i')\right).
$$
В этих координатах левоинвариантный репер на $G$ есть
\begin{align*}
&X_i = \dfrac{\partial}{\partial x_i} + \frac 12 b_i y_i\dfrac{\partial}{\partial z}, \qquad i = 1, \dots, k,\\
&Y_i = \dfrac{\partial}{\partial y_i} - \frac 12 b_i x_i\dfrac{\partial}{\partial z}, \qquad i = 1, \dots, k,\\
&Z = \dfrac{\partial}{\partial z}.
\end{align*}

\subsubsection{Постановка задачи}
Рассмотрим субриманову структуру на $G$ с ортонормированным репером $(X_1, \dots, X_k, Y_1, \dots, Y_k)$. Соответствующая задача оптимального управления имеет вид
\begin{align*}
&\dx_i = u_i, \qquad i =  1, \dots k,\\
&\dy_i = v_i, \qquad i =  1, \dots k,\\
&\dz = \frac 12 \sum_{i=1}^k b_i (u_i y_i - v_i x_i),\\
&g(0) = (0, 0, 0), \qquad g(t_1) = (x^1, y^1 , z^1),\\
&\int_0^{t_1} \left( \sum_{i=1}^k u_i^2 + v_i^2   \right)^{1/2}\, dt \to \min.
\end{align*}
Существование решений следует из теорем Рашевского-Чжоу и Филиппова.

\subsubsection{Экстремали}
Анормальные экстремальные траектории постоянны.

Введем линейные на слоях $T^*G$ гамильтонианы
$$
h_{X_i}(\lam) = \lan\lam, X_i(g)  \ran, \qquad h_{Y_i}(\lam) = \lan\lam, Y_i(g)\ran, \qquad h_{Z}(\lam) = \lan\lam, Z(g)  \ran.
$$
Вдоль нормальных экстремалей имеем
$$
u_i(t) = h_{X_i}(\lam(t)), \quad v_i(t) = h_{Y_i}(\lam(t)), \quad w(t) = h_Z(\lam(t)).
$$

Переходя к натуральной параметризации геодезических, можно считать, что
$$
u_1^2(t) + \dots + u_k^2(t) + v_1^2(t)  + \dots + v_k^2(t) \equiv 1.
$$

Если $w = 0$, то
\begin{align*}
&(u_i(t), v_i(t)) \equiv (u_i^0, v_i^0) = \const, \\ 
&x_i(t) = u_i^0 t, \quad y_i(t) = v_i^0 t, \quad z(t) \equiv 0.
\end{align*}

Если $w \neq 0$, то (обозначая $a_i = b_i w$)
\begin{align*}
&u_i(t) = u_i^0 \cos a_i t - v_i^0 \sin a_i t,\\
&v_i(t) = u_i^0 \sin a_i t + v_i^0\cos a_i t,\\
&w(t) = w
\end{align*}
и
\begin{align}
&x_i(t) = \frac{1}{a_i}(u_i^0 \sin a_i t + v_i^0\cos a_i t -v_i^0), \nonumber\\
&y_i(t) =  \frac{1}{a_i}( - u_i^0 \cos a_i t + v_i^0 \sin a_i t +u_i^0),\label{carank1} \\
&z(t) = \frac{1}{2w^2}\left(w t - \sum_{i}\frac{1}{b_i}((u_i^0)^2 + (v_i^0)^2)\sin a_i t\right).\nonumber
\end{align}
В полярных координатах 
\begin{align*}
&u_i = r_i\cos \t_i, \qquad v_i = r_i \sin \t_i, \qquad i = 1, \dots, k,\\
&r_1^2 + \dots + r_k^2 \equiv 1,
\end{align*}
формулы \eq{carank1} переписываются в виде 
\begin{align*}
&x_i(t) = \frac{r_i}{a_i}(\cos(a_i t + \t_i) - \cos \t_i),\\
&y_i(t) = \frac{r_i}{a_i}(\sin(a_i t + \t_i) - \sin\t_i),\\
&z(t) = \frac{1}{2w^2}\left(w t - \sum_{i}\frac{r_i^2}{b_i}\sin a_i t\right).
\end{align*}
Поэтому проекция геодезической на любую плоскость $(x_i, y_i)$ есть окружность периода $T_i$, радиуса $\r_i$ с центром $C_i$, где
$$
T_i = \frac{2\pi}{b_i w}, \quad \r_i = \frac{r_i}{b_i w}, \quad C_i = -\frac{r_i}{b_i w}(\cos \t_i, \, \sin \t_i), \quad i = 1, \dots, k.
$$
Компонента $z(t)$ геодезических есть взвешенная сумма (с коэффициентами $b_i$) площадей, заметенных радиус-векторами $(x_i(t), y_i(t))$ на плоскостях $\R^2_{x_iy_i}$.

\subsubsection{Время разреза}
\begin{theorem}
Пусть $g(t)$ есть натурально параметризованная геодезическая, выходящая из начала координат. Тогда время разреза  вдоль нее совпадает с первым сопряженным временем и равно
\begin{align*}
&t_{\cut} =  \frac{2 \pi}{w\max_i b_i} \quad \text{при} \quad  w \neq 0,\\
&t_{\cut} = +\infty\quad \text{при} \quad w = 0.
\end{align*}
\end{theorem}

\subsubsection{Библиографические комментарии}
Результаты этого раздела получены в работе \cite{corank1}.

В более ранней работе \cite{mon_anz_heis} для этой же субримановой задачи получена параметризация геодезических и описано первое сопряженное время.

\subsection{Двухступенная субриманова задача коранга 2}\label{subsec:corank2}
\subsubsection{Алгебра Ли и группа Ли}
Рассмотрим алгебру Ли 
\begin{align*}
&\gg = \gg^{(1)} \oplus \gg^{(2)},\\
&[\gg^{(1)}, \gg^{(1)}] = \gg^{(2)}, \quad [\gg^{(1)}, \gg^{(2)}] =[\gg^{(2)}, \gg^{(2)}] = \{ 0 \},\\
&\dim \gg^{(1)} = k \geq 2, \quad \dim \gg^{(2)} =  2.
\end{align*}
Существует базис
\begin{align*}
&\gg = \spann(X_1, \dots, X_k, Y_1, Y_2),\\
&\gg^{(1)} = \spann(X_1, \dots, X_k), \quad \gg^{(2)} = \spann(Y_1,  Y_2),
\end{align*}
в котором
\begin{align*}
&[X_i, X_j] = \sum_{h=1}^2b_{ij}^h \, Y_h, \quad i, j = 1, \dots, k,\\
&[X_i, Y_j] = [Y_1, Y_2]=0, \quad i = 1, \dots, k, \quad j = 1, 2,\\
&L_h = (b_{ij}^h) \in \so(k), \quad h = 1, 2.
\end{align*}
Пусть $G$ --- связная односвязная группа Ли с алгеброй Ли  $\gg$. Тогда на группе $G\cong \R^{k+2}$ существуют координаты $(x_1, \dots, x_k, y_1, y_2)$, в которых 
\begin{align*}
&X_i = \frac{\partial}{\partial x_i } + \frac 12 \sum_{j, h}b_{ij}^h \,x_j\frac{\partial}{\partial y_h }, \quad i = 1, \dots, k,\\
&Y_h = \frac{\partial}{\partial y_h }, \quad h = 1, 2.
\end{align*}

\subsubsection{Постановка задачи}
Рассмотрим субриманову задачу на $G$ с ортонормированным репером $(X_1, \dots, X_k)$. В координатах соответствующая задача оптимального управления имеет вид:
\begin{align*}
&\dx_i = u_i, \quad i = 1, \dots, k,\\
&\dy_h = \frac 12 x^TL_h u, \quad h = 1, 2,\\
&g = (x_1, \dots, x_k, y_1, y_2) \in G\cong \R^{k + 2}, \quad (u_1, \dots, u_k) \in \R^k,\\
&g(0) = \Id = (0, \dots, 0), \quad g(t_1) = g^1,\\
&\int_0^{t_1} \left(\sum_{i = 1}^k u_i^2 \right)^{1/2}\, dt \to \min.
\end{align*}

Существование решений следует из теорем Рашевского-Чжоу и Филиппова.

\subsubsection{Экстремальные управления и траектории}
В силу условия Гоха, локально оптимальные анормальные траектории нормальны.

Экстремальные управления, соответствующие натурально параметризованным нормальным траекториям, имеют вид
$$
u(t) = e^{t(r_1L_1 + r_2L_2)}\,u_0, \qquad \|u_0 \| = 1, \quad r_1, r_2 \in \R,\\
$$
а сами эти траектории суть
\begin{align*}
&x(t) = \int_0^te^{s(r_1L_1 + r_2L_2)}u_0\, ds,\\
&y_i(t) = \frac 12 \int_0^t(x(s))^TL_i u(s)ds, \quad i = 1, 2.
\end{align*}

\subsubsection{Время разреза}
\begin{theorem}
Натурально параметризованные геодезические, соответствующие начальному ковектору  $\lam_0 = (u_0, r) \in S^{k-1} \times \R^2$, имеют время разреза 
\begin{align*}
& t_{\cut}(\lam_0) =\frac{2\pi}{\max \s(r_1L_1 + r_2L_2)}, \quad r\neq 0,\\
&t_{\cut}(\lam_0) = + \infty, \quad r = 0,
\end{align*}
где $\max \s(A)$ обозначает максимальный модуль собственного значения матрицы $A$.
Вообще говоря, время разреза отлично от первого сопряженного времени.
\end{theorem}
\begin{theorem}
В случае $k = 4$ время разреза совпадает с первым сопряженным временем.
\end{theorem}
\subsubsection{Библиографические комментарии}
Результаты этого раздела получены в работе \cite{corank2}.

 \subsection{Субримановы $\ddd\oplus\sss$ задачи}  \label{subsec:s+d}

Левоинвариантная субриманова структура $(\D, \lan\cdot,\cdot\ran)$ на группе Ли $G$  с алгеброй Ли $\gg$ называется \ddef{$\ddd\oplus\sss$  структурой},  если выполнены следующие условия:
\begin{itemize}
\item[(1)]
на группе Ли $G$ имеется биинвариантное (левоинвариантное и правоинвариантное) скалярное произведение $\widetilde{\lan\cdot,\cdot\ran}$,
\item[(2)]
$\gg = \ddd \oplus \sss$,
\item[(3)]
$\D = \ddd$, $\lan\cdot,\cdot\ran = \widetilde{\lan\cdot,\cdot\ran}|_{\ddd}$,
\item[(4)]
$\sss = \ddd^{\perp}$,  где ортогональность понимается в смысле $\widetilde{\lan\cdot,\cdot\ran}$,
\item[(5)]
$[\sss, \sss] \subset \sss$.
\end{itemize}

\begin{theorem}\label{th:d+s}
Геодезические субримановой $\ddd\oplus\sss$  структуры на группе Ли $G$, начинающиеся в точке $\Id$, суть произведения двух однопараметрических подгрупп:
$$
g(t) = e^{t(x_0 + y_0)} e^{-ty_0}, \qquad x_0 \in \ddd, \quad y_0 \in \sss.
$$
\end{theorem}

Частные случаи $\ddd\oplus\sss$  структур  рассматриваются далее:
\begin{itemize}
\item[(1)]
на группе $\SU(2)$, см. раздел \ref{subsec:su2},
\item[(2)]
на группе $\SO(3)$, см. раздел \ref{subsec:sr_so3},
\item[(3)]
на группе $\SL(2)$, см. раздел \ref{subsec:sr_sl2}.
\end{itemize}

\subsubsection{Библиографические комментарии}
 Теорема \ref{th:d+s} была получена независимо А.А. Аграчевым \cite{agr95} и Р. Брокеттом \cite{brock99} (эта работа базируется на предыдущей работе \cite{brock73}).  Затем этот результат интенсивно исследовался в \cite{jurd99,jurd01,jurd16,BCG02a,boscain_rossi},   а также в \cite[С. 200]{montgomery_book}. Его красивая геометрическая интерпретация получена в \cite{BZ15b, BZ15a,BZ16}.

\subsection{Осесимметричная субриманова задача на группе $\SU(2)$}\label{subsec:su2}
\subsubsection{Группа Ли и алгебра Ли}\label{subsubsec:su2_gr}
Группа Ли $\SU(2)$ есть группа унитарных унимодулярных  $2\times 2$ комплексных матриц 
$$
\SU(2) = \left\{
\begin{pmatrix}
\a       & \b     \\
-\bar \b & \bar \a
\end{pmatrix} \in \GL(2, \C) \mid |\a|^2 +|\b|^2 = 1
      \right\}.
$$
Эта группа компактна, связна и односвязна. Группа Ли $\SU(2)$ диффеоморфна трехмерной сфере
$$
S^3 = \left\{
\begin{pmatrix}
\a \\
 \b
\end{pmatrix} \in \C^2 \mid |\a|^2 +|\b|^2 = 1
  \right\}
$$
в силу диффеоморфизма 
$$
\F: \SU(2) \to S^3, \qquad \begin{pmatrix}
                            \a       & \b     \\
                           -\bar \b & \bar \a
                           \end{pmatrix} 
\mapsto        
\begin{pmatrix}
\a \\
 \b
\end{pmatrix}.
$$
Поэтому будем далее записывать элементы группы $\SU(2)$ как пары комплексных чисел $(\a, \b)$.

Алгебра Ли группы $\SU(2)$ есть алгебра косоэрмитовых бесследовых $2\times 2$ комплексных матриц
$$
\su(2) = \left\{
\begin{pmatrix}
i\a       & \b     \\
-\bar \b & -i \a
\end{pmatrix} \in \gl(2, \C) \mid \a \in \R, \quad \b \in \C
      \right\}.
$$
В этой алгебре можно выбрать базис
$$
p_1 = \frac 12 \begin{pmatrix}
                 0 & 1\\
								-1 & 0
                \end{pmatrix}, \quad
p_2 = \frac 12 \begin{pmatrix}
                 0 & i\\
								 i & 0
                \end{pmatrix},\quad
k = \frac 12 \begin{pmatrix}
                 i & 0\\
								 0 & -i
                \end{pmatrix} 																								
$$
с таблицей умножения
$$
[p_1, p_2] = k, \quad [p_2, k] = p_1, \quad [k, p_1] = p_2.
$$
Форма Киллинга  для $\su(2)$ есть
$\Kil(X, Y) = 4\Tr(XY)$, поэтому $\Kil(p_i, p_j) = -2\d_{ij}$.

Подпространства ${\mathbf \k} = \R k$,   ${\mathbf \ppp} = \spann(p_1, p_2)$ образуют картановское разложение для $\su(2)$. Более того, $(p_1, p_2)$ есть ортонормированный репер для скалярного произведения 
$\lan \cdot, \cdot\ran = - \frac 12 \Kil(\cdot, \cdot)$, суженного на подпространство ${\mathbf \ppp}$.

\subsubsection{Субриманова ${\k}\oplus {\ppp}$ структура}
Рассмотрим левоинвариантную  субриманову структуру на $\SU(2)$ с ортонормированным репером 
$$
X_1(g) = L_{g*}p_1, \qquad X_2(g) = L_{g*}p_2,\qquad g \in \SU(2),
$$
то есть распределение  $\D_g =  L_{g*}{\mathbf \ppp}$ со скалярным произведением $\lan v_1, v_2 \ran_g = \lan L_{g^{-1}*}v_1, L_{g^{-1}*}v_2 \ran$.

Такая субриманова  структура называется {${\k}\oplus {\ppp}$ структурой}.
Эти структуры определяются следующим образом. Пусть $G$ есть простая группа Ли с алгеброй Ли $\gg$.
Пусть $\gg = \k \oplus \ppp$  есть картановское разложение алгебры $\gg$:
$$
[\k, \k] \subset \k, \quad 
[\ppp, \ppp] \subset \k, \quad
[\k, \ppp] \subset \ppp. 
$$
Рассмотрим на $G$  распределение $\D_g = L_{g*} \ppp$  с метрикой 
$$
\lan v_1, v_2\ran_g = \lan L_{g^{-1}*}v_1, L_{g^{-1}*}v_2\ran, 
\qquad g \in G,
$$
  где $\lan\cdot, \cdot\ran = \a \Kil|_{\ppp} (\cdot, \cdot)$,  и $\a < 0$ (соотв. $\a > 0$)  если $G$  компактна (соотв. некомпактна). Тогда $(\D, \lan\cdot, \cdot\ran)$  называется \ddef{субримановой $\k \oplus \ppp$  структурой}  на $G$.

Все $\k \oplus\ppp$  субримановы структуры на $\SU(2)$  эквивалентны между собой.

\subsubsection{Геодезические и симметрии}\label{subsubsec:su2_geod}
Задача трехмерная контактная, поэтому анормальные траектории постоянны.

Введем линейные на слоях $T^*G$  гамильтонианы $h_i(\lam) = \lan\lam, X_i(g)\ran$, $i = 1, 2, 3$, $X_3 = [X_1, X_2]$, и максимизированный нормальный гамильтониан ПМП $H = \frac 12 (h_1^2+h_2^2)$.  Натурально параметризованные нормальные экстремали задаются точками цилиндра $C = \gg^* \cap \{H = \frac 12\}$:
\begin{align*}
&\lam = (h_1, h_2, h_3, \Id) \in C, \\
&h_1 = \cos \t, \quad h_2 = \sin \t, \quad h_3 = c.
\end{align*}
Тогда экспоненциальное отображение
$$
\map{\Exp}{C \times \R_+}{G}, \quad (\lam, t) \mapsto g(t),
$$
имеет следующую параметризацию:
\begin{align*}
&\Exp(\lam, t) = \Exp(\t, c, t) = e^{(\cos \t p_1 + \sin \t p_2 + c k)} e^{-ckt} = 
\left(\begin{array}{c} \a \\ \b \end{array}\right),
\intertext{где}
&\a = \frac{c \sin(\frac{ct}{2}) \sin(\sqrt{1+c^2} \frac t2)}{\sqrt{1+c^2}} + 
\cos\left(\frac{ct}{2}\right) \cos\left(\sqrt{1+c^2} \frac t2\right)+ \\
&\qquad + i \left(
\frac{c \cos(\frac{ct}{2}) \sin(\sqrt{1+c^2} \frac t2)}{\sqrt{1+c^2}} - 
\sin\left(\frac{ct}{2}\right) \cos\left(\sqrt{1+c^2} \frac t2\right)
\right),\\
&\b = 
 \frac{\sin(\sqrt{1+c^2} \frac t2)}{\sqrt{1+c^2}}
\left(\cos\left(\frac{ct}{2} + \t\right) + i \sin \left(\frac{ct}{2} + \t\right) \right).
\end{align*}
Экспоненциальное отображение имеет следующие симметрии:
\begin{itemize}
\item
вращения
$$
\Exp(\t, c, t) = 
\left(\begin{array}{cc}
1 & 0 \\
0 & e^{i \t}
\end{array}\right)  \Exp(0, c, t),
$$
\item
центральную симметрию: если $\Exp(\t, c, t) = \left(\begin{array}{c} \a \\ \b\end{array}\right)$, то
$$
\Exp(\t, -c, t) = 
\begin{cases}
 \left(\begin{array}{c} 
\overline{\a} \\ e^{2 i(\t - \arg \b)}\b\end{array}\right)
&\text{  если } \b \neq 0, \\
 \left(\begin{array}{c} 
\overline{\a} \\ 0 \end{array}\right)
&\text{  если } \b = 0.
\end{cases}
$$
\end{itemize}

\subsubsection{Сопряженные точки}\label{subsubsec:su2_conj}
Момент времени $t > 0$  является сопряженным временем вдоль геодезической $g(t) = \Exp(\t, c, t)$  тогда и только тогда, когда
$$
\sin\left(\sqrt{1+c^2} \frac t2\right) \left( 2 \sin \left( \sqrt{1+c^2} \frac t2\right) - \sqrt{1+c^2} t  \cos \left(\sqrt{1+c^2} \frac t2\right)\right) = 0.
$$
Поэтому $n$-е сопряженное время $t_n^{\conj}(\lam)$  вдоль $g(t)$  имеет вид:
\begin{align*}
&t_{2m-1}^{\conj}(\lam) = \frac{2 \pi m}{\sqrt{1+c^2}}, \\
&t_{2m}^{\conj}(\lam) = \frac{2 x_m}{\sqrt{1+c^2}},
\end{align*}
где $\{x_1, x_2, \dots\}$  суть упорядоченные по возрастанию положительные корни уравнения $x = \tg x$.

Соответствующая $n$-я каустика
$$
\Conj^n = \{\Exp(\lam, t) \mid t = t_n^{\conj}(\lam), \ \lam \in C \}
$$
есть
\begin{align*}
&\Conj^{2m-1} = e^{\k} \setminus \{\Id\}, \\
&\Conj^{2m} = 
\left\{
\left(\begin{array}{c}
\frac{c \sin x_m}{\sqrt{1+c^2}} e^{i(\frac{\pi}{2} - y_n)} + \cos x_n e^{-i y_n} \\
\frac{\sin x_n}{\sqrt{1+c^2}} e^{i \t}
\end{array}\right) 
\mid c \in \R, \ \t \in \R/(2 \pi \Z)
\right\},
\end{align*}
где $ y_n = \frac{cx_n}{\sqrt{1+c^2}}$.

\subsubsection{Множество разреза}\label{subsubsec:su2_cut}
\begin{theorem}
Множество разреза есть
$$
\Cut = \Conj^1 = e^{\k} \setminus \{\Id\} = \{ e^{c k} \mid c \in (0, 4 \pi)\}.
$$
\end{theorem}
Топологически, $\Cut$  есть интервал (большая окружность $e^{\k}$  с выколотой точкой $\Id$). Так как множество разреза совпадает с первой каустикой, то геодезические одновременно теряют локальную и глобальную оптимальность (так же, как га группе Гейзенберга, см. раздел \ref{subsec:heis}).

\subsubsection{Субриманово расстояние}\label{subsubsec:su2_dist}

\begin{theorem}
Пусть $g = (\a, \b) \in \SU(2)$, тогда субриманово расстояние $d_0(g) = d(g, \Id)$  имеет следующее представление.
\begin{itemize}
\item[$(1)$]
Если $\a = 0$,  то $d_0 = \pi$.
\item[$(2)$]
Если $|\a| = 1$,  то $d_0 =  2 \sqrt{|\arg \a|(2 \pi - |\arg \a|)}$, где $\arg \a \in [-\pi, \pi]$.
\item[$(3)$]
Если $0 < |\a| < 1$ и $\Re \a = |\a| \sin(\frac{\pi}{2} |\a|)$,  то $d_0 = \pi \sqrt{1 - |\a|^2}$.
\item[$(4)$]
Если $0 < |\a| < 1$ и $\Re \a > |\a| \sin(\frac{\pi}{2} |\a|)$,  то 
$$
d_0 = \frac{2}{ \sqrt{1 + \b^2}} \arcsin \sqrt{(1-|\a|^2)(1+\b^2)} \in \left(0, \frac{\pi}{\sqrt{1+\b^2}}\right),
$$
где $\b$  --- единственное решение системы уравнений
$$
\begin{cases}
\cos \left( - \frac{\b}{\sqrt{1+\b^2}} \arcsin  \sqrt{(1-|\a|^2)(1+\b^2)} 
+ \arcsin \frac{\b \sqrt{1-|\a|^2}}{|\a|}\right) = \frac{\Re \a}{|\a|},\\
\sin \left( - \frac{\b}{\sqrt{1+\b^2}} \arcsin  \sqrt{(1-|\a|^2)(1+\b^2)} 
+ \arcsin \frac{\b \sqrt{1-|\a|^2}}{|\a|}\right) = \frac{\Im \a}{|\a|}.
\end{cases}
$$
\item[$(5)$]
Если $0 < |\a| < 1$ и $\Re \a < |\a| \sin(\frac{\pi}{2} |\a|)$,  то 
$$
d_0 = \frac{2}{ \sqrt{1 + \b^2}} \left( \pi - \arcsin \sqrt{(1-|\a|^2)(1+\b^2)}\right) \in \left(\frac{\pi}{\sqrt{1+\b^2}}, \frac{2 \pi}{\sqrt{1+\b^2}}\right),
$$
где $\b$  --- единственное решение системы уравнений
$$
\begin{cases}
\cos \left(  \frac{\b}{\sqrt{1+\b^2}} \left( \pi - \arcsin \sqrt{(1-|\a|^2)(1+\b^2)}\right) 
+ \arcsin \frac{\b \sqrt{1-|\a|^2}}{|\a|}\right) = - \frac{\Re \a}{|\a|},\\
\sin \left(  \frac{\b}{\sqrt{1+\b^2}} \left( \pi - \arcsin \sqrt{(1-|\a|^2)(1+\b^2)}\right)
+ \arcsin \frac{\b \sqrt{1-|\a|^2}}{|\a|}\right) = \frac{\Im \a}{|\a|}.
\end{cases}
$$
\end{itemize}
\end{theorem}

\subsubsection{Геодезические со специальными граничными условиями}\label{subsubsec:su2_spec}
\begin{theorem}
Если точка $g\in \SU(2)$  принадлежит интервалу
$$
A = \{(\cos \f + i \sin \f, 0) \in \SU(2) \mid \f \in (0, 2 \pi)\},
$$
 то существует счетное число геометрически разных геодезических $\g_n$, соединяющих $\Id$  и $g$:
\begin{align*}
&\g_n = \{\a_n(t), \b_n(t)) \in \SU(2) \mid t \in [0, t_n]\},\\
&\a_n(t) = \left( \cos\left( t \frac{\pi n}{\sqrt{\pi^2n^2 - \f^2}}\right) - i \frac{\f}{\pi n} 
\sin \left( t \frac{\pi n}{\sqrt{\pi^2n^2 - \f^2}}\right)\right) e^{\frac{it\f}{\sqrt{\pi^2n^2 - \f^2}}}, \\
&\b_n(t) = \dot\b_n(0) \frac{\sqrt{\pi^2n^2 - \f^2}}{\pi n}
\sin \left( t \frac{\pi n}{\sqrt{\pi^2n^2 - \f^2}}\right) e^{\frac{it\f}{\sqrt{\pi^2n^2 - \f^2}}}, 
\end{align*}
где $n \in \Z \setminus \{0, \ \pm 1\}$,  а $l_n = \frac{1}{\sqrt 2} t_n = \frac{1}{\sqrt 2} \sqrt{\pi^2n^2 - \f^2}$  есть длина геодезической $\g_n$.
\end{theorem}

\begin{theorem}
Если точка $g\in \SU(2)$  не  принадлежит ни интервалу
$
A
$, ни сфере $S^2 = \{(\a, \b) \in \SU(2)  \mid \Im \a = 0 \}$,
 то существует конечное число геометрически разных геодезических, соединяющих $\Id$  и $g$.
\end{theorem}

\subsubsection{Библиографические комментарии}
Изложение разделов \ref{subsubsec:su2_gr}--\ref{subsubsec:su2_cut}  опирается на работу \cite{boscain_rossi}, раздела \ref{subsubsec:su2_dist} ---  на работу \cite{BZ15b}, раздела \ref{subsubsec:su2_spec} ---  на работу \cite{chang_markina_vasiliev}.

\subsection{Осесимметричная субриманова  задача   на группе $\SO(3)$}  \label{subsec:sr_so3}

\subsubsection{Группа Ли и алгебра Ли}\label{subsubsec:so3_gr}
 Группа Ли \ddef{$\SO(3)$}   есть группа унимодулярных ортогональных $3 \times 3$  вещественных матриц
$$
\SO(3) = \{ g \in \GL(3,\R) \mid g g^T = \Id, \ \det g = 1\}.
$$
Эта группа компактна, связна и неодносвязна: ее фундаментальная группа есть $\Z_2$.
Алгебра Ли этой группы есть алгебра кососимметрических $3 \times 3 $  вещественных матриц:
$$
\so(3) = \{X \in \gl(3,\R) \mid X^T = - X\}.
$$
В этой алгебре Ли можно выбрать базис $\gg = \so(3) = \spann(p_1, p_2, k)$,
\be{so3_tab}
p_1 = 
\left(\begin{array}{ccc}
0 & 0 & 0 \\
0 & 0 & - 1 \\
0 & 1 & 0
\end{array}\right), 
\quad 
p_2 = 
\left(\begin{array}{ccc}
0 & 0 & 1 \\
0 & 0 & 0 \\
-1 & 0 & 0
\end{array}\right), 
\quad 
k = 
\left(\begin{array}{ccc}
0 & -1 & 0 \\
1 & 0 & 0 \\
0 & 0 & 0
\end{array}\right), 
\ee
с таблицей умножения
$
[p_1, p_2] = k$, $[p_2, k] = p_1$, $[k, p_1] = p_2$.
Алгебры Ли $\so(3)$  и $\su(2)$ изоморфны, а группа Ли $\SU(2)$  есть 
односвязная накрывающая группы $\SO(3)$, см. раздел \ref{subsec:su2}.
Двулистное накрытие $\map{\Pi}{\SU(2)}{\SO(3)}$  можно задать следующим образом:
\be{so3Pi}
 \Pi \ : \ 
\left(\begin{array}{c}
a + i b \\ c + i d
\end{array}\right)
\mapsto 
\left(\begin{array}{ccc}
1 - 2b^2 - 2 d^2 & 2 cd - 2 ab & 2bc + 2 ad \\
2cd + 2ab & 1-2b^2-2c^2 & 2bd+2ac \\
2bc - 2ad & 2bd-2ac & 1-2c^2-2d^2
\end{array}\right).
\ee

Форма Киллинга в алгебре Ли $\so(3)$  есть $\Kil(X,Y) = \Tr(XY)$,  поэтому $\Kil(p_i, p_j) = - 2 \d_{ij}$.
Подпространства $\k = \R k$, $\ppp = \spann(p_1, p_2)$  образуют картановское разложение алгебры Ли $\so(3)$.  Векторы $p_1, p_2$  образуют ортонормированный репер для скалярного произведения $\langle\cdot, \cdot\rangle = - \frac 12 \Kil(\cdot, \cdot)$,  суженного на $\ppp$.

\subsubsection{Субриманова $\k \oplus \ppp$  структура}\label{subsubsec:so3_struct}
Рассмотрим левоинвариантную субриманову структуру на $\SO(3)$  с ортонормированным репером
$$
X_1(g) = L_{g*} p_1, \quad X_2(g) = L_{g*}p_2, \qquad g \in \SO(3),
$$
то есть распределение $\D_g = L_{g*} \ppp$  со скалярным произведением 
$\langle v_1, v_2\rangle_g = \langle L_{g^{-1}*}v_1, L_{g^{-1}*}v_2\rangle$, она является $\k \oplus \ppp$  структурой на $\SO(3)$.

\subsubsection{Геодезические и симметрии}\label{subsubsec:so3_geod}
Анормальные траектории постоянны.

Пусть $h_i(\lam) = \langle \lam, X_i(g)\rangle$, $i = 1, 2, 3$, $X_3 = [X_1, X_2]$, $H = \frac 12 (h_1^2 + h_2^2)$, и $C = \gg^* \cap \{H=\frac12\}$. Далее, пусть
\begin{align*}
&\lam = (h_1, h_2, h_3,\Id) \in C, \\
&h_1 = \cos \t, \quad h_2 = \sin \t, \quad h_3 = c.
\end{align*}
Тогда экспоненциальное отображение 
$$
\map{\Exp}{C \times \R_+}{G}, \qquad (\lam,t) \mapsto g(t),
$$
параметризуется следующим образом:
\begin{align*}
&\Exp(\lam,t) = \Exp(\t, c, t) = e^{(\cos \t p_1 + \sin \t p_2 + ck)t} e^{-ckt} = \\
&= 
\left(\begin{array}{ccc}
K_1 \cos ct + K_2 \cos(2\t + ct) + K_3 c \sin ct & K_1 \sin ct + K_2 \sin(2\t + ct) - K_3 c \cos ct & K_4 \cos \t + K_3 \sin \t \\
-K_1 \sin ct + K_2 \sin(2\t + ct) + K_3 c \cos ct & K_1 \cos ct - K_2 \cos(2\t + ct)+K_3 c \sin ct & - K_3 \cos \t + K_4 \sin \t \\
K_4 \cos(\t + ct) - K_3 \sin(\t + ct) & K_3 \cos(\t + ct) + K_4 \sin(\t + ct) & \frac{\cos(\sqrt{1+c^2} t) + c^2}{1+c^2}
\end{array}\right),
\end{align*}
где 
$K_1 = \frac{1+(1+2c^2)\cos(\sqrt{1+c^2}t)}{2(1+c^2)}$, $K_2 = \frac{1-\cos(\sqrt{1+c^2}t)}{2(1+c^2)}$, $K_3 = \frac{\sin(\sqrt{1+c^2}t)}{\sqrt{1+c^2}}$, $K_4 = \frac{c(1-\cos(\sqrt{1+c^2}t))}{1+c^2}$.

Семейство геодезических имеет симметрии, аналогичные случаю $\SU(2)$, см. п. \ref{subsubsec:su2_geod}.
\subsubsection{Каустика}
Каустика (множество сопряженных точек) в $\SO(3)$  получается из случая $\SU(2)$ (п. \ref{subsubsec:su2_conj})  с помощью канонической проекции $\map{\Pi}{\SU(2)}{\SO(3)}$, см. \eq{so3Pi}.  Так же как в случае $\SU(2)$, все геодезические в $\SO(3)$  имеют счетное число сопряженных точек.

\subsubsection{Множество разреза}\label{subsubsec:so3_cut}
\begin{theorem}
Множество разреза на $\SO(3)$  имеет стратификацию
\begin{align*}
&\Cut = \Cut^{\loc} \cup \Cut^{\glob}, \\
&\Cut^{\loc} = e^{\k} \setminus \{\Id\}, \\
&\Cut^{\glob} = \left\{\Pi\left(\begin{array}{c} \a \\ \b \end{array}\right) \mid \a, \b \in \C, \ \Re \a = 0, \ \Im^2 \a + |\b|^2 = 1 \right\}.
\end{align*}
\end{theorem}

Топологически $\Cut^{\loc}$  есть интервал (окружность $e^{\k}$  с выколотой точкой $\Id$), а $\Cut^{\glob}$  есть проективная плоскость $\R P^2$. Начальная точка $\Id$  находится в замыкании локальной компоненты $\Cut^{\loc}$  и изолирована от глобальной компоненты $\Cut^{\glob}$. Эти компоненты пересекаются в единственной точке $e^{\pi k}$, поэтому $\Cut$  есть стратифицированное пространство.

\subsubsection{Расстояние}\label{subsubsec:so3_dist}

\begin{theorem}
Пусть $3\times 3$  матрица $g = (g_{ij}) \in \SO(3)$.  Тогда субриманово расстояние $d_0(g) = d(\Id, g)$  задается следующим образом.
\begin{itemize}
\item[$(1)$]
Если $g_{33} = - 1$, то $d_0 = \pi$.
\item[$(2)$]
Если $g_{33} = 1$  и $g \neq \Id$, то $d_0 = \frac{2\pi}{\sqrt{1+\b^2}}$,  где $\b$  --- единственное решение системы уравнений
$$
\begin{cases}
\cos \frac{\pi \b}{\sqrt{1+\b^2}} = - \frac 12 \sqrt{1 + g_{11} + g_{22} + g_{33}}, \\
\sin \frac{\pi \b}{\sqrt{1+\b^2}} = \frac 12 \sgn(g_{21}-g_{12}) \sqrt{1 - g_{11}- g_{22}+ g_{33}}.
\end{cases}
$$
\item[$(3)$]
Если $-1< g_{33} < - 1$ и $\cos\left(\pi\sqrt{\frac{1+g_{33}}{2}}\right) = - \frac{g_{11}+g_{22}}{1+g_{33}}$, то $d_0 = \pi \sqrt{\frac 12 (1 - g_{33})}$.
\item[$(4)$]
Если $-1< g_{33} < - 1$ и $\cos\left(\pi\sqrt{\frac{1+g_{33}}{2}}\right) > - \frac{g_{11}+g_{22}}{1+g_{33}}$, то $d_0 = \frac{2}{\sqrt{1+\b^2}} \arcsin \sqrt{\frac 12 (1- g_{33})(1 + \b^2)}$, где $\b$  --- единственное решение системы уравнений
$$
\begin{cases}
\cos \left(- \frac{\b}{\sqrt{1+\b^2}} \arcsin \sqrt{\frac 12 (1- g_{33})(1 + \b^2)} + \arcsin\left(\b \sqrt{\frac{1-g_{33}}{1+g_{33}}}\right)\right) = 
\sqrt{\frac{1 + g_{11} + g_{22} + g_{33}}{2(1 + g_{33})}}, \\
\sin \left(- \frac{\b}{\sqrt{1+\b^2}} \arcsin \sqrt{\frac 12 (1- g_{33})(1 + \b^2)} + \arcsin\left(\b \sqrt{\frac{1-g_{33}}{1+g_{33}}}\right)\right) = 
\sgn(g_{21} - g_{12}) \sqrt{\frac{1-g_{11}-g_{22}+ g_{33}}{2(1 + g_{33})}}.
\end{cases}
$$
\item[$(5)$]
Если $-1< g_{33} < - 1$ и 
$\cos\left(\pi \sqrt{\frac{1 + g_{33}}{2}}\right) < - \frac{g_{11}+ g_{22}}{1+ g_{33}}$, то
$d_0 = \frac{2}{\sqrt{1+\b^2}} \left( \pi - \arcsin \sqrt{\frac 12 (1 - g_{33})(1 + \b^2)} \right)$,
 где $\b$  --- единственное решение системы уравнений
$$
\begin{cases}
\cos\left( \frac{\b}{\sqrt{1+\b^2}} \left( \pi - 
\arcsin \sqrt{ \frac 12 (1 - g_{33})(1 + \b^2)}\right) + 
\arcsin \left( \b \sqrt{\frac{1 - g_{33}}{1+ g_{33}}}\right)\right) = 
- \sqrt{\frac{1+ g_{11}+g_{22} + g_{33}}{2(1 + g_{33})}},\\
\sin\left( \frac{\b}{\sqrt{1+\b^2}} \left( \pi - 
\arcsin \sqrt{ \frac 12 (1 - g_{33})(1 + \b^2)}\right) + 
\arcsin \left( \b \sqrt{\frac{1 - g_{33}}{1+ g_{33}}}\right)\right) = 
\sgn(g_{21} - g_{12}) \sqrt{\frac{1 - g_{11} - g_{22}+ g_{33}}{2(1 + g_{33})}}.
\end{cases}
$$

\end{itemize}
\end{theorem}

\subsubsection{Сферы}\label{subsubsec:so3_sphere}
Рассмотрим другую модель группы $\SO(3)$. Пусть $M_1$  --- ориентированная двумерная сфера гауссовой кривизны 1, заданная в объемлющем пространстве $\R^3$  равенством $x^2+y^2+z^2 = 1$;  риманова метрика $d_R$ на $M_1$  индуцирована евклидовой метрикой на $\R^3$;  ориентация сферы $M_1$  задана ее внешней нормалью. Пусть $V_1$  есть расслоение единичных касательных векторов к $M_1$.

Пусть распределение $\D$  на $V_1$  есть горизонтальное распределение связности Леви-Чивита, а расстояние на $V_1$  определяется как
$$
d(x,y) = \inf d_R(l),
$$
 где нижняя грань берется по всем кривым в $M_1$, горизонтальные лифты которых в $V_1$  соединяют $x$  и $y$.  Здесь $d_R(l)$  есть длина кривой $l$  в метрике $d_R$.

Многообразие $V_1$  с метрикой $d$  изометрично группе $\SO(3)$  с субримановой $\k \oplus \ppp$  структурой, определенной в п. \ref{subsubsec:so3_struct}.

Введем в $V_1$  систему координат
$
\{(r, \a, \b) \mid r \geq 0, \ -\pi \leq \a \leq \pi, \ -\pi \leq \b \leq \pi\}
$ 
 с началом в некотором элементе $v_O \in V_1$.  
Тогда точка $A \in M_1$  имеет декартовы координаты $(\cos \b \sin r, \sin \b \sin r, \cos r)$,  а вектор $v_A \in V_1$  имеет декартовы компоненты 
$$
(\cos r \cos \b \cos(\a-\b) - \sin \b \sin(\a-\b), \cos r \sin \b \cos(\a-\b) + \cos \b \sin(\a-\b), - \sin r \cos(\a-\b)).
$$
Поэтому координатное отображение
$$
\map{f}{\{(r,\a,\b) \mid 0 \leq r \leq \pi, \ -\pi \leq \a \leq \pi, \ -\pi \leq \b \leq \pi\}}{V_1}
$$
можно продолжить по непрерывности до отображения $g$  с включением случая $r = \pi$. В этом случае точка $A$  имеет декартовы координаты $(0, 0, -1)$, вектор $V_A$ --- декартовы компоненты $(-\cos(\a-2\b), \sin(\a-2\b), 0)$.
Тогда 
$$
\map{g}{\{(r,\a.\b) \mid 0 \leq r \leq \pi, \ -\pi \leq \a \leq \pi, \ -\pi \leq \b \leq \pi\}}{V_1}
$$
 есть отображение отождествления замкнутого полнотория на пространство $V_1$.  При этом кривая $\b = \frac 12 \a + c$, $c \in \R$,  на граничном торе $\T^2$  переходит под действием отображения $g$  в один элемент пространства~$V_1$.

\begin{theorem}
Диаметр пространства $V_1$  с метрикой $d$  равен $\sqrt 3 \pi$.

Введем в $V_1$  систему координат $\{r,\a,\b\}$  с началом в некотором элементе $v_O$. В этой системе координат сфера с центром $v_O$  радиуса $\sqrt 3 \pi $  есть точка $(0, \pi, 0)$. Сфера с центром $v_O$  радиуса $0 < R < \sqrt 3 \pi$  --- поверхность вращения вокруг оси $\a$  той части кривой $S_R$, определяемой параметрическими уравнениями $r = \pm r(t)$, $\a = \pm \a(t)$, $1 \leq t \leq \frac{2\pi}{R}$,  которая расположена в полосе $-\pi \leq \a \leq \pi$  на плоскости $\b = 0$, причем
\begin{itemize}
\item[$(1)$]
Если $0 < R < \pi$, то
\begin{align*}
&r(t) = 2 \arcsin\left(\frac 1t \sin \frac{Rt}{2}\right), \qquad 1 \leq t \leq \frac{2\pi}{R},\\
&\a(t) = 
\begin{cases}
2 \left| \frac R2 \sqrt{t^2 - 1} - \arcsin \frac{\sqrt{t^2-1} \sin \frac {Rt}{2}}{\sqrt{t^2 - \sin^2 \frac {Rt}{2}}}\right|, & \text{ если } 1 \leq t < \frac{\pi}{R}, \\
2 \left( \pi - \frac R2 \sqrt{t^2 - 1} - \arcsin \frac{\sqrt{t^2-1} \sin \frac {Rt}{2}}{\sqrt{t^2 - \sin^2 \frac {Rt}{2}}}\right), & \text{ если } \frac{\pi}{R} \leq t < \frac{2\pi}{R}.
\end{cases}
\end{align*}
\item[$(2)$]
Если $\pi \leq R < \sqrt 3 \pi$, то
\begin{align*}
&r(t) = 2 \arcsin\left(\frac 1t \sin \frac{Rt}{2}\right), \qquad 1 \leq t \leq \frac{2\pi}{R},\\
&\a(t) = 
2 \left(\pi -  \frac R2 \sqrt{t^2 - 1} - \arcsin \frac{\sqrt{t^2-1} \sin \frac {Rt}{2}}{\sqrt{t^2 - \sin^2 \frac {Rt}{2}}}\right),  \qquad 1 \leq t < \frac{\pi}{R}.
\end{align*}
\end{itemize}
При $R = \pi$  значение $\a(1)$  определяется по непрерывности и равно $\pi$.
\end{theorem}

Сферы пространства $V_1$  радиуса $R \in (0, \sqrt 3 \pi)$,  $R \neq \pi$,  гомеоморфны $S^2$.  
Сфера  радиуса  $R = \pi$ гомеоморфна сфере $S^2$, у которой диаметрально противоположные точки отождествлены. 
Сфера  радиуса  $R = \sqrt 3\pi$  есть точка. При $R \in (0, \pi)$  сфера имеет две конические особенности. Сфера радиуса $R = \pi$  диффеоморфна двум пересекающимся по окружности одинарным конусам, направленным в разные стороны, на общей окружности которых отождествлены диаметрально противоположные точки. Сфера радиуса $R \in (\pi, \sqrt 3 \pi)$  диффеоморфна двум пересекающимся по окружности одинарным конусам, направленным в разные стороны.

\subsubsection{Библиографические комментарии}
Разделы \ref{subsubsec:so3_geod},  \ref{subsubsec:so3_cut}  опираются на статью \cite{boscain_rossi}, а разделы \ref{subsubsec:so3_dist} и \ref{subsubsec:so3_sphere} ---  на статьи \cite{BZ15b}  и \cite{ber01}  соответственно. Рассматриваемой субримановой задаче на $\SO(3)$  посвящена также работа \cite{BZ15a}.

Насколько нам известно, впервые геодезические и сферы для осесимметричной субримановой задачи на $\SO(3)$  описаны в \cite{gersh84, versh_gersh2}.

\subsection{Осесимметричная субриманова  задача   на группе $\SL(2)$}  \label{subsec:sr_sl2}

\subsubsection{Группа Ли и алгебра Ли}\label{subsubsec:sl2_gr}
 Группа Ли \ddef{$\SL(2)$}   есть группа унимодулярных   $2\times 2$  вещественных матриц
$$
\SL(2)= \{ g \in \GL(2,\R) \mid   \det g = 1\}.
$$
Эта группа некомпактна, связна и неодносвязна: ее фундаментальная группа есть $\Z$.
Алгебра Ли  группы Ли $G = \SL(2)$ есть алгебра бесследовых $2 \times 2 $  вещественных матриц:
$$
\sla(2) = \{X \in \gl(2,\R) \mid \Tr  X = 0\}.
$$
Матрицы
\be{sl2p12k}
p_1 = \frac 12 
\left(\begin{array}{cc} 1 & 0 \\
0 & -1
\end{array}\right), 
\quad
p_2 = \frac 12 
\left(\begin{array}{cc} 0 & 1 \\
1 & 0
\end{array}\right), 
\quad
k = \frac 12 
\left(\begin{array}{cc} 0 & -1 \\
1 & 0
\end{array}\right)
\ee
образуют базис в алгебре Ли $\gg = \sla(2)$  с таблицей умножения
$$
[p_1, p_2] = - k, \quad [p_2, k] = l_1, \quad [k, p_1] = p_2.
$$
Форма Киллинга для $\sla(2)$  есть $\Kil(X,Y) = 4 \Tr(XY)$,  поэтому  $\Kil(p_i, p_j) = 2 \d_{ij}$.

Подпространства 
$$
\k = \R k, \quad \ppp = \spann(p_1, p_2)
$$
 образуют картановское разложение в $\sla(2)$;  оно единственно т.к. $\k$ должно быть максимальной компактной подалгеброй.

Векторы $p_1, p_2$  образуют ортонормированный репер для скалярного произведения $\langle\cdot, \cdot\rangle = \frac 12 \Kil(\cdot, \cdot)$, суженного на $\ppp$. 

\subsubsection{Субриманова $\k \oplus \ppp$  структура}\label{subsubsec:sl2_struct}
Рассмотрим левоинвариантную субриманову структуру на $\SL(2)$  с ортонормированным репером
$$
X_1(g) = L_{g*} p_1, \quad 
X_2(g) = L_{g*} p_2, \qquad
g \in \SL(2),
$$
то есть распределение $\D_g = L_{g*} \ppp$  со скалярным произведением
$\langle v_1, v_2\rangle_g = \langle L_{g^{-1}*} v_1,  L_{g^{-1}*} v_2 \rangle$. Эта структура есть  $\k \oplus\ppp$  структура на $\SL(2)$.

\subsubsection{Геодезические и симметрии}\label{subsubsec:sl2_geod}

Анормальные траектории постоянны.

Пусть $h_i(\lam) = \langle \lam, X_i(g)\rangle$, $i = 1, 2, 3$, $X_3 = [X_1, X_2]$, $H = \frac 12 (h_1^2 + h_2^2)$,  и  $C = \gg^* \cap \{ H = \frac 12\}$. Далее, пусть 
\begin{align*}
&\lam = (h_1, h_2, h_3, \Id) \in C, \\
&h_1 = \cos \t, \quad h_2 = \sin \t, \quad h_3 = c.
\end{align*}
Тогда экспоненциальное отображение 
$$
\map{\Exp}{C \times \R_+}{G}, \qquad
(\lam, t) \mapsto g(t),
$$
 параметризуется следующим образом:
\begin{align*}
&\Exp(\lam, t) = \Exp(\t, c, t) = e^{(\cos \t p_1 + \sin \t p_2 + c k)t} e^{- c k t} = \\
&=
\left(\begin{array}{cc}
K_1 \cos(c \frac t2) + K_2(\cos(\t + c \frac t2) + c \sin(c \frac t2)) & 
K_1 \sin(c \frac t2) + K_2(\sin(\t + c \frac t2) - c \cos(c \frac t2))\\
-K_1 \sin(c \frac t2) + K_2(\sin(\t + c \frac t2) + c \cos(c \frac t2))&
K_1 \cos(c \frac t2) + K_2(-\cos(\t + c \frac t2) + c \sin(c \frac t2))
\end{array}\right),
\end{align*}
где
\begin{align*}
&K_1 = 
\begin{cases}
\ch(\sqrt{1-c^2} \frac t2) &\text{ при } c \in [-1, 1], \\
\cos(\sqrt{c^2-1} \frac t2) &\text{ при } c \in (-\infty,-1) \cup (1, + \infty), 
\end{cases}\\
&K_2 = 
\begin{cases}
\sh(\sqrt{1-c^2} \frac t2) &\text{ при } c \in (-1, 1), \\
\frac t2 &\text{ при } c \in \{-1, 1\}, \\
\frac{\sin(\sqrt{c^2-1} \frac t2)}{\sqrt{c^2-1}} &\text{ при } c \in (-\infty,-1) \cup (1, + \infty).
\end{cases}
\end{align*}

Семейство геодезических имеет следующие симметрии:
\begin{itemize}
\item
 вращения
$\Exp(\t, c, t) = e^{z_0 k}e^{x p_1 + y p_2}$,  где 
$$
\left(\begin{array}{c} x \\ y \end{array}\right) 
= 
\left(\begin{array}{cc} \cos \t  &  -\sin \t \\ \sin \t & \cos \t  \end{array}\right) 
\left(\begin{array}{c} x_0 \\ y_0 \end{array}\right), 
$$
 а $(x_0, y_0, z_0)$  определяются условием $\Exp(0, c, t) = e^{z_0 k}e^{x_0 p_1 + y_0 p_2}$,
\item
 отражения
$\Exp(\t, -c, t) = e^{-z_0 k}e^{x p_1 + y p_2}$,  где 
$$
\left(\begin{array}{c} x \\ y \end{array}\right) 
= 
\left(\begin{array}{cc} \cos 2 \t  &  \sin 2 \t \\ \sin 2 \t & - \cos 2 \t  \end{array}\right) 
\left(\begin{array}{c} x_0 \\ y_0 \end{array}\right), 
$$
 а $(x_0, y_0, z_0)$  определяются условием $\Exp(\t c, t) = e^{z_0 k}e^{x_0 p_1 + y_0 p_2}$.
\end{itemize}

\subsubsection{Сопряженные точки}
Геодезические $\Exp(\t, c, t)$, $| c| \leq 1$,  не содержат сопряженных точек.

Если $|c| > 1$, то сопряженные времена вдоль геодезической $\Exp(\t, c, t)$  следующие:
\begin{align*}
&t_{2n-1} = \frac{2 \pi n}{\sqrt{c^2-1}}, \\
&t_{2n } = \frac{2 x_n}{\sqrt{c^2-1}}, \qquad n \in \N,
\end{align*}
где $\{x_1, x_2, \dots, \}$  суть упорядоченные по возрастанию положительные корни уравнения $x = \tg x$.

Соответствующая $m$-я каустика
$$
\Conj^m = \{ \Exp(\lam, t) \mid t = t_m^{\conj}(\lam), \ \lam \in C \}
$$
есть:
\begin{align*}
&\Conj^{2n-1} = e^{\k} \setminus \{\Id\}, \\
&\Conj^{2k} = 
\left\{
\left(\begin{array}{cc}
\cos x_n \cos y_n + \frac{\sin x_n}{\sqrt{c^2-1}}(\cos \t + c \sin y_n)) &
\cos x_n \sin y_n + \frac{\sin x_n}{\sqrt{c^2-1}}(\sin \t - c \cos y_n)) \\
- \cos x_n \sin y_n + \frac{\sin x_n}{\sqrt{c^2-1}}(\sin \t + c \cos y_n)) &
\cos x_n \cos y_n + \frac{\sin x_n}{\sqrt{c^2-1}}(-\cos \t + c \sin y_n))
\end{array}\right)\right.
\\
& \qquad\qquad\qquad\qquad\qquad
\left.
\vphantom{
\left(\begin{array}{cc}
\cos x_n \cos y_n + \frac{\sin x_n}{\sqrt{c^2-1}}(\cos \t + c \sin y_n)) &
\cos x_n \sin y_n + \frac{\sin x_n}{\sqrt{c^2-1}}(\sin \t - c \cos y_n)) \\
- \cos x_n \sin y_n + \frac{\sin x_n}{\sqrt{c^2-1}}(\sin \t + c \cos y_n)) &
\cos x_n \cos y_n + \frac{\sin x_n}{\sqrt{c^2-1}}(-\cos \t + c \sin y_n))
\end{array}\right)
}
\mid 
c \in \R, \ \t \in \R /(2 \pi \Z)
\right\},
\end{align*}
где $y_n = \frac{c x_n}{\sqrt{c^2-1}}$.

\subsubsection{Множество разреза}\label{subsubsec:sl2_cutloc}

\begin{theorem}
Множество разреза есть стратифицированное пространство
\begin{align*}
&\Cut = \Cut^{\loc} \cup \Cut^{\glob}, \\
&\Cut^{\loc} = e^{\k} \setminus \{\Id\} = 
\left\{\left(
\begin{array}{cc} \cos \a & - \sin \a \\ \sin \a & \cos \a
\end{array}\right) \mid \a \in (0, 2 \pi)\right\},\\
&\Cut^{\glob} = e^{2 \pi k} e^{\ppp}  = 
\left\{g \in \SL(2)
\mid g = g^T, \ \Tr g < - 2 \right\}.
\end{align*}
\end{theorem}

Топологически, $\Cut^{\loc}$  есть интервал (окружность $e^{\k}$  с выколотой точкой $\Id$), а $\Cut^{\glob}$  есть плоскость~$\R^2$.

\subsubsection{Геодезическая орбитальность}\label{subsubsec:sl2_orb}

\begin{theorem}
Группа $\SL(2)$  с рассматриваемой субримановой структурой геодезически орбитальна.
\end{theorem}

\subsubsection{Время разреза}\label{subsubsec:sl2_cuttime}
 Пусть $\lam = (\t, c) \in C$. Опишем время разреза $\tcut(\lam)$  вдоль соответствующей геодезической.

\begin{proposition}
$\tcut(\t, 0) = + \infty$.
\end{proposition}

\begin{theorem}
Пусть $c \neq 0$, тогда число $T = \tcut(\t, c) \in (0, + \infty)$  выражается следующим  образом.
\begin{itemize}
\item[$(1)$]
Если $|c| \geq \frac{2}{\sqrt 3}$, то $T = \frac{2\pi}{\sqrt{c^2-1}}$.
\item[$(2)$]
Если $|c| =1$, то $T$  принадлежит интервалу $(2\pi, 3 \pi)$  и удовлетворяет системе уравнений
$$
\cos \frac T2 = - \frac{1}{\sqrt{1+T^2/4}}, \qquad
\sin \frac T2 = - \frac{T/2}{\sqrt{1+T^2/4}}.
$$
\item[$(3)$]
Если $c^2 < 1$, то $T$  принадлежит интервалу $(2\pi/|c|, 3 \pi/|c|)$  и удовлетворяет системе уравнений
$$
\cos kx = - \frac{1}{\sqrt{1+k^2\th^2x}}, \qquad
\sin kx = - \frac{k\th x}{\sqrt{1+k^2\th^2x}},
$$
где
$$
k = \frac{|c|}{\sqrt{1-c^2}}, \quad
x = \frac{T \sqrt{1-c^2}}{2} = \frac{T}{2 \sqrt{1+k^2}}.
$$
\item[$(4)$]
Если $|c| = \frac{3}{2 \sqrt 2}$, то $T =  2 \sqrt 2 \pi$.
\item[$(5)$]
Если $\frac{3}{2 \sqrt 2} < |c| <  \frac{2}{\sqrt 3}$, то $\frac{3\pi}{|c|} < T < 2 \pi(|c| + \sqrt{c^2-1}) < \frac{4\pi}{|c|}$  и $T$   удовлетворяет системе уравнений
$$
\cos kx = \frac{1}{\sqrt{1+k^2\tg^2x}}, \qquad
\sin kx =  \frac{k\tg x}{\sqrt{1+k^2\tg^2x}} < 0,
$$
где
\be{sl2_kx}
k = \frac{|c|}{\sqrt{c^2-1}}, \quad
x = \frac{T \sqrt{1-c^2}}{2} = \frac{T}{2 \sqrt{k^2-1}}.
\ee
\item[$(6)$]
Если $1 < |c| <  \frac{3}{2 \sqrt 2}$, то $\frac{2\pi}{|c|} < 2 \pi(|c| + \sqrt{c^2-1}) < T  < \frac{3\pi}{|c|}$  и $T$   удовлетворяет системе уравнений
$$
\cos kx = -\frac{1}{\sqrt{1+k^2\tg^2x}}, \qquad
\sin kx = - \frac{k\tg x}{\sqrt{1+k^2\tg^2x}} < 0,
$$
где $k$ и $x$  определяются формулами \eq{sl2_kx}.
\end{itemize}
\end{theorem}

В следующей теореме описаны свойства монотонности и регулярности времени разреза.

\begin{theorem}
Функция $T(c) = \tcut(\t, c)$  имеет следующие свойства:
\begin{itemize}
\item[$(1)$]
$T(c)$  строго убывает на промежутках $(0, \frac{3}{2 \sqrt 2}]$, $[\frac{2}{\sqrt 3}, + \infty)$  и строго возрастает на отрезке $[\frac{3}{2 \sqrt 2}, \frac{2}{\sqrt 3}]$.
\item[$(2)$]
$T(c)$  непрерывна, кусочно вещественно аналитична и $T(0, + \infty) = (0, + \infty)$.
\item[$(3)$]
$T(c)$  имеет локальный минимум $2 \sqrt 2 \pi$  при $c = \frac{3}{2 \sqrt 2}$  и локальный максимум $2 \sqrt 3 \pi$  при $c = \frac{2}{\sqrt 3}$.
\end{itemize}
\end{theorem}

\subsubsection{Библиографические комментарии}
Разделы \ref{subsubsec:sl2_struct}--\ref{subsubsec:sl2_cutloc}  опираются на статью \cite{boscain_rossi}, а разделы \ref{subsubsec:sl2_orb}--\ref{subsubsec:sl2_cuttime} ---  на работу \cite{BZ16}.

  Рассматриваемой субримановой задаче на $\SL(2)$  частично посвящена также работа \cite{ber01}.
	
	В работе \cite{versh_gersh86}  описаны динамические свойства геодезического потока для рассматриваемой в данном разделе субримановой структуры на $\SL(2)$.

\subsection{Осесимметричные римановы задачи на группах $\PSL(2; \R)$ и $\SL(2; \R)$}  \label{subsec:riemsl2}

\subsubsection{Постановка задачи на $\PSL(2; \R)$}
Пусть  $\SL(2; \R)$ есть группа вещественных $2\times 2$ матриц с единичным определителем, а $\PSL(2; \R) = \SL(2; \R)/ \{ \pm \Id\}$.
Обозначим группу Ли $G = \PSL(2; \R)$ и ее алгебру Ли $\gg = \sll(2;\R)$.

Рассмотрим левоинвариантную риманову структуру на группе Ли $G$, она задается квадратичной формой на $\gg$ с собственными значениями $I_1$, $I_2$, $I_3 > 0$. Выберем базис $e_1, e_2, e_3 \in \gg$, в котором форма Киллинга и риманова метрика имеют матрицы $\diag(1, 1, -1)$ и $\diag(I_1, I_2, I_3)$ соответственно.

Отождествим $\gg$ с $\gg^*$ с помощью формы Киллинга, тогда базис $e_1, e_2, e_3 \in \gg$ перейдет в базис $\eps_1, \eps_2, \eps_3 \in \gg^*$. 

Пусть $p = p_1 \eps_1 + p_2\eps_2 + p_3\eps_3 \in\gg^*$. Введем обозначение:
$$
\Kil(p) = p_1^2 + p_2^2 - p_3^2, \qquad |p| = \sqrt{|\Kil(p)|}, \qquad \type(p) = \sgn(-\Kil(p)),
$$
где $\Kil(p)$ есть значение квадратичной формы Киллинга на ковекторе $p$. Этот ковектор называется времениподобным,  светоподобным  или пространственноподобным, если $\type(p)$ равно 1, 0 или $-1$
 соответственно.

Римановы кратчайшие суть решения задачи оптимального управления
\begin{equation}\label{riemsl21}
\begin{split}
&\dot{Q} = Q \Omega,  \quad Q \in G, \quad  \Omega = u_1 e_1 + u_2 e_2 + u_3 e_3 \in \mathfrak{g},   \quad (u_1, u_2, u_3) \in \R^3, \\
&Q(0) = \Id, \quad  Q(t_1) = Q_1, \\
&\frac{1}{2} \int_0^{t_1}{(I_1 u_1^2 + I_2 u_2^2 + I_3 u_3^2) \ dt} \rightarrow \min.
\end{split}
\end{equation}

Далее рассматривается случай осесимметричной метрики: $I_1 = I_2$. Обозначим через
$$
\eta = -\frac{I_1}{I_3} - 1 < -1
$$
параметр римановой метрики, измеряющий вытянутость малых сфер. Для $p \in \gg^*$, $|p| \neq 0$, обозначим
$$
\bar{p} = \frac{p}{|p|}, \qquad \tau(p) = \frac{t |p|}{2 I_1}.
$$

Через $R_{v, \f}$ обозначим поворот трехмерного ориентированного евклидова пространства вокруг оси $\R v $ на угол $\f$ в положительном направлении.

\subsubsection{Геодезические на $\PSL(2; \R)$}
Натурально параметризованные геодезические соответствуют начальным импульсам
$$
p \in C = \gg^* \cap \left\{ H = \dfrac 12 \right\},
$$
где $H(p) = \dfrac{1}{2} \left(\dfrac{p_1^2}{I_1}+  \dfrac{p_2^2}{I_2}+\dfrac{p_3^2}{I_3} \right)$ есть максимизированный гамильтониан принципа максимума Понтрягина.
\begin{theorem}
Геодезическая $Q(t)$, начинающаяся в единице и имеющая начальный импульс $p \in C$, есть произведение двух однопараметрических групп:
\begin{equation}\nonumber
Q(t) = \exp \left(\frac{tp}{I_1}\right) \exp \left(\frac{t \eta p_3 e_3}{I_1} \right).
\end{equation}
\end{theorem}

Форма Киллинга  есть функция Казимира на $\gg^*$. Поэтому вдоль экстремалей $\type(p) \equiv \const$, то есть тип ковектора есть интеграл гамильтоновой системы.

\subsubsection{Модель группы $\PSL(2; \R)$}
Рассмотрим группу $\SU(1, 1)$, реализованную как группа  сплит-кватернионов единичной длины
$$
\SU(1, 1) = \{ q_0 + q_1 i + q_2 j + q_3 k \ | \ q_0^2 - q_1^2 - q_2^2 + q_3^2 = 1, \ q_0, q_1, q_2, q_3 \in \R\}.
$$
Умножение сплит-кватернионов дистрибутивно и удовлетворяет соотношениям 
\begin{equation*}
i^2 = j^2 = 1, \qquad k^2 = -1, \qquad
ij = -k, \qquad
jk = i, \qquad
ki = j.
\end{equation*}
Существует изоморфизм 
$$
\psi: \SL(2; \R) \rightarrow \SU(1, 1), \quad
\psi
\left(
\begin{array}{ll}
a & b \\
c & d
\end{array}
\right)
= \frac{a+d}{2} + \frac{a-d}{2} i + \frac{b+c}{2} j + \frac{c-b}{2} k, 
$$
где $ad - bc =1$, $\quad $ $a, b, c, d \in \R$.

Рассмотрим проекцию группы $\SU(1, 1)$ на трехмерное вещественное пространство с координатами $q_1, q_2, q_3$. Условие 
$$
q_3^2 - q_1^2 - q_2^2 = 1 - q_0^2 \leqslant 1
$$
означает, что образ группы $\SU(1, 1)$ есть область между двумя полостями гиперболоида $q_3^2 - q_1^2 - q_2^2 = 1$. Для фиксированных $q_1, q_2, q_3 $ (таких, что $q_3^2 - q_1^2 - q_2^2 \neq 1$)
значение $q_0$ можно выбрать двумя разными способами. Поэтому группа $\SU(1, 1)$ есть объединение двух таких областей с отождествленными граничными точками (соответствующими условию $q_0 = 0$). Группа 
$\SU(1, 1)$ гомеоморфна открытому полноторию.

Группу $\PSL(2; \R) \cong  \SU(1, 1) / \{ \pm \Id\}$ можно представить как область между полостями гиперболоида  с  отождествленными противоположными точками  на полостях гиперболоида:
$(q_1, q_2, q_3) \sim (-q_1, -q_2, -q_3)$.

\subsubsection{Геодезические на $\SU(1, 1)$}
Рассмотрим левоинвариантную риманову задачу на группе $\SU(1, 1)$, являющуюся лифтом задачи \eq{riemsl21} на $\PSL(2; \R)$. 
\begin{theorem}
Геодезические на группе $\SU(1, 1)$, выходящие из единицы с начальным ковектором $p = p_1 \frac i2 + p_2\frac j2 + p_3\frac k2 \in C$, имеют  следующую параметризацию:
\begin{itemize}
\item[$(1)$]
для времениподобного ковектора $p$ ($p_3^2 - p_1^2 - p_2^2 > 0$)
\begin{equation}\label{riemsl2q0e}
\begin{aligned}
\begin{array}{ccl}
q_0^e(\tau) & = & \ct \ce - \pp \st \se, \\
\left(
  \begin{array}{l}
     q_1^e(\tau )\\
     q_2^e(\tau )
  \end{array}
\right) & = &  \st R_{e_3, -\tau \eta \pp}
\left(
  \begin{array}{l}
     \bar{p}_1\\
     \bar{p}_2
  \end{array}
\right), \\
q_3^e(\tau ) & = & \ct \se + \pp \st \ce,
\end{array}
\end{aligned}
\end{equation}
\item[$(2)$]
для светоподобного ковектора $p$ ($p_3^2 - p_1^2 - p_2^2 = 0$)
\begin{equation}\label{riemsl2q0p}
\begin{aligned}
\begin{array}{ccl}
q_0^p(t) & = & \ctl - \frac{t}{2 I_1} p_3 \stl, \\
\left(
  \begin{array}{l}
     q_1^p(t)\\
     q_2^p(t)
  \end{array}
\right) & = &  \frac{t}{2 I_1} R_{e_3, -\argl}
\left(
  \begin{array}{l}
     p_1\\
     p_2
  \end{array}
\right), \\
q_3^p(t) & = & \stl + \frac{t}{2 I_1}p_3 \ctl,
\end{array}
\end{aligned}
\end{equation}
\item[$(3)$]
для пространственноподобного ковектора   $p$ ($p_3^2 - p_1^2 - p_2^2 < 0$)
\begin{equation}\label{riemsl2q0h}
\begin{aligned}
\begin{array}{ccl}
q_0^h(\tau) & = & \cht \ce - \pp \sht \se, \\
\left(
  \begin{array}{l}
     q_1^h(\tau )\\
     q_2^h(\tau )
  \end{array}
\right) & = & \sht R_{e_3, -\tau \eta \pp}
\left(
  \begin{array}{l}
     \bar{p}_1\\
     \bar{p}_2
  \end{array}
\right), \\
q_3^h(\tau ) & = & \cht \se + \pp \sht \ce.
\end{array}
\end{aligned}
\end{equation}
\end{itemize}
\end{theorem}

\subsubsection{Сопряженные времена}\label{subsubsec:riemsl2_conj}
\begin{theorem}
Рассмотрим геодезическую на $G = \PSL(2; \R)$ или $G=\SL(2; \R)$, начинающуюся в единице, с начальным ковектором $p \in C$. 

Для времениподобного  начального импульса $p$ с $\pp \neq \pm 1$ есть две серии сопряженных времен:
$$
t_{2k-1} = \frac{2 I_1 \pi k}{|p|}, \qquad t_{2k} = \frac{2 I_1 \tau_k(p)}{|p|}, \qquad k \in \mathbb{N},
$$
где $\tau_k(p)$ есть $k$-ый положительный корень уравнения 
$$
\tan{\tau} = -\tau\eta\frac{1-\pp^2}{1+\eta\pp^2}.
$$
В случае $p_3 = \pm 1$ эти две серии сливаются в одну серию:
$$
t_k = \frac{2 I_1 \pi k}{|p|}, \qquad k \in \N.
$$

Для свето- и пространственноподобных начальных ковекторов $p$ соответствующие геодезические не имеют сопряженных точек.
\end{theorem}
\begin{corollary}
Первое сопряженное время для геодезической, соответствующей ковектору $p$, есть 
$$
t_{\conj}^1(p) = \left\{
\begin{array}{lcl}
\frac{2\pi I_1}{|p|}, & \text{при} & \type{(p)} = 1,\\
+\infty, & \text{при} & \type{(p)} \leqslant 0.\\
\end{array}
\right.
$$
\end{corollary}

\subsubsection{Время разреза и множество разреза}
Обозначим первые положительные  нули функций $q_0^e(\tau)$, $q_0^p(t)$, $q_0^h(\tau)$, см.~\eq{riemsl2q0e}--\eq{riemsl2q0h}:
$$
\begin{array}{lclclcl}
\tau^e_0(\pp) & = & \min\{\tau \in \R_+ \ | \ q_0^e(\tau, \pp) = 0 \}, \\
t^p_0(p) & = & \min\{t \in \R_+ \ | \ q_0^p(t, p) = 0 \}, \\
\tau^h_0(\pp) & = & \min\{\tau \in \R_+ \ | \ q_0^h(\tau, \pp) = 0 \}.
\end{array}
$$

Обозначим через 
$C^e$, $C^p$, $C^h$ 
времени-, свето- и
пространственноподобные части поверхности уровня гамильтониана $C$.

\begin{theorem}
\begin{itemize}
\item[$(1)$]
Если $\eta \leq - \frac 32$, то
$$
t_{\cut}(p) =
\left\{
\begin{array}{rcll}
\frac{2 I_1 \tau^e_0(\bar{p}_3)}{|p|}, & \text{при} & p \in C^e, & \\
t^p_0(p_3),                            & \text{при} & p \in C^p, & \\
\frac{2 I_1 \tau^h_0(\bar{p}_3)}{|p|}, & \text{при} & p \in C^h, & \pp \neq 0, \\
+\infty, & \text{при} & \pp = 0. & \\
\end{array}
\right.
$$
\item[$(2)$]
Если $\eta > - \frac 32$, то
$$
t_{\cut}(p) =
\left\{
\begin{array}{rcll}
\frac{2 I_1 \tau^e_0(\bar{p}_3)}{|p|}, & \text{при} & p \in C^e, & |\pp| > -\frac{3}{2\eta}, \\
\frac{2 I_1 \pi}{|p|},                 & \text{при} & p \in C^e, & |\pp| \leqslant -\frac{3}{2\eta}, \\
t^p_0(p_3),                            & \text{при} & p \in C^p, & \\
\frac{2 I_1 \tau^h_0(\bar{p}_3)}{|p|}, & \text{при} & p \in C^h, & \pp \neq 0,\\
+\infty, & \text{при} & \pp = 0. & \\
\end{array}
\right.
$$
\end{itemize}
\end{theorem}

Группу $\PSL(2; \R)$ можно
 интерпретировать
как группу собственных движений плоскости Лобачевского.

\begin{theorem}
\begin{itemize}
\item[$(1)$]
Если $\eta \leq - \frac 32$, то множество разреза есть плоскость, состоящая из центральных симметрий
\begin{align*}
&Z = \Pi(\{q \in \SU(1, 1) \mid q_0 = 0  \}),\\
& \Pi: \SU(1, 1) \to \SU(1, 1)/\{\pm \Id\} \cong \PSL(2; \R).
\end{align*}
\item[$(2)$]
Если $\eta > - \frac 32$, то множество разреза есть стратифицированное многообразие $Z\cup { R_{\eta}}$, где
$$
{R_{\eta}} = \{ R_{0, \pm \varphi} \in \PSL_2(\R) \ | \ \varphi \in [-2\pi(1+\eta), \pi] \}
$$
есть отрезок, состоящий из некоторых поворотов вокруг центра модели Пуанкаре гиперболической плоскости.
\end{itemize}
\end{theorem}

\subsubsection{Радиус инъективности }
\begin{theorem}
Радиус   инъективности осесимметричной римановой метрики  на группе $\PSL(2;\R)$ равен
\begin{itemize}
\item[$(1)$]
$\pi \sqrt{I_1} \sqrt{-\frac{1}{1+\eta}}\quad \text{при}\quad\eta \leqslant -2$; 
\item[$(2)$]
$\pi \sqrt{I_1} \sqrt{-\frac{\eta+4}{\eta}}\quad\text{при}\quad-2 < \eta \leqslant \frac{-3-\sqrt{73}}{8}$; 
\item[$(3)$]
 $2\pi \sqrt{I_1}\sqrt{-(1+\eta)}\quad\text{при}\quad\frac{-3-\sqrt{73}}{8} < \eta < -1$.
\end{itemize} 
\end{theorem}

\subsubsection{Осесимметричная левоинвариантная риманова задача на $\SL(2;\R)$}
Рассмотрим риманову задачу на группе $\SL(2;\R)$, являющуюся  лифтом задачи \eq{riemsl21} на $\PSL(2;\R)$. Геодезические для задачи на $\SL(2;\R)$ задаются прежними формулами \eq{riemsl2q0e}--\eq{riemsl2q0h}. Сопряженные времена для $\SL(2;\R)$ выражаются так же, как в п.~\ref{subsubsec:riemsl2_conj}.

Время разреза и множество разреза для задачи на $\SL(2;\R)$ описываются следующим образом.
\begin{theorem}
Пусть $p \in C$. Время разреза для соответствующей  геодезической на $\SL(2;\R)$ есть
$$
t_{\cut}(p) = \dfrac{2I_1}{|p|}\tau_{\cut}(\pp),
$$
где 
\begin{align*}
&\tau_{\cut}(\pp) = \tau_3^e(\pp) \quad \text{при} \quad \eta \leq - \dfrac 32,\\
&\tau_{\cut}(\pp) =
\begin{cases}
&\pi, \quad \pp \in [1, - {2}/{\eta}]\\
&\tau_3^e(\pp), \quad \pp > -{2}/{\eta}
\end{cases}
\quad \text{при} \quad \eta > - \dfrac 32,
\end{align*}
а $\tau_3^e(\pp)$ есть первый положительный корень функции $q_3^e(\tau, \pp)$ \eq{riemsl2q0e}.
\end{theorem}

\begin{theorem}
\begin{itemize}
\item[$(1)$]
Если $\eta \leq - \dfrac 32$, то множество разреза есть плоскость
$$
H := \{q \in \SU(1, 1) \ | \ q_3 = 0 \},
$$
которая представляется плоскостью гиперболических изометрий, соответствующих пучкам ультрапараллельных прямых, симметричных в диаметрах модели Пуанкаре гиперболической плоскости.
\item[$(2)$]
Если $\eta > - \dfrac 32$, то множество разреза есть стратифицированное многообразие 
$$
H \cup T_{\eta},
$$
где 
$T_{\eta} = \{ q = \pm (\cos(2\pi\pp) + \sin(2\pi\pp)k) \ | \ \pp \in [1, -\frac{2}{\eta}]  \}$
есть отрезок, состоящий из некоторых поворотов вокруг центра модели Пуанкаре гиперболической плоскости.
\end{itemize}
\end{theorem}
\subsubsection{Связь с субримановой задачей}
Пусть $G = \PSL(2;\R)$ или $\SL(2;\R)$. Отождествим алгебру Ли $\gg$ с пространством чисто мнимых сплит-кватернионов, и рассмотрим разложение 
$$
\gg = \k \oplus \ppp,
$$
где $\k = \R k$ и $\ppp = \R i \oplus \R j$.

Рассмотрим левоинвариантное распределение $\D$ на $TG$,  полученное левыми сдвигами подпространства $\ppp \subset \gg$. Снабдим распределение $\D$ левоинвариантной римановой структурой, полученной левыми сдвигами из формы Киллинга. Полученная субриманова структура  есть субриманова $\k \oplus \ppp$ структура на группе  $\PSL(2;\R)$ или $\SL(2;\R)$, см.~раздел~\ref{subsec:sr_sl2}.
\begin{theorem}
Для указанной субримановой $\k \oplus \ppp$-задачи на группе $\PSL(2;\R)$ (или $\SL(2;\R)$) следующие объекты:
\begin{itemize}
\item[$(1)$]
параметризация геодезических,
\item[$(2)$]
сопряженные времена,
\item[$(3)$]
каустика,
\item[$(4)$]
время разреза,
\item[$(5)$]
множество разреза
\end{itemize}
получаются из тех же объектов для осесимметричной левоинвариантной римановой задачи на $\PSL(2;\R)$ (или $\SL(2;\R)$ соответственно) при переходе к пределу $I_3 \to \infty$.
\end{theorem}
\subsubsection{Библиографические комментарии}
Результаты этого раздела получены в работе~\cite{sl2}.

\subsection{Осесимметричные римановы задачи на группах $\SO(3)$ и $\SU(2)$}  \label{subsec:riemso3}

\subsubsection{Постановка задачи на $\SO(3)$ }
Любая левоинвариантная риманова метрика на группе Ли $G = \SO(3)$ задается положительно определенной квадратичной формой $J$ на касательном пространстве $T_{\Id}\SO(3) = \so(3)$. Пусть $e_1, e_2, e_3$ есть ортонормированный базис в алгебре Ли $\gg = \so(3)$ относительно формы Киллинга, в которой $J$ диагональна. Пусть $I_1, I_2, I_3 $ суть соответствующие собственные значения $J$.

Римановы кратчайшие для рассматриваемой метрики суть решения задачи оптимального управления
\begin{align}
&\dot{Q} = Q \Omega, \quad \quad \Omega = u_1 e_1 + u_2 e_2 + u_3 e_3 \in \gg, \label{riemso31}\\
&Q \in G, \quad  (u_1, u_2, u_3) \in \R^3, \label{riemso32}\\
&Q(0) = \Id, \quad \quad Q(t_1) = Q_1, \label{riemso33}\\
&\frac{1}{2} \int_0^{t_1}{(I_1 u_1^2 + I_2 u_2^2 + I_3 u_3^2) \ dt} \rightarrow \min.\label{riemso34}
\end{align}

Если существует треугольник со сторонами $I_1, I_2, I_3$, то задача имеет механическую интерпретацию: она описывает вращения твердого тела вокруг неподвижной точки по инерции. Числа $I_1, I_2, I_3 $ суть моменты инерции этого твердого тела.

Существование римановых кратчайших следует из теоремы Филиппова.

Далее рассматривается только \ddef{случай Лагранжа} $I_1 = I_2$ (случай $I_1 = I_2 = I_3$ называется \ddef{случаем Эйлера}).

\subsubsection{Параметризация экстремалей для $\SO(3)$}\label{subsubsec:riemso3_par}
Экстремали задачи \eq{riemso31}--\eq{riemso34} в случае Лагранжа имеют вид
\begin{align}
&Q(t) = R_{p, \frac{t}{I_1}|p|} \ R_{e_3,  \frac{t}{I_1} \eta p_3},\label{riemso35}\\
&p(t) = R_{e_3, -\frac{t}{I_1} \eta p_3} \ p, \nonumber
\end{align}
где $p(0) = p_1 \eps_1 + p_2 \eps_2 + p_3 \eps_3 \in \gg^*$, $\{\eps_i\}$ есть базис в $\gg^*$, двойственный к $\{ e_i\}$ относительно формы Киллинга, $Q(0) = \Id$, а $R_{v, \f}$  обозначает поворот  пространства $\R^3$ на угол $\f$  вокруг вектора $v \in \gg^*$ (направление поворота должно быть таким, чтобы для любого вектора $w \not \in \spann(v)$ репер $(w, R_{v, \f}w, v)$ был положительно ориентированным).

Параметр
$$
\eta = \frac{I_1}{I_3} - 1 > -1
$$
задает сплюснутость твердого тела. Элементы группы $\SO(3)$ отождествляются с ортогональными преобразованиями коалгебры $\gg^*$ с помощью коприсоединенного представления. 

Представим геодезические с помощью кватернионов. Рассмотрим двулистное накрытие
\be{riemso36}
\Pi : \{q \in \mathbb{H} \ | \ |q| = 1\} \simeq S^3 \rightarrow \SO(3).
\ee
Любой кватернион единичной нормы может быть записан в форме
$$
q = \cos \left(\frac{\varphi}{2}\right) + \sin \left(\frac{\varphi}{2}\right)(a_1i + a_2j + a_3k),
$$
где $a_1, a_2, a_3 \in \R$, $ a_1^2 + a_2^2 + a_3^2 = 1$. По определению
 $\Pi(\pm q) = R_{a, \varphi}$ есть поворот на угол $\varphi$ вокруг вектора
$a = a_1 e_1 + a_2 e_2 + a_3 e_3$.

Пусть $\pm (q_0(\tau ) + q_1(\tau ) i + q_2(\tau )j + q_3(\tau )k) \in \mathbb{H}$ есть лифт на  $S^3$ геодезической \eq{riemso35}, где использовано новое время $\tau = \frac{t}{2I_1}|p|$. 
Будем рассматривать натурально параметризованные геодезические, это соответствует начальному ковектору $p \in C = 
\left\{ p \in  \gg^* \mid \frac{p_1^2}{I_1} + \frac{p_2^2}{I_2} + \frac{p_3^2}{I_3} = 1\right\}$.
Тогда
\begin{equation}\label{riemso37}
\left\{
\begin{array}{ccl}
q_0(\tau ) & = & \cos (\tau) \cos (\tau \eta \bar{p}_3) - \bar{p}_3 \sin (\tau) \sin (\tau \eta \bar{p}_3), \\
\left(\begin{array}{l}q_1(\tau )\\q_2(\tau )\end{array}\right) & = & \sin (\tau) R_{e_3, -\tau
\eta \bar{p}_3} \left(\begin{array}{l}\bar{p}_1\\ \bar{p}_2\end{array}\right), \\
q_3(\tau ) & = & \cos (\tau) \sin (\tau \eta \bar{p}_3) + \bar{p}_3\sin (\tau) \cos(\tau \eta \bar{p}_3),
\end{array}
\right.
\end{equation}
где  $\bar{p} = \frac{p}{|p|}$, и где ограничение поворота $R_{e_3, \alpha }$ на плоскость $\spann\{e_1, e_2\}$ обозначено тем же символом.

\subsubsection{Сопряженное время}\label{subsubsec:riemso3_conj}
Обозначим через $t^1_{\conj}(p)$ первое сопряженное время для геодезической \eq{riemso35}, соответствующей  начальному  ковектору $p = p_1 \eps_1 + p_2 \eps_2 + p_3 \eps_3 \in C$. 
Пусть $\tau^1_{\conj}(p) = \dfrac{|p|}{2I_1}t^1_{\conj}(p)$.
\begin{theorem}
\begin{itemize}
\item[$(0)$]
Функция $\tau^1_{\conj} $ зависит только от $\pp \in [-1, 1]$.
\item[$(1)$]
Если $\eta \in (-1, 0]$, то $\tau^1_{\conj}(\pp)  = \pi$ для всех $\pp \in [-1, 1]$.
\item[$(2)$]
Если $\eta > 0$, то $\tau^1_{\conj}(\pp)$ есть первый  положительный корень уравнения
\begin{align*}
&\tan \tau = - \eta \dfrac{1 - \pp^2}{1 + \eta\pp^2}\,\tau,
\end{align*}
причем выполняется включение $\tau^1_{\conj}(\pp) \in (\frac{\pi}{2}, \pi]$. Равенство выполняется только при $\pp = \pm 1$.
\item[$(3)$]
Функция $\tau^1_{\conj}: [-1, 1] \to \R $ гладкая и возрастающая.
\end{itemize}

\end{theorem}

\subsubsection{Время разреза и множество  разреза  в $\SO(3)$}\label{subsubsec:riemso3_cut}
Обозначим первые положительные корни уравнений $q_0(\tau) = 0$ и $q_3(\tau) = 0$ через $\tau_0(\pp)$ и $\tau_3(\pp)$ соответственно. Функции $\tau_0(\pp)$ и $\tau_3(\pp)$ зависят от параметра $\eta$.
Если $\pp = 0$, то  $q_3(\tau) \equiv 0$, и значение $\tau_3(0)$ не определено. Таким образом,
\begin{align*}
&\tau_0: [-1, 1] \rightarrow (0, +\infty],\\
&\tau_3: [-1, 1] \setminus \{0\} \rightarrow (0, +\infty].
\end{align*}
Обозначим через $t_{\cut}(p)$ время разреза для геодезической, соответствующей начальному ковектору $p \in \gg^*$.
\begin{theorem}
\begin{itemize}
\item[$(1)$]
Если $\eta \geq - \dfrac 12$, то  $t_{\cut}(p) = \dfrac{2 I_1 \tau_0(\pp)}{|p|}$.
\item[$(2)$]
Если $\eta < - \dfrac 12$, то
\begin{align*}
&t_{\cut}(p) = 
\begin{cases}
&\dfrac{2\pi I_1}{|p|} \quad \text{при} \quad \dfrac{1}{2|\eta|}\leq |\pp| < 1,\\
&\dfrac{2I_1\tau_0(\pp)}{|p|} \quad \text{при} \quad  |\pp| < \dfrac{1}{2|\eta|}.
\end{cases}
\end{align*}
\end{itemize}
\end{theorem}

\begin{theorem}
\begin{itemize}
\item[$(1)$]
Если $\eta \geq - \dfrac 12$, то множество разреза есть проективная плоскость центральных
симметрий сферы
$$
P = \{R_{v, \pi} \in \SO(3) \ | \ v \in \R^3, \ v \neq 0 \}.
$$
\item[$(2)$]
Если $\eta < - \dfrac 12$, то множество разреза есть стратифицированное множество  $P\cup L_{\eta}$, где
$$
L_{\eta} = \{R_{e_3, \pm \f} \in \SO(3)\mid \f \in [2\pi(1+ \eta), \pi]   \}
$$
 есть отрезок, состоящий  из некоторых вращений вокруг оси $e_3$, соответствующей собственному значению метрики, отличному от двух других.
\end{itemize}
\end{theorem}

\subsubsection{Диаметр группы $\SO(3)$ в случае Лагранжа}\label{subsubsec:riemso3_diam}
\begin{theorem}
\begin{itemize}
\item[$(1)$]
Диаметр группы $\SO(3)$ в рассматриваемой римановой метрике равен 
$$
\begin{array}{lll}
2 \pi \sqrt{I_1} \sqrt{1 + \frac{1}{4 \eta}}, &  \text{при} & \eta \in (-1, -\frac{1}{2}), \\
\pi \sqrt{I_3}, &  \text{при} & \eta \in [-\frac{1}{2}, 0], \\
\pi \sqrt{I_1}, & \text{при} & \eta \in (0, +\infty).
\end{array}
$$
\item[$(2)$]
Множество точек,  наиболее удаленных от единицы, есть
$$
\begin{array}{lll}
\{ R_{\pm e_3, \pi} \}, &  \text{при} & \eta \in (-1, 0), \\
P, &  \text{при} & \eta = 0, \\
\{ R_{e, \pi} \ | \ e \in \spann\{e_1, e_2\} \}, & \text{при} & \eta \in (0, +\infty).
\end{array}
$$
\end{itemize}
\end{theorem}

\subsubsection{Осесимметричная  риманова задача на  $\SU(2)$ }\label{subsubsec:riemso3_su2}
Группа Ли $\SU(2)$ есть односвязная двулистная накрывающая группы $\SO(3)$,  см.~\eq{riemso36}.

Рассмотрим осесимметричную левоинвариантную риманову метрику на $\SU(2)$,  являющуюся поднятием метрики на $\SO(3)$, рассмотренной в предыдущих пунктах. Геодезические для нее задаются формулами 
\eq{riemso37}. Сопряженное время для задачи на $\SU(2)$ совпадает с сопряженным временем для задачи на $\SO(3)$, см.~п.~\ref{subsubsec:riemso3_conj}.
\begin{theorem}
 Время разреза для задачи на $\SU(2)$ есть
$$
t_{\cut}(\pp) = \dfrac{2I_1}{|p|}\,\tau_{\cut}(\pp),
$$
где 
\begin{align*}
&\tau_{\cut}(\pp) =
\begin{cases}
&\pi \quad \quad \quad \text{при} \quad  \eta \leq 0,\\
&\tau_3(\pp) \quad \text{при} \quad  \eta > 0.
\end{cases}
\end{align*}
Здесь $\tau_3(\pp)$ есть первый положительный корень уравнения~\eq{riemso37},  дополненный по непрерывности равенством $\tau_3(0) = \tau_{\conj}^1(0)$.
\end{theorem}

\begin{theorem}\label{th:su2_cut}
Множество разреза для задачи на $\SU(2)$ есть:
\begin{itemize}
\item[$(1)$]
отрезок
$$
T_{\eta} := \{-\cos (\pi \eta \bar{p}_3) - \sin (\pi \eta \bar{p}_3) k \in \mathbb{H} \ | \ \bar{p}_3 \in [-1, 1] \}
$$
при $-1 < \eta \leqslant 0$ (если $\eta = 0$, то $T_{\eta}$ есть точка $\{ -1\}$),
\item[$(2)$]
диск
$$
\{q_0(p, \tau_3(\bar{p}_3)) + q_1(p, \tau_3(\bar{p}_3)) i + q_2(p, \tau_3(\bar{p}_3)) j \in \mathbb{H}
\ | \ p \in C \},
$$
ограниченный окружностью из сопряженных точек 
$$
 \{\cos \tau_{\conj}(0) + \sin \tau_{\conj}(0) (i \cos \varphi  + j \sin \varphi) \in \mathbb{H}
\ | \ \varphi \in [0, 2 \pi] \},
$$
при $\eta > 0$.
\end{itemize}
\end{theorem}

\begin{theorem}
Диаметр группы $\SU(2)$ для левоинвариантной  римановой  метрики с  собственными значениями $I_1 = I_2$, $I_3 > 0$ равен
\begin{align*}
&\diam_{g(I_1, I_2, I_3)}\SU(2) =
\begin{cases}
&2\pi \sqrt{I_1}       \quad  \quad \text{при} \quad I_1 \leq I_3,\\
& 2\pi \sqrt{I_3} \quad   \quad \text{при}  \quad I_3 < I_1 \leq 2I_3,\\
& \dfrac{\pi I_1}{\sqrt{I_1 - I_3}}\quad \text{при} \quad  2 I_3 < I_1.
\end{cases}
\end{align*}
\end{theorem}

\subsubsection{Связь с субримановой задачей на $\SO(3)$ }\label{subsubsec:riemso3_connect}
Рассмотрим наряду с осесимметричной левоинвариантной римановой задачей на группе $\SO(3)$  также левоинвариантную субриманову $\k \oplus \ppp$-задачу на группе $\SO(3)$, см.~раздел~\ref{subsec:sr_so3}. 
\begin{theorem}
Для осесимметричной левоинвариантной римановой задачи на $\SO(3)$ следующие объекты сходятся к соответствующим объектам левоинвариантной субримановой  $\k \oplus \ppp$-задачи при $I_3 \to \infty$:
\begin{itemize}
\item[$(1)$]
параметризация геодезических,
\item[$(2)$]
сопряженное время,
\item[$(3)$]
первая каустика,
\item[$(4)$]
время разреза,
\item[$(5)$]
множество разреза.
\end{itemize}
\end{theorem}

\subsubsection{Библиографические комментарии}

Параметризация геодезических левоинвариантной римановой метрики на $\SO(3)$ (п.~\ref{subsubsec:riemso3_par}) есть классический результат Л.~Эйлера~\cite{lan_liv}. Первое  сопряженное время для осесимметричной задачи на $\SO(3)$ (п.~\ref{subsubsec:riemso3_conj}) было описано в работе~\cite{bates-fasso}. Результаты п.~\ref{subsubsec:riemso3_cut}, \ref{subsubsec:riemso3_diam}, 
\ref{subsubsec:riemso3_connect} получены в работе~\cite{so3}. Множество разреза для задачи на $\SU(2)$ при $I_1 > I_3$ (п.~(2) теоремы~\ref{th:su2_cut}) описано в работе~\cite{sakai}. Остальные результаты 
п.~\ref{subsubsec:riemso3_su2} получены в работах~\cite{so3, diamsu2}.

\subsection[Задача о качении сферы с прокручиванием, без проскальзывания]{Задача о качении сферы с прокручиванием, без проскальзывания}  \label{subsec:roll_besch}

\subsubsection{Постановка задачи}
Рассматривается механическая система, состоящая из сферы, катящейся по плоскости с прокручиванием, но без проскальзывания. Состояние такой системы в каждый момент времени характеризуется точкой на плоскости и ориентацией сферы в пространстве. Требуется перекатить сферу из заданного начального состояния в заданное конечное так, чтобы достигался минимум действия. Отсутствие проскальзывания означает, что точка контакта сферы и плоскости имеет нулевую мгновенную скорость, наличие прокручивания означает, что вектор угловой скорости сферы может быть направлен в произвольном направлении.

Эта задача является естественной модификацией задачи об оптимальном качении сферы по плоскости без прокручивания и проскальзывания, которая рассмотрена в разделе~\ref{subsec:roll}.

Выберем в  пространстве $\mathbb{R}^3$ такой неподвижный правый ортонормированный репер $(e_1, e_2, e_3)$, чтобы плоскость, по которой катится сфера, была натянута на $(e_1, e_2)$, а $e_3$ направлен в верхнее полупространство. Выберем также подвижный правый ортонормированный репер $(e'_1, e'_2, e'_3)$, закрепленный в центре сферы. Тогда ориентация сферы задается матрицей поворота $\SO (3) \ni R: (e'_1, e'_2, e'_3)\mapsto (e_1, e_2, e_3) $, а её положение --- координатами центра $(x, y) \in \mathbb{R}^2$ в базисе $(e_1, e_2)$. В качестве управляющих параметров возьмем компоненты вектора угловой скорости сферы в неподвижном репере $\vec{\Omega} = (u_2, -u_1, u_3) \in \mathbb{R}^3$. Тогда кинематика системы задается уравнениями
\begin{equation}\label{controlEquations}
\dot{x} = u_1,\qquad  \dot{y} = u_2,\qquad  \dot{R} = R \begin{pmatrix} 0 & -u_3 & -u_1 \\ u_3 & 0 & -u_2 \\ u_1 & u_2 & 0 \end{pmatrix},
\end{equation}
$$
(x,y)\in\mathbb{R}^2, \quad R \in \SO(3), \quad (u_1, u_2, u_3)\in \mathbb{R}^3.
$$
В качестве минимизируемого функционала   рассмотрим квадратичный функционал типа действия
\begin{equation}
\label{functional}
J = \frac{1}{2}\int_0^{t_1} (u_1^2 + u_2^2 + u_3^2) dt,
\end{equation} 
который с точностью до постоянного множителя представляет собой интеграл от вращательной энергии сферы. Требуется перекатить сферу из начального состояния $Q_0$ в конечное $Q_1$ так, чтобы достигался минимум функционала $J$. 

Эта задача формулируется естественным образом, как субриманова левоинвариантная задача на группе Ли $G = \mathbb{R}^2 \times \SO(3)$. 

Группу $G$ можно представить, как подгруппу группы $\GL(6)$ с помощью матриц
$$
Q =  \begin{pmatrix}
& & & \\ 
& R & & 0 \\
& & & \\
& 0 & & \begin{matrix}
1 & 0 & x \\
0 & 1 & y\\
0 & 0 & 1 
\end{matrix}
\end{pmatrix}.
$$

Введем левоинвариантный репер на $G$:
$$
e_1 = \frac{\partial}{\partial x}, \quad e_2 = \frac{\partial}{\partial y}, \quad V_i(R) = R\tilde{A}_i, \quad i=1,2,3;
$$
где $\tilde{A}_i$ --- базис алгебры Ли $\so(3)$:
$$
\tilde{A}_1 = \begin{pmatrix}
0 & 0 & 0 \\
0 & 0 & -1 \\
0 & 1 & 0 
\end{pmatrix}, \quad
\tilde{A}_2 = \begin{pmatrix}
0 & 0 & 1 \\
0 & 0 & 0 \\
-1 & 0 & 0 
\end{pmatrix}, \quad
\tilde{A}_3 = \begin{pmatrix}
0 & -1 & 0 \\
1 & 0 & 0 \\
0 & 0 & 0 
\end{pmatrix}.
$$
Тогда уравнения (\ref{controlEquations}) задают управляемую систему на группе $G$, и могут быть записаны в виде
\begin{equation}
\label{group_control}
\dot{Q} = u_1X_1(Q) + u_2X_2(Q) + u_3X_3(Q), \quad Q \in G, \quad (u_1, u_2 ,u_3) \in \mathbb{R}^3,
\end{equation}  
где
$$
X_1 = e_1 - V_2,\quad X_2 = e_2 + V_1, \quad X_3 = V_3. 
$$ 
Векторные поля $X_i$ задают распределение $\Delta = \spann\{X_1, X_2, X_3\} \subset T G$. Если $u(t)$ есть измеримое локально ограниченное отображение, то решение системы (\ref{group_control}) является допустимой кривой.

На распределении $\Delta$ можно задать скалярное произведение $\langle \cdot,\cdot \rangle$ следующим образом:
$$
\langle QA, QB \rangle_Q = \langle A,B \rangle_{\Id} = -\frac{1}{2}\tr(AB), \quad Q\in G,\quad A,B \in \Delta,
$$
где $\tr A$ --- след матрицы $A$. При этом длина допустимой кривой $Q(t), t\in[0,t_1]$ выражается стандартным образом:
\begin{equation}
\label{group_functional}
l(Q) = \int_0^{t_1} \sqrt{\langle \dot{Q}, \dot{Q} \rangle} dt = \int_0^{t_1} \sqrt{u_1^2 + u_2^2 + u_3^2} dt.
\end{equation}

Из неравенства Коши-Буняковского следует, что функционал длины и его минимизация (\ref{group_functional}) эквивалентны действию (\ref{functional}) и его минимизации.  

Поскольку распределение $\Delta$ и метрика $g$ являются левоинвариантными, то можно, не ограничивая общности, левыми сдвигами перевести $Q_0$ в единичный элемент, т.е. принять $Q_0 = (0,0,\Id)$.

Таким образом, получаем следующую левоинвариантную субриманову задачу оптимального управления на группе $G = \mathbb{R}^2 \times \SO(3)$:
\begin{align}
&\dot{Q} = u_1X_1(Q) + u_2X_2(Q) + u_3X_3(Q), \label{roll_besch1}\\
&Q = (x,y,R) \in G = \mathbb{R}^2 \times \SO(3), \quad (u_1,u_2,u_3)\in \mathbb{R}^3,\label{roll_besch2}\\
&Q(0)=Q_0=(0,0,\Id), \quad Q(t_1) = Q_1 = (x_1, y_1, R_1),\label{roll_besch3}\\
&l(Q) = \int_0^{t_1} \sqrt{u_1^2 + u_2^2 + u_3^2}\; dt \rightarrow \min.\label{roll_besch4}
\end{align}
Значение минимизирующего функционала длины (\ref{roll_besch4}) не зависит от параметризации кривой $Q(t)$, поэтому можно считать, что она имеет постоянную скорость, т.е. $u_1^2+u_2^2+u_3^2 \equiv \const$. Более того, из (\ref{roll_besch1}) и (\ref{functional}) видно, что, если управление $u(t)$, $t\in[0,t_1]$ переводит сферу из состояния $Q_0$ в $Q_1$ за время $t_1$, то управление $u'(t) = ku(kt)$, где $k$ --- некоторое положительное число, переводит $Q_0$ в $Q_1$ за время $t_1/k$. При этом $Q(t)$ переходит в $Q(kt)$. Это позволяет, не ограничивая общность, считать что $u_1^2+u_2^2+u_3^2 \equiv 1$.

Так как таблица умножения в алгебре Ли $L = \spann(e_1, e_2, V_1, V_2, V_3)$ имеет вид:
$$
\ad e_i = 0, \quad [V_1, V_2] = V_3, \quad [V_2, V_3] = V_1, \quad [V_3, V_1] = V_2,
$$
то для векторных полей $X_1, X_2, X_3$ имеем
$$
[X_1, X_2] = X_3, \quad [X_1, X_3] = - V_1, \quad  [X_2, X_3] = - V_2.
$$
Тогда видно, что $\spann(X_1, X_2, X_3, [X_1, X_3], [X_2, X_3]) = L$, и из теоремы Ра\-шевс\-ко\-го-Чжоу следует, что система является вполне управляемой. Существование оптимальных траекторий в задаче (\ref{roll_besch1})--(\ref{roll_besch4}) следует из теоремы Филиппова.

В дальнейших вычислениях будет использоваться изоморфизм между $\mathbb{R}^3$ и алгеброй Ли $\so(3)$. А именно, каждому вектору $\vec{A} \in \mathbb{R}^3$ можно поставить в соответствие матрицу $\tilde{A}\in \so(3)$ по следующему правилу:
$$
\vec{A} =
\begin{pmatrix}
a_1 \\ a_2 \\a_3
\end{pmatrix}, \quad
\tilde{A} = \begin{pmatrix}
0 & -a_3 & a_2 \\
a_3 & 0 & -a_1\\
-a_2 & a_1 & 0 
\end{pmatrix}.
$$

\subsubsection{Анормальные траектории}

\begin{proposition}\label{abnormal_theorem}
\begin{itemize}
\item[$(1)$]
 Все анормальные экстремальные траектории постоянной скорости имеют вид  
\begin{equation}
\label{abnormal_geodesics}
x = -\Omega_2 t, \quad y=\Omega_1 t, \quad R = e^{\tilde{\Omega}t},
\end{equation}
где $\tilde{\Omega}$ --- кососимметрическая матрица, соответствующая вектору угловой скорости $\vec{\Omega}$ с компонентой $\Omega_3 = 0$:
$$
\tilde{\Omega} = \begin{pmatrix}
0 & 0 & \Omega_2 \\
0 & 0 & -\Omega_1 \\
-\Omega_2 & \Omega_1 & 0
\end{pmatrix}.
$$
\item[$(2)$] 
Любая анормальная экстремальная траектория $Q(t)$, $t\in[0,t_1]$, оптимальна для любого $t_1>0$.
\end{itemize}
\end{proposition}
Следовательно, в анормальном случае вектор угловой скорости  является постоянным горизонтальным вектором, и сфера равномерно катится по прямой без прокручивания.

\subsubsection{Нормальные экстремали}
Гамильтонова система ПМП в нормальном случае имеет вид
\begin{eqnarray*}
&\dot{x} &= -\Omega_2, \\
&\dot{y} &= \Omega_1, \\
&\dot{R} &= R\tilde{\Omega}, \\
&\dot{\Omega}_1 &= \omega_2\Omega_3, \\
&\dot{\Omega}_2 &= -\omega_1\Omega_3, \\
&\dot{\Omega}_3 &= \omega_1\Omega_2 - \omega_2\Omega_1, \\
&\dot{\omega}_1 &= 0, \; \dot{\omega}_2 = 0. 
\end{eqnarray*}

\begin{theorem}
Если $\omega \neq 0$ и $\vec{\Omega} \neq \lambda \vec{\omega}$, то параметризованные длиной дуги нормальные экстремали описываются уравнениями
\begin{equation}\nonumber
\begin{pmatrix}
\Omega_1 \\
\Omega_2 \\
\Omega_3
\end{pmatrix}=
\begin{pmatrix}
\dfrac{\omega_1^2+\omega_2^2\cos\omega t}{\omega^2} & \dfrac{\omega_1\omega_2}{\omega^2}\left( 1 - \cos\omega t\right) & \dfrac{\omega_2}{\omega}\sin \omega t\\
\dfrac{\omega_1\omega_2}{\omega^2}\left( 1 - \cos\omega t\right)&
\dfrac{\omega_2^2+\omega_1^2\cos\omega t}{\omega^2} & -\dfrac{\omega_1}{\omega}\sin \omega t \\
-\dfrac{\omega_2}{\omega}\sin \omega t & \dfrac{\omega_1}{\omega}\sin \omega t & \cos\omega t
\end{pmatrix}
\begin{pmatrix}
\Omega_1^0 \\
\Omega_2^0 \\
\Omega_3^0
\end{pmatrix},
\end{equation}
\begin{equation}\nonumber
\begin{array}{l}
x = \dfrac{\Omega_3^0 \omega_1}{\omega^2}(1 - \cos\omega t) + \dfrac{(\Omega_1^0\omega_2 - \Omega_2^0\omega_1)\omega_1}{\omega^3}\sin\omega t - \dfrac{\omega_2(\Omega_2^0\omega_2 + \Omega_1^0\omega_1)}{\omega^2}t, \\
y = \dfrac{\Omega_3^0 \omega_2}{\omega^2}(1 - \cos\omega t) + \dfrac{(\Omega_1^0\omega_2 - \Omega_2^0\omega_1)\omega_2}{\omega^3}\sin\omega t + \dfrac{\omega_1(\Omega_2^0\omega_2 + \Omega_1^0\omega_1)}{\omega^2}t,
\end{array}
\end{equation}
\begin{equation}\nonumber
\label{general_geodesic_R}
R(t) = e^{t(\tilde{\omega} + \tilde{\Omega}_0)}e^{-t\tilde{\omega}}.
\end{equation}
В оставшихся случаях нормальные экстремальные траектории описываются уравнениями
\begin{equation}\nonumber
\label{line_geodesic}
x = -\Omega_2 t, \quad y=\Omega_1 t, \quad R = e^{\tilde{\Omega}t}.
\end{equation}
где $\tilde{\Omega}$ --- кососимметрическая матрица, соответствующая произвольному единичному вектору $\vec{\Omega}$.
\end{theorem}

\subsubsection{Диаметр субримановой метрики}
Рассмотрим субриманову метрику $d$ на группе $\SO(3)$, соответствующую задаче~\eq{roll_besch1}--\eq{roll_besch4}.
\begin{theorem}
Для метрики $d$ на $\SO(3)$ наиболее удаленными точками являются $\Id$ и $e^{\pi(a_1 A_1 + a_2 A_2)}$, $a_1^2 + a_2^2 = 1$, с расстоянием между ними $d(\Id,\,\, e^{\pi(a_1 A_1 + a_2 A_2)}) = \pi \sqrt 3$. Состояниям $e^{\pi(a_1 A_1 + a_2 A_2)}$, $a_1^2 + a_2^2 = 1$, соответствует сфера, перевернутая на противоположный полюс.
\end{theorem}

\subsubsection{Библиографические комментарии}
Результаты этого раздела получены в работе~\cite{roll_besch}.

\section[Задачи, интегрируемые в эллиптических функциях и интегралах]{Задачи, интегрируемые в эллиптических функциях \\и интегралах}\label{sec:elliptic}
\subsection{Эллиптические интегралы и функции}
Стандартные  источники по эллиптическим интегралам и функциям --- книги \cite{akhiezer, whit_watson, lawden}. Мы приведем ниже минимальные сведения о них, необходимые для изложения в последующих разделах.
\paragraph{Эллиптические интегралы в форме Якоби}

Эллиптические интегралы Лежандра первого рода:
$$
F(\f, k) = \int _0^{\f} \frac{dt}{\sqrt{1 - k^2 \sin^2 t}},
$$
второго рода:
$$
E(\f, k) = \int_0^{\f} \sqrt{1 - k^2 \sin^2 t} \, dt,
$$
третьего рода:
$$
\Pi(m; \f, k) = \int_0^{\f} \frac{dt}{(1 + m \sin^2 t)\sqrt{1 - k^2 \sin^2 t}},
$$
здесь и далее эллиптический модуль $k \in (0, 1)$.
Дополнительный модуль есть $k' = \sqrt{1-k^2}$.

Полные эллиптические интегралы:
\begin{align*}
&K(k) = F\left(\frac{\pi}{2}, k\right), \\
&E(k) = E\left(\frac{\pi}{2}, k\right).
\end{align*}

\paragraph{Эллиптические функции Якоби:}
\begin{align*}
&\f = \am (u,k) \quad \Leftrightarrow \quad  u = F(\f, k), \\
&\sn (u, k) = \sin \am (u, k), \\
&\cn (u, k) = \cos \am (u, k), \\
&\dn (u, k) = \sqrt{1 - k^2 \sn^2 (u, k)}, \\
&\E(u, k) = E(\am u, k).
\end{align*}
При записи эллиптических функций модуль $k$  часто опускается.

\paragraph{Стандартные формулы}
Производные и интегралы:
\begin{align*}
&\am' u = \dn u, \\
&\sn'u = \cn u \dn u, \\
&\cn' u = - \sn u \dn u, \\
&\dn'u = - k^2 \sn u \cn u,  \\
&\int_0^u \dn^2 t \, dt  = \E(u).
\end{align*}

 Вырождение:
\begin{align*}
k \to +0 \quad &\Rightarrow \quad \sn u  \to \sin u, \quad  \cn u \to \cos u,
\quad  \dn u \to 1,
\quad  \E(u) \to u,\\
k \to 1-0 \quad &\Rightarrow \quad \sn u  \to \tanh u, \quad  \cn u,  \ \dn u
\to \frac{1}{\cosh u},
\quad  \E(u) \to \tanh u.
\end{align*}

\subsection{Математический маятник}
Во всех субримановых задачах разделов \ref{subsec:martinet}--\ref{subsec:cartan} вертикальная подсистема гамильтоновой системы принципа максимума Понтрягина загадочным образом сводится к уравнению маятника, поэтому все они интегрируются в эллиптических функциях и интегралах.

\subsubsection{Уравнение маятника и его решение}\label{subsubsec:pend_eq}
Рассмотрим \ddef{математический маятник} --- материальную точку, закрепленную на невесомом  нерастяжимом стержне длины $L$, который может  свободно вращаться в вертикальной плоскости вокруг точки подвеса. Пусть 
$\t$ обозначает угол отклонения маятника от нижнего вертикального  положения. Тогда движение маятника удовлетворяет уравнениям 
\be{pend0}
\dot \t = c, \quad \dc = - r \sin \t,
\ee
где $r = \dfrac{g}{L} > 0$  и $g$ есть ускорение силы тяжести. Полная энергия маятника (первый интеграл уравнений \eq{pend0}) есть
$$
E = \dfrac{c^2}{2} - r\cos \t \in [-r, +\infty).
$$
Характер движения маятника определяется значением энергии $E$:
\begin{itemize}
\item
если $E = -r$, то $(\t, c) \equiv (0,0)$, и маятник покоится в устойчивом положении равновесия;
\item
если $E \in (-r, r)$, то маятник колеблется  вокруг устойчивого положения равновесия, он совершает периодические движения с периодом $T = \frac{4}{\sqrt r}K(k)$, $k = \sqrt{\frac{E + r}{2r}} \in (0, 1)$ по 
закону
$$
\sin\dfrac{\t}{2} = k \sn(\sqrt r \,t, k);
$$
\item
если $E = r$, $c = 0$, то маятник покоится в неустойчивом положении равновесия $(\t, c) \equiv (\pm\pi, 0)$;
\item
если $E = r$, $c \neq 0$, то маятник совершает непериодическое движение вдоль сепаратрисы, стремясь к неустойчивым положениям равновесия при $t \to \pm \infty$ по закону
$$
\sin \frac{\t}{2} = \th(\sqrt r \,t);
$$
\item
если $E > r$, то маятник неравномерно вращается по часовой $(c<0)$ или против часовой $(c > 0)$ стрелки, он совершает периодические движения с периодом $T = \frac{2}{\sqrt r} k K(k)$, 
$k = \sqrt{\frac{2r}{E+r}} \in (0, 1)$ по закону
$$
\sin\dfrac{\t}{2} =\pm \sn \left( \frac{\sqrt r}{k} \,t, k \right), \qquad \pm = \sgn c.
$$
\end{itemize}

Выше указан характер движений маятника \eq{pend0} при $r = \frac{g}{L} >0$. Если же $r=0$ (что можно истолковать как отсутствие силы тяжести), то:
\begin{itemize}
\item
при $c \neq 0$ маятник равномерно  вращается по часовой $(c< 0)$ или против часовой $(c > 0)$ стрелки;
\item
при $c = 0$ маятник покоится в неустойчивом положении равновесия.
\end{itemize}

Случай $r < 0$ (сила тяжести направлена вверх) сводится к случаю $r > 0$ заменой  переменных $(\t, c, r) \mapsto (\t + \pi, c, -r)$.

\subsubsection{Выпрямляющие координаты}\label{subsubsec:phik}
При $r > 0$
фазовый цилиндр маятника \eq{pend0}, 
$$
C = \{(\t, c) \mid \t\in S^1, \quad c \in \R  \}, \qquad S^1 = \R/{2\pi}\Z,
$$
стратифицируется в зависимости от типа движения маятника:
\begin{align*}
&C = \sqcup_{i = 1}^5 C_i,\\
&C_1 = \{(\t, c) \in C \mid E \in (-r, r)  \},\\
&C_2 = \{(\t, c) \in C \mid E >r  \},\\
&C_3 = \{(\t, c) \in C \mid E = r, \quad c \neq 0  \},\\
&C_4 = \{(\t, c) \in C \mid  c = 0, \quad \t = 0 \},\\
&C_5 = \{(\t, c) \in C \mid  c = 0, \quad \t = \pi \}.
\end{align*}
В областях $C_1, C_2, C_3$ можно ввести координаты $(\f, k)$, выпрямляющие  уравнение маятника.

Если $(\t, c) \in C_1$, то
\begin{align*}
&k = \sqrt{\frac{E + r}{2r}} \in (0, 1), \quad \sqrt r \f (\mod 4K(k)) \in [0, 4 K(k)],\\
& \sin \frac{\t}{2} = k \sn(\sqrt r \f, k), \quad \cos \frac{\t}{2} = \dn(\sqrt r \f, k), \\
&c = 2 k \sqrt r \cn(\sqrt r \f, k).
\end{align*}
Если $(\t, c) \in C_2$, то
\begin{align*}
&k = \sqrt{\frac{2r}{E + r}} \in (0, 1), \qquad \sqrt r \f (\mod 2 k K(k)) \in [0, 2 k K(k)],\\
& \sin \frac{\t}{2} = \pm\sn\left(\frac{\sqrt r \f}{k}, k\right), \quad \cos \frac{\t}{2} = \cn\left(\frac{\sqrt r \f}{k}, k\right), \\
&c = \pm 2 \frac{\sqrt r}{k}\dn\left(\frac{\sqrt r \f}{k}, k\right), \quad \pm = \sgn c.
\end{align*}
Если $(\t, c) \in C_3$, то
\begin{align*}
&k = 1, \qquad \f \in \R,\\
&\sin \frac{\t}{2} = \pm \th(\sqrt r \f), \quad \cos \frac{\t}{2} = \frac{1}{\ch(\sqrt r \f)},\\
&c = \pm \dfrac{2 \sqrt r}{\ch(\sqrt r \f)}, \qquad \quad \pm = \sgn c.
\end{align*}

В координатах $(\f, k)$ уравнение маятника \eq{pend0} выпрямляется:
$$
\dot\f = 1, \quad \dot k = 0,
$$
поэтому оно имеет решение
$$
\f_t = \f + t, \quad k \equiv \const.
$$
Эти выпрямляющие координаты и их модификации используются для параметризации экстремальных траекторий в разделах \ref{subsec:se2}--\ref{subsec:cartan}.

\subsubsection{Библиографические комментарии}
Раздел \ref{subsubsec:pend_eq} опирается на \cite{akhiezer}, а раздел \ref{subsubsec:phik} --- на \cite{dido_exp}.

\subsection{Плоская субриманова задача  Мартине}\label{subsec:martinet}
\subsubsection{Постановка задачи}
\ddef{Плоская субриманова структура  Мартине} задается  метрикой $ds^2 = dx^2 + dy^2$ на распределении Мартине $\D = \{dz - \frac 12 y^2 dx = 0  \}$ в пространстве $M = \R^3_{x, y, z}$. Ортонормированный репер может быть выбран в форме 
$$
X_1 = 
\frac{\partial}{\partial x } + \frac{y^2}{2}\frac{\partial}{\partial z }, \qquad X_2 = \frac{\partial}{\partial y }.
$$
Пусть $X_3 = \dfrac{\partial}{\partial z }$, тогда алгебра Ли, порожденная полями $X_1$, $X_2$, имеет таблицу умножения 
\begin{align*}
&[X_1, X_2] = -yX_3, \qquad [X_2, [X_1, X_2]] = -X_3,\\
&[X_1, [X_1, X_2]] = 0, \qquad\ad X_3 = 0,
\end{align*}
то есть это алгебра Энгеля (см.~раздел~\ref{subsec:engel}).

Плоская субриманова структура  Мартине {\em не левоинвариантна}, но мы включаем ее в данный обзор из-за ее особой роли в субримановой геометрии:
\begin{itemize}
\item
это простейшая субриманова структура с анормальными кратчайшими,
\item
это простейшая субриманова структура, в которой сфера не субаналитична,
\item
эта структура является нильпотентной аппроксимацией общих субримановых структур на распределении Мартине,
\item
это простейшая субриманова структура, интегрируемая в эллиптических функциях и интегралах.
\end{itemize}
Кроме того, плоская субриманова структура  Мартине есть фактор-структура левоинвариантной субримановой структуры на группе Энгеля (см. раздел \ref{subsec:engel}), поэтому гамильтонова система для экстремалей Мартине сводится к уравнению маятника, а сами эти экстремали проецируются на плоскость $(x, y)$ в эйлеровы эластики (см. раздел \ref{subsec:elastica}).

Задача оптимального управления  для плоской субримановой структуры  Мартине имеет вид
\begin{align*}
& \dq = u_1X_1 + u_2X_2, \quad q = (x, y, z) \in\R^3, \quad u = (u_1, u_2)\in\R^2,\\
&q(0) = q_0, \quad q(t_1) = q_1,\\
&J = \frac{1}{2}\int_0^{t_1}(u_1^2 + u_2^2) dt \to \min.
\end{align*}

\subsubsection{Принцип максимума Понтрягина}
\begin{proposition}
Анормальные траектории суть $\{y = 0, z = z_0  \}$. Они нестрого анормальны.
\end{proposition}
Нормальные экстремали суть траектории гамильтонова поля с гамильтонианом
$$
H = \frac 12(h_1^2 + h_2^2) = \frac 12[(p_x + \frac{y^2}{2}p_z)^2 + p_y^2],
$$
где $(p_x, p_y, p_z)$ --- канонические координаты ковектора $\lam \in T^*M$, и $h_i(\lam) = \lan \lam, X_i(q) \ran$, $i = 1, 2, 3$. Соответствующая гамильтонова система $\dot \lambda = \vH(\lambda)$ имеет вид 
\begin{align*}
&\dx = p_x + \frac{y^2}{2}p_z, \quad  &&\dot p_x = 0,\\
&\dy = p_y, \quad &&\dot p_y = -(p_x + \frac{y^2}{2}p_z)p_z y,\\
&\dz = (p_x + \frac{y^2}{2}p_z)\frac{y^2}{2},\quad  && \dot p_z = 0,
\end{align*}
или
\begin{align}\label{mart_ham2}
&\dx = h_1,  && \dot h_1 = yh_2h_3,\nonumber  \\
&\dy = h_2, && \dot h_2 = -yh_1h_3,    \\
&\dz  = \frac{y^2}{2}h_1, &&\dot h_3 = 0.\nonumber
\end{align}

Будем рассматривать экстремали на поверхности уровня $\{H = \frac 12  \}$, на которой  введем координаты
$$
h_1 = \cos \t, \quad h_2 = \sin \t, \quad h_3 = c.
$$

\subsubsection{Симметрии}
\paragraph{Отражения}
Субриманова структура $(\D, ds^2)$ сохраняется группой отражений
\begin{align*}
&\Sym = \{ \Id, \eps^1, \eps^2, \eps^3  \} \cong \Z_2 \times \Z_2,\\
&\eps^1: (x, y, z)\mapsto (x, -y, z), && (\t, c)\mapsto (\pi - \t, c),\\
&\eps^2: (x, y, z)\mapsto (- x, y, -z), && (\t, c)\mapsto (-\t, - c),\\
&\eps^3: (x, y, z)\mapsto (- x, -y, -z), && (\t, c)\mapsto (\t- \pi, - c).
\end{align*}
\paragraph{Дилатации}
Гамильтонова система \eq{mart_ham2} сохраняется однопараметрической группой дилатаций	
\begin{align*}
&(x, y, z) \mapsto (\d^{-1}x, \,\d^{-1}y, \,\d^{-3}z),\\
&(h_1, h_2, h_3)\mapsto (\d^{-1}h_1, \,\d^{-1}h_2, \,\d h_3).
\end{align*}

\subsubsection{Параметризация геодезических}
Далее предполагается, что $q_0 = 0$.
\begin{proposition}
Натурально параметризованные геодезические, выходящие из $q_0 = 0$, суть кривые
\begin{align*}
&x_t = -t + \frac{2}{\sqrt c}(\E(u) - E(k)),\\
&y_t = - \frac{2 k}{\sqrt c} \cn u,\\
&z_t = \frac{2}{3c^{3/2}}[(2k^2 - 1)(\E(u) -E(k)) + k'^2 t \sqrt c + 2 k^2 \sn u \cn u \dn u],
\end{align*}
где $u = K + t \sqrt c$, $k = \sin(\frac{\pi}{4} - \frac{\t}{2})$, $\t \in (-\frac{\pi}{2}, \frac{\pi}{2})$, $c > 0$, а также
$$
x_t = t\sin\t, \quad y_t = t\cos \t, \quad z_t = \frac{t^3}{6}\sin\t\cos^2\t,
$$
где $\t \in (-\frac{\pi}{2}, \frac{\pi}{2}]$, и кривые, получающиеся из указанных с помощью симметрий $\eps^1$,  $\eps^2$.
\end{proposition}

Обозначим экспоненциальное отображение 
\begin{align*}
&\Exp: C \times \R_+ \to M, \quad (\lam, t)\mapsto q_t = \pi \circ e^{t\vec H}(\lam),\\
&C = T_{q_0}^* M\cap\left\{H = \dfrac 12 \right\}.
\end{align*}

\subsubsection{Сопряженное время}
Если геодезическая проецируется на плоскость $(x, y)$ в прямую и строго нормальна, то она оптимальна, потому свободна от сопряженных точек.  В анормальном случае  геодезическая оптимальна и состоит из сопряженных точек.

Пусть $\lam = (\t, c) \in C$, и пусть геодезическая $q_t = \Exp(\lam, t)$ проецируется на плоскость $(x, y)$ не в прямую. Благодаря симметриям $\eps^1$ и $\eps^2$  можно считать, что $c > 0$ и $\t \in (-\frac{\pi}{2}, \frac{\pi}{2})$. 
Тогда первое сопряженное время есть
\begin{align*}
&\tconj(\lam) = \min\{t > 0 \mid v^2 c_1(v) + vc_2(v) + c_3(v) = 0 \},\\
&c_1(v) = k'^2\frac{\cn v}{\dn v},\\
&c_2(v) = k'^2 \sn v - 2k'^2 \E(v)\frac{\cn v}{\dn v},\\
&c_3(v) =\E^2(v)\frac{\cn v}{\dn v} - \E(v)\sn v,\\
&v = t \sqrt c.
\end{align*}

\begin{theorem}
Пусть $q_t = \Exp(\lam, t)$, $\lam \in C$, $t > 0$, есть геодезическая, которая проецируется  на плоскость $(x, y)$ не в прямую. Тогда
$$
\tconj(\lam )\in \left(\frac{2 K}{\sqrt{|c|}}, \frac{3 K}{\sqrt{|c|}}    \right).
$$
\end{theorem}
Приближенные вычисления показывают, что отношение $\dfrac{\tconj\cdot\sqrt{|c|}}{3K}$ есть приближенно константа $0, 97$.

\subsubsection{Время разреза и множество разреза}
\begin{theorem}
Геодезические, проецирующиеся на плоскость $(x, y)$ в прямую, суть кратчайшие. Геодезическая $q_t = \Exp(\lam, t)$, $\lam \in C$, $t > 0$, проецирующаяся на плоскость $(x, y)$ не в прямую, имеет время разреза $t_{\cut}(\lam)  = \dfrac{2 K}{\sqrt{|c|}}$, соответствующее ее первому пересечению с плоскостью Мартине $\{y = 0 \}$.

Множество разреза есть
$$
\Cut = \{q \in M \mid y = 0,\quad z \neq 0\}.
$$
Это множество не пересекается с первой каустикой.
\end{theorem}

\subsubsection{Сфера и фронт}
Разные сферы с центром $q_0 = 0$ переводятся друг в друга дилатациями, поэтому достаточно рассмотреть  единичную сферу
$$
S = \{q \in M \mid d(q_0, q)  = 1 \}.
$$
Сфера $S$  изображена на
 Рис. \ref{fig:martinet_sphere} в координатах $(x,y,v)$, $v = z - xy^2/6$. .

\figout{
\onefiglabelsize{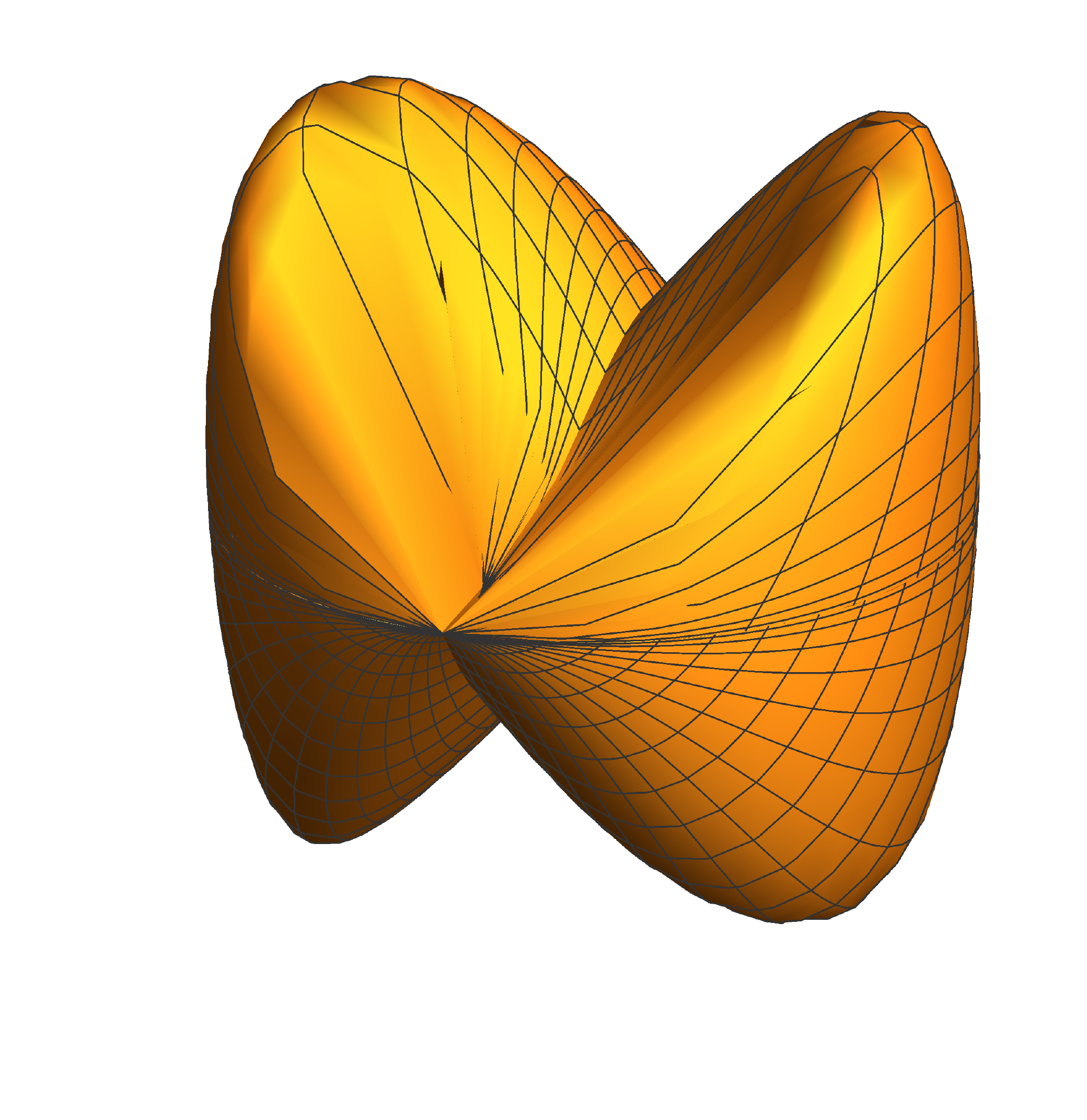}{Сфера в плоском случае Мартине}{fig:martinet_sphere}{0.4}
}

\begin{theorem}\label{th:martinet_sphere}
Пересечение сферы $S$ со множеством разреза (см. Рис. {\em\ref{fig:martinet_sphere_cut}}) есть кривая $k \mapsto \g(k)$, содержащаяся в плоскости Мартине $\{ y = 0 \}$  и заданная параметрическими уравнениями 
\begin{align}
&x(k) = -1 + 2\frac{E(k)}{K(k)}, \label{mart_x}  \\
&z(k) = \frac{1}{6K^3(k)}[(2 k^2 - 1)E(k) + k'^2K(k)]\label{mart_z},
\end{align}
где $k \in (0, 1)$, и кривая, полученная из $\g$ симметрией $\eps^2|_{\{y = 0\}}:  (x, z) \mapsto (-x, -z)$.

Если $k  \to +0$, то кривая $\g$ есть сужение на полуплоскость $\{ z > 0 \}$ графика аналитической функции $z = - \dfrac{2}{3 \pi^2} (x -1) + o(x-1)$, $x \to 1 - 0$.

Если $k  \to 1-0$, то кривая $\g$ есть график гладкой неаналитической функции
$$
z = \frac{X^3}{6} + F\left(X\right), \qquad X = \frac{x+1}{2},
$$
где $F$ есть плоская функция
$$
F(X) = - 4 X^3 e^{-\frac{2}{X}} + o\left(X^3e^{-\frac{2}{X}}\right), \quad  X\to +0.
$$
\end{theorem}

\figout{
\onefiglabelsize{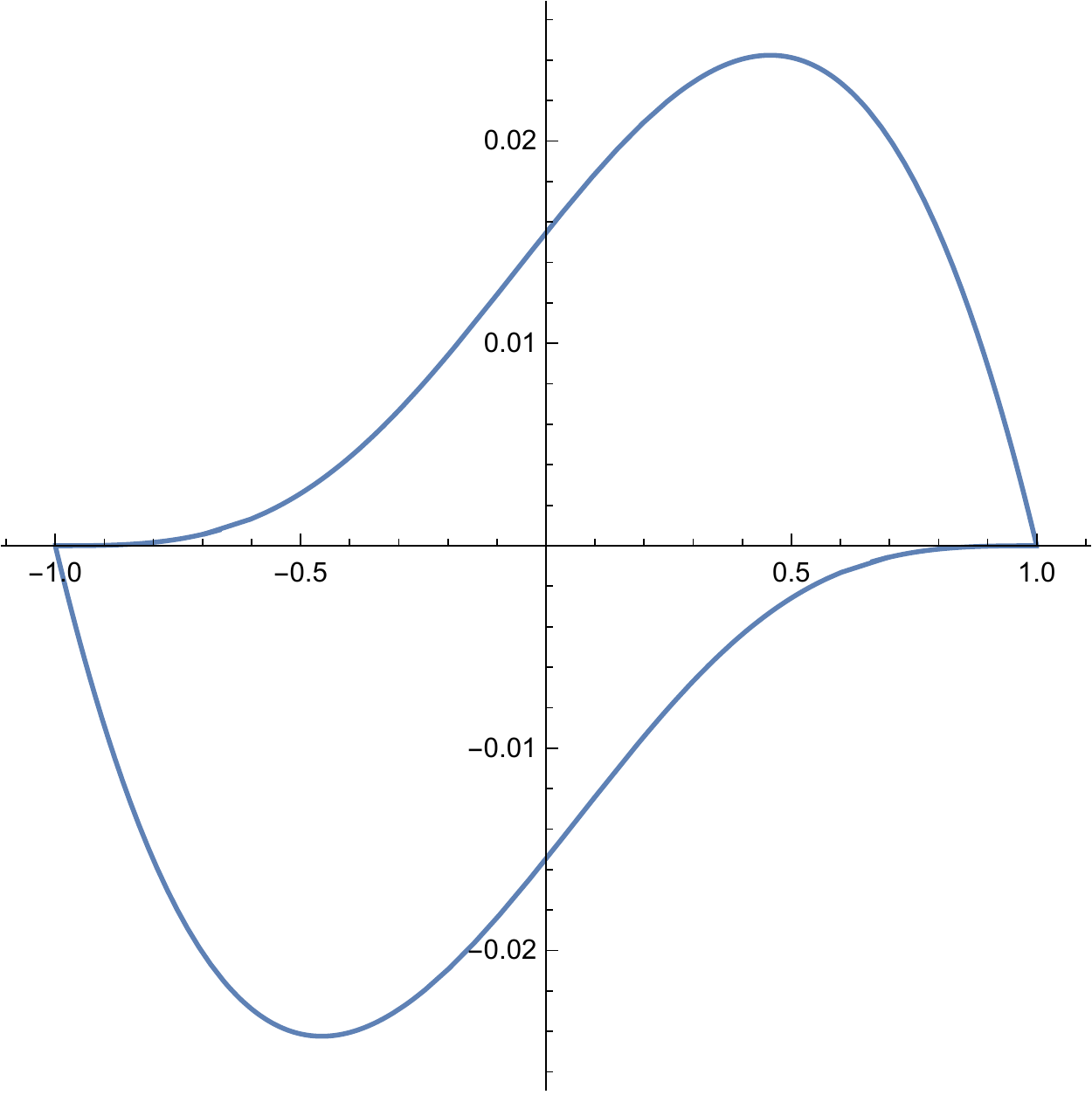}{Пересечение сферы с плоскостью Мартине $\{y=0\}$}{fig:martinet_sphere_cut}{0.4}
}

\begin{theorem}
Пересечение сферы $S$ с плоскостью Мартине не субаналитично, поэтому сфера $S$ не субаналитична.
\end{theorem}

Рассмотрим  волновой фронт из точки $q_0$ за единичное время:
$$
W = \{q \in M \mid q = \Exp(\lam, 1), \quad \lam \in C  \},
$$
остальные фронты из точки $q_0$ переводятся в этот фронт дилатациями.
\begin{theorem}
Пересечение волнового фронта $W$ с плоскостью Мартине $\{ y = 0 \}$ и полупространством $\{z > 0 \}$ есть объединение кривых $\g_n$, $n \in \N$, замыкание которых имеет две точки ветвления $x = \pm 1$,
$z = 0$. Кривая $\g_n$ задается параметрическими уравнениями
\begin{align*}
&x_n(k) = -1 + 2 \frac{E(k)}{K(k)},\\
&z_n(k) = \frac{1}{6n^2K^3(k)}[(2k^2 -1)E(k) + k'^2K(k)].
\end{align*}
Эта кривая вблизи точки $x = -1$, $z= 0$ есть график  функции
$$
z = \frac{1}{6n^2}X^3 +F \left(X\right),
$$
где $F(X) =  \a X^3 e^{-\frac{2}{X}} + o\left(X^3e^{-\frac{2}{X}}\right)$, $\a \neq 0$, а вблизи точки $x = 1$, $z = 0$ есть график функции
$$
z = - \frac {2}{3n^2\pi^2}(x-1)+ o(x-1).
$$
Внешняя кривая $\g_1$ есть пересечение $\g$ сферы с плоскостью Мартине $\{y = 0\}$ и полупространством $\{z > 0 \}$, см. теорему {\em \ref{th:martinet_sphere}}.
\end{theorem}

Пересечение сферы $S$ с плоскостью Мартине  и полупространством $\{ z > 0 \}$ есть параметрически заданная кривая $k \mapsto (x(k), z(k))$, $k \in (0, 1)$, см.~\eq{mart_x}, \eq{mart_z}.  Эта кривая  продолжается по непрерывности в полуплоскость $\{z \geq 0\}$ условием $k \in [0, 1]$. Полученная кривая полуаналитична при $k \neq 1$. Однако при $k = 1$ эта кривая не полуаналитична, поэтому не субаналитична.
\begin{theorem}
Пересечение сферы $S$ с плоскостью Мартине $\{y = 0 \}$ и полуплоскостью $\{ z \geq 0 \}$ вблизи точки $X=0$, где $X = \dfrac{x+1}{2}$, является  графиком функции вида
$$
z = F\left(X, \frac{e^{-\frac{ 1}{ X}}}{X^2}\right),
$$
где $X \geq 0$, и $F$ есть аналитическое отображение из окрестности точки $(0, 0) \in \R^2$ в $\R$.

Поэтому пересечение сферы $S$ с плоскостью Мартине  принадлежит $\exp$-$\log$ категории{\em\cite{dries, 
lion}}.
\end{theorem}

\subsubsection{Библиографические комментарии}
Этот раздел опирается на работу \cite{martinet}.

\subsection{Субриманова задача на группе   $\SE(2)$ евклидовых движений плоскости }\label{subsec:se2}
\subsubsection{Постановка задачи}\label{subsubsec:SE2_state}
\paragraph{Механическая постановка}
Рассмотрим задачу об оптимальном  движении для кинематической модели мобильного робота на плоскости.  Состояние робота задается его положением на плоскости $(x, y) \in \R^2$ и углом ориентации $\t \in S^1 = \R^2/2\pi \Z$ относительно положительного направления оси абсцисс. Робот может двигаться  с произвольной линейной скоростью $u_1 \in \R$ и при этом поворачиваться с произвольной угловой скоростью $u_2\in\R$. Требуется перевести робот из начального состояния $g_0 = (x_0, y_0, \t_0)$ в конечное состояние $g_1 = (x_1, y_1, \t_1)$ вдоль кратчайшего пути  в пространстве состояний. Длина пути в пространстве состояний $\R^2_{x, y}\times S^1_{\t}$ измеряется интегралом
$\int_0^{t_1} (\dx^2 + \dy^2 + \a^2\dot\t^2)^{1/2}\, dt$, где  $\a > 0$ --- некоторое заданное число, определяющее компромисс между линейной и угловой скоростью.
\paragraph{Задача оптимального управления  и ее нормализация}
Описанная задача для мобильного робота формализуется как задача оптимального управления:
\begin{align*}
&\dx = u_1 \cos\t, \quad \dy = u_1 \sin \t, \quad \dot\t = u_2,\\
&g = (x, y, \t) \in  \R^2_{x,y}\times S^1_{\t}, \quad u = (u_1, u_2)\in\R^2,\\
&g(0) = g_0, \quad g(t_1) = g_1,\\
&l = \int_0^{t_1}\sqrt{u_1^2 + \a^2u_2^2}\, dt \to \min.
\end{align*}
Заменой масштаба в плоскости $(x, y)$: 
$$
(x, y, \t) \mapsto \left(\frac{x}{\a}, \frac{y}{\a}, \t \right), \quad (u_1, u_2)\mapsto \left(\frac{u_1}{\a}, u_2\right) 
$$
можно свести эту задачу к случаю $\a =1$.

Параллельными переносами и поворотами плоскости $(x, y)$ можно добиться равенства $g_0 = (0, 0, 0)$.

В итоге получаем задачу оптимального управления:
\begin{align}
&\dx = u_1 \cos \t, \quad \dy = u_1 \sin \t, \quad \dot \t = u_2, \label{se2p21} \\
& g = (x, y, \t) \in  \R^2_{x, y}\times S^1_{\t}, \quad u = (u_1, u_2)\in \R^2,\label{se2p22} \\
&g(0) = g_0 = (0, 0, 0), \quad g(t_1)= g_1 = (x_1, y_1, \t_1),\label{se2p223} \\
&l =  \int_0^{t_1}\sqrt{u_1^2 + u_2^2}\, dt \to \min.\label{se2p24} 
\end{align}
Это  субриманова задача, заданная ортонормированным репером
\begin{align}
X_1 = \cos\t\frac{\partial}{\partial x} + \sin \t\frac{\partial}{\partial y}, \quad X_2 = \frac{\partial}{\partial \t}. \label{X12se2}
\end{align}

\paragraph{Группа движений плоскости}
\ddef{Группа собственных евклидовых движений плоскости} $G = \SE(2)$ есть полупрямое произведение  группы параллельных переносов $\R^2$ и группы вращений $\SO(2)$: 
$$\SE(2) = \R^2 \ltimes \SO(2).$$
Эта группа имеет линейное представление
$$
\SE(2) = 
\left\{\left(
\begin{array}{ccc}
\cos \t & - \sin \t & x \\
\sin \t & \cos \t & y \\
0 & 0 & 1
\end{array} 
\right)
\mid
\t \in S^1 = \R / (2 \pi \Z), \ x, y \in \R
\right\}.
$$
Действие движения $g = (x, y, \t)$ на вектор $(a, b) \in \R^2$ вычисляется с помощью матричного произведения:
$$
\begin{pmatrix}
\cos \t & -\sin \t & x\\
\sin \t & \cos \t  &  y\\
0  & 0 & 1
\end{pmatrix}\cdot
\begin{pmatrix}
a\\b\\1
\end{pmatrix}=
\begin{pmatrix}
a \cos \t - b \sin\t+x\\
a\sin\t + b\cos \t + y\\
1
\end{pmatrix},
$$
то есть $$
g : (a, b) \mapsto (a\cos\t - b\sin\t+x, \quad a\sin\t + b \cos\t +y).
$$

Алгебра Ли группы Ли $\SE(2)$ есть
$$
\gg = \sea(2) = \spann(E_{21} - E_{12},\, E_{13},\, E_{23}),
$$
где $E_{ij}$ есть $3\times 3$ матрица с единственным ненулевым элементом --- единицей в строке $i$ и столбце $j$.
Базисные левоинвариантные векторные поля на группе $\SE(2)$ суть
\begin{align*}
&X_1 =g E_{13 } =\cos\t \frac{\partial}{\partial x} + \sin \t\frac{\partial}{\partial y},\\
&X_2 = g(E_{21} - E_{12}) = \frac{\partial}{\partial \t},\\
&X_3 = - g E_{23} = \sin\t\frac{\partial}{\partial x} - \cos \t \frac{\partial}{\partial y},
\end{align*}
с таблицей умножения
\begin{align}
&[X_1, X_2] = X_3, \quad [X_2, X_3] = X_1, \quad [X_1, X_3] = 0.\label{se2_tab}
\end{align}

Ортонормированный репер \eq{X12se2} для субримановой задачи \eq{se2p21}--\eq{se2p24} состоит из левоинвариантных полей, поэтому эта задача --- левоинвариантная субриманова задача на группе $G = \SE(2)$. 

Согласно теореме \ref{th:class3}, это единственная, с точностью до локальных изометрий, вполне неголономная субриманова задача на $\SE(2)$, ей соответствуют инварианты $\chi = \kappa = 1$.

 Существование оптимальных управлений в задаче \eq{se2p21}--\eq{se2p24} следует из теорем Рашевского-Чжоу и Филиппова: система имеет полный ранг, так как
$$
\gg = \spann(X_1, X_2, X_3), \quad X_3 = [X_1, X_2].
$$

\subsubsection{Принцип максимума Понтрягина}\label{subsubsec:se2_PMP}
Анормальные  траектории постоянны.

Нормальные экстремали суть траектории гамильтоновой системы $\dot \lam = \vec H(\lam)$, $\lam \in T^*G$, где $H = (h_1^2 + h_2^2)/2$, $h_i(\lam) = \lan \lam, X_i \ran$, $i = 1, 2, 3$.
В координатах эта система записывается как 
\begin{align}
&\dot h_1 = - h_2 h_3, \quad \dot h_2 = h_1 h_3, \quad \dot h_3 = h_1 h_2,     \label{dh123} \\
&\dot x = h_1 \cos \t, \quad \dot y = h_1 \sin \t, \quad \dot \t = h_2.  \nonumber
\end{align}
На поверхности уровня $\{H=1/2  \}$ в  координатах $(\g, c)$, где 
$$
h_1 = \sin \frac{\g}{2}, \quad h_2 = -\cos \frac{\g}{2}, \quad c=2h_3,
$$
вертикальная подсистема \eq{dh123} гамильтоновой системы принимает форму двулистного накрытия маятника:
\begin{align}\label{se2_pend}
&\dot\g = c, \quad \dc = -\sin\g, \quad (\g, c) \in C = \gg^*\cap\left\{H = \frac 12 \right \} \cong (2S^1_{\g})\times \R_c, \quad 2S^1 = \R/(4\pi\Z).
\end{align}
Первый интеграл этого уравнения --- энергия маятника
\be{E}
E = \frac{c^2}{2} - \cos \g \in [-1, + \infty).
\ee

\paragraph{Симплектическое слоение}
На коалгебре Ли $\gg^*$ имеется функция Казимира $F = h_1^2 + h_3^2$. Симплектическое слоение состоит из круговых цилиндров $\{ h_1^2 + h_3^2 = \const > 0\}$ и точек $\{h_1 = h_3 = 0, \quad h_2 = \const  \}$.

Энергия маятника есть линейная комбинация функции Казимира и гамильтониана:
$$
E = 2F - 2H.
$$
\paragraph{Стратификация цилиндра $C$ и выпрямляющие координаты}\label{par:sl2_stratif}
Цилиндр $C$ разбивается на инвариантные  множества маятника \eq{se2_pend} критическими линиями уровня энергии $E$:
\begin{align}
&C = \sqcup_{i=1}^5 C_i, \label{decompC} \\
&C_1 = \{ \lam \in C \mid E \in (-1, 1) \}, \nonumber \\
&C_2 = \{ \lam \in C \mid E \in (1, + \infty) \}, \nonumber \\
&C_3 = \{ \lam \in C \mid E =1, \ c \neq 0 \}, \nonumber \\
&C_4 = \{ \lam \in C \mid E = - 1 \} = \{ (\g, c) \in C \mid \g = 2 \pi n, \ c = 0 \}, \nonumber \\
&C_5 = \{ \lam \in C \mid E = 1, \ c = 0 \} = \{ (\g, c) \in C \mid \g = \pi + 2 \pi n, \ c = 0 \}, \qquad n \in \Z. \nonumber 
\end{align}

Для регулярного интегрирования уравнения маятника \eq{se2_pend}  на стратах $C_1, C_2, C_3$ вводятся  координаты $(\f, k)$, выпрямляющие это уравнение.

Если $\lam = (\g, c) \in C_1$, то
\begin{align*}
&k = \sqrt{\frac{E+1}{2}} = \sqrt{\sin^2 \frac{\g}{2} + \frac{c^2}{4}} \in (0,1),\\
&\sin \frac{\g}{2} = s_1 k \sn(\f,k), \qquad s_1 = \sgn \cos(\g/2),\\
&\cos \frac{\g}{2} = s_1 \dn(\f,k), \\
&\frac{c}{2} = k \cn(\f,k), \qquad \f \in [0, 4 K(k)].
\end{align*}

Если $\lam = (\g, c) \in C_2$, то
\begin{align*}
&k = \sqrt{\frac{2}{E+1}} = \frac{1}{\sqrt{\sin^2 \frac{\g}{2} + \frac{c^2}{4}}} \in (0,1),\\
&\sin \frac{\g}{2} = s_2  \sn(\f/k,k), \qquad s_2 = \sgn c, \\
&\cos \frac{\g}{2} = \cn(\f/k,k), \\
&\frac{c}{2} = (s_2/ k) \dn(\f/k,k), \qquad \f \in [0, 4 k K(k)].
\end{align*}

Если $\lam = (\g,c) \in C_3$, то
\begin{align*}
&k = 1,\\
&\sin \frac{\g}{2} = s_1 s_2 \tanh \f,   \qquad s_1 = \sgn \cos(\g/2), \quad s_2 = \sgn c, \\
&\cos \frac{\g}{2} = s_1 / \cosh \f, \\
&\frac{c}{2} = s_2/ \cosh \f, \qquad \f \in (-\infty, +\infty).
\end{align*}

В координатах $(\f, k)$ поток маятника \eq{se2_pend} выпрямляется:
$$
\dot \f = 1, \quad \dot k = 0, \qquad \lam = (\f, k) \in \cup_{i=1}^3 C_i.
$$

\paragraph{Параметризация геодезических}
Если $\lam = (\f, k) \in C_1$, то $\f_t = \f+t$ и
\begin{align*}
&\cos \t_t = \cn \f \cn \f_t + \sn \f \sn \f_t, \\ 
&\sin \t_t = s_1(\sn \f \cn \f_t - \cn \f \sn \f_t), \\
&\t_t = s_1(\am \f - \am  \f_t) \pmod {2 \pi}, \\
&x_t = (s_1/k) [ \cn \f (\dn \f - \dn \f_t) + \sn \f (t + \E(\f) - \E(\f_t))], \\  
&y_t = (1/k) [ \sn \f (\dn \f - \dn \f_t) - \cn \f (t + \E(\f) - \E(\f_t))].
\end{align*}

Если $\lam \in C_2$, то
\begin{align*}
&\cos \t_t = k^2 \sn \psi \sn \psi_t + \dn \psi \dn \psi_t, \\ 
&\sin \t_t = k(\sn \psi \dn \psi_t - \dn \psi \sn \psi_t), \\
&x_t = s_2 k [\dn \psi(\cn \p - \cn \p_t) + \sn \p (t/k + \E(\p) - \E(\p_t))], \\
&y_t = s_2  [k^2 \sn \psi (\cn \p - \cn \p_t) - \dn \p (t/k + \E(\p) - \E(\p_t))],
\end{align*}
где
$$
\psi = \f/k, \qquad \psi_t = \f_t/k = \psi + t/k.
$$

Если $\lam = (\f, k) \in C_3$, $k=1$, то $\f_t = \f+t$ и
\begin{align*}
&\cos \t_t = 1/ (\cosh \f \cosh \f_t)  + \tanh \f  \tanh \f_t, \\
&\sin \t_t = s_1 (\tanh \f /\cosh \f_t - \tanh \f_t /\cosh \f), \\
&x_t = s_1 s_2 [(1/\cosh \f)(1/\cosh \f - 1/\cosh \f_t) + \tanh \f(t + \tanh \f - \tanh \f_t)],\\
&y_t = s_2 [\tanh \f (1/\cosh \f - 1/\cosh \f_t) -(1/\cosh \f) (t + \tanh \f - \tanh \f_t)].
\end{align*}

Если $\lam \in C_4$, то
$$
\t_t = -s_1 t, \qquad x_t = 0, \qquad y_t = 0.
$$

Если $\lam \in C_5$, то
$$
\t_t = 0, \qquad
x_t = t \, \sgn \sin (\g/2), \qquad  y_t = 0.
$$

Проекции геодезических на плоскость $(x,y)$  в случаях $C_1$, $C_2$, $C_3$  изображены соответственно на Рис.~\ref{fig:xyC1}, \ref{fig:xyC2},  \ref{fig:xyC3}. 

\figout{
\begin{figure}[htbp]
\includegraphics[width=0.32\textwidth]{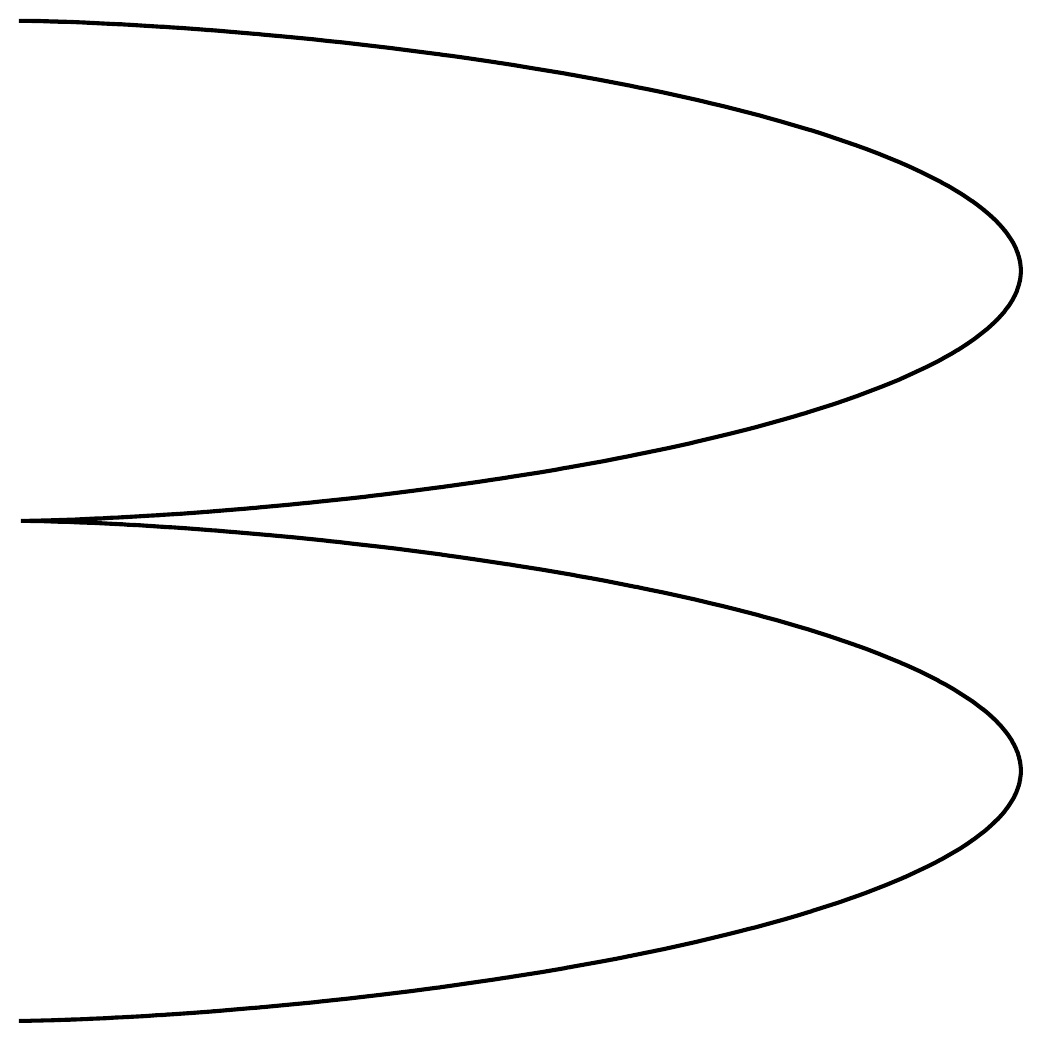}
\hfill
\includegraphics[width=0.32\textwidth]{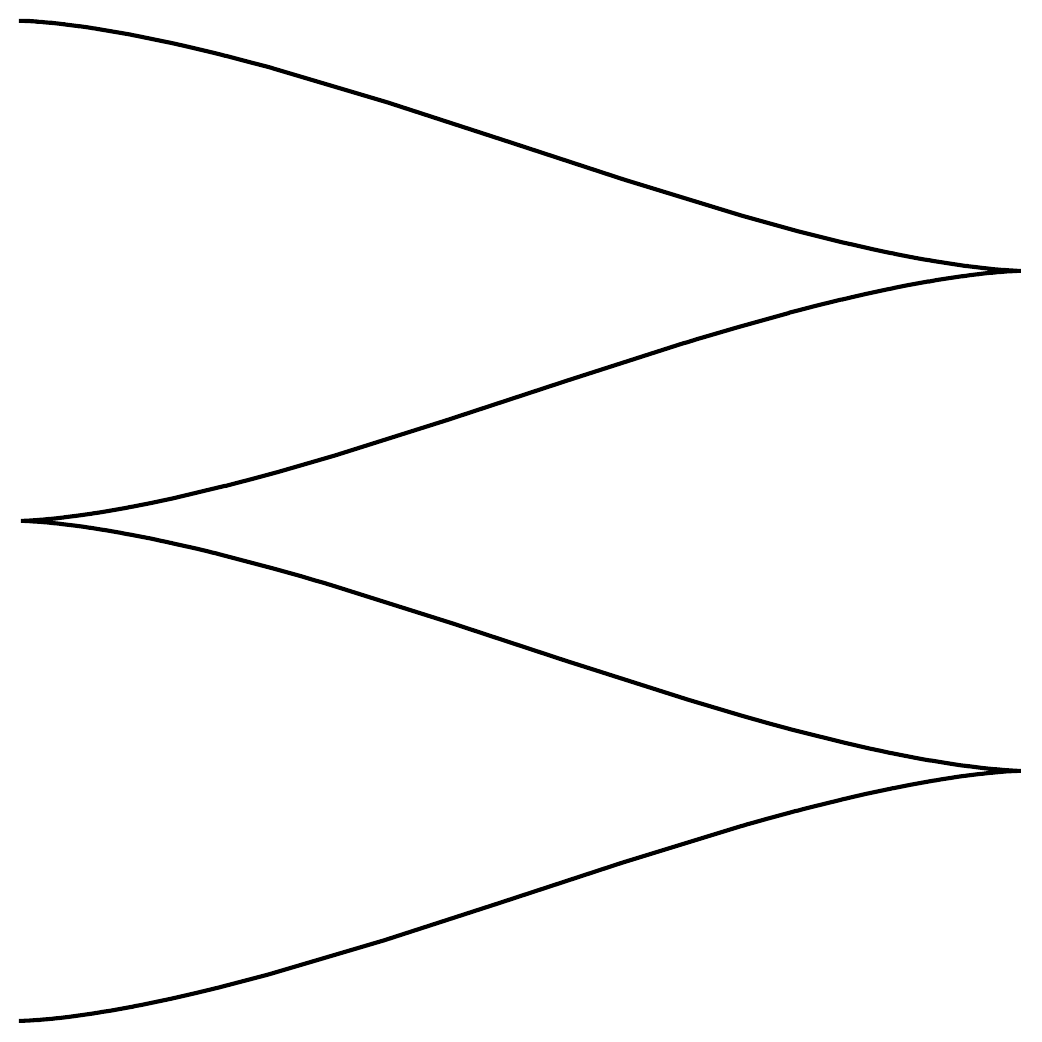}
\hfill
\includegraphics[width=0.32\textwidth]{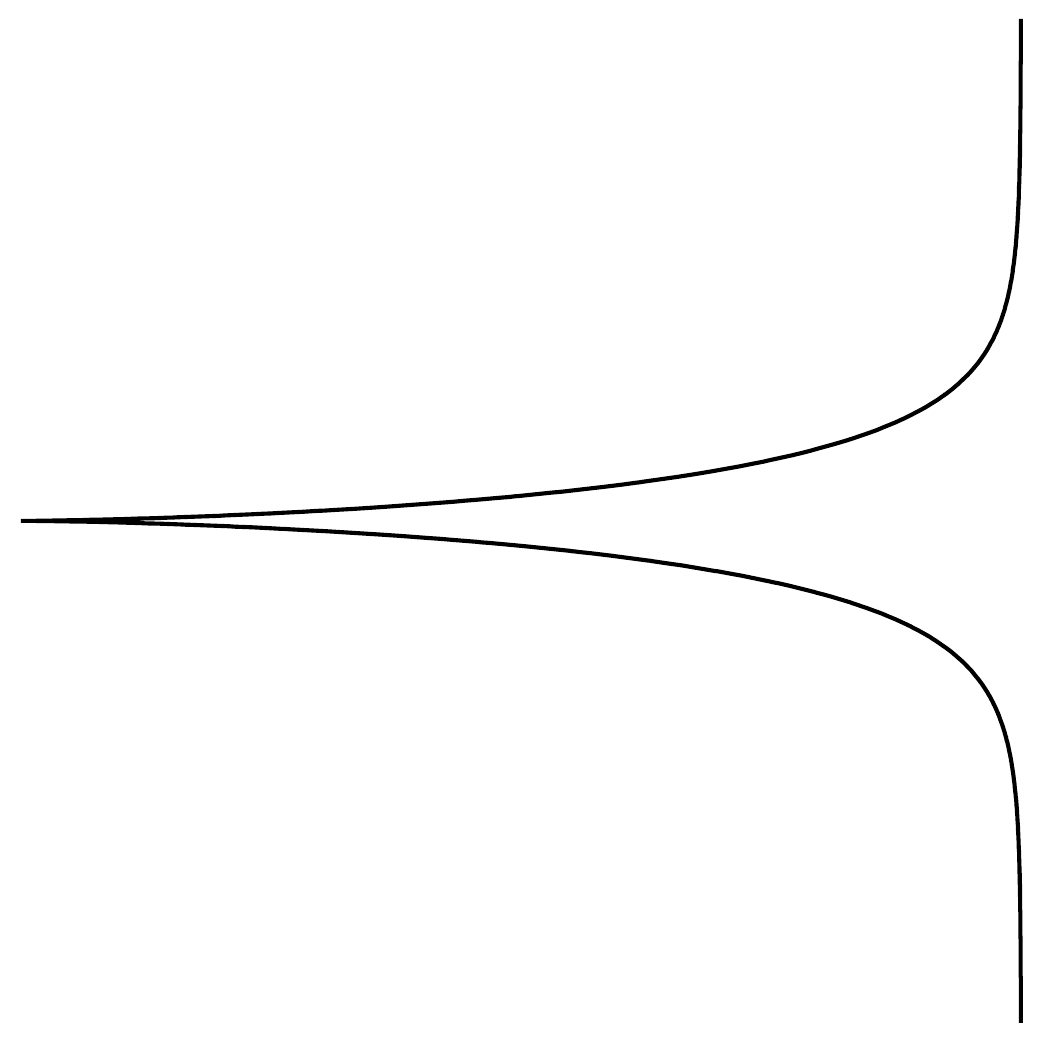}
\\
\parbox[t]{0.32\textwidth}
{\caption{Неинфлексионная кривая $(x_t, y_t)$: $\lam \in C_1$}\label{fig:xyC1}}
\hfill
\parbox[t]{0.32\textwidth}
{\caption{Инфлексионная кривая $(x_t, y_t)$: $\lam \in C_2$}\label{fig:xyC2}}
\hfill
\parbox[t]{0.32\textwidth}
{\caption{Трактриса $(x_t, y_t)$: $\lam \in C_3$}\label{fig:xyC3}}
\end{figure}
}

\subsubsection{Симметрии и страты Максвелла}\label{subsubsec:se2_max}
Фазовый портрет маятника \eq{se2_pend} сохраняется  группой симметрий $\Sym$,  порожденной отражениями цилиндра $C$ в осях координат $\g$, $c$, в начале координат $(\g, c) = (0, 0)$,  и поворотом на угол 
$2\pi$:
$$
\Sym = \{ \Id, \eps^1, \dots \eps^7  \}\cong \Z_2\times\Z_2\times\Z_2,
$$
где
\begin{align*}
&\map{\eps^1}{(\g,c)}{(\g,-c)},\\
&\map{\eps^2}{(\g,c)}{(-\g,c)},\\
&\map{\eps^3}{(\g,c)}{(-\g,-c)},\\
&\map{\eps^4}{(\g,c)}{(\g+2 \pi,c)},\\
&\map{\eps^5}{(\g,c)}{(\g + 2 \pi, -c)},\\
&\map{\eps^6}{(\g,c)}{(-\g + 2 \pi,c)},\\
&\map{\eps^7}{(\g,c)}{(-\g + 2 \pi,-c)}.
\end{align*}

Эти симметрии естественно продолжаются на прообраз и образ экспоненциального отображения.

Если $\nu = (\lam, t) = (\g, c, t) \in N = C \times \R_+ $, то $\eps^i(\nu) = \nu^i = (\lam^i, t) = (\g^i, c^i, t) \in N$, 
\begin{align*}
&(\g^1,c^1) = (\g_t, -c_t), \\
&(\g^2,c^2) = (-\g_t, c_t), \\
&(\g^3,c^3) = (-\g, -c), \\
&(\g^4,c^4) = (\g + 2 \pi, c), \\
&(\g^5,c^5) = (\g_t + 2 \pi, -c_t), \\
&(\g^6,c^6) = (-\g_t+ 2 \pi, c_t), \\
&(\g^7,c^7) = (-\g, -c).
\end{align*}

Если $g = (x, y, \t) \in G$, то $g^i = \eps^i(g) = (x^i, y^i, \t^i) \in G$, где
\begin{align*}
&(x^1,y^1,\t^1) = (x \cos \t + y \sin \t, x \sin \t - y \cos \t, \t), \\
&(x^2,y^2,\t^2) = (-x \cos \t - y \sin \t, -x \sin \t + y \cos \t, \t), \\
&(x^3,y^3,\t^3) = (-x, -y, \t), \\
&(x^4,y^4,\t^4) = (-x, y, -\t), \\
&(x^5,y^5,\t^5) = (-x \cos \t - y \sin \t, x \sin \t - y \cos \t, -\t), \\
&(x^6,y^6,\t^6) = (x \cos \t + y \sin \t, -x \sin \t + y \cos \t, -\t), \\
&(x^7,y^7,\t^7) = (x, -y, -\t).
\end{align*}

\begin{proposition}\label{propos:se2_sym}
Группа $\Sym = \{\Id, \eps^1, \dots, \eps^7 \}$ есть подгруппа группы симметрий экспоненциального отображения.
\end{proposition}

\begin{theorem}\label{th:se2_max}
Первое время Максвелла, соответствующее группе симметрий	$\Sym$, для почти всех геодезических выражается следующим образом:
\begin{align*}
&\lam \in C_1 \then \tmax(\lam) = 2 K(k),  \\
&\lam \in C_2 \then \tmax(\lam) = 2 k p_1^1(k), \\
&\lam \in C_3 \then \tmax(\lam) = +\infty, \\
&\lam \in C_4 \then \tmax(\lam) = \pi, \\
&\lam \in C_5 \then \tmax(\lam) = +\infty, 
\end{align*}
где $p = p_1^1(k) \in (K(k), 2K(k))$ есть первый положительный корень функции 
$$
f_1(p, k) = \cn p (\E(p) - p) - \dn p \sn p.
$$
\end{theorem}

\begin{remark}
Для тех геодезических, для которых первое время Максвелла, соответствующее группе $\Sym$, не равно $\tmax$, оно больше этого значения, а $\tmax$  есть первое сопряженное время.
\end{remark}

\begin{theorem}\label{th:se2_max_inv}
Функция $\tmax: C \to (0, +\infty]$ имеет следующие свойства инвариантности:
\begin{itemize} 
\item[$(1)$]
$\tmax(\lam)$ зависит только от $E$,
\item[$(2)$]
$\tmax(\lam)$ есть первый интеграл поля $\vec H_v$,
\item[$(3)$]
$\tmax(\lam)$ инвариантно относительно отражений $\eps^i \in \Sym$: если $(\lam, t) \in C\times\R_+$, $(\lam^i, t) = \eps^i(\lam, t)$, то $\tmax(\lam^i) = \tmax(\lam)$.
\end{itemize} 
\end{theorem}

\subsubsection{Оценки сопряженного времени}\label{subsubsec:se2_conj}
\begin{theorem}
\begin{itemize} 
\item[$(1)$] 
Если $\lam \in C_1\cup C_3\cup C_4\cup C_5$, то $\tconj(\lam) = +\infty$.
\item[$(2)$] 
Если $\lam \in C_2$, то $\tconj(\lam)\in [2kp_1^1, 4kK] $.
\item[$(3)$]
Следовательно, $\tconj(\lam) \geq \tmax(\lam)$ для всех $\lam \in C$.
\end{itemize} 
\end{theorem}

\subsubsection{Диффеоморфная структура экспоненциального отображения} \label{subsubsec:se2_diff}
Рассмотрим подмножество в пространстве состояний, не содержащее неподвижных точек отражений~$\eps^i$:
\begin{align*}
&\widetilde G = \{g \in G \mid \eps^i(g) \neq g, \quad i = 1, \dots, 7\} = \{g \in G \mid R_1(g) R_2(g)\sin\t\neq 0  \},\\
&R_1 = y\cos \frac{\t}{2} - x \sin\frac{\t}{2},  \quad R_2 = x \cos \frac{\t}{2}  + y\sin\frac{\t}{2},
\end{align*}
и его разбиение на компоненты связности 
$$
\widetilde G = \sqcup_{i=1}^8 G_i,
$$
где каждое множество $G_i$ характеризуется постоянными знаками функций $\sin \t$, $R_1$, $R_2$, описанными в таблице \ref{tab:se2_Gi}.
\begin{table}[htbp]
$$
\begin{array}{|c|c|c|c|c|c|c|c|c|}
\hline 
G_i      & G_1      & G_2 & G_3 &G_4 & G_5 & G_6 & G_7 & G_8 \\
\hline 
\sgn(\sin \t) & - & - & - & - & + & + & + & + \\
\hline 
\sgn(R_1) & + & + & - & - & - & - & + & + \\
\hline 
\sgn(R_2) & + & - & - & + & + & - & - & + \\
\hline 
\end{array}
$$
\caption{Определение областей $G_i$}\label{tab:se2_Gi}
\end{table}

Также рассмотрим открытое плотное подмножество в пространстве всех потенциально оптимальных геодезических:
$$
\widetilde N = \{(\lam, t) \in N  \mid t < \tmax(\lam), \quad c_{t/2}\sin\g_{t/2} \neq 0  \},
$$
и его связные компоненты
\begin{align*}
&D_1 = \{ (\lam, t) \in N  \mid t < \tmax(\lam), \quad c_{t/2} > 0, \quad \g_{t/2}\in(-\pi, 0) \},\\
&D_2 = \{ (\lam, t) \in N  \mid t < \tmax(\lam), \quad c_{t/2} > 0, \quad \g_{t/2}\in(0, \pi) \},\\
&D_3 = \{ (\lam, t) \in N  \mid t < \tmax(\lam), \quad c_{t/2} < 0, \quad \g_{t/2}\in(0, \pi) \},\\
&D_4 = \{ (\lam, t) \in N  \mid t < \tmax(\lam), \quad c_{t/2} < 0, \quad \g_{t/2}\in(-\pi, 0) \},\\
&D_5 = \{ (\lam, t) \in N  \mid t < \tmax(\lam), \quad c_{t/2} > 0, \quad \g_{t/2}\in(\pi, 2\pi) \},\\
&D_6 = \{ (\lam, t) \in N  \mid t < \tmax(\lam), \quad c_{t/2} > 0, \quad \g_{t/2}\in(2\pi, 3\pi) \},\\
&D_7 = \{ (\lam, t) \in N  \mid t < \tmax(\lam), \quad c_{t/2} < 0, \quad \g_{t/2}\in(2\pi, 3\pi) \},\\
&D_8 = \{ (\lam, t) \in N  \mid t < \tmax(\lam), \quad c_{t/2} < 0, \quad \g_{t/2}\in(\pi, 2\pi) \},\\
&\widetilde N = \sqcup_{i=1}^8D_i.
\end{align*}

\begin{theorem}
Следующие отображения являются диффеоморфизмами:
\begin{align*}
&\Exp : D_i \to G_i, \quad i = 1, \dots, 8,\\
&\Exp : \widetilde N \to \widetilde G.
\end{align*}
\end{theorem}

\subsubsection{Время разреза}\label{subsubsec:se2_cut}

\begin{theorem}
Для любого $\lam \in C$
$$
t_{\cut}(\lam) = \tmax(\lam).
$$
\end{theorem}

Время разреза инвариантно относительно вертикальной компоненты гамильтонова поля $\vec H_v$, поэтому субриманова структура на группе $\SE(2)$ эквиоптимальна.

\subsubsection{Множество разреза и его стратификация}\label{subsubsec:se2_strat}
\begin{theorem}
Множество разреза есть $2$-мерное стратифицированное многообразие со
стратификацией
\begin{align*}
&\Cut = \Cut_{\glob} \sqcup \Cut_{\loc}^+ \sqcup \Cut_{\loc}^-,  \\
&\Cut_{\glob} = \{ q \in M \mid \t = \pi\}, \\
&\Cut_{\loc}^+ = \{ q \in M \mid \t \in (- \pi, \pi), \ R_2 = 0, \ R_1 \geq R_1^1(|\t|)\}, \\
&\Cut_{\loc}^- = \{ q \in M \mid \t\in (- \pi, \pi), \ R_2 = 0, \ R_1 \leq -  R_1^1(|\t|)\}, 
\end{align*}
где
\begin{align*}
&R_1 = R_1^1(\t), \qquad \t \in [0,\pi], \\
&R_1^1(\t) = 2 (F(v_1^1(k),k) - E(v_1^1(k),k)), \qquad k = k_1^1(\t),\\
& v_1^1(k) = \am(p_1^1(k),k), \qquad k \in [0,1),
\end{align*}
а функция $k = k_1^1(\t)$, $\t \in [0, \pi]$, есть обратная функция к убывающей функции
$$
\t(k) = 2\arcsin(k\sin v_1^1(k)), \quad k \in [0, 1].
$$
Начальная точка $g_0 = \Id$ содержится в замыкании каждой компоненты $\Cut_{\loc}^+$, $\Cut_{\loc}^-$,  и отделена от компоненты $\Cut_{\glob}$.
\end{theorem}

Множество разреза $\Cut \subset \SE(2)$  изображено на Рис. \ref{fig:se2_cut1} (в выпрямляющих координатах 
$R_1 = y \cos \frac{\theta}{2} - x \sin \frac{\theta}{2}$, 
$R_2 = x \cos \frac{\theta}{2} + y \sin \frac{\theta}{2}$)  и на Рис. \ref{fig:se2_cut2} (при вложении в полноторий --- модель группы $\SE(2)$).

\figout{
\begin{figure}[htbp]
\begin{center}
\includegraphics[height=6cm]{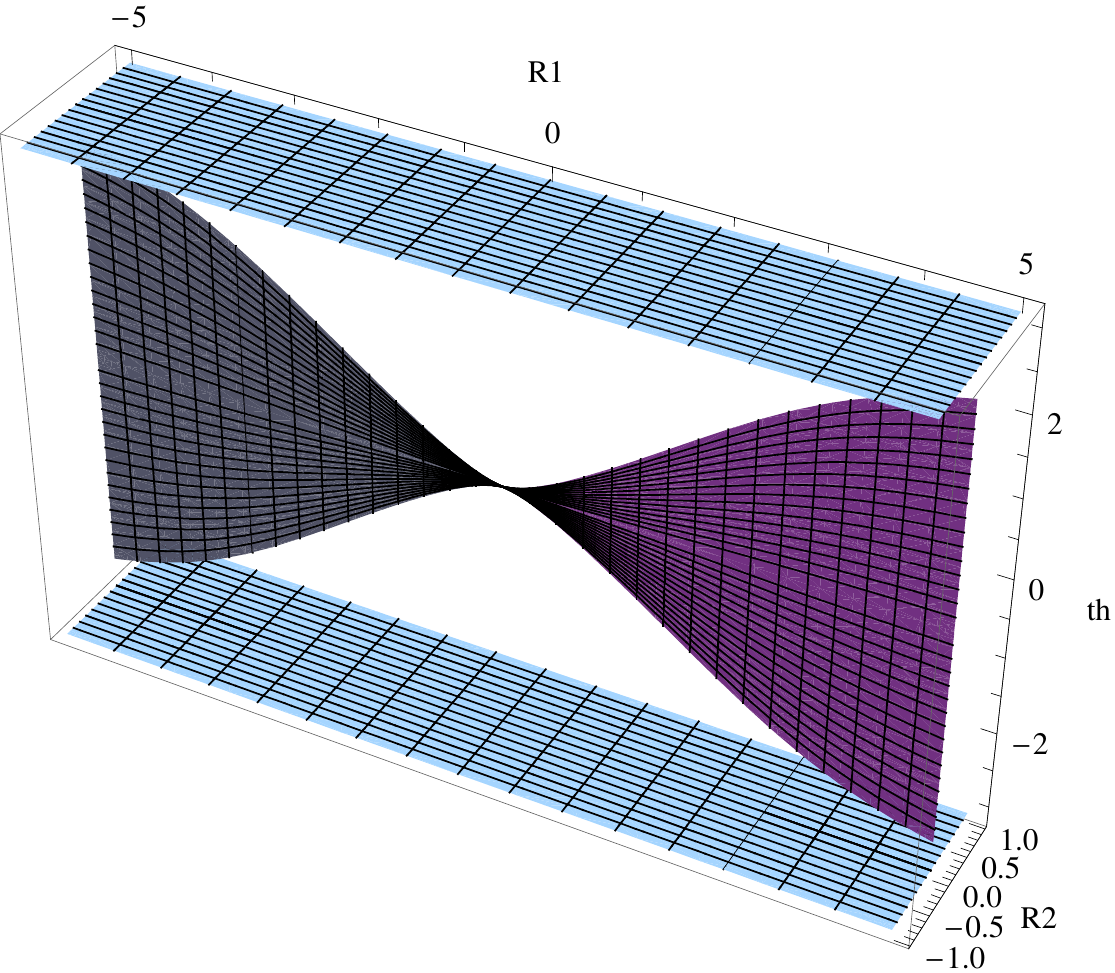}
\\
{\caption{Множество разреза в выпрямляющих координатах $(R_1, R_2, \theta)$}\label{fig:se2_cut1}}
\end{center}
\end{figure}


\begin{figure}[htbp]
\begin{center}
\includegraphics[height=4cm]{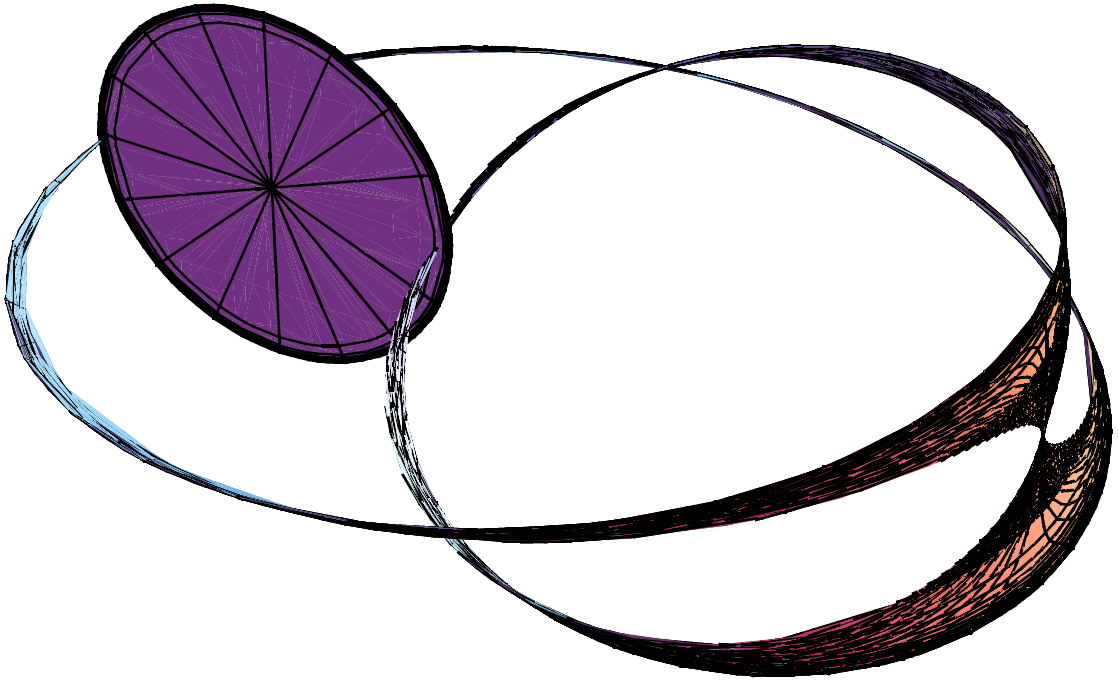}
\\
{\caption{Множество разреза в $\SE(2)$}\label{fig:se2_cut2}}
\end{center}
\end{figure}
}

\subsubsection{Сферы}\label{subsubsec:se2_sphere}
Субримановы сферы $S_R$ гомеоморфны (но не диффеоморфны):
\begin{itemize}
\item
евклидовой сфере $S^2$ при $R\in (0, \pi)$, 
\item
сфере с отождествленными полюсами $S$ и $N$: $S^2/\{S\sim N\}$ при $R=\pi$,
\item
тору $\T^2$ при $R > \pi$,
\end{itemize}

На Рис. \ref{fig:sphereRpi2torus},  \ref{fig:sphereRpitorus},  \ref{fig:sphereR32pitorus} изображены субримановы сферы радиусов $\pi/2$, $\pi$, $3\pi/2$  соответственно, вложенные в полноторий --- модель группы $\SE(2)$. 

\figout{
\onefiglabelsizen
{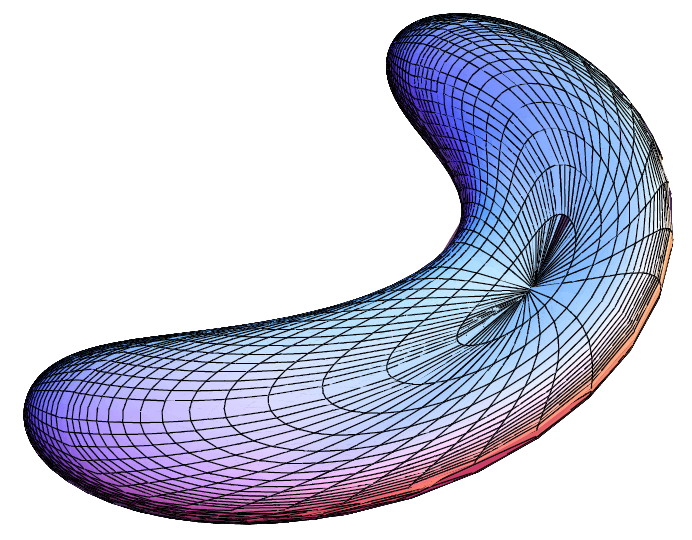}{Субриманова сфера $S_{\pi/2} \subset \SE(2)$}{fig:sphereRpi2torus}{6}

\onefiglabelsizen
{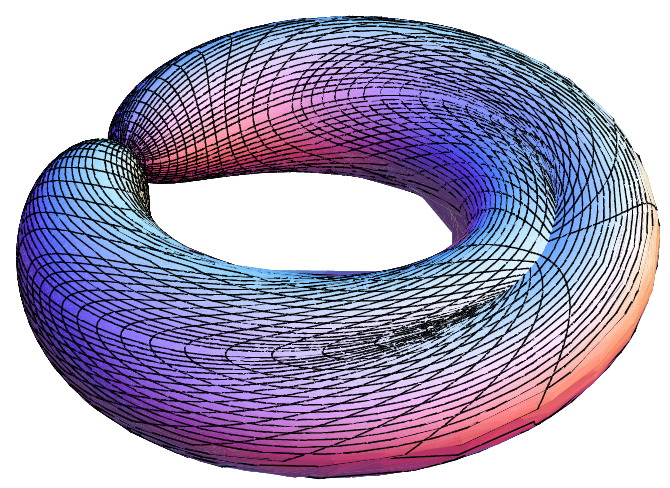}{Субриманова сфера $S_{\pi}\subset \SE(2)$}{fig:sphereRpitorus}{6} 

\onefiglabelsizen
{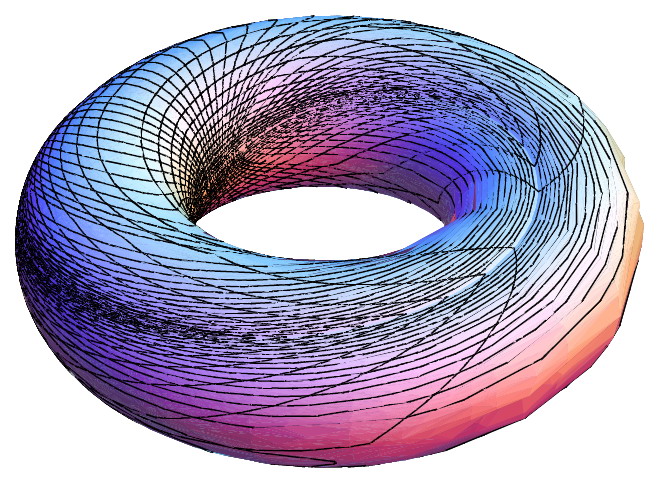}{Субриманова сфера $S_{3\pi/2}\subset \SE(2)$}{fig:sphereR32pitorus}{6}
}

\subsubsection{Метрические прямые}\label{subsubsec:se2_line}
Метрические прямые, проходящие через единичный элемент $g_0 = \Id$, суть $g(t) = \Exp(\lam, t)$, $t \in \R$, где $\lam \in C_3\cup C_5$.
Геодезические $\Exp(\lam, t)$, $\lam \in C_3$, проецируются на плоскость $(x, y)$ в трактрисы, а геодезические $\Exp(\lam, t)$, $\lam \in C_5$ --- в прямые $(x, y) =(\pm t, 0)$.

\subsubsection{Модель велосипеда}\label{subsubsec:se2_bicycle}
Субриманову задачу на группе $\SE(2)$ можно рассматривать как задачу об оптимальном движении модели велосипеда.

Пусть переднее и заднее колеса велосипеда касаются земли в точках $\mathbf {f}$ и $\mathbf {b}$ соответственно, а расстояние между этими точками (длина рамы велосипеда) постоянно и равно $\ell$. При движении велосипеда точки $\mathbf {f}$ и $\mathbf {b}$ пробегают две кривые ---  передний и задний пути. При этом отрезок $\mathbf {f} - \mathbf {b}$ в каждый момент времени касается заднего пути. Назовем движение велосипеда оптимальным, если оно минимизирует длину переднего пути.
Тогда задача об оптимальном  движении велосипеда есть в точности субриманова задача на группе
 $\SE(2)$ \eq{se2p21}--\eq{se2p24}.

Будем говорить, что две кривые на плоскости \ddef{имеют одинаковую форму}, если одну из них можно перевести в другую композицией движений и растяжений. \ddef{Ширина } плоской кривой есть нижняя грань расстояний  между двумя параллельными прямыми, ограничивающими полосу, содержащую эту кривую.

\begin{theorem}
Оптимальная траектория переднего колеса велосипеда $\mathbf {b}(t)$ есть либо прямая, либо дуга неинфлексионной эластики ширины не больше $2 \ell$. Таким образом возникает любая форма неинфлексионной эластики.
\end{theorem}

\begin{theorem}
Бесконечное движение велосипеда является оптимальным на каждом своем отрезке тогда и только тогда, когда оно имеет один из следующих двух типов:
\begin{itemize} 
\item[$(1)$] 
передний путь $\mathbf {b}(t)$ есть прямая, а задний путь $\mathbf {f}(t)$ есть трактриса или прямая,
\item[$(2)$] 
передний путь $\mathbf {b}(t)$ есть солитон Эйлера (критическая эластика) ширины $2\ell$, а задний путь $\mathbf {f}(t)$ есть трактриса.
\end{itemize}
\end{theorem}

\subsubsection{Группа изометрий и однородные геодезические}\label{subsubsec:se2_isom}
\begin{theorem}
Группа изометрий субримановой структуры на $\SE(2)$ есть
$\Isom(\SE(2)) = \SE(2) \rtimes (\Z_2 \times \Z_2)$,
где справа первый сомножитель $\SE(2)$ действует на себе левыми сдвигами, второй сомножитель $\Z_2$ действует на пару $(\mathbf {b}, \mathbf {f})$ как отражение плоскости в какой-нибудь оси, а третий сомножитель $\Z_2$ действует как отражение $(\mathbf {b}, \mathbf {f}) \mapsto (\mathbf {b}, 2 \mathbf {b} - \mathbf {f})$.
\end{theorem}

Геодезическая $\gamma$  на субримановом многообразии $M$  называется {\em однородной}, если она является однородным пространством некоторой однопараметрической подгруппы в группе изометрий $\Isom(M)$, т.е. существует однопараметрическая подгруппа $\{\f_s \mid s \in \R\} \subset \Isom(M)$  такая, что:
\begin{enumerate}
\item
$\forall \ s \in \R \quad \f_s(\gamma) \subset \gamma$, 
\item
$\forall \ g_1, \, g_2 \in \gamma \quad \exists \ s \in \R \ : \quad \f_s(g_1) = g_2$.
\end{enumerate}

Субриманово многообразие называется {\em геодезически орбитальным}, если все его геодезические однородны.

\begin{theorem}
Однородные геодезические на $\SE(2)$ есть $g(t) = \Exp(\lam, t)$, $\lam \in C_4\cup C_5$. Это однопараметрические подгруппы $e^{t X_2}$ и $e^{t X_1}$, они проецируются на плоскость $(x, y)$ соответственно в точку $(0,0)$ и прямую $y=0$.

Поэтому $\SE(2)$ не является геодезически орбитальным пространством.
\end{theorem}

\subsubsection{Библиографические комментарии}
Разделы \ref{subsubsec:SE2_state}--\ref{subsubsec:se2_max} опираются на \cite{max_sre}, разделы \ref{subsubsec:se2_conj}--\ref{subsubsec:se2_cut}, \ref{subsubsec:se2_sphere}, \ref{subsubsec:se2_line} --- на \cite{cut_sre1}, раздел \ref{subsubsec:se2_strat} --- на \cite{cut_sre2}, разделы \ref{subsubsec:se2_bicycle} и \ref{subsubsec:se2_isom} --- на \cite{bicycle, se2_GO}.

\subsection[Субриманова задача на группе   $\SH(2)$ движений псевдоевклидовой плоскости]{Субриманова задача на группе   $\SH(2)$ \\движений псевдоевклидовой плоскости}\label{subsec:sh2}
\subsubsection{Группа $\SH(2)$ движений псевдоевклидовой плоскости}\label{subsubsec:SH2}
\paragraph{Псевдоевклидова плоскость}
\ddef{Псевдоевклидовой плоскостью} называется двумерное вещественное линейное пространство, в котором задана знакопеременная билинейная форма
$$
(\x, \y) = x_1 y_1 - x_2 y_2, \qquad \x = (x_1, x_2), \quad \y = (y_1, y_2).
$$ 
\ddef{Расстояние} $r$  между точками $\x = (x_1, x_2)$ и  $\y = (y_1, y_2)$  определяется формулами
\begin{align*}
&r^2 = (\x-\y, \x-\y) = (x_1 - y_1)^2 - (x_2 - y_2)^2,\\
&r = \begin{cases}
|r| &\text{ при } r^2 \geq 0, \\ 
i |r| &\text{ при } r^2 < 0.
\end{cases}
\end{align*}
Множество точек $\x = (x_1, x_2)$, находящихся на нулевом расстоянии от начала координат $(x_1^2 - x_2^2 = 0$),  называется \ddef{световым конусом}. Дополнение псевдоевклидовой плоскости до светового конуса распадается на 4 связные компоненты ---  \ddef{квадранты} ($\sgn(x_1 - x_2) = \pm 1$, $\sgn(x_1 + x_2 ) = \pm 1$).

\paragraph{Группа Ли $\SH(2)$ и алгебра Ли $\sha(2)$}
\ddef{Движением} псевдоевклидовой плоскости называется ее линейное преобразование, сохраняющее ориентацию, квадранты, и расстояние между точками этой плоскости. \ddef{Группа движений псевдоевклидовой плоскости}  обозначается $\SH(2)$.  Эта группа имеет линейное представление
$$
\SH(2) = \left\{\left(
\begin{array}{ccc}
\ch z & \sh z & x \\
\sh z & \ch z & y \\
0 & 0 & 1 
\end{array}
\right)\mid x, y, z, \in \R \right\}.
$$
Действие движения $g = (x,y,z)$  на точку $\ab = (a_1, a_2)$ псевдоевклидовой плоскости вычисляется с помощью матричного произведения:
$$
\left(
\begin{array}{ccc}
\ch z & \sh z & x \\
\sh z & \ch z & y \\
0 & 0 & 1 
\end{array}
\right)
\left(\begin{array}{c} a_1 \\ a_2 \\ 1\end{array}\right) = 
\left(
\begin{array}{ccc}
a_1 \ch z + a_2 \sh z + x \\
a_1 \sh z + a_2 \ch z + y \\
1 
\end{array}
\right),
$$
т.е. $g \ : \ (a_1, a_2) \mapsto (a_1 \ch z + a_2 \sh z + x, a_1 \sh z + a_2 \ch z + y)$.

$G  = \SH(2)$  есть группа Ли с алгеброй Ли $\gg = \sha(2) = \spann(E_{21} + E_{12}, E_{13}, E_{23})$. Базисные левоинвариантные векторные поля на группе $\SH(2)$  суть
\begin{align*} 
&X_1 = L_{g*} E_{13} = \ch z \pder{}{x} + \sh z \pder{}{y},\\
&X_2 = L_{g*} (E_{21} + E_{12}) = \pder{}{z}, \\
&X_3 = L_{g*} E_{23} = \sh z \pder{}{x} + \ch z \pder{}{y},
\end{align*}
 с таблицей умножения
\be{sh2tab}
[X_1, X_2] = - X_3, \quad [X_2. X_3] = X_1, \quad [X_1, X_3] = 0.
\ee

\subsubsection{Субриманова задача на $\SH(2)$}\label{subsubsec:sh2_state}
Рассмотрим субриманову задачу на группе $\SH(2)$  с ортонормированным репером $(X_1, X_2)$:
\begin{align}
&\dot g = u_1 X_1 + u_2 X_2, \qquad g \in G = \SH(2), \quad u = (u_1, u_2) \in \R^2, \label{sh2pr1}\\
&g(0) = g_0 = \Id, \quad g(t_1) = g_1, \label{shepr2}\\
&l = \int_0^{t_1} \sqrt{u_1^2 +u_2^2} dt \to \min. \label{sh2pr3}
\end{align}

Согласно теореме Аграчева-Барилари (см. раздел \ref{subsec:class3}), это единственная, с точностью до локальных изометрий, неинтегрируемая субриманова задача ранга 2 на группе $\SH(2)$, ей соответствуют инварианты $\chi = - \kappa = 1$.

\subsubsection{Геодезические}\label{subsubsec:sh2_geod}
Существование оптимальных управлений в задаче \eq{sh2pr1}--\eq{sh2pr3}  следует из теорем Рашевского-Чжоу и Филиппова.

\paragraph{Принцип максимума Понтрягина}
Анормальные траектории постоянны.

Нормальные экстремали суть проекции траекторий гамильтоновой системы $\dot \lam = \vec{H}(\lam)$, $\lam \in T^* G$, где $H = (h_1^2 + h_2^2)/2$, $h_i(\lam) = \langle \lam, X_i\rangle$, $i = 1, 2, 3$.  В координатах эта система записывается как 
\begin{align}
&\dot{h}_{1} = h_{2}h_{0}, \label{sh2dh1}\\
&\dot{h}_{2} = -h_{1}h_{0},  \label{sh2dh2}\\
&\dot{h}_{0}=  h_{1}h_{2},  \label{sh2dh3}\\
&\dot{x} = h_{1}\cosh z, \nonumber \\
&\dot{y} = h_{1}\sinh z,  \nonumber\\
&\dot{z} = h_{2}. \nonumber
\end{align}
На поверхности уровня
$\{H=1/2\}$  в координатах $(\g, c)$,  где 
$$
h_1 = \cos \frac{\gamma}{2}, \quad h_2 = \sin \frac{\gamma}{2}, \quad c = - 2 h_3,
$$
вертикальная подсистема \eq{sh2dh1}--\eq{sh2dh3}  принимает форму двулистного накрытия маятника
\be{sh2_pend}
\dot\g = c, \quad \dot c = - \sin \g, \qquad (\g, c) \in \gg^* \cap \{H = 1/2\} \simeq (2 S^1_{\g}) \times \R_c.
\ee 
Первый интеграл этого уравнения --- энергия маятника
$$
E = \frac{c^2}{2} - \cos \g = 2 h_3^2 - h_1^2 + h_2^2 \in [-1, + \infty).
$$

\paragraph{Симплектическое слоение}
На коалгебре Ли $\gg^*$  имеется функция Казимира $F = h_1^2 - h_3^2$.  Симплектическое слоение состоит из:
\begin{itemize}
\item
гиперболических цилиндров (компонент связности поверхностей $\{h_1^2 - h_3^2 = \const \neq 0\}$),
\item
полуплоскостей (компонент связности поверхности $\{h_1^2 - h_3^2 = 0, \ h_1^2 + h_3^2 \neq 0\}$),
\item
точек $\{h_1 =  h_3 = 0, \ h_2 = \const\}$.
\end{itemize}
Энергия маятника есть линейная комбинация функции Казимира и гамильтониана: $E = 2 H - 2 F$.

\paragraph{Стратификация цилиндра $C$  и выпрямляющие координаты}
Так как вертикальная подсистема гамильтоновой системы для задачи на $\SH(2)$ --- маятник \eq{sh2_pend} ---  совпадает с таковой системой \eq{se2_pend} для задачи на $\SE(2)$, то стратификация цилиндра $C$  и выпрямляющие координаты $(\f, k)$ для задачи на $\SH(2)$  совпадают с таковыми для задачи на $\SE(2)$, см. п. \ref{subsubsec:se2_PMP}.

\paragraph{Параметризация геодезических}
Если $\lam = (\f, k) \in C_1$, то $\f_t = \f + t$  и 
$$
\left(\begin{array}{c}
x\\
y\\
z
\end{array}\right)=\left(\begin{array}{c}
\frac{s_{1}}{2}\left[\left(w+\frac{1}{w\left(1-k^{2}\right)}\right)\left[\mathrm{E}(\varphi)-\mathrm{E}(\varphi_{0})\right]+\left(\frac{k}{w(1-k^{2})}-kw\right)\left[\mathrm{sn}\,\varphi-\mathrm{sn}\,\varphi_{0}\right]\right]\\
\frac{1}{2}\left[\left(w-\frac{1}{w\left(1-k^{2}\right)}\right)\left[\mathrm{E}(\varphi)-\mathrm{E}(\varphi_{0})\right]-\left(\frac{k}{w\left(1-k^{2}\right)}+kw\right)\left[\mathrm{sn}\,\varphi-\mathrm{sn}\,\varphi_{0}\right]\right]\\
s_{1}\ln\left[(\mathrm{dn}\,\varphi-k\mathrm{cn}\,\varphi)w\right]
\end{array}\right)\label{eq:42}
$$
где $w=\frac{1}{\mathrm{dn}\varphi_{0}-k\mathrm{cn}\varphi_{0}}$.

Если $\lam = (\f, k) \in C_2$, то $\p = \frac{\f}{k}$, $\p_t = \frac{\f_t}{k} = \p + \frac tk$  и 
\begin{eqnarray*}
x  &=&  \frac{1}{2}\left(\frac{1}{w(1-k^{2})}-w\right)\left[\mathrm{E}(\psi)-\mathrm{E}(\psi_{0})-k^{\prime2}(\psi-\psi_{0})\right]\nonumber \\
 && \qquad +  \frac{1}{2}\left(kw+\frac{k}{w(1-k^{2})}\right)\left[\mathrm{sn}\,\psi-\mathrm{sn}\,\psi_{0}\right],\nonumber \\
y  &=&  -\frac{s_{2}}{2}\left(\frac{1}{w(1-k^{2})}+w\right)\left[\mathrm{E}(\psi)-\mathrm{E}(\psi_{0})-k^{\prime2}(\psi-\psi_{0})\right]\nonumber \\
 && \qquad +  \frac{s_{2}}{2}\left(kw-\frac{k}{w(1-k^{2})}\right)\left[\mathrm{sn}\,\psi-\mathrm{sn}\,\psi_{0}\right],\nonumber \\
z  &=&  s_{2}\ln[(\mathrm{dn}\,\psi-k\mathrm{cn}\,\psi)w],
\end{eqnarray*}
где  $w=\frac{1}{\mathrm{dn}\,\psi_{0}-k\mathrm{cn}\,\psi_{0}}$. 

Если $\lam = (\f, k) \in C_3$, $k = 1$, то ${\f_t} = \f +  t$  и 
$$
\left(\begin{array}{c}
x\\
y\\
z
\end{array}\right)=\left(\begin{array}{c}
\frac{s_{1}}{2}\left[\frac{1}{w}\left(\varphi-\varphi_{0}\right)+w\left(\tanh\varphi-\tanh\varphi_{0}\right)\right]\\
\frac{s_{2}}{2}\left[\frac{1}{w}\left(\varphi-\varphi_{0}\right)-w\left(\tanh\varphi-\tanh\varphi_{0}\right)\right]\\
-s_{1}s_{2}\ln[w\,\textrm{sech}\,\varphi]
\end{array}\right),\label{eq:52}
$$
где  $w=\cosh\varphi_{0}$.

Если $\lam =  (\g, c) \in C_4$,  то  
$$
\left(\begin{array}{c}
x\\
y\\
z
\end{array}\right)=\left(\begin{array}{c}
\sgn\left(\cos\frac{\gamma}{2}\right)t\\
0\\
0
\end{array}\right). 
$$

Если $\lam =  (\g, c) \in C_5$,  то 
$$
\left(\begin{array}{c}
x\\
y\\
z
\end{array}\right)=\left(\begin{array}{c}
0\\
0\\
\sgn\left(\sin\frac{\gamma}{2}\right)t
\end{array}\right).
$$

Проекция геодезической на плоскость $(x,y)$  имеет кривизну $\ds\frac{\tg \frac{\g}{2}}{(\ch 2 z)^{3/2}}$.  Она имеет точки перегиба при $\sin \frac{\g}{2} = 0$ (если $\lam \in C_1 \cup C_2 \cup C_3$)    и точки возврата при $\cos \frac{\g}{2} = 0$ (если $\lam \in C_2$).

\subsubsection{Симметрии и страты Максвелла}\label{subsubsec:sh2_max}
Фазовый портрет маятника \eq{sh2_pend}  имеет группу симметрий $\Sym = \{\Id, \eps^1, \dots, \eps^7\}$, описанную в разделе \ref{subsubsec:se2_max}. Продолжение этой группы симметрий на прообраз экспоненциального отображения $N = C \times \R_+$  описано в том же разделе. Продолжение этой группы симметрий на образ экспоненциального отображения имеет вид
$$
\eps^i \ : \ g = (x,y,z) \mapsto g^i = \eps^i (g) = (x^i, y^i, z^i),
$$
где
\begin{align}
(x^{1},y^{1},z^{1}) & =(x\cosh z-y\sinh z,\, x\sinh z-y\cosh z,\, z),\nonumber \\
(x^{2},y^{2},z^{2}) & =(x\cosh z-y\sinh z,\,-x\sinh z+y\cosh z,\,-z),\nonumber \\
(x^{3},y^{3},z^{3}) & =(x,\,-y,\,-z),\nonumber \\
(x^{4},y^{4},z^{4}) & =(-x,\, y,\,-z),\label{eq:6-6}\\
(x^{5},y^{5},z^{5}) & =(-x\cosh z+y\sinh z,\, x\sinh z-y\cosh z,\, -z),\nonumber \\
(x^{6},y^{6},z^{6}) & =(-x\cosh z+y\sinh z,\,-x\sinh z+y\cosh z,\, z),\nonumber \\
(x^{7},y^{7},z^{7}) & =(-x,\,-y,\, z).\nonumber 
\end{align}

Имеет место предложение, аналогичное предложению \ref{propos:se2_sym}.

\begin{theorem}
Первое время Максвелла, соответствующее группе симметрий $\Sym$, для почти всех геодезических выражается следующим образом:
\begin{eqnarray*}
\lambda\in C_{1} & \implies & \tmax(\lambda)=4K(k),\\
\lambda\in C_{2} & \implies &  \tmax(\lambda)=4kK(k),\\
\lambda\in C_{3}\cup C_{4}\cup C_{5} & \implies &  \tmax(\lambda)=+\infty.
\end{eqnarray*}
\end{theorem}

Имеет место 
\begin{corollary}
Для любого $\lam \in C$  первое время Максвелла $\tmax$   равно периоду колебаний маятника \eq{sh2_pend}.
\end{corollary}

Имеет место теорема, аналогичная теореме \ref{th:se2_max_inv}.

\subsubsection{Оценки сопряженного времени}\label{subsubsec:sh2_conj}
Обозначим через $p_1^1(k) \in (2 K, 3 K)$  первый положительный корень уравнения $\cn p \E(p) - \sn p \dn p = 0$.

\begin{theorem}\label{th:sh2_conjC1}
Если $\lam \in C_1$,  то $4 K(k) \leq \tconj(\lam) \leq 2 p_1^1(k)$. Более того, 
$$\lim_{k \to + 0} \tconj(\lam) = 2 \pi, 
\qquad
\lim_{k \to 1 - 0} \tconj(\lam) = + \infty.
$$
\end{theorem}

\begin{theorem}\label{th:sh2_conjC2}
Если $\lam \in C_2$,  то $4 k K(k) \leq \tconj(\lam) \leq 2 k p_1^1(k)$. Более того, 
$$\lim_{k \to + 0} \tconj(\lam) = 0, 
\qquad
\lim_{k \to 1 - 0} \tconj(\lam) = + \infty.
$$
\end{theorem}

\begin{theorem}\label{th:sh2_conjC4}
Если $\lam \in C_4$,  то $  \tconj(\lam) = 2 \pi$. 

Если $\lam \in C_3 \cup C_5$,  то $  \tconj(\lam) = + \infty$. 
\end{theorem}

\begin{theorem}
Нижние оценки для $\tconj(\lam)$ при $\lam \in C_1 \cup C_2$, приведенные в теоремах {\em \ref{th:sh2_conjC1}} и {\em \ref{th:sh2_conjC2}}, точны:
\begin{itemize}
\item[$(1)$]
если $\lam = (\f, k) \in C_1$  и $\sn \f = 0$,  то $\tconj(\lam) = 4 K(k)$,
\item[$(2)$]
если $\lam = (\f, k) \in C_2$  и $\sn \frac{\f}{k} = 0$,  то $\tconj(\lam) = 4 k K(k)$.
\end{itemize}
\end{theorem}

\subsubsection{Время разреза}\label{subsubsec:sh2_cut}
\begin{theorem}
Для любого $\lam \in C$
$$
\tcut(\lam) = 
\min(\tmax(\lam), \, \tconj(\lam)) = 
\begin{cases}
4 K(k), & \lam \in C_1,\\
4 k K(k), & \lam \in C_2,\\
2 \pi , & \lam \in C_4,\\
+ \infty , & \lam \in C_3 \cup C_5.
\end{cases}
$$
\end{theorem}

\begin{theorem}
\begin{itemize}
\item[$(1)$]
Функция $\map{\tcut}{C}{(0, + \infty]}$  зависит только от энергии $E$  маятника \eq{sh2_pend}.
\item[$(2)$]
Функция $\tcut$  инвариантна относительно вертикальной компоненты гамильтонова поля $\vec{H}_v$  и симметрий $\eps^i \in \Sym$.
\item[$(3)$]
Функция $\tcut$ является непрерывной на $C$  и гладкой на $C_1 \cup C_2$. 
\item[$(4)$]
$\lim_{E \to - 1} \tcut = 2 \pi$, $\lim_{E \to 1} \tcut = + \infty$, $\lim_{E \to + \infty} \tcut = 0$. 
\end{itemize}
\end{theorem}

\subsubsection{Диффеоморфная структура экспоненциального отображения}
Рассмотрим открытое всюду плотное подмножество в $G$, не содержащее первых точек Максвелла:
$$
\tG = \{ g \in G \mid z \neq 0\}
$$
 и его разбиение на компоненты связности
$$
\tG = G_1 \sqcup G_2, \qquad 
G_1 = \{ g \in G \mid z > 0\}, \quad  G_2 = \{ g \in G \mid z < 0\}.
$$
Также рассмотрим открытое плотное подмножество в пространстве всех потенциально оптимальных геодезических
$$
\widetilde N = \left\{(\lam, t) \in \cup_{i=1}^3 N_1 \cup N_5 \mid t < \tcut(\lam), \ \sin\left(\frac{\g_{t/2}}{2}\right) \neq 0 \right\}
$$
и его разбиение на компоненты связности
\begin{align*}
&\widetilde N = D_1 \sqcup D_2, \\
&D_1 = \left\{(\lam,t) \in \widetilde N \mid \sin\left(\frac{\g_{t/2}}{2}\right) > 0 \right\},\\
&D_2 = \left\{(\lam,t) \in \widetilde N \mid \sin\left(\frac{\g_{t/2}}{2}\right) < 0\right \}.
\end{align*}

\begin{theorem}
Отображения
\begin{align*}
&\map{\Exp}{D_i}{G_i}, \qquad i = 1, 2,\\
&\map{\Exp}{\widetilde N}{\tG}
\end{align*}
 суть диффеоморфизмы.
\end{theorem}

\subsubsection{Множество разреза}

\begin{theorem}
Множество разреза $\Cut$  содержится в плоскости $\{z = 0\}$. Имеет место разбиение на связные компоненты:
$$
\Cut = \Cut_{\loc}^+ \sqcup \Cut_{\loc}^- \sqcup \Cut_{\glob}^+ \sqcup \Cut_{\glob}^-,  
$$
где 
\begin{itemize}
\item
$\Cut_{\loc}^+$  есть часть плоскости $\{z=0\}$, ограниченная кривой 
\begin{align*}
&x = \pm \frac{4 k a(k)}{1-k^2}, \quad y = \frac{4a(k)}{1-k^2}, \qquad k \in [0, 1), \\
&a(k) = E(k) - (1-k^2)K(k),
\end{align*}
содержащая луч $\{z = x = 0, \ y > 0\}$  за вычетом начальной точки $\Id = \{x = y =  z = 0\}$,
\item
$\Cut_{\loc}^- $  получается из $\Cut_{\loc}^+$ отражением $(x,y) \mapsto (x, - y)$,
\item
$\Cut_{\glob}^+$  есть часть плоскости $\{z=0\}$, ограниченная кривой 
\begin{align*}
&x = \frac{4 E(k)}{1-k^2}, \quad y = \pm \frac{4kE(k)}{1-k^2}, \qquad k \in [0, 1), \\
&a(k) = E(k) - (1-k^2)K(k),
\end{align*}
содержащаяся в полуплоскости $\{z =  0, \ x > 0\}$,
\item
$\Cut_{\glob}^- $  получается из $\Cut_{\glob}^+$ отражением $(x,y) \mapsto (-x, -y)$.
\end{itemize}
Компоненты связности $\Cut_{loc}^{\pm}$  содержат в своем замыкании начальную точку $\Id$, а компоненты $\Cut_{glob}^{\pm}$  нет.
\end{theorem}

Множество разреза изображено на Рис. \ref{fig:sh2_cut}. На Рис. \ref{fig:sh2_cut_conj}  изображено множество разреза и первая каустика $\Conj^1$.

\figout{
\onefiglabelsize{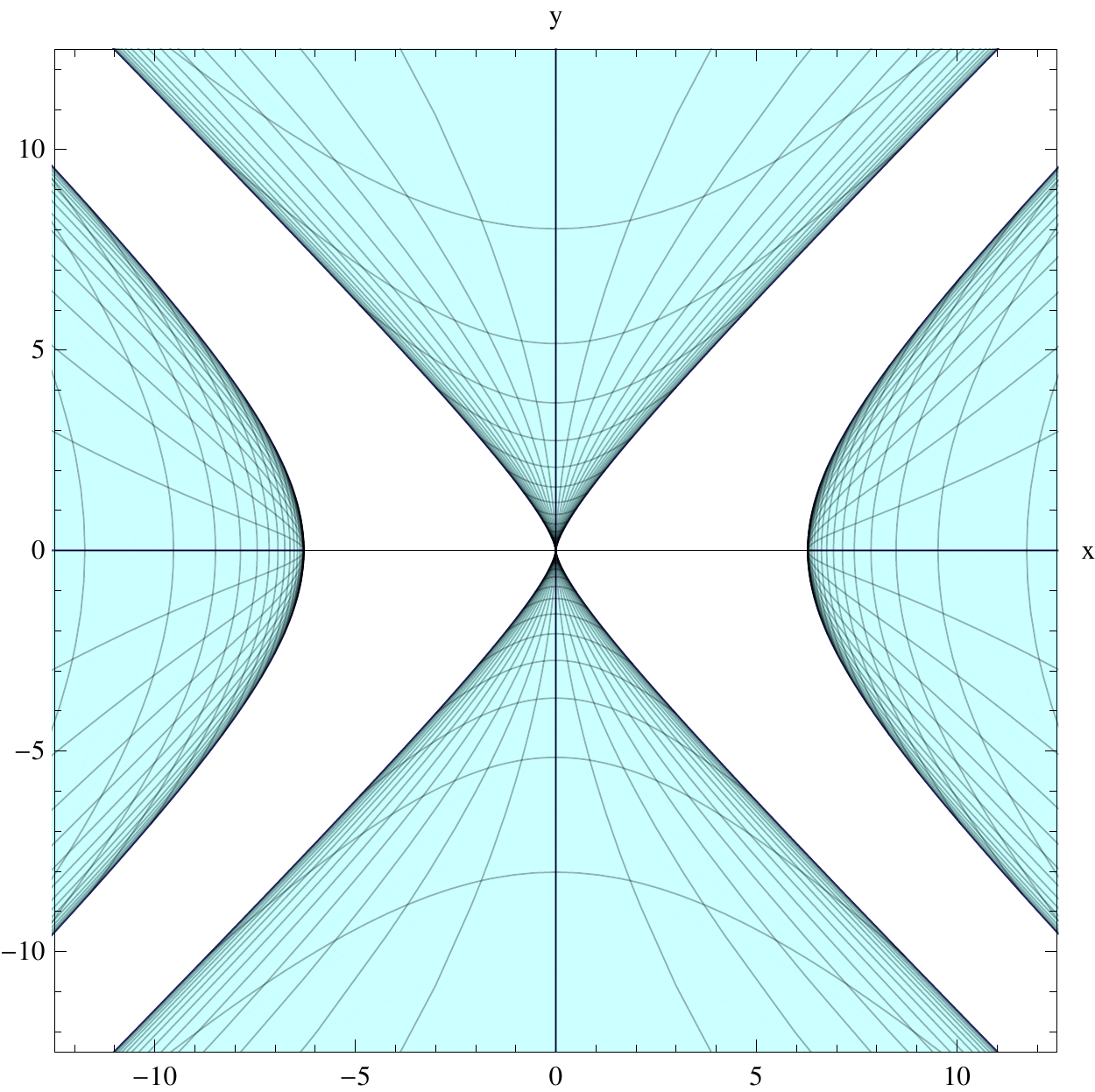}{Множество разреза на $\SH(2)$}{fig:sh2_cut}{0.5}
}

\figout{
\onefiglabelsize{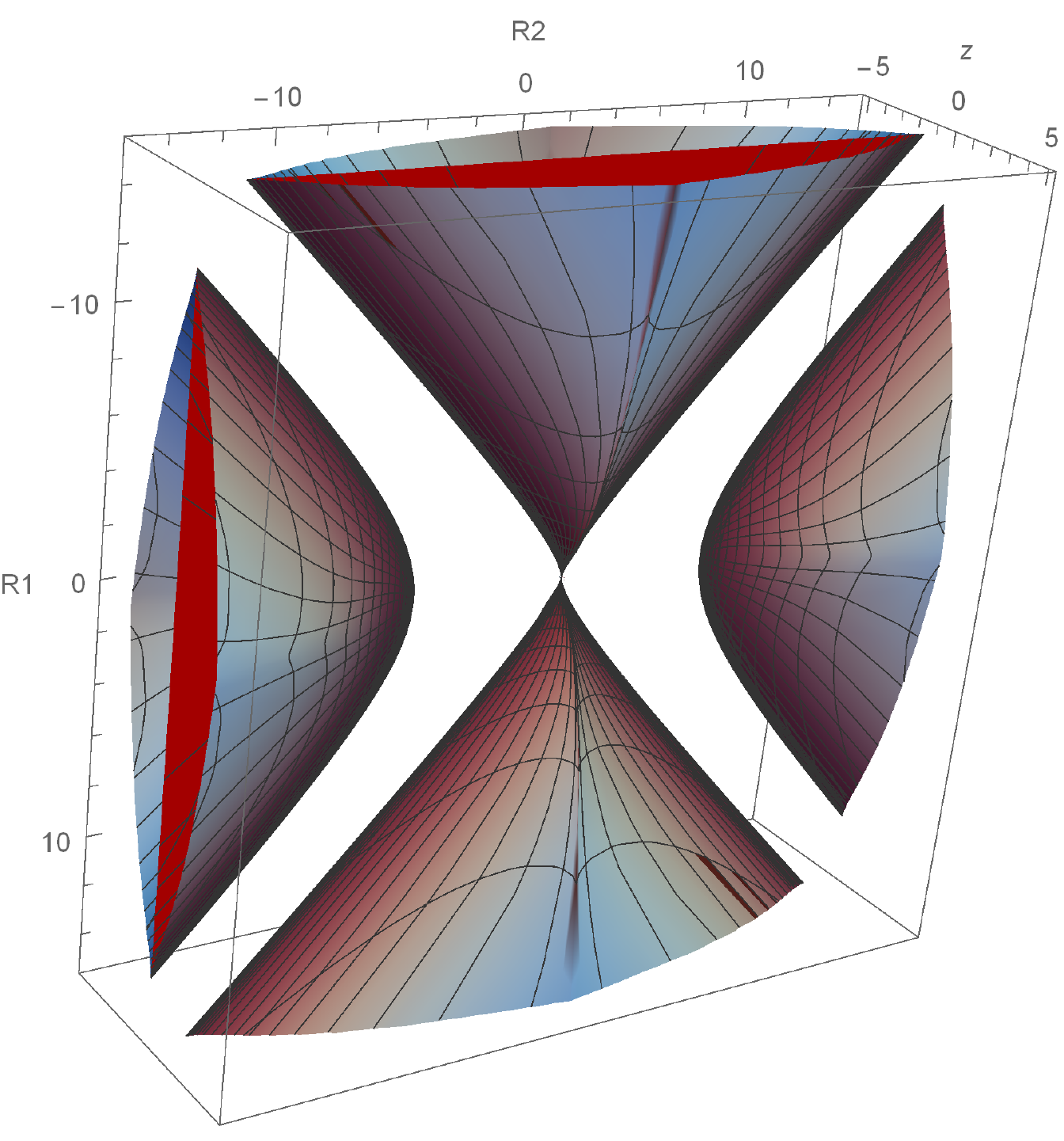}{Первая каустика и множество разреза на $\SH(2)$}{fig:sh2_cut_conj}{0.5}
}

\subsubsection{Сферы}
Субримановы сферы $S_R$, $R > 0$,  гомеоморфны двумерной евклидовой сфере, см. сферу $S_{\pi}$  на Рис. \ref{fig:sh2_Spi}  и сферу $S_{2 \pi}$  на Рис. \ref{fig:sh2_S2pi}.

\figout{
\twofiglabelsize
{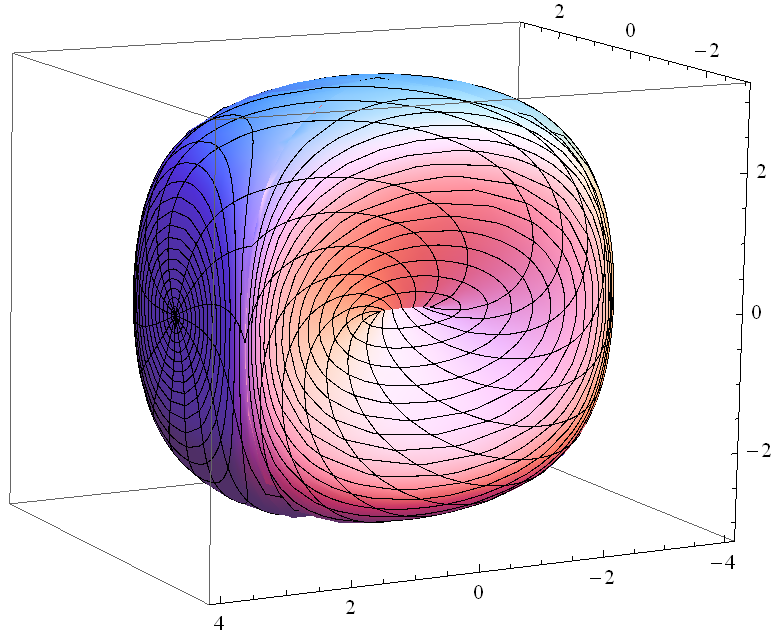}{Сфера $S_{\pi}  \subset \SH(2)$}{fig:sh2_Spi}{0.35}
{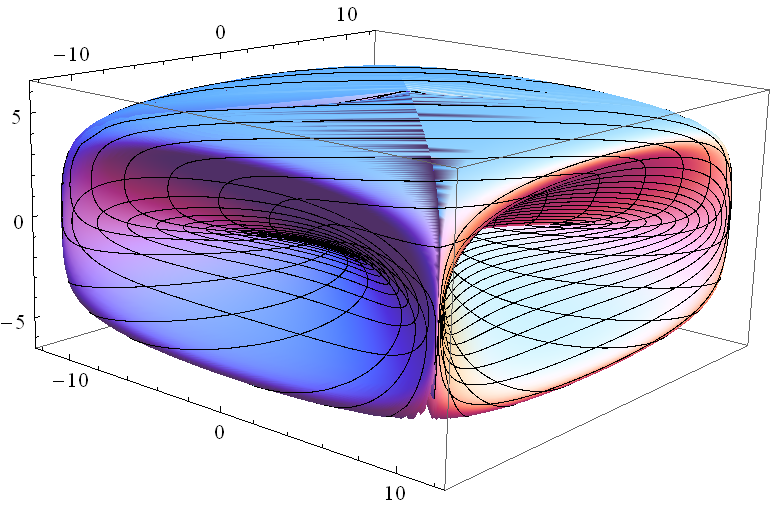}{Сфера $S_{2\pi}  \subset \SH(2)$}{fig:sh2_S2pi}{0.35}
}

Сферы имеют особенности при пересечении со множеством разреза, см. пересечение $\Cut$  и $S_{\pi} \cap \{z < 0\}$  на Рис. \ref{fig:sh2_Spi_cut}  и  пересечение $\Cut$  и $S_{2\pi} \cap \{z < 0\}$  на Рис. \ref{fig:sh2_S2pi_cut}.

\figout{
\twofiglabelsize
{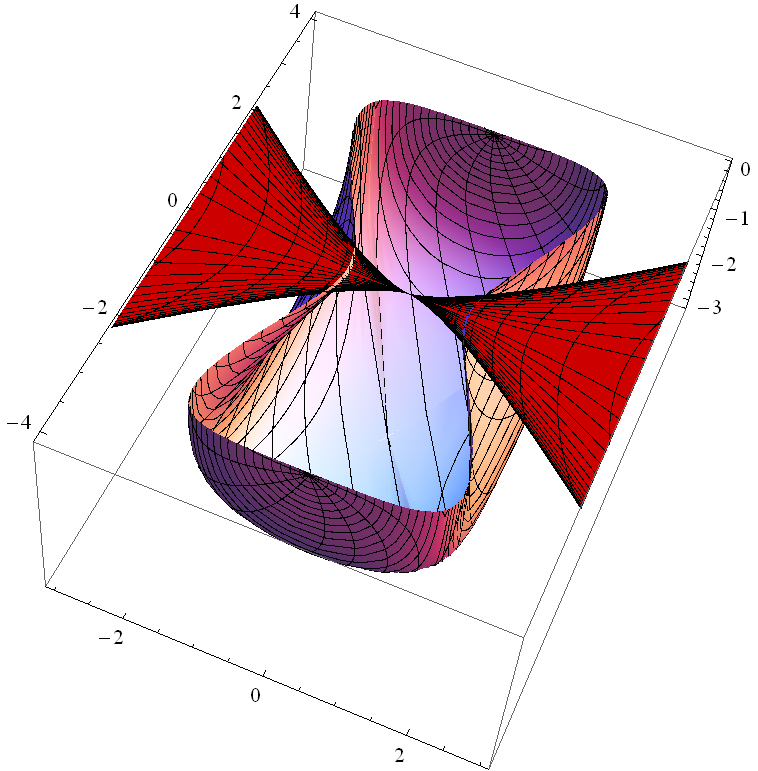}{Пересечение полусферы $S_{\pi}  \cap \{z < 0\}$  со множеством разреза}{fig:sh2_Spi_cut}{0.45}
{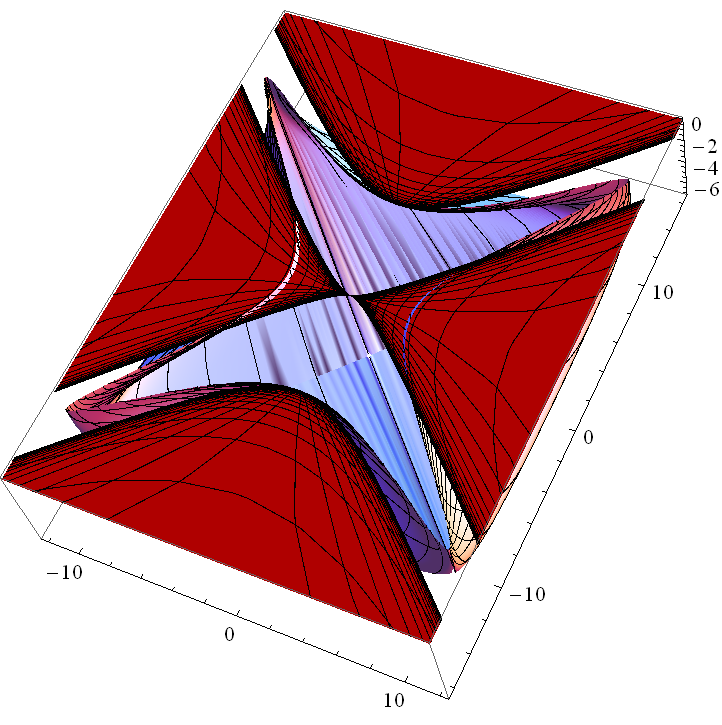}{Пересечение полусферы $S_{2 \pi}  \cap \{z < 0\}$  со множеством разреза}{fig:sh2_S2pi_cut}{0.45}
}

\subsubsection{Структура оптимального синтеза}
\begin{theorem}
\begin{itemize}
\item[$(1)$]
 Для любой точки $g_1 \in \Cut \setminus \Conj^1 = \intt_{\{z=0\}} \Cut$  существуют ровно две кратчайшие, соединяющие точки $\Id$  и $g_1$, причем для этих кратчайших  $g_1$ есть точка разреза и точка Максвелла, но не сопряженная точка.
\item[$(2)$]
 Для любой точки $g_1 \in \Cut \cap \Conj^1 = (\partial_{\{z=0\}} \Cut) \setminus \{\Id\}$  существует единственная кратчайшая, соединяющая точки $\Id$  и $g_1$, причем для этой кратчайшей  $g_1$ есть точка разреза и сопряженная точка, но не точка Максвелла.
\item[$(3)$]
 Для любой точки $g_1 \in G \setminus(\Cut \cup  \Id)$  существует единственная кратчайшая, соединяющая точки $\Id$  и $g_1$, причем для этой кратчайшей  $g_1$ не является ни точка разреза, ни сопряженной точкой, ни точкой Максвелла.
\end{itemize}
\end{theorem}

\subsubsection{Метрические прямые}\label{subsubsec:sh2_line}
Метрические прямые, проходящие через единичный элемент $\Id$, суть 
$$
g(t) = \Exp(\lam, t), \qquad t \in \R, \quad \lam \in C_3 \cup C_5.
$$

\subsubsection{Библиографические комментарии}
Раздел \ref{subsubsec:SH2}  опирается на книгу \cite{vilenkin},  разделы \ref{subsubsec:sh2_state},  \ref{subsubsec:sh2_geod} --- на \cite{sh2_1}, 
разделы \ref{subsubsec:sh2_max},  \ref{subsubsec:sh2_conj} --- на \cite{sh2_2}, 
разделы \ref{subsubsec:sh2_cut}--\ref{subsubsec:sh2_line} --- на \cite{sh2_3}.

\subsection{Задача Эйлера об эластиках}\label{subsec:elastica}
\subsubsection{История задачи}\label{subsubsec:el_history}
В 1691 году Я.~Бернулли рассмотрел задачу о форме однородного плоского упругого стержня, сжимаемого внешней силой. Он вывел уравнения для упругого стержня, закрепленного вертикально в горизонтальной стене и согнутого силой, направляющей его верхний конец горизонтально (\ddef{прямоугольная эластика}):
$$
dy = \dfrac{x^2 d x}{\sqrt{1 - x^4}}, \qquad ds = \dfrac{d x}{\sqrt{1 - x^4}}, \qquad x \in [0, 1),
$$ 
где $(x, y)$ есть упругий стержень, а $s$ --- его параметр длины (стержень отклоняется по горизонтали на расстояние 1). Я.~Бернулли  проинтегрировал эти дифференциальные  уравнения в рядах и получил двусторонние оценки их решения в конечной точке $x = 1$ \cite{JBernoulli}.

В 1742 году Д.~Бернулли  в своем письме \cite{DBernoulli} к Эйлеру написал, что упругая энергия стержня пропорциональна величине $\ds J = \int\frac{ds}{R^2}$, где $R$ --- радиус кривизны стержня, и предложил
отыскивать форму упругого стержня из вариационного принципа $J \to \min$. В это время Эйлер писал свой трактат по вариационному исчислению <<Methodus inveniendi \dots  >> \cite{euler_rus}, опубликованный  в 1744 году, и снабдил свою книгу приложением <<De curvis elesticis>>, в котором он применил только что разработанные  методы к задаче об упругих стержнях. Эйлер  рассмотрел тонкую однородную упругую пластину, прямолинейную в естественном (не напряженном) состоянии. Он поставил следующую  задачу для профиля  пластины:

<<{\it \dots среди всех кривых одной и той же длины, которые не только проходят через $A$ и $B$,  но и касаются в этих точках прямых, заданных по положению, определить ту, для которой значение выражения 
$\ds\int_A^B\frac{ds}{R^2}$ будет наименьшим}>>. 

Эйлер написал уравнение, известное сейчас как уравнение Эйлера-Лагранжа, для соответствующей вариационной задачи и свел его к уравнениям
\begin{equation*}
 dy = \frac{(\alpha + \beta x + \gamma x^2)\, dx}
 {\sqrt{a^4 -(\alpha + \beta x + \gamma x^2)^2}},
 \quad
 ds = \frac{a^2\,dx}{\sqrt{a^4 - (\alpha + \beta x + \gamma x^2)^2}},
\end{equation*}
параметры которых выражаются через упругие характеристики и длину стержня, а также величину нагрузки. Говоря современным языком, Эйлер исследовал качественное поведение эллиптических функций, параметризующих упругие кривые с помощью качественного анализа определяющих их уравнений. После работы Леонарда Эйлера кривые, представляющие форму однородного плоского стрежня, называются 
\ddef {эластиками Эйлера}.
Эйлер описал все типы эластик и указал значения параметров, для которых эти типы реализуются. Эйлер разделил все эластики на 9 классов, изображенных на рисунках:
\begin{enumerate}
\item
прямая линия,
\item
синусообразная кривая, Рис. \ref{fig:el2},
\item
прямоугольная эластика, Рис. \ref{fig:el3},
\item
Рис. \ref{fig:el4},
\item
замкнутая эластика в форме восьмерки, Рис. \ref{fig:el5},
\item
Рис. \ref{fig:el6},
\item
непериодическая эластика с одной петлей, <<солитон Эйлера>>, Рис. \ref{fig:el7},
\item
Рис. \ref{fig:el8},
\item
окружность.
\end{enumerate}
Эластики типов 2--6, имеющие точки перегиба, называются \ddef{инфлексионными}, эластика типа 7 называется \ddef{критической},  а эластики типа 8 без точек перегиба называются \ddef{неинфлексионными}.
Семейство всех эластик изображено на Рис.  \ref{fig:elastica_all}.

\figout{
\begin{figure}
\includegraphics{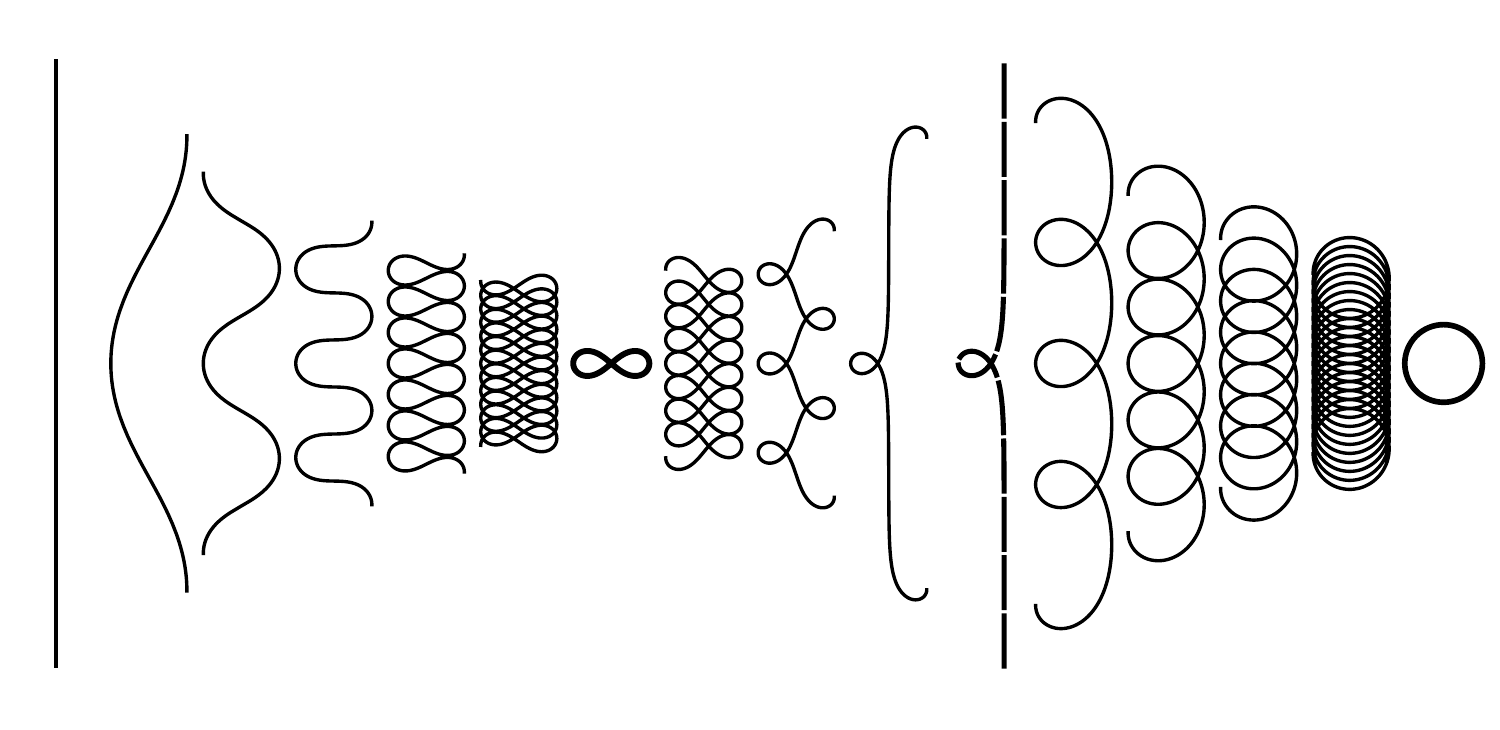}
\caption{Эластики Эйлера\label{fig:elastica_all}}
\end{figure}
}

Первую явную параметризацию эластик Эйлера получил Л. Заалшютц в 1880 г. \cite{saalchutz}.

 В 1906 г. будущий нобелевский лауреат Макс Борн защитил диссертацию <<Устойчивость упругих кривых на плоскости и в пространстве>> \cite{born}. Он рассмотрел задачу об эластиках методами вариационного исчисления и вывел из уравнения Эйлера-Лагранжа уравнения
 \begin{align*}
 &\dot x = \cos \theta, \quad \dot y = \sin \theta,
 \\
 &A \ddot \theta + R \sin (\theta - \gamma) = 0,
 \quad A, R, \gamma = \const.
\end{align*}
То есть угол $\t$  наклона эластик удовлетворяет уравнению математического маятника. Далее, Борн изучил устойчивость эластик с закрепленными концами и касательными на концах. Он доказал, что дуга эластики без точек перегиба устойчива (в этом случае угол $\t$ монотонен и может быть выбран параметром на эластике; Борн показал, что вторая вариация функционала упругой энергии $J = \frac 12 \int \dot\t^2 dt$  положительна).  В общем случае Борн записал якобиан, обращающийся в нуль в сопряженных точках. В силу сложности функций, входящих в якобиан, Борн ограничился численным исследованием сопряженных точек. Он первым численно построил чертежи эластик и проверил теоретические результаты с помощью экспериментов с упругими стержнями. Более того, Борн исследовал устойчивость эластик с различными другими граничными условиями и получил некоторые результаты для трехмерных упругих кривых.

В 1993 г. В. Джурджевич \cite{jurd_ball}  обнаружил эластики Эйлера в задаче о качении шара по плоскости без прокручивания и проскальзывания (см. раздел \ref{subsec:roll}), а Р. Брокетт и Л. Даи \cite{brock_dai}  --- в субримановой задаче на группе Картана (см. раздел \ref{subsec:cartan}). Эластики Эйлера также  удивительным образом появляются в плоской субримановой задаче Мартине (см. раздел \ref{subsec:martinet}), субримановых задачах на группах $\SE(2)$ (см. раздел \ref{subsec:se2})  и на группе Энгеля (см. раздел \ref{subsec:engel}). Было бы интересно понять, почему эластики Эйлера появляются в стольких задачах оптимального управления.

Далее задача Эйлера об эластиках исследовалась в работах \cite{el_max, el_conj, el_dan, el_AiT, el_stable, el_cut, el_closed, el_SE2}, на которые опирается изложение в этом разделе.

\subsubsection{Постановка задачи}\label{subsubsec:el-state}
\paragraph{Механическая постановка}
Пусть однородный упругий стержень на плоскости $\R^2$  имеет длину $l > 0$. Выберем любые точки $a_0, a_1 \in \R^2$  и произвольные единичные касательные вектора $v_i \in T_{a_i} \R^2$, $|v_i| = 1$, $i = 0, 1$. Задача заключается в том, чтобы найти профиль стержня $\map{\g}{[0,l]}{\R^2}$, $|\dot \g(s)| \equiv 1$, выходящего из точки $a_0$  и приходящего в точку $a_1$  с соответствующими касательными векторами $v_0$ и $v_1$:
\begin{align*}
&\g(0) = a_0, \quad \g(l) = a_1, \\ 
&\dot\g(0) = v_0, \quad \dot\g(l) = v_1, 
\end{align*}
 с минимальной упругой энергией
$$
J = \frac 12 \int_0^l k^2(s) ds \to \min,
$$
 где $k(s)$ ---  кривизна кривой $\g(s)$. 

\paragraph{Задача оптимального управления}
Выберем на плоскости $\R^2$  декартовы координаты $(x,y)$. Будем обозначать параметр длины $s$  на кривой $\g$  через $t$, и пусть $t_1 = l$. Искомая кривая имеет параметризацию $\g(t) = (x(t), y(t))$, $t \in [0, t_1]$,  а ее граничные точки имеют координаты $a_i = (x_i, y_i)$, $i = 0, 1$.  Обозначим через  $\t(t)$  угол между касательным вектором $\dot \g(t)$  и положительным направлением оси $x$. Наконец, пусть касательные векторы в граничных точках кривой $\g$ имеют координаты $v_i = (\cos \t_i, \sin \t_i)$, $i = 0, 1$, см. Рис. \ref{fig:el_problem}.
\figout{
\onefiglabelsize{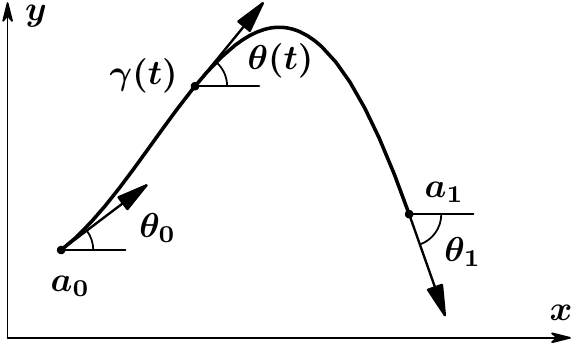}{Постановка задачи об эластиках Эйлера}{fig:el_problem}{0.4}
}

Тогда искомая кривая $\g(t) = (x(t), y(t))$  есть проекция траектории следующей управляемой системы:
\begin{align}
 &\dot x = \cos \theta, \label{el_pr1} \\
 &\dot y = \sin \theta, \label{el_pr2} \\
 &\dot \theta = u, \label{el_pr3} \\
 &g=(x,y,\theta) \in M= \R^2_{x,y} \times S^1_{\theta},
 \quad u \in \R,  \label{el_pr4} \\
 &g(0) = g_0 = (x_0, y_0, \theta_0), \quad
 g(t_1) = g_1 = (x_1, y_1, \theta_1), \quad t_1 \text{ фиксировано}.
 \label{el_pr5}
\end{align}
Для натурально параметризованной кривой $\g$  кривизна равна угловой скорости: $k = \dot \t = u$, откуда получаем функционал качества
\begin{equation}\label{el_pr6}
 J = \frac 12 \int_0^{t_1} u^2(t) \, dt\to \min.
\end{equation}
Естественный класс допустимых управлений для задачи \eq{el_pr1}--\eq{el_pr6} есть $u(\cdot) \in L^2[0, t_1]$, поэтому допустимая траектория есть $g(\cdot) \in W^{1,2}([0, t_1], M)$. 

В векторных обозначениях задача принимает форму:
\begin{align}
& \dot g = X_1(g) + u X_2(g),
 \quad g \in M = \R^2 \times S^1,
 \quad u \in \R, \label{el_sys}
 \\
 \nonumber
 &g(0) = g_0, \quad g(t_1) = g_1, \quad t_1 \text{ фиксировано},
 \\
 \nonumber
 &J = \frac 12 \int_0^{t_1} u^2 dt \to \min, \quad
 u \in L^2[0, t_1],
\end{align}
где векторные поля в правой части системы \eq{el_sys}  суть
$$
 X_1 = \cos \theta \pder{}{x} + \sin \theta \pder{}{y}, \quad
 X_2 = \pder{}{\theta}.
$$
Пространство состояний $M = \R^2 \times S^1$ имеет естественную структуру группы движений плоскости $G = \R^2 \ltimes \SO(2)$, см. раздел \ref{subsec:se2}. При этом векторные поля $X_1$, $X_2$  становятся левоинвариантными полями на группе Ли $G$.
Таблица умножения в алгебре Ли $\gg = \sea(2)$ приведена в \eq{se2_tab}.   

Таким образом, задача Эйлера об эластиках \eq{el_pr1}--\eq{el_pr6}  есть левоинвариантная задача оптимального управления на группе $\SE(2)$. Поэтому можно считать, что $g_0 = \Id = (0,0,0)$.

\subsubsection{Множество достижимости}
\begin{theorem}
Множество достижимости системы \eq{el_sys}  из точки $\Id = (0,0,0)$  за время $t_1 > 0$  есть
$$
\A(t_1) = \{(x,y,\t) \in G \mid x^2 + y^2 < t_1^2 \text{ или } (x,y,\t) = (t_1, 0, 0)\}.
$$
\end{theorem}

Топологически множество достижимости $\A(t_1)$  есть открытый полноторий (внутренность тора) с одной точкой на границе. Будем далее рассматривать задачу об эластиках при естественном условии управляемости: $g_1 \in \A(t_1)$.

\subsubsection{Существование и ограниченность оптимальных управлений}
\begin{theorem}
Пусть $g_1 \in \A(t_1)$.  Тогда существует оптимальное управление $u \in L^2[0, t_1]$.  Более того, $u \in L^{\infty}[0, t_1]$.  Поэтому оптимальное управление удовлетворяет принципу максимума Понтрягина.
\end{theorem}

\subsubsection{Экстремали}
\paragraph{Анормальные траектории}
Проходящая через точку $\Id$  натурально параметризованная анормальная траектория есть $(x,y,\t) = (t, 0, 0)$, $t \in [0, t_1]$. Она проецируется на плоскость $(x,y)$  в отрезок --- это упругий стержень в отсутствие внешних сил. Упругая энергия в этом случае достигает абсолютного минимума $J = 0$, поэтому анормальная траектория оптимальна. Именно эта траектория приходит в единственную точку $(t_1, 0, 0)$ на границе  множества достижимости $\A(t_1)$. Анормальная траектория одновременно нормальна.

\paragraph{Нормальные экстремали}
Нормальные экстремали удовлетворяют гамильтоновой системе $\dot \lam = \vec{H}(\lam)$, $\lam \in T^*G$,  где $H = h_1 + \frac 12 h_2^2$, $h_i(\lam) = \langle \lam, X_i\rangle$, $i = 1, 2, 3$. В координатах эта система имеет вид
\begin{align}
&\dot h_1 = - h_2h_3, \quad \dot h_2 = h_3, \quad \dot h_3 = h_1h_2, \label{el_vert}\\
&\dot g = X_1 + h_2 X_2. \label{el_hor}
\end{align}
Вертикальная подсистема \eq{el_vert}  имеет интеграл --- функцию Казимира $F = h_1^2 + h_3^2$. 

Введем координаты
$$
c = h_2, \quad h_1 = - r \cos \g, \quad h_2 = - r \sin \g,
$$
 в которых вертикальная подсистема \eq{el_vert}  принимает форму \ddef{математического маятника}
\be{el_pend}
\dot \g = c, \quad \dot c = - r \sin \g, \qquad c \in \R, \quad \g \in S^1, \quad r \equiv \const \geq 0,
\ee
 известного как \ddef{кинетический аналог Кирхгофа}  для эластик. Полная энергия маятника есть 
$$
E = H = \frac{c^2}{2} - r \cos \g \in [-r, + \infty).
$$

\paragraph{Стратификация прообраза экспоненциального отображения и выпрямляющие координаты}
\ddef{Экспоненциальное отображение} за время $t_1 > 0$  в задаче об эластиках есть
$$
\map{\Exp_{t_1}}{N = \gg^*}{G}, \qquad \lam \mapsto \pi \circ e^{t_1 \vec{H}}(\lam),
$$
 где $\map{\pi}{T^*G}{G}$  есть каноническая проекция.

Прообраз экспоненциального отображения $N = \gg^*$  разбивается на инвариантные многообразия гамильтонова поля $\vec{H}$ критическими множествами энергии $E = H$:
\begin{align*}
 &N = \sqcup_{i=1}^7 N_i,  \\
 &N_1 = \{ \lambda \in N \mid r \neq 0, \ E \in (-r, r)\}, \\
 &N_2 = \{ \lambda \in N \mid r \neq 0, \ E \in (r, + \infty)\},  \\
 &N_3 = \{ \lambda \in N \mid r \neq 0,
 \ E  = r, \ \g \neq \pi \},  \\
 &N_4 = \{ \lambda \in N \mid r \neq 0, \ E =-r \},  \\
 &N_5 = \{ \lambda \in N \mid r \neq 0, \ E =r, \ \g = \pi \}, \\
 &N_6 = \{ \lambda \in N \mid r = 0, \ c \neq 0 \}, \\
 &N_7 = \{ \lambda \in N \mid r = c =  0\}. 
\end{align*}
На множествах $N_1$, $N_2$, $N_3$  введем координаты $(\f, k, r)$  следующим образом:
\begin{gather*}
 \lambda = (\g,c, r) \in N_1 \then
 \begin{cases}
 \sin \frac{\g}{2} = k \sn(\sqrt{r} \varphi, k), \\
 \frac{c}{2} = k \sqrt{r} \cn(\sqrt{r} \varphi, k), \\
 \cos \frac{\g}{2} = \dn(\sqrt{r} \varphi, k),
 \end{cases}
 \\
 k = \sqrt{\frac{E+r}{2r}} \in (0,1), \quad
 \sqrt{r} \varphi  \pmod{ 4 K(k)} \in [0, 4 K(k)],
 \\
 \lambda = (\g,c,r) \in N_2 \then
 \begin{cases}
 \sin \frac{\g}{2} = \pm \sn\left(\frac{\sqrt{r} \varphi}{k}, k\right),
 \\
 \frac{c}{2} = \pm \frac{\sqrt{r}}{k}
 \dn\left(\frac{\sqrt{r} \varphi}{k}, k\right),
 \\
 \cos \frac{\g}{2} = \cn\left(\frac{\sqrt{r} \varphi}{k}, k\right),
 \end{cases}
 \\
 k = \sqrt{\frac{2 r}{E+r}} \in (0,1),
 \ \sqrt{r} \varphi \pmod{ 2 K(k)k} \in [0, 2 K(k)k],
 \ \pm = \sgn c,
 \\
 \lambda = (\g,c, r) \in N_3 \then
 \begin{cases}
 \sin \frac{\g}{2} = \pm \tanh (\sqrt{r} \varphi), \\
 \frac{c}{2} = \pm \frac{\sqrt{r}}{\cosh(\sqrt{r} \varphi)}, \\
 \cos \frac{\g}{2} = \frac{1}{\cosh(\sqrt{r} \varphi)},
 \end{cases}
 \\
 k = 1, \quad \varphi \in \R, \quad \pm = \sgn c.
\end{gather*}

\paragraph{Параметризация экстремалей}
В области $N_1 \cap N_2 \cup N_3$  уравнение маятника выпрямляется:
$$
\dot \f = 1, \quad \dot k = \dot r = 0,
$$
поэтому имеет решения
$$
\f_t = \f + t, \quad k, \ r \equiv \const.
$$
 В исходных координатах $(\g, c)$  уравнение маятника \eq{el_pend}  имеет решения:
\begin{align*}
 \lambda  \in N_1 &\then
 \begin{cases}
 \sin \frac{\g_t}{2} = k_1 \sn(\sqrt{r} \varphi_t), \\
 \cos \frac{\g_t}{2} = \dn(\sqrt{r} \varphi_t), \\
 \frac{c_t}{2} = k \sqrt{r} \cn(\sqrt{r} \varphi_t),
 \end{cases}
 \\
 \lambda  \in N_2 &\then
 \begin{cases}
 \sin \frac{\g_t}{2} = \pm \sn\left(\frac{\sqrt{r} \varphi_t}{k}\right), \\
 \cos \frac{\g_t}{2} = \cn\left(\frac{\sqrt{r} \varphi_t}{k}\right), \\
 \frac{c_t}{2} = \pm \frac{\sqrt{r}}{k} \dn\left(\frac{\sqrt{r} \varphi_t}{k}\right), \quad \pm = \sgn c,
 \end{cases}
 \\
 \lambda  \in N_3 &\then
 \begin{cases}
 \sin \frac{\g_t}{2} = \pm \tanh (\sqrt{r} \varphi_t), \\
 \cos \frac{\g_t}{2} = \frac{1}{\cosh(\sqrt{r} \varphi_t)}, \\
 \frac{c_t}{2} = \pm \frac{\sqrt{r}}{\cosh(\sqrt{r} \varphi_t)}, \quad \pm = \sgn c.
 \end{cases}
\end{align*}
В вырожденных случаях $\cup_{i=4}^7 N_i$  уравнение маятника \eq{el_pend}  интегрируется в элементарных функциях:
\begin{align*}
\lambda  \in N_4 &\then
 \g_t \equiv 0, \quad c_t \equiv 0,
 \\
 \lambda  \in N_5 &\then
 \g_t \equiv \pi, \quad c_t \equiv 0,
 \\
 \lambda  \in N_6 & \then \g_t  = c t + \g,
 \quad c_t \equiv  c,
 \\
 \lambda  \in N_7 &\then  c_t \equiv 0, \quad r \equiv 0.
\end{align*}

Параметризация решений горизонтальной подсистемы \eq{el_hor}  имеет следующий вид.

Если $\lam \in N_1$, то
\begin{align*}
 \sin \frac{\theta_t}{2} &= k \dn (\sqrt{r} \varphi)
 \sn (\sqrt{r} \varphi_t)  - k \sn (\sqrt{r} \varphi)
 \dn (\sqrt{r} \varphi_t),
 \\
 \cos \frac{\theta_t}{2} &= \dn (\sqrt{r} \varphi)
 \dn (\sqrt{r} \varphi_t)  + k^2 \sn (\sqrt{r} \varphi)
 \sn (\sqrt{r} \varphi_t),
 \\
 x_t &= \frac{2}{\sqrt r} \dn^2 (\sqrt r \varphi)
 (\E(\sqrt r \varphi_t) - \E(\sqrt r \varphi))
 \\
 &\qquad + \frac{4k^2}{\sqrt r}  \dn (\sqrt{r} \varphi) \sn
 (\sqrt{r} \varphi)
 (\cn \sqrt{r} \varphi) - \cn (\sqrt{r} \varphi_t))
 \\
 &\qquad + \frac{2k^2 }{\sqrt r}  \sn^2 (\sqrt{r} \varphi)
 (\sqrt r t + \E(\sqrt r \varphi) - \E(\sqrt r \varphi_t)) - t,
 \\
 y_t &= \frac{2k}{\sqrt r}(2 \dn^2 (\sqrt{r} \varphi) -1)
 (\cn (\sqrt{r} \varphi) - \cn (\sqrt{r} \varphi_t))
 \\
 &\qquad - \frac{2k}{\sqrt r} \sn (\sqrt{r} \varphi) \dn (\sqrt{r}
 \varphi)(2(\E(\sqrt r \varphi_t) - \E(\sqrt r \varphi))
 - \sqrt r t).
\end{align*}

Если $\lam \in N_2$, то
\begin{align*}
 \sin \frac{\theta_t}{2} &=
 \pm( \cn(\sqrt r \psi) \sn(\sqrt r \psi_t)
 - \sn(\sqrt r \psi) \cn(\sqrt r \psi_t)),
 \\
 \cos \frac{\theta_t}{2} &=
 \cn(\sqrt r \psi) \cn(\sqrt r \psi_t)
 + \sn(\sqrt r \psi) \sn(\sqrt r \psi_t),
 \\
 x_t &= \frac{1}{\sqrt r} (1 - 2 \sn^2(\sqrt r \psi))
 \left( \frac 2k (\E(\sqrt r \psi_t) - \E(\sqrt r \psi))
 - \frac{2-k^2}{k^2} \sqrt r t \right)
 \\
 &+ \frac{4}{k \sqrt r} \cn(\sqrt r \psi)
 \sn(\sqrt r \psi) (\dn(\sqrt r \psi) - \dn(\sqrt r \psi_t)),
 \\
 y_t &= \pm \left(
 \frac{2}{k \sqrt r}(2 \cn^2(\sqrt r \psi) - 1)(\dn(\sqrt r \psi)
 - \dn(\sqrt r \psi_t)) \right.
 \\
 &\left. -
 \frac{2}{\sqrt r} \sn(\sqrt r \psi) \cn(\sqrt r \psi)
 \left( \frac{2}{k}(\E(\sqrt r \psi_t) - \E(\sqrt r \psi)) -
 \frac{2-k^2}{k^2} \sqrt r t\right)\right).
\end{align*}
где $\pm = \sgn c$, $\p_t = \frac{\f_t}{k} = \frac{\f+t}{k}$.

Если $\lam \in N_3$, то
\begin{align*}
 \sin \frac{\theta_t}{2} &= \pm \left( \frac{\tanh (\sqrt r
 \varphi_t)}{\cosh (\sqrt r \varphi)} -
 \frac{\tanh \sqrt r \varphi)}{\cosh (\sqrt r \varphi_t)}\right),
 \\
 \cos \frac{\theta_t}{2} &=
 \frac{1}{\cosh(\sqrt r \varphi) \cosh (\sqrt r \varphi_t)}
 + \tanh (\sqrt r \varphi) \tanh (\sqrt r \varphi_t),
 \\
 x_t &=  (1 - 2 \tanh^2 (\sqrt r \varphi)) t
 \\
 &+ \frac{4 \tanh (\sqrt r \varphi)}{\sqrt r \cosh (\sqrt r \varphi)}
 \left( \frac{1}{\cosh (\sqrt r \varphi)}
 -  \frac{1}{\cosh (\sqrt r \varphi_t)}\right),
 \\
 y_t &= \pm \left( \frac{2}{\sqrt r}
 \left( \frac{2}{\cosh^2\sqrt r \varphi)}-1\right)
 \left( \frac{1}{\cosh (\sqrt r \varphi)}
 - \frac{1}{\cosh (\sqrt r \varphi_t)}\right) \right.
 \\
 &\qquad\qquad \left. - 2 \frac{\tanh (\sqrt r \varphi)}
 {\cosh (\sqrt r \varphi)} t \right).
\end{align*}
где $\pm = \sgn c$.

Если $\lam \in N_4 \cup N_5 \cup N_7$, то
$$
 \theta_t = 0, \quad x_t = t, \quad y_t = 0.
$$

Если $\lam \in N_6$, то
$$
 \theta_t = ct, \quad x_t = \frac{\sin ct}{c},
 \quad y_t = \frac{1-\cos ct}{c}.
$$

\paragraph{Эластики Эйлера}
Проекции экстремальных траекторий на плоскость $(x,y)$  суть эйлеровы эластики. Эти кривые удовлетворяют уравнениям
\begin{align}
&\dot x = \cos \t, \quad \dot y = \sin \t, \nonumber\\
&\ddot \t = - r \sin(\t - \g), \qquad r, \ \g \equiv \const. \label{el_pend2}
\end{align}
В зависимости от значения энергии маятника $E = \frac{\dot \t^2}{2} - r \cos (\t-\g) \in [-r, + \infty)$  и функции Казимира $r \geq 0$, эластики имеют разные качественные типы, открытые Эйлером.

Если энергия $E$  принимает минимальное значение $-r < 0$, т.е. $\lam \in N_4$,  то эластика $(x_t,y_t)$  есть прямая. Соответствующее движение маятника \eq{el_pend2} (кинетический аналог Кирхгофа)  есть устойчивое положение равновесия.

Если $E \in (-r, r)$, $r > 0$, т.е. $\lam \in N_1$,  то маятник \eq{el_pend2} колеблется между экстремальными значениями угла, и угловая скорость $\dot \t$  меняет знак. Соответствующие эластики имеют точки перегиба при $\dot \t = 0$  и вершины при $|\dot \t| = \max$, т.к. $\dot \t$  есть кривизна эластики. Такие эластики называются \ddef{инфлексионными}, см. Рис. \ref{fig:el2}--\ref{fig:el6}.  Разные случаи на этих рисунках определяются значениями модуля эллиптических функций $k = \frac{\sqrt{E+r}}{2r} \in (0, 1)$:
\begin{alignat*}{2}
 &k \in \left(0, \frac{1}{\sqrt 2}\right)
 &&\then \text{Рис.~\ref{fig:el2}}, \\
 &k = \frac{1}{\sqrt 2}
 &&\then \text{Рис.~\ref{fig:el3}}, \\
 &k \in \left(\frac{1}{\sqrt 2}, k_0\right)
 &&\then \text{Рис.~\ref{fig:el4}}, \\
 &k = k_0
 &&\then \text{Рис.~\ref{fig:el5}}, \\
 &k \in \left(k_0, 1\right)
 &&\then \text{Рис.~\ref{fig:el6}}.
\end{alignat*}
Значение $k = 1/\sqrt 2$ соответствует \ddef{прямоугольной эластике}, исследованной Я. Бернулли (см. раздел \ref{subsubsec:el_history}), Рис. \ref{fig:el3}. Значение $k \approx 0,909$  соответствует периодической эластике в форме восьмерки, см. Рис.  \ref{fig:el5}. Как отмечал Эйлер, при $k \to 0$  инфлексионные эластики похожи на синусоиды, что соответствует гармоническому осциллятору $\ddot \t = - r (\t - \g)$  как кинетическому аналогу Кирхгофа, см. Рис.  \ref{fig:el2}.

Если $E = r > 0$  и $\t - \g \neq \pi$, т.е. $\lam \in N_3$,  то маятник \eq{el_pend2}  стремится к неустойчивому положению равновесия ($\t - \g = \pi$, $\dot \t = 0$) вдоль сепаратрисы седла, а соответствующая критическая эластика (<<солитон Эйлера>>) имеет одну петлю, см. Рис. \ref{fig:el7}.

Если $E = r > 0$  и $\t - \g = \pi$, т.е. $\lam \in N_5$,  то маятник \eq{el_pend2}  находится в   неустойчивом положении равновесия ($\t - \g = \pi$, $\dot \t = 0$) и эластика есть прямая.

Если $E > r > 0$, т.е.  $\lam \in N_2$, то кинетический аналог Кирхгофа есть маятник \eq{el_pend2}, вращающийся против часовой стрелки ($\dot \t > 0$)  или по часовой стрелке ($\dot \t < 0$). Соответствующие эластики имеют ненулевую кривизну $\dot \t$, не имеют точек перегиба и называются \ddef{неинфлексионными}, см. Рис. \ref{fig:el8}.

Если $r = 0$  и $\dot \t \neq 0$, т.е. $\lam \in N_6$, то маятник \eq{el_pend2} равномерно вращается в невесомости, и соответствующая эластика есть окружность.

Наконец, если $r = 0$  и $\dot \t = 0$, т.е. $\lam \in N_7$, то маятник \eq{el_pend2}  неподвижен в невесомости (положение равновесия неустойчиво), и эластика есть прямая.

Изображения эластик на Рис. \ref{fig:el2}--\ref{fig:el8}  не всегда передают отношение $x/y$  для экономии места.

\figout{
\twofiglabelsize
{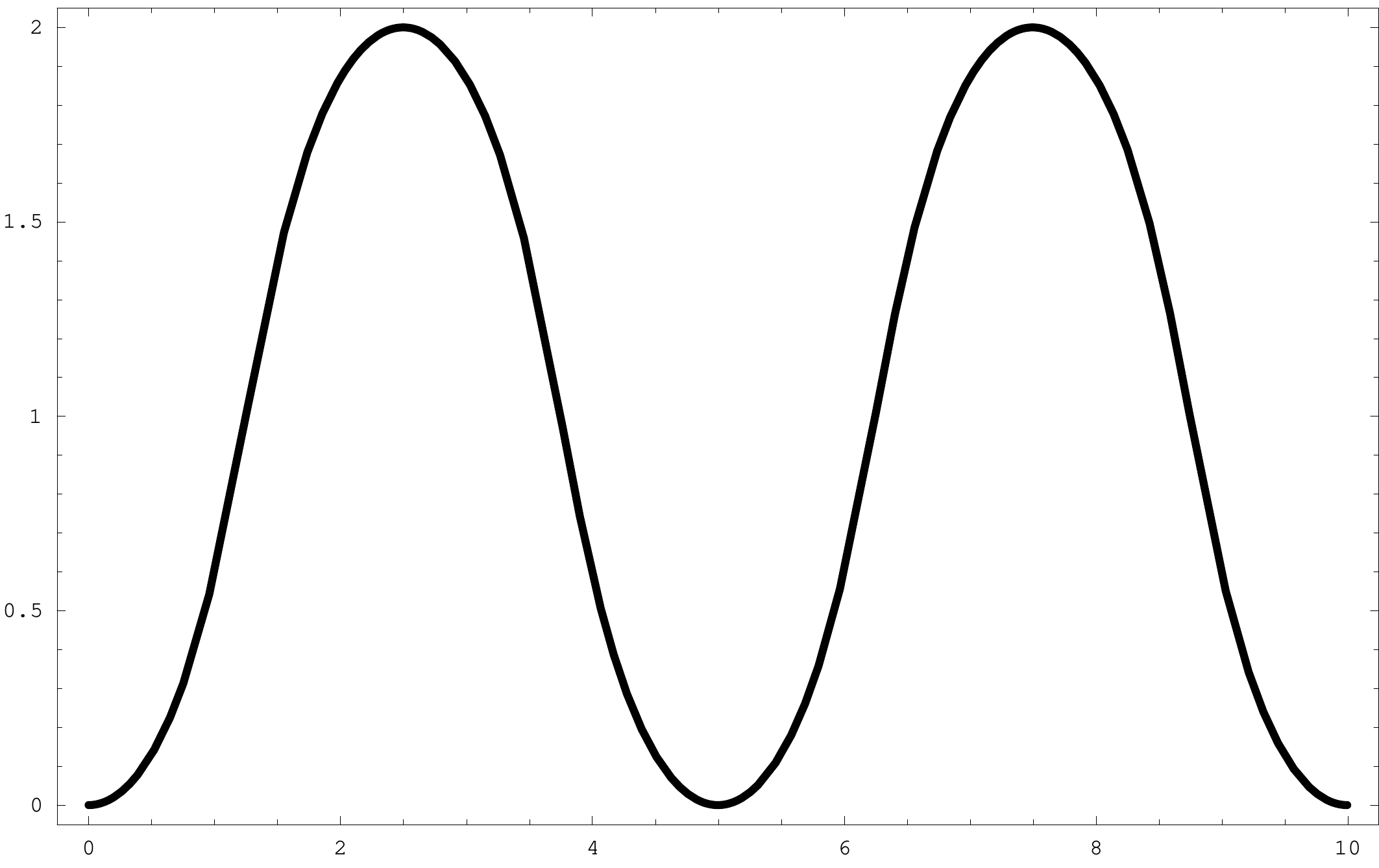}{Инфлексионная эластика}{fig:el2}{0.25}
{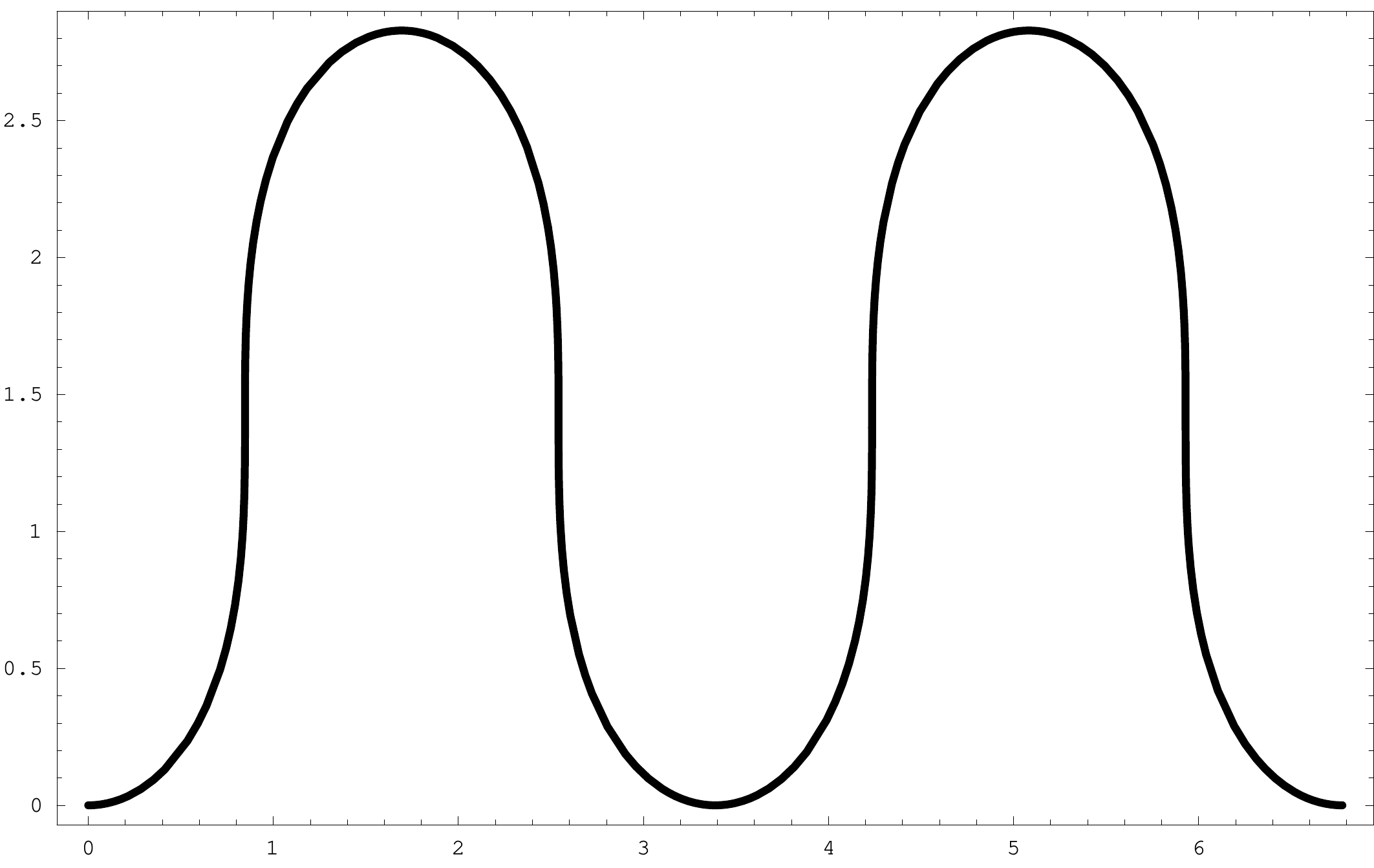}{Прямоугольная эластика}{fig:el3}{0.25}
\twofiglabelsize
{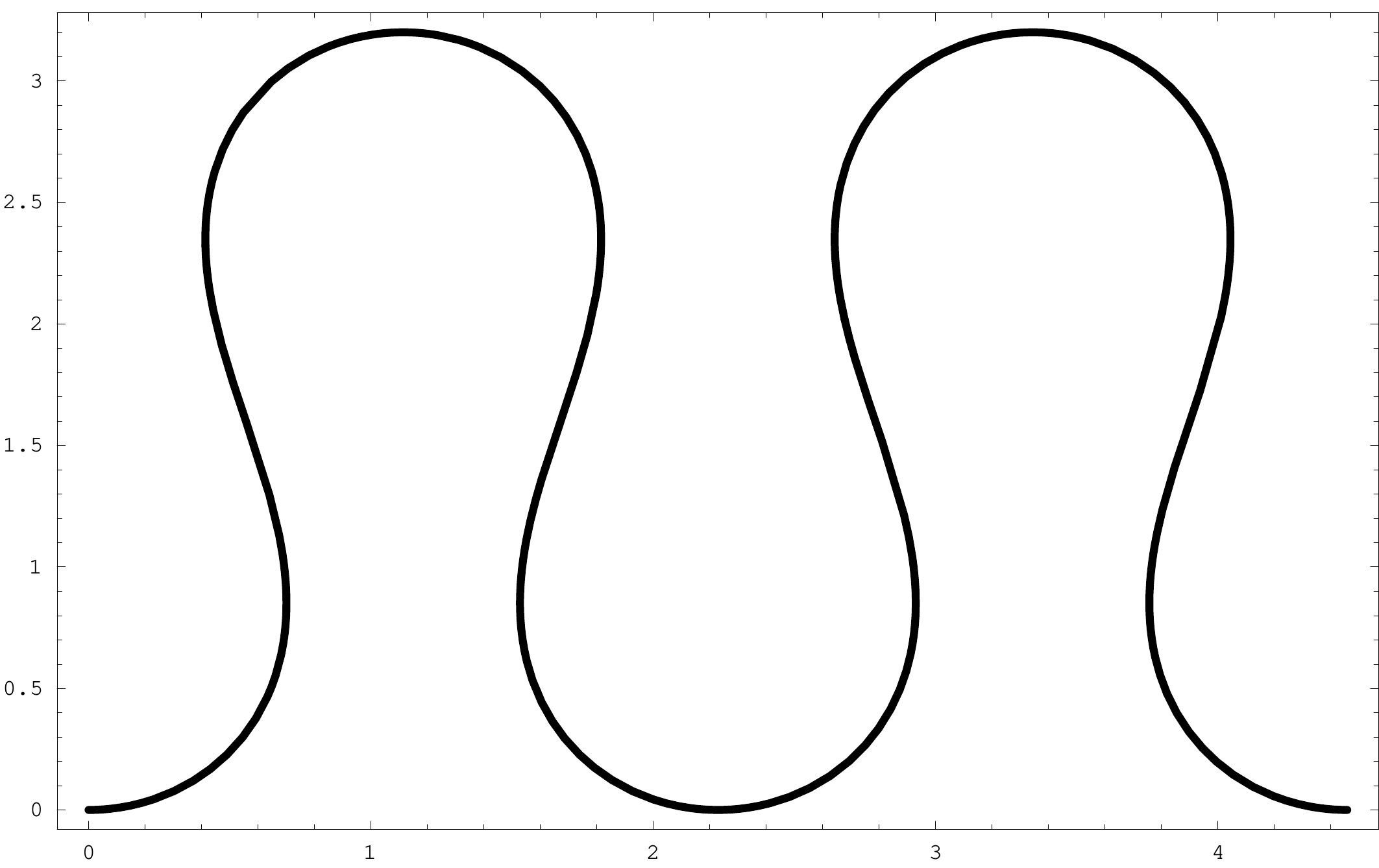}{Инфлексионная эластика}{fig:el4}{0.25}
{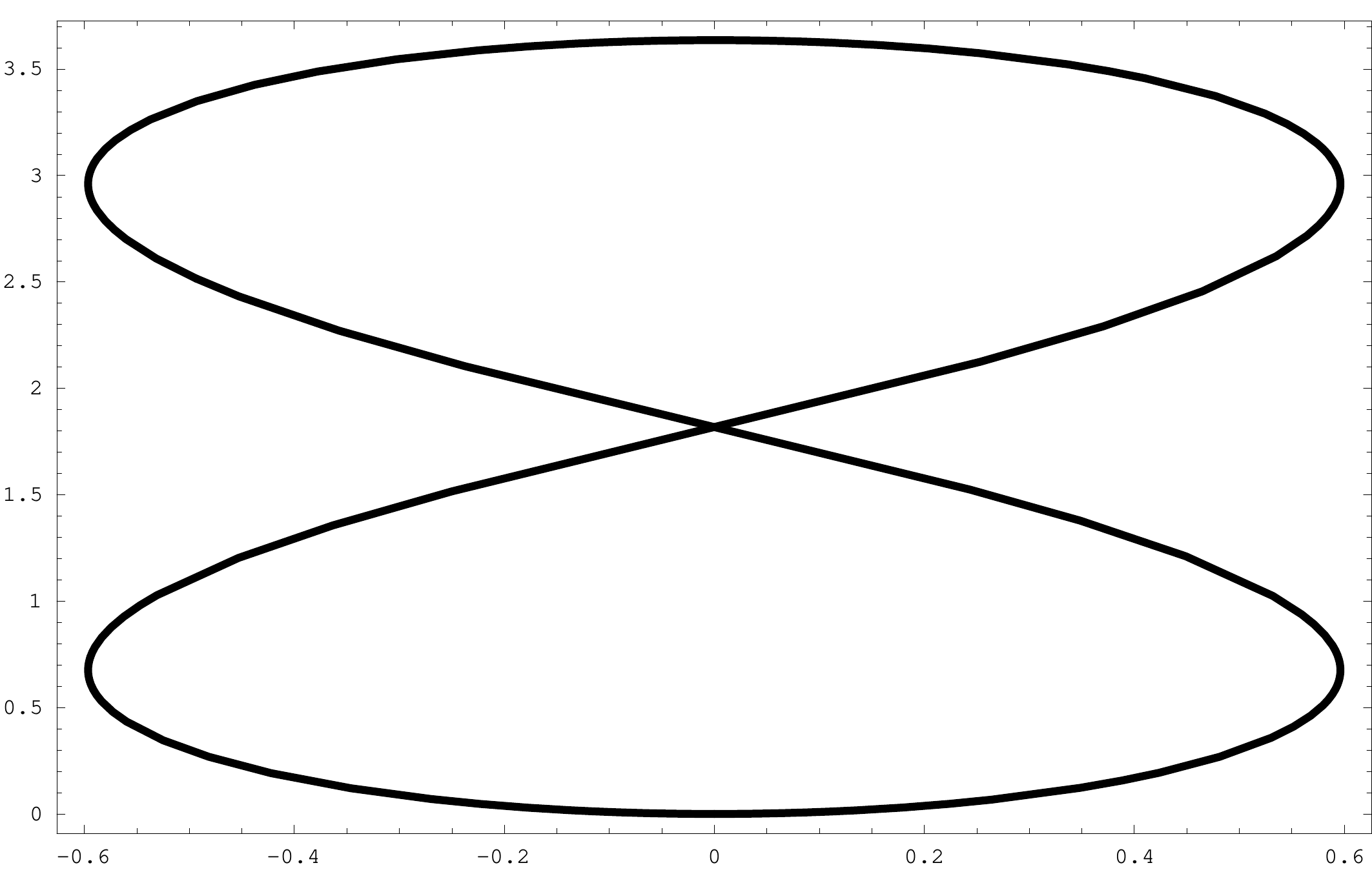}{Эластика-восьмерка}{fig:el5}{0.25}
\twofiglabelsize
{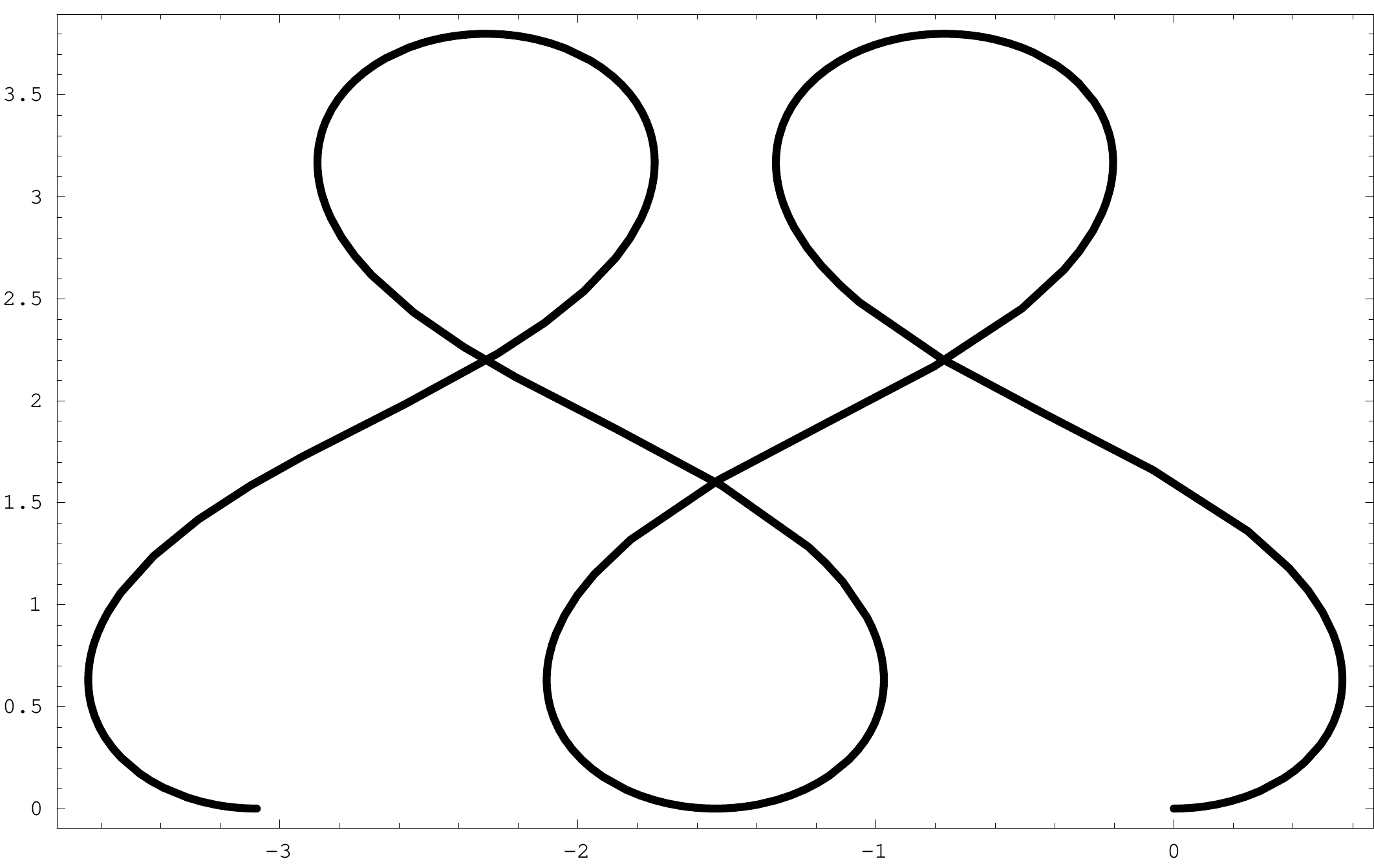}{Инфлексионная эластика}{fig:el6}{0.25}
{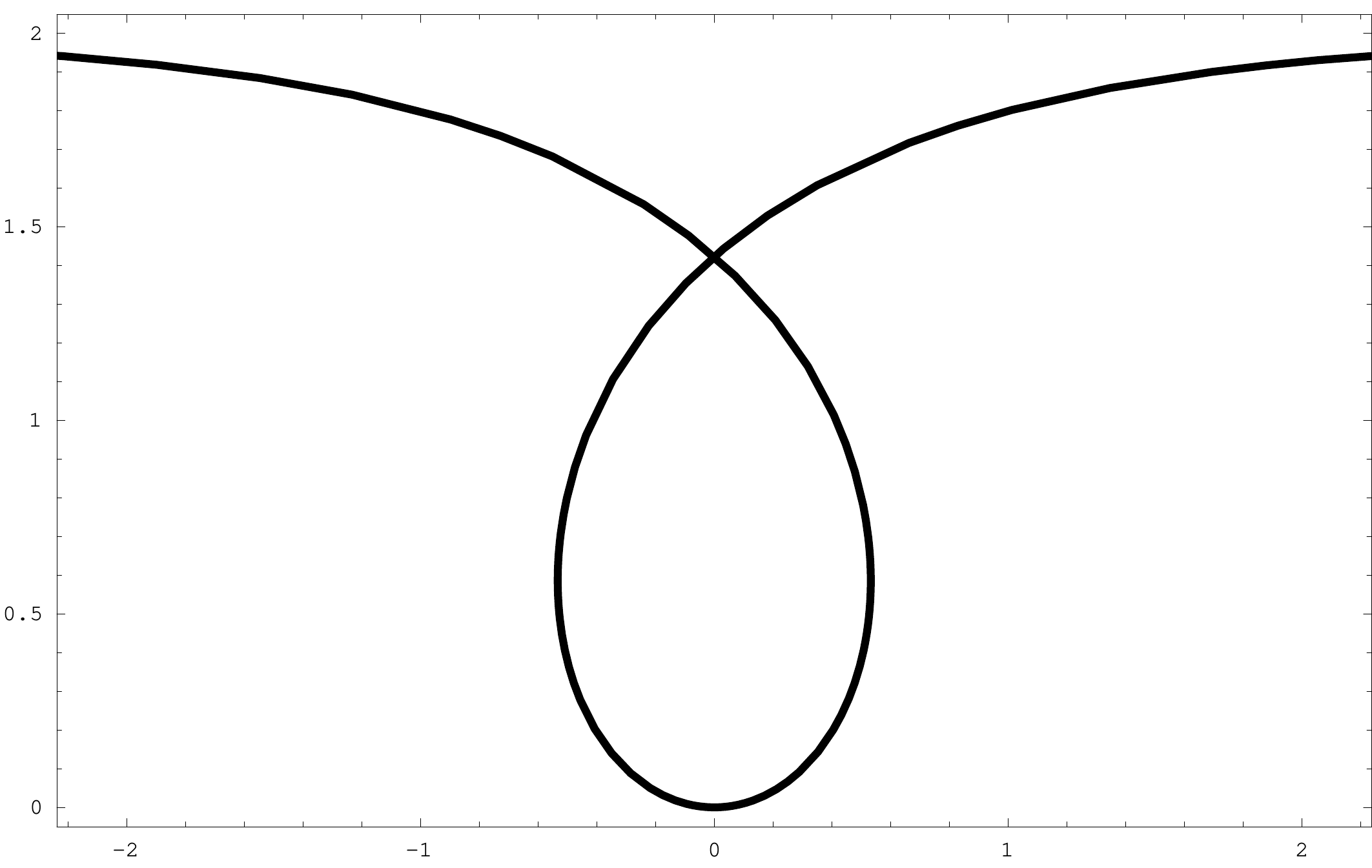}{Критическая эластика}{fig:el7}{0.25}
\onefiglabelsize
{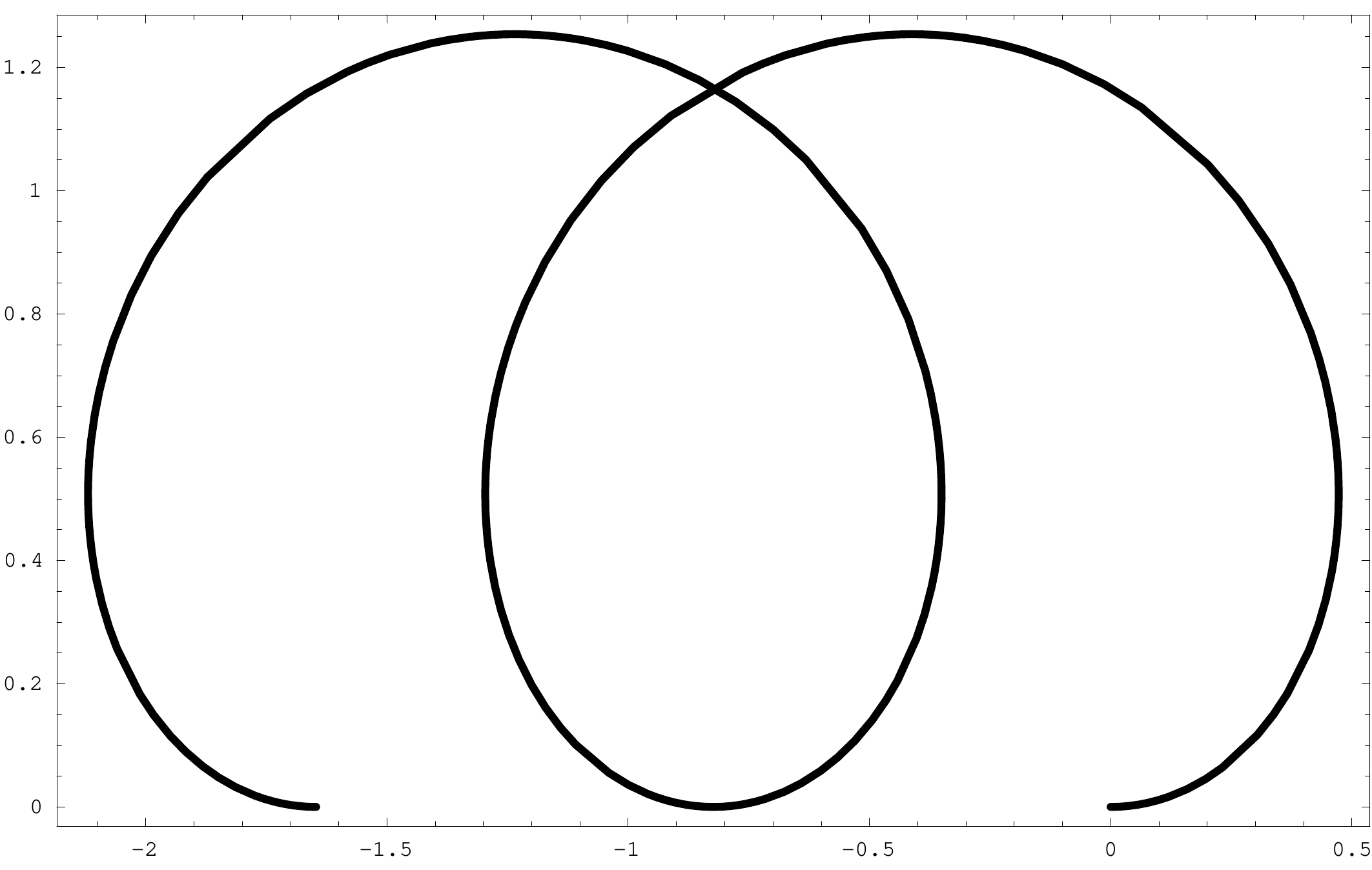}{Неинфлексионная эластика}{fig:el8}{0.35}
}


Периодические движения маятника \eq{el_pend}, \eq{el_pend2}  имеют период
$$
T = \begin{cases}
4 \frac{K(k)}{\sqrt r}, & \lam \in N_1, \\
2 \frac{k K(k)}{\sqrt r}, & \lam \in N_2, \\
\frac{2 \pi}{|c|}, & \lam \in N_6.
\end{cases}
$$

\subsubsection{Симметрии и страты Максвелла}\label{subsubsec:el_max}
Фазовый портрет маятника \eq{el_pend}  сохраняется группой симметрий $\Sym$, порожденной отражением $\eps^1$  в оси $\g$, отражением $\eps^2$  в оси $c$, и отражением $\eps^3$ в начале координат $(\g,c) = (0,0)$:
$$
\Sym = \{ \Id, \eps^1, \eps^2, \eps^3\} \simeq \Z_2 \times \Z_2.
$$
 Эти симметрии естественно продолжаются на прообраз $N = \gg^*$   и образ $G$  экспоненциального отображения $\Exp_t$.   Если $\nu = (\g, c, r) \in N$, то 
$$
\eps^i(\nu) = \nu^i = (\g^i, c^i, r) \in N,
$$
где
\begin{align*}
&(\g^1, c^1) = (\g_t, - c_t), \\
&(\g^2, c^2) = (-\g_t,  c_t), \\
&(\g^3, c^3) = (-\g, - c).
\end{align*}
Если $g = (x,y,\t) \in G$, то $\eps^i(g) = (x^i, y^i, \t^i) \in G$,  где
\begin{align*}
&(x^1, y^1, \t^1) = (x \cos \t + y \sin \t, - x \sin \t + y \cos \t, - \t), \\
&(x^2, y^2, \t^2) = (x \cos \t + y \sin \t, x \sin \t - y \cos \t, \t), \\
&(x^3, y^3, \t^3) = (x, - y, - \t).
\end{align*}

\begin{proposition}
Группа $\Sym = \{ \Id, \eps^1, \eps^2, \eps^3\} $  состоит из симметрий экспоненциального отображения.
\end{proposition}

\begin{theorem}\label{th:el_max}
Первое время Максвелла, соответствующее группе симметрий $\Sym$, для почти всех экстремальных траекторий $g_t = \Exp_t(\lam)$, $\lam \in N$,  выражается следующим образом:
\begin{align*}
&\lam  \in N_1 \then \tmax = \frac{2}{\sqrt r} p_1(k),\\
&p_1(k) = \min(2K(k), p_z^1(k)) = \begin{cases}
2 K(k), & k \in (0, k_0], \\
p_z^1(k), & k \in (k_0, 1), 
\end{cases}\\
&\lam  \in N_2 \then \tmax = \frac{2}{\sqrt r} kK(k),\\
&\lam  \in N_6 \then \tmax = \frac{2\pi}{|c|},\\
&\lam  \in N_3 \cup N_4 \cup N_5 \cup N_7 \then \tmax = + \infty.
\end{align*}
Здесь $p = p_z^1(k) \in (K, 3 K)$  есть первый положительный корень уравнения $\sn p \dn p - (2 \E(p) - p) \cn p = 0$,  а $ k_0 \approx 0,909$  есть корень уравнения $2E(k) - K(k) = 0$. 
\end{theorem}

Имеет место замечание, аналогичное замечанию после теоремы \ref{th:se2_max}, и теорема об инвариантных свойствах функции $\map{\tmax}{N}{(0, + \infty]}$,   аналогичная теореме \ref{th:se2_max_inv}.

\subsubsection{Оценки сопряженного времени}\label{subsubsec:el_conj}
Для эластик Эйлера вопрос локальной оптимальности очень важен с прикладной точки зрения, т.к. локальная оптимальность эластики означает ее устойчивость относительно малых возмущений профиля при закрепленных концах и касательных на концах. С теоретической точки зрения решение этого вопроса важно как шаг в направлении исследования глобальной оптимальности эластик.

\begin{theorem}
Пусть $\lam = (k, \f, r) \in N_1$.  Тогда первое сопряженное время $\tconj(\lam)$  на траектории $\Exp_t(\lam)$  принадлежит отрезку с концами $\ds \frac{4 K(k)}{\sqrt r}$ и  $\ds \frac{2 p_1(k)}{\sqrt r}$, а именно:
\begin{itemize}
\item[$(1)$]
$k \in (0, k_0) \then \tconj \in \left[\ds\frac{4K(k)}{\sqrt r}, \frac{2 p_1^1(k)}{\sqrt r}\right]$,
\item[$(2)$]
$k = k_0 \then \tconj = \ds\frac{4K(k)}{\sqrt r} = \frac{2 p_1^1(k)}{\sqrt r}$,
\item[$(3)$]
$k \in (k_0,1) \then \tconj \in \left[\ds\frac{2 p_1^1(k)}{\sqrt r}, \frac{4K(k)}{\sqrt r}\right]$,
\end{itemize}
где функция $p_1(k)$  определена в теореме {\em \ref{th:el_max}}.
\end{theorem}

\begin{corollary}
Пусть $\lam = (k, \f, r) \in N_1$.  Тогда 
\begin{itemize}
\item[$(1)$]
$k \in (0, k_0) \then \tconj \in [T, t_1^1] \subset [T, 3 T/2 ), \quad t_1^1 = 2 p_1^1/\sqrt r \in (T, 3 T/2)$,
\item[$(2)$]
$k = k_0 \then \tconj = T$, 
\item[$(3)$]
$k \in (k_0, 1) \then \tconj \in [t_1^1, T] \subset (T/2 , T], \quad t_1^1 = 2 p_1^1/\sqrt r \in (T/2, T)$,
\end{itemize}
где $\ds T = \frac{4 K(k)}{\sqrt r}$  есть период колебаний маятников \eq{el_pend}, \eq{el_pend2}.
\end{corollary}

\begin{corollary}\label{cor:el_stab}
Пусть $\lam = (k, \f, r) \in N_1$, $t_1 > 0$, и пусть 
\be{elG}
\G = \{(x_t, y_t) \mid t \in [0, t_1]\}, 
\qquad g(t) = (x_t, y_t, \t_t) = \Exp_t(\lam),
\ee
  есть дуга соответствующей эластики.
\begin{itemize}
\item[$(1)$]
Если дуга $\G$  не содержит точек перегиба, то она локально оптимальна.
\item[$(2)$]
Если $k \in (0, k_0]$  и дуга $\G$  содержит ровно одну точку перегиба, то она локально оптимальна.
\item[$(3)$]
Если дуга $\G$  содержит не менее трех точек перегиба внутри себя, то она не является локально оптимальной.
\end{itemize}
\end{corollary}

Рассмотрим дуги инфлексионных эластик \eq{elG}, центрированные в вершине, т.е. пусть в точке $(x_{t_1/2}, y_{t_1/2})$  достигается локальный экстремум кривизны эластики. 
Примеры таких дуг см. на Рис. \ref{fig:el_cent_v}, \ref{fig:el_cent_v2}.

\figout{
\twofiglabelsize
{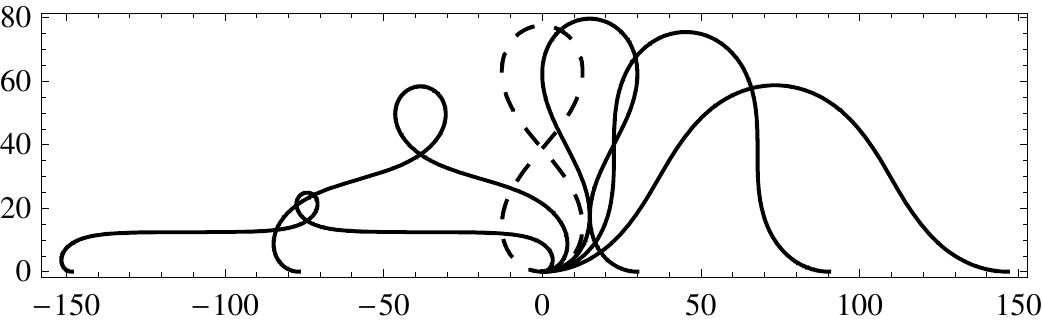}{Эластики, центрированные в вершине}{fig:el_cent_v}{0.18}
{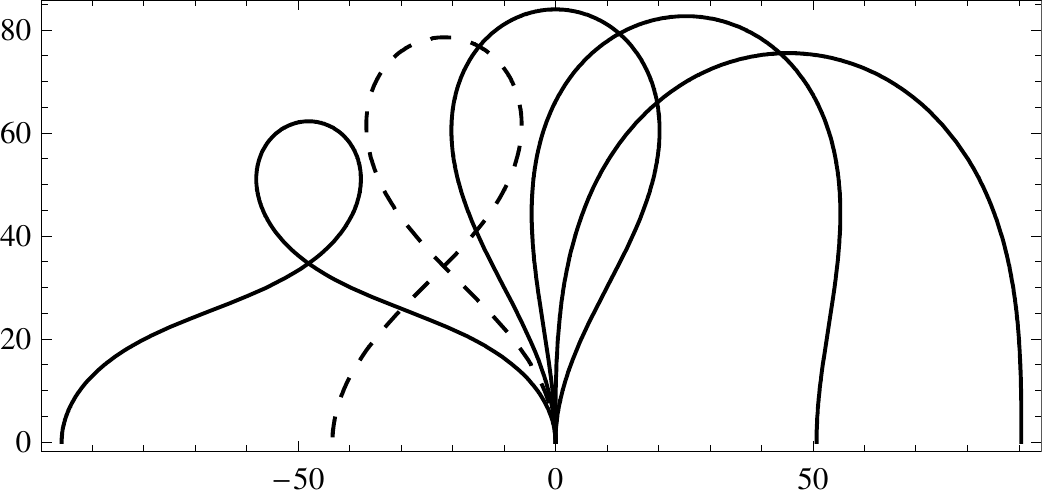}{Эластики, центрированные в вершине}{fig:el_cent_v2}{0.18}
}

Обозначим $t_1^1 = \frac{2}{\sqrt r} p_1(k)$, где функция  $p_1(k)$  определена в теореме \ref{th:el_max}.

\begin{theorem}
\label{th:vertex}
Пусть инфлексионная эластика $\G$   центрирована в вершине.
\begin{itemize}
\item[$(1)$]
Если $t < t_1^1$, то эластика $\G$  устойчива.
\item[$(2)$]
Если $t = t_1^1$, то  конец эластики $\G$  является первой сопряженной точкой.
\item[$(3)$]
Если $t > t_1^1$, то эластика $\G$  неустойчива.
\end{itemize}
\end{theorem}

Рассмотрим дуги инфлексионных эластик \eq{elG}, центрированные в точке перегиба, т.е. пусть в точке $(x_{t_1/2}, y_{t_1/2})$    эластика имеет нулевую кривизну. 
Примеры таких дуг см. на Рис. \ref{fig:el_cent_i}.
 
\figout{
\onefiglabelsize
{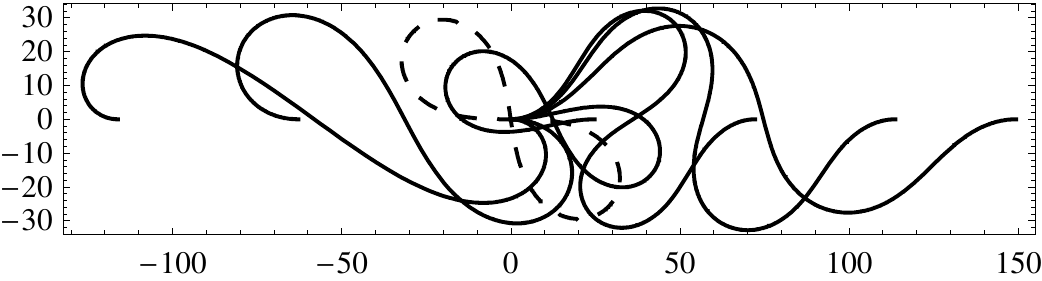}{Эластики, центрированные в точке перегиба}{fig:el_cent_i}{0.55}
}

\begin{theorem}
\label{th:inflex}
Пусть   эластика $\G$  центрирована в точке перегиба. Пусть 
также $k \in (0, k_0]$.
\begin{itemize}
\item[$(1)$]
Если $t < T$, то эластика $\G$  устойчива.
\item[$(2)$]
Если $t = T$, то    конец эластики $\G$ является первой сопряженной точкой.
\item[$(3)$]
Если $t > T$, то эластика $\G$  неустойчива.
\end{itemize}
\end{theorem}

\begin{theorem}
Пусть $\lam \in N_2 \cup N_3 \cup N_6$.  Тогда экстремальная траектория $g(t) = \Exp_t(\lam)$  не содержит сопряженных точек при $t > 0$.
\end{theorem}

Итак, если дуга эластики не содержит точек перегиба, то она устойчива; если она содержит не менее трех точек перегиба внутри себя, то она неустойчива. Если есть одна или две точки перегиба, то эластика может быть устойчивой или неустойчивой.

\subsubsection{Диффеоморфная структура экспоненциального отображения}\label{subsubsec:el_diff}
Пусть $t_1 = 1$, $\Exp = \Exp_1$, 
$$
\A = \A_1 = \{(x,y,\t) \in G \mid x^2 + y^2 < 1  \text{  или }  (x,y,\t) = (1,0,0)\}.
$$
Случай общего $t_1 > 0$  сводится к частному случаю $t_1 = 0$  гомотетиями плоскости $(x,y)$:
$$
(x,y,\t,t,u,t_1,J) \mapsto (\tilde x, \tilde y, \tilde \t, \tilde t, \tilde u, \tilde t_1, \tilde J) = (e^sx, e^sy, \t, e^st, e^{-s}u, e^st_1, e^{-s}J).
$$

Рассмотрим подмножество в $\A$, не содержащее неподвижных точек отражений $\eps^1$, $\eps^2$:
\begin{align*}
&\tG = \{g \in \A \mid \eps^i (g) \neq g, \ i = 1, 2\} = 
\left\{g \in \A \mid \sin \frac{\t}{2} P(g) \neq 0 \right\}, \\
&P(g) = x \sin \frac{\t}{2} - y \cos \frac{\t}{2},
\end{align*}
и его разбиение на компоненты связности
\begin{align*}
&\tG = G_+ \sqcup G_-, \\
&G_{\pm} = \{ g \in G \mid \t \in (0, 2 \pi), \ x^2 + y^2 < 1, \ \sgn P(g) = \pm 1 \}.
\end{align*}
Также рассмотрим открытое плотное подмножество в пространстве всех потенциально оптимальных экстремальных траекторий:
$$
\tN = \{ \lam \in \cup_{i=1}^3 N_i \mid t_1 < \tmax(\lam), \ c_{t_1/2} \sin \g_{t_1/2} \neq 0 \},
$$
 и его связные компоненты
\begin{align*}
&\tN = \sqcup_{i=1}^4 D_i, \\
&D_1 = \{ \lam \in \cup_{i=1}^3 N_i \mid t_1 < \tmax(\lam), \ c_{t_1/2} > 0, \  \sin \g_{t_1/2} > 0 \}, \\
&D_2 = \{ \lam \in \cup_{i=1}^3 N_i \mid t_1 < \tmax(\lam), \ c_{t_1/2} < 0, \  \sin \g_{t_1/2} > 0 \}, \\
&D_3 = \{ \lam \in \cup_{i=1}^3 N_i \mid t_1 < \tmax(\lam), \ c_{t_1/2} < 0, \  \sin \g_{t_1/2} < 0 \}, \\
&D_4 = \{ \lam \in \cup_{i=1}^3 N_i \mid t_1 < \tmax(\lam), \ c_{t_1/2} > 0, \  \sin \g_{t_1/2} < 0 \}.
\end{align*}

\begin{theorem}
Следующие отображения являются диффеоморфизмами:
$$
\map{\Exp}{D_1}{G_+}, \quad
\map{\Exp}{D_2}{G_-}, \quad
\map{\Exp}{D_3}{G_+}, \quad
\map{\Exp}{D_4}{G_-}.
$$
\end{theorem}

\begin{corollary}
Отображение $\map{\Exp}{\tN}{\tG}$ есть двулистное накрытие.
\end{corollary}

\subsubsection{Оптимальные эластики для различных граничных условий}\label{subsubsec:el_opt}
\paragraph{Граничные условия общего положения}
Если $ g_1 \in G_+$, то существует единственная пара $(\lam_1, \lam_3) \in D_1 \times D_3$, для которой $\Exp(\lam_1) = \Exp(\lam_3) = g_1$.  Оптимальная траектория находится среди траекторий $q^1(t) = \Exp_t(\lam_1)$ и $q^3(t) = \Exp_t(\lam_3)$, $t \in [0, 1]$.  Для отыскания оптимальной траектории необходимо взять ту из них, для которой функционал качества $J[q^i(\cdot)] = \frac 12 \int_0^1 (c_t^i)^2 dt $  принимает меньшее значение. Если $J[q^1(\cdot)] = J[q^3(\cdot)]$,  то оптимальны обе траектории, этот случай изображен на Рис.  \ref{fig:el_opt2}. 

\figout{
\onefiglabelsize{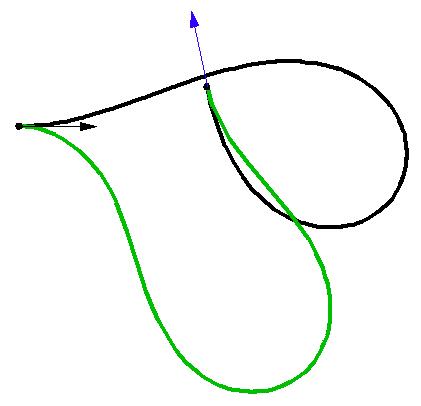}{Две оптимальные несимметричные эластики с одинаковыми граничными условиями}{fig:el_opt2}{0.25}
}

Если $g_1 \in G_-$, то оптимальные траектории выбираются аналогично среди соответствующих ковекторам $\lam_2 \in D_2$  и $\lam_4 \in D_4$, для которых $\Exp(\lam_2) = \Exp(\lam_4) = g_1$.

\paragraph{Случай $(x_1, y_1, \t_1) = (1, 0, 0)$}
Оптимальная эластика есть отрезок $(x,y) = (t, 0)$, $t \in [0, 1]$.

\paragraph{Случай $x_1>0$, $y_1=0$, $\t_1=\pi$}
В этом случае $g_1 \in G_+$  и уравнение $\Exp(\lam) = g_1$, $\lam \in \tG$,  имеет два корня $\lam_1 \in D_1$  и  $\lam_3 \in D_3$.  Траектории $q^1(t) = \Exp_t(\lam_1)$ и $q^3(t) = \Exp_t(\lam_3)$  имеют одинаковое значение функционала $J$, поэтому оптимальны. Соответствующие оптимальные инфлексионные эластики симметричны относительно оси $x$,  см. Рис.  \ref{fig:elx+0pi}.

\figout{
\twofiglabelsize
{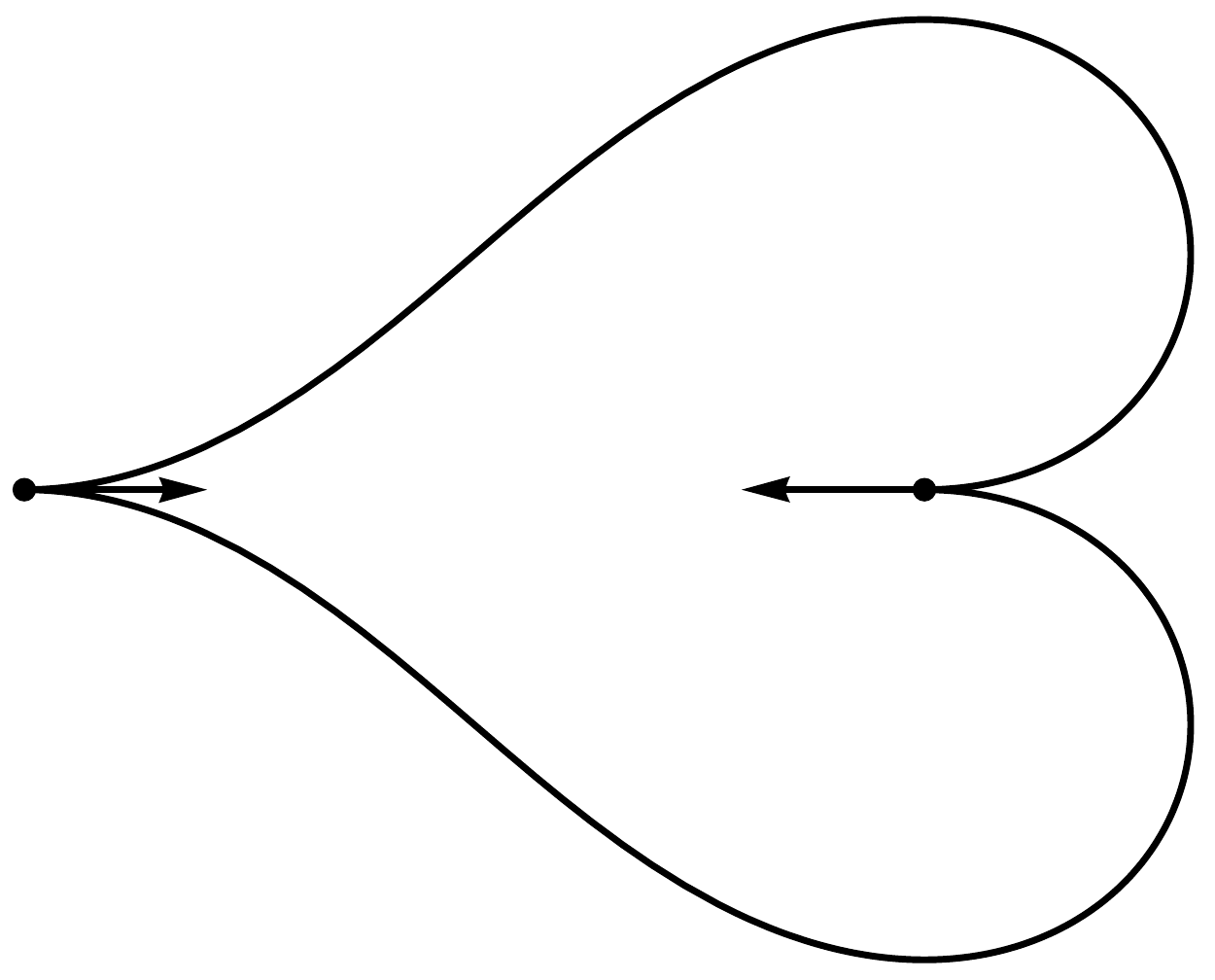}{Оптимальные эластики для $x_1>0$, $y_1=0$, $\t_1=\pi$}{fig:elx+0pi}{0.3}
{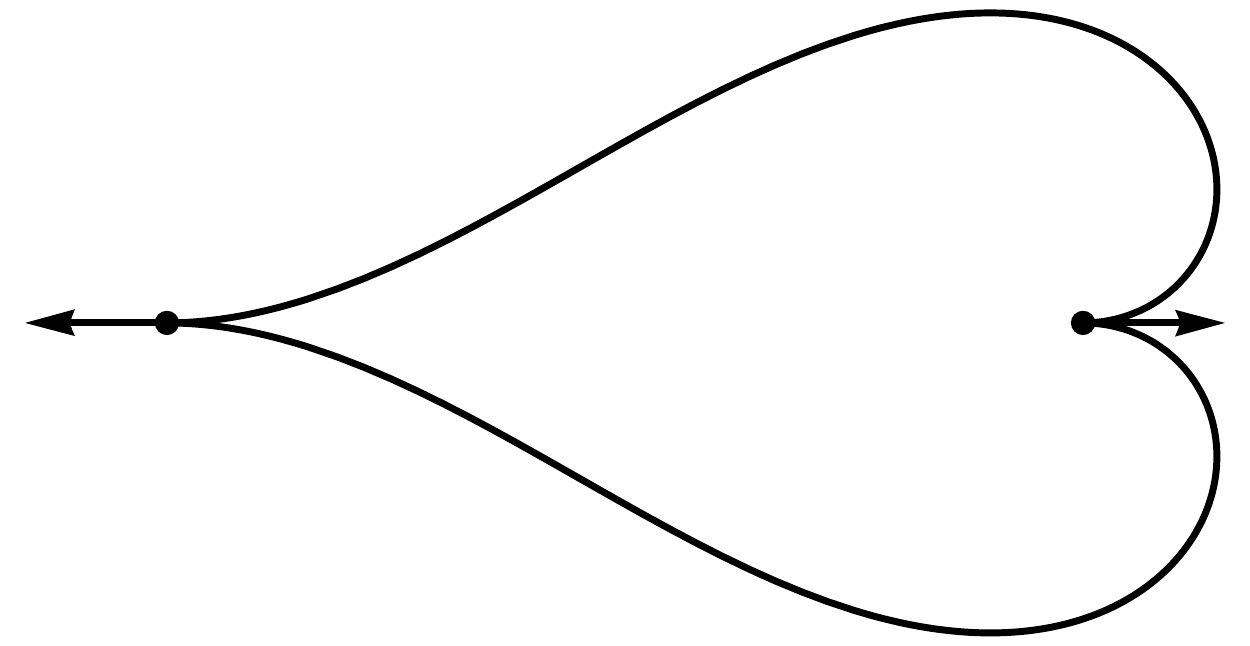}{Оптимальные эластики для $x_1<0$, $y_1=0$, $\t_1=\pi$}{fig:elx-0pi}{0.2}
}

\paragraph{Случай $x_1<0$, $y_1=0$, $\t_1=\pi$}
Этот случай аналогичен предыдущему случаю, 
см. Рис.  \ref{fig:elx-0pi}.

\paragraph{Случай $x_1=0$, $y_1=0$, $\t_1=\pi$}
Единственная оптимальная эластика -<<капля>> определяется параметрами $\lam = (\f, k, r) \in N_1$, $\f = \frac{\tau}{2p}
- \frac 12$, $r = 4 p^2$, $\sn \tau = 0$, $1 - 2 k^2 \sn^2 p - 0$, $2 \E(p) - p = 0$,  см. Рис.  \ref{fig:elx00pi}.

\figout{
\onefiglabelsize
{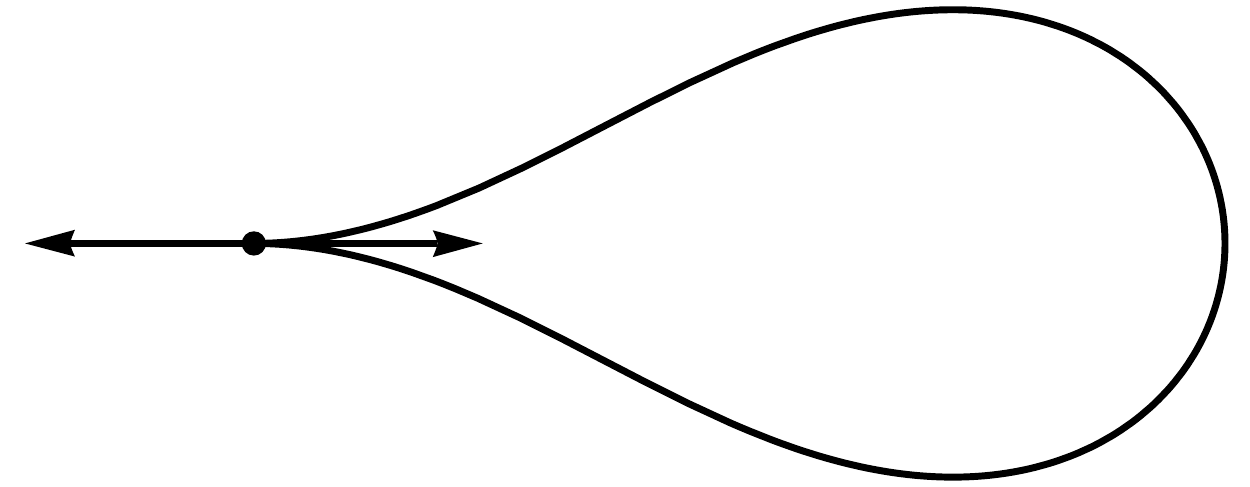}{Оптимальная эластика -<<капля>> для $x_1=0$, $y_1=0$, $\t_1=\pi$}{fig:elx00pi}{0.4}
}

\paragraph{Случай $x_1=0$, $y_1=0$, $\t_1=0$}
Существуют две оптимальные эластики --- окружности, симметричные относительно оси $x$.

\paragraph{Случай $x_1>0$, $y_1=0$, $\t_1=0 $}
Имеются две или четыре оптимальных эластики; существует такое $x_* \in (0,4, \, 0,5)$,  что:
\begin{itemize}
\item
если $x_1 \in (0, x_*)$, то имеются две оптимальные неинфлексионные эластики, см. Рис. \ref{fig:el_x<x*},
\item
если $x_1 \in (x_*, 1)$, то имеются две оптимальные инфлексионные эластики, см. Рис. \ref{fig:el_x>x*},
\item
если $x_1 = x_* $, то существуют четыре оптимальные эластики (две инфлексионные  и две неинфлексионные), см. Рис. \ref{fig:el_x=x*},
\end{itemize}

\figout{
\twofiglabelsize
{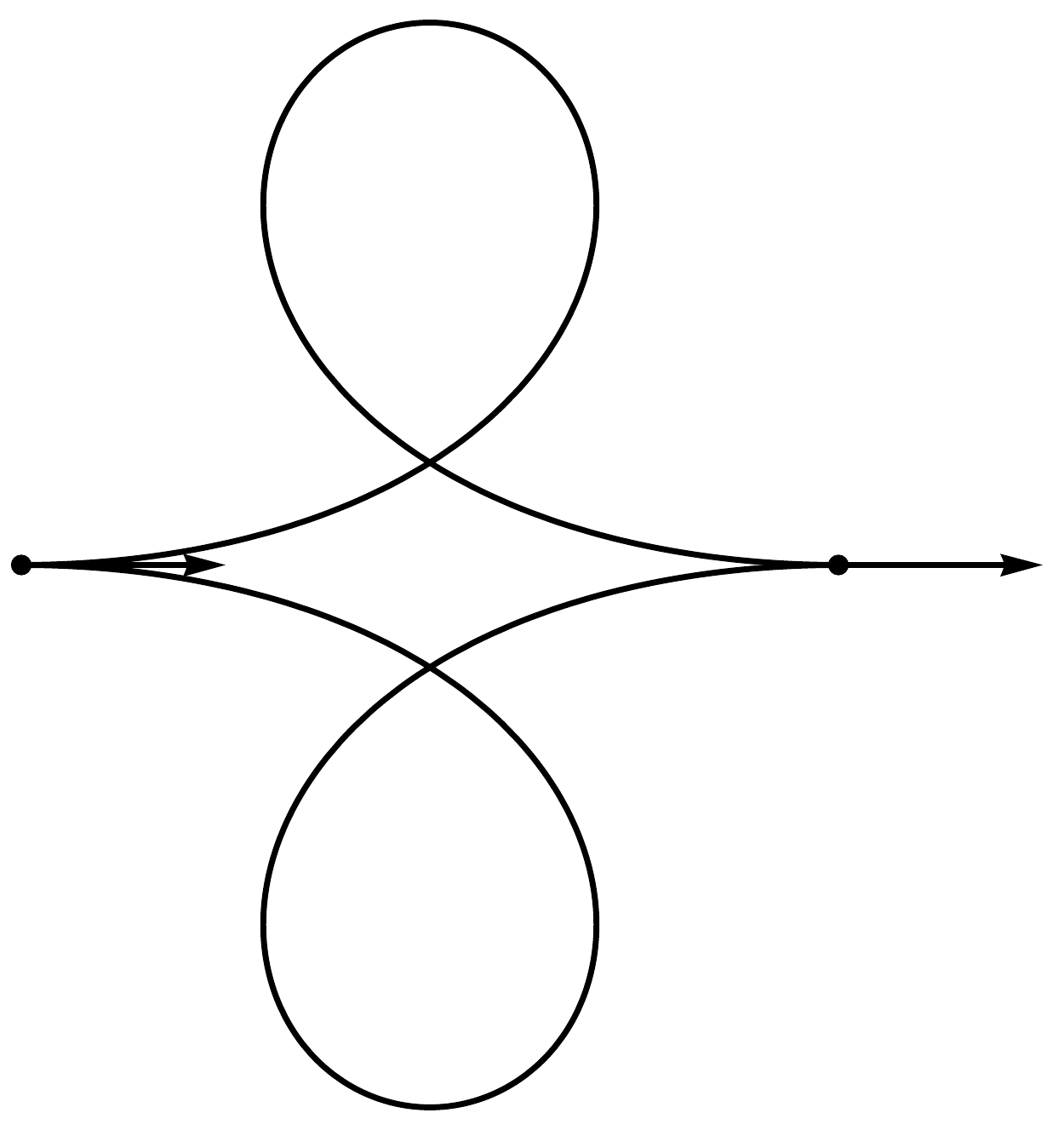}{Оптимальные эластики для $x_1>0$, $y_1=0$, $\t_1=0 $,  $x_1 \in (0, x_*)$}{fig:el_x<x*}{0.3}
{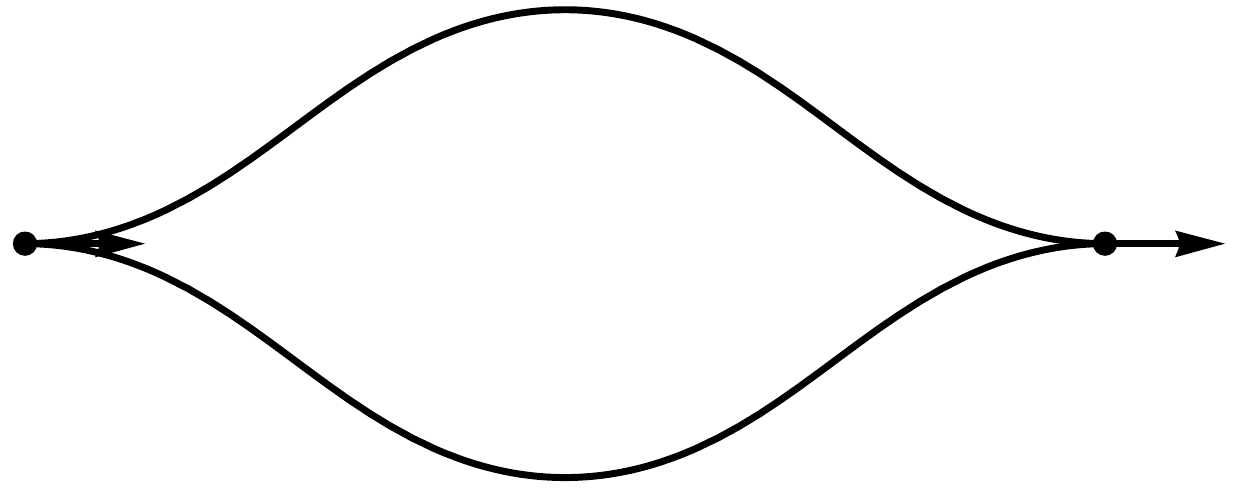}{Оптимальные эластики для $x_1>0$, $y_1=0$, $\t_1=0 $,  $x_1 \in (x_*, 1)$}{fig:el_x>x*}{0.15}

\onefiglabelsize
{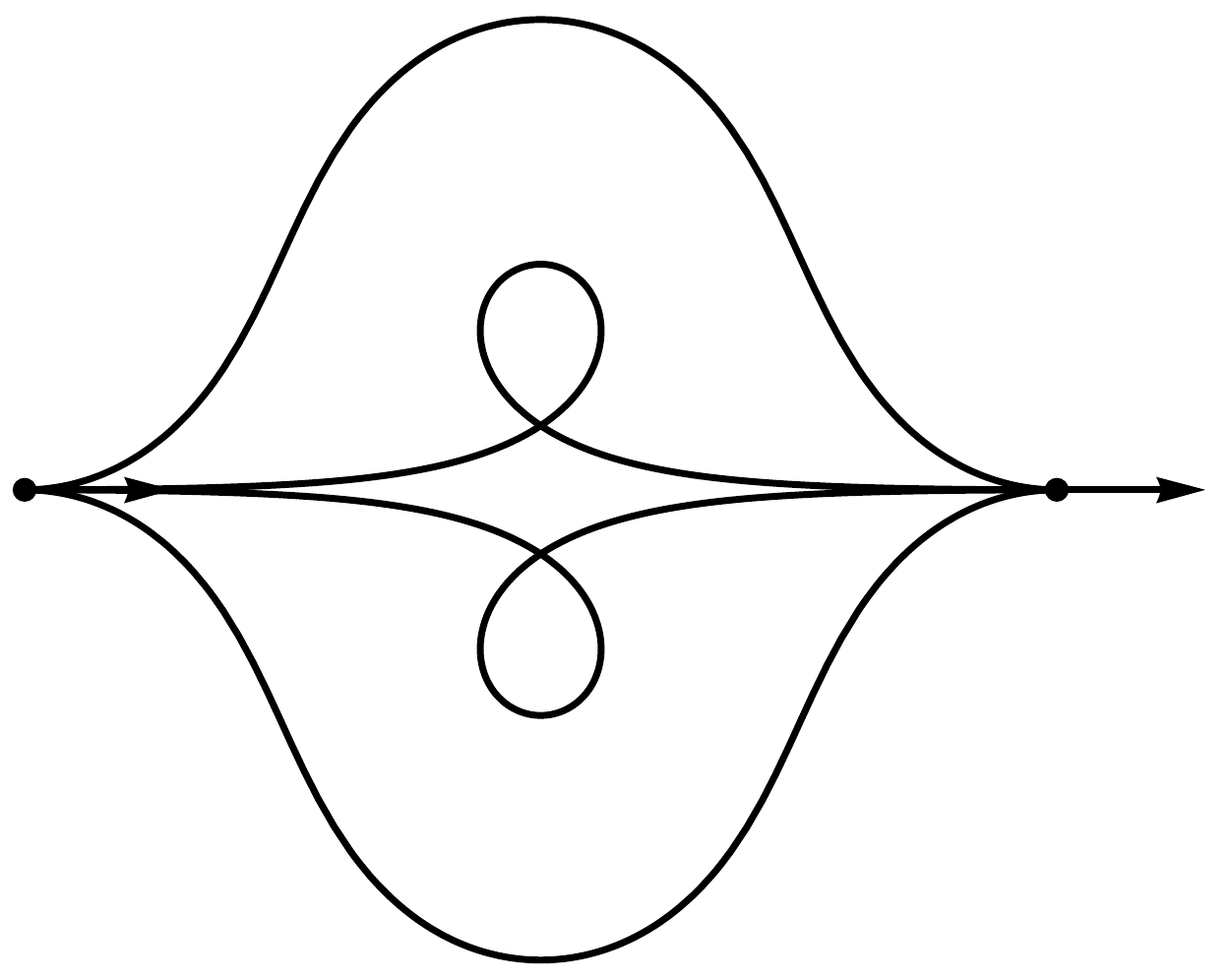}{Оптимальные эластики для $x_1>0$, $y_1=0$, $\t_1=0 $,  $x_1 = x_*$}{fig:el_x=x*}{0.3}
}

\paragraph{Случай $x_1<0$, $y_1=0$, $\t_1=0 $}
Существуют две оптимальные неинфлексионные эластики, см. Рис.  \ref{fig:el_x<000}.

\figout{
\onefiglabelsize
{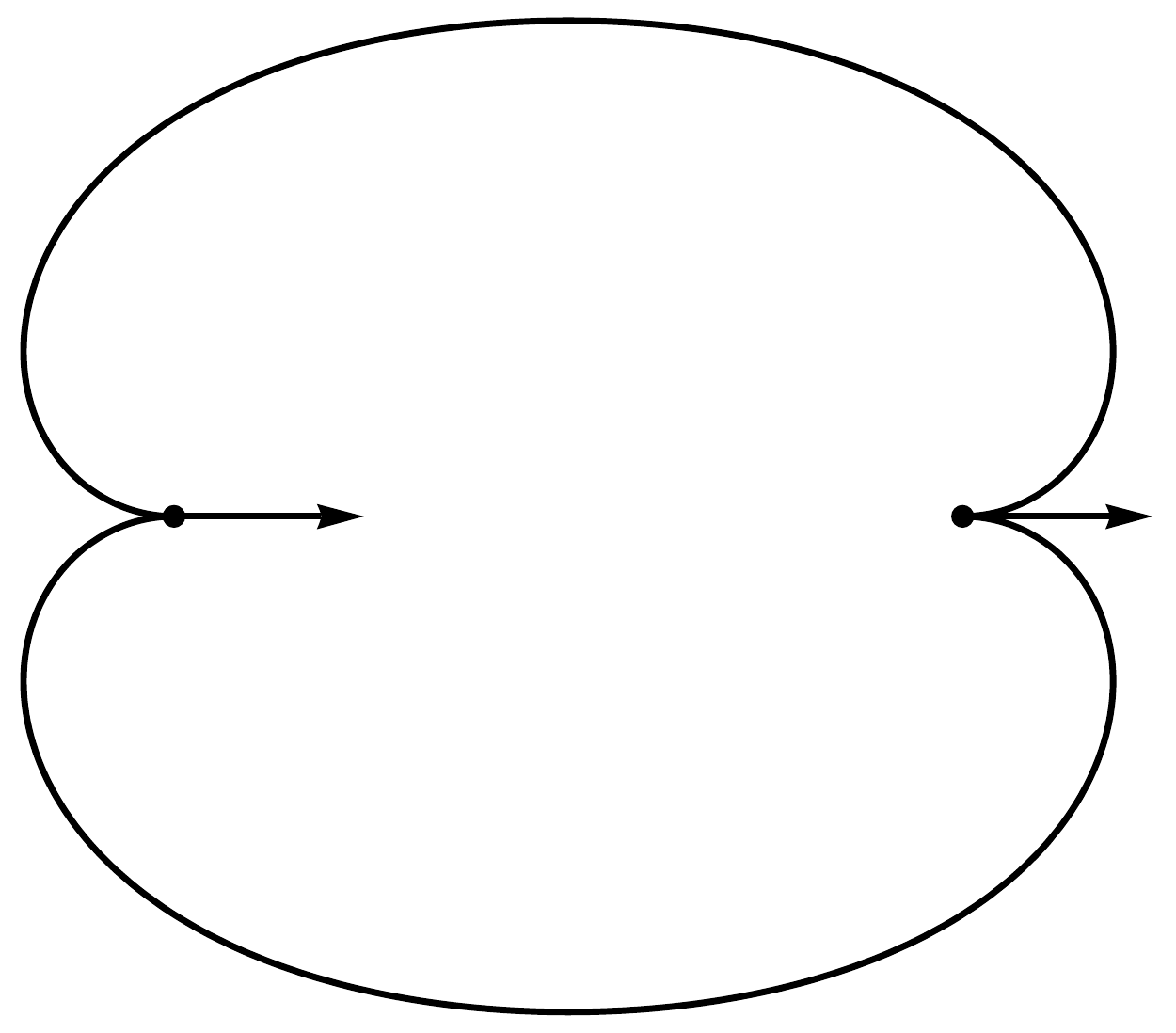}{Оптимальные эластики для $x_1<0$, $y_1=0$, $\t_1=0 $}{fig:el_x<000}{0.3}
}

\subsubsection{Библиографические  комментарии}
Раздел  \ref{subsubsec:el_history}  по истории задачи об эластиках опирается на классические источники \cite{truesdell, love, timoshenko}. Имеется также замечательное описание \cite{levien}  этой истории.

Заметим, что задача об эластиках долгое время представляла лишь теоретический интерес и служила одним из примеров приложения теории эллиптических функций (см., например, \cite{Greenhill, love}). В связи с широким внедрением стали в практику проектирования и появлением гибких тонкостенных конструкций, стимулировавшим развитие теории устойчивости деформируемых систем, решение задачи об эластиках стало приобретать практическое значение. Возникли, в частности, важные для инженерных приложений вопросы: каково поведение сжатой стойки при нагрузках, превышающих эйлерово критическое значение, какова при этом форма стойки, единственна ли эта форма и устойчива ли она? 
Решению этих вопросов посвящены многочисленные исследования~\cite{Frisch, Bisshopp, Lardner, Wang, Panayotounakos, Seide, Naschie, Stampouloglou}, где рассматривались различные условия опирания и нагружения гибких нерастяжимых стержней. В последние десятилетия интерес к эластикам возрос в связи с применением теории гибких стержней к анализу микро- и наноструктур в биологии и нанотехнологиях \cite{Glassmaker, Tang, Mikata, Heijden}.  Подтверждено существование множественных форм равновесия при фиксированной нагрузке. 

Разделы \ref{subsubsec:el-state}--\ref{subsubsec:el_max} опираются на работу \cite{el_max}, 
разделы \ref{subsubsec:el_conj} --- на работы \cite{el_conj} и \cite{el_stable}, 
разделы \ref{subsubsec:el_diff} и  \ref{subsubsec:el_opt} --- на работы \cite{el_cut} и \cite{el_closed}.

\subsection{Левоинвариантная субриманова  задача общего вида на группе $\SO(3)$}  \label{subsec:so3_gen}

\subsubsection{Постановка задачи }
Из классификации контактных левоинвариантных субримановых структур на трехмерных группах Ли (см. раздел \ref{subsec:class3}) следует, что для произвольной такой структуры на группе $G = \SO(3)$  можно выбрать ортонормированный репер $(X_1, X_2)$ с таблицей умножения
\be{so3_gen_tab}
[X_2, X_1] = X_3, \quad [X_1, X_3] = (\kappa + \chi) X_2,  \quad [X_2, X_3] = (\chi - \kappa) X_1,  
\ee
 где $\kappa \geq \chi \geq 0$  суть дифференциальные инварианты субримановой структуры. Равномерное растяжение полей $(X_1, X_2)$  пропорционально изменяет функцию расстояния и оба инварианта $\kappa$  и $\chi$. В разделе \ref{subsec:class3}  использована нормализация $\kappa^2 + \chi^2 = 1$.  В этом разделе удобнее принять $\kappa + \chi = 1$ и использовать инвариант $ a = \sqrt{2 \chi} \in [0, 1)$. Случай $a = 0$  соответствует осесимметричной субримановой структуре, рассмотренной в разделе \ref{subsec:sr_so3}.

Следующие векторные поля удовлетворяют таблице умножения \eq{so3_gen_tab}:
$$
X_1(g)  = L_{g*} A_2, \quad  X_2(g)  =  \sqrt{1-a^2} L_{g*} A_1, \quad X_3(g)  =  \sqrt{1-a^2} L_{g*} A_3,
$$
где базис $A_1, A_2, A_3$  алгебры Ли $\gg = \so(3)$  имеет вид
\be{s03_tab}
A_1 = \begin{pmatrix}
0 & 0 & 0 \\
0 & 0 & -1 \\
0 & 1 & 0
\end{pmatrix},
\quad
A_2 = \begin{pmatrix}
0 & 0 & 1 \\
0 & 0 & 0 \\
-1 & 0 & 0
\end{pmatrix},
\quad
A_3                                                = \begin{pmatrix}
0 & -1 & 0 \\
1 & 0 & 0 \\
0 & 0 & 0
\end{pmatrix}.
\ee

\subsubsection{Параметризация геодезических}
Анормальные экстремальные траектории постоянны.

Для  параметризации нормальных геодезических введем гамильтонианы 
$h_i(\lam) = \langle \lam, X_i(g)\rangle$, $i = 1, 2, 3$, $H = \frac 12 (h_1^2 + h_2^2)$. 
Натурально параметризованные экстремали параметризуются точками цилиндра $C = \gg^* \cap \{H = \frac 12\}$. Введем на этом цилиндре координаты $(\p, c)$:
$$
h_1 = \cos \p, \quad h_2 = - \sin \p, \quad h_3 = c.
$$

Нормальная гамильтонова система ПМП имеет вид
\begin{align}
&\dot h_1 = h_2 h_3, \quad \dot h_2 = - h_1 h_3, \quad \dot h_3 = a^2 h_1 h_2, \label{so3gen_ham} \\
&\dot g = h_1 X_1(g) + h_2 X_2(g). \label{so3gen_hor}
\end{align}
Вертикальная подсистема \eq{so3gen_ham}  задает на цилиндре $C$  уравнение маятника
\be{so3gen_pend}
\dot \p = c, \quad \dot c = - \frac{a^2}{2} \sin 2 \p.
\ee
Этот цилиндр имеет стратификацию
$$
C = \sqcup_{i=1}^5 C_i
$$
 на инвариантные множества системы \eq{so3gen_pend}, которые определяются значением полной энергии маятника $E = 2 c^2 - a^2 \cos 2 \p$:
\begin{align*}
C_1 &= \{\lambda\in C: E\in(-a^2,a^2)\},& &(\text{область внутри сепаратрис}), \\
C_2 &= \{\lambda\in C: E\in(a^2, +\infty)\}, & &(\text{область вне сепаратрис}),\\
C_3 &= \{\lambda\in C: E = a^2, c \neq 0\},& &(\text{сепаратрисы}),\\
C_4 &= \{\lambda\in C: E = -a^2\}, &  &(\text{устойчивое положение равновесия}), \\
C_5 &= \{\lambda\in C: E = a^2, c = 0\},& &(\text{неустойчивое положение равновесия}).
\end{align*}
Введем на множествах $C_1$, $C_2$, $C_3$  координаты $(\t, k)$, выпрямляющие уравнение маятника  \eq{so3gen_pend}.
В области~$C_1$:
\begin{align*}
\sin\psi &= s_1k\sn\left(a\theta,k \right), & &s_1 = \sgn\left( \cos\psi \right),\\
\cos\psi &= s_1\dn\left(a\theta,k \right), & &k = \sqrt{\frac{E + a^2}{2a^2}}\in(0,1),\\
c &= ak\cn\left(a\theta,k \right), & &\theta\in[0,4K(k)/a].
\end{align*}
В области $C_2$:
\begin{align*}
\sin\psi &= s_2\sn\left(\frac{a\theta}{k},k\right), & &s_2 = \sgn(c),\\
\cos\psi &= \cn\left(\frac{a\theta}{k},k\right), & &k = \sqrt{\frac{2a^2}{E + a^2}}\in (0,1),\\
c &= \frac{s_2a}{k}\dn\left(\frac{a\theta}{k},k\right), & &\theta \in[0,4kK(k)/a].
\end{align*}
На множестве $C_3$:
\begin{align*}
\sin\psi &= s_1s_2\th a\theta,\\
\cos\psi &= \frac{s_1}{\ch a\theta},\\
c &= \frac{s_2a}{\ch a\theta}, \\
\theta &\in(-\infty,+\infty), \quad
k = 1.
\end{align*}

Тогда при $(\p_0, c_0) \in C_1 \cup C_2 \cup C_3$  решение уравнения маятника есть $\t(t) = t + \t_0$, $k \equiv \const$. При $(\p_0, c_0) \in C_4$  имеем $\p \equiv \pi n$, $n \in \Z$, $c = 0$,  а при $(\p_0, c_0) \in C_5$  имеем $\p \equiv - \frac{\pi}{2} + \pi n$, $n \in \Z$, $c = 0$.

Для параметризации решений горизонтальной подсистемы \eq{so3gen_hor}  представим их с помощью углов Эйлера
$$
g_t = \exp(-\f_1(0) A_3) \exp(-\f_2(0) A_1) \exp(\f_3(t) A_3)  \exp(\f_2(t) A_1) \exp(\f_1(t) A_3).  
$$
Тогда
\begin{align}
&\cos \f_2 = \frac{c}{\sqrt M}, & \sin \f_2 = \sqrt{\frac{M - c^2}{M}}, \label{so3gen_f2}\\
&\cos \f_1 = \frac{h_1 \sqrt{1-a^2}}{\sqrt {M- c^2}}, & \sin \f_1 = \sqrt{\frac{h_2}{M - c}}, \label{so3gen_f1}
\end{align}
где $M = h_2^2 + (1-a^2) h_1^2 + c^2$  есть первый интеграл подсистемы \eq{so3gen_ham}.

Угол $\f_3$ удовлетворяет уравнению
\be{so3gen_f3}
\dot\f_3 = \frac{\sqrt{M (1-a^2)}}{M - c^2} = \frac{\sqrt{M(1-a^2)}}{1 - a^2 h_1^2}
\ee
 и является монотонной функцией времени т.к. $0 < \sqrt{M(1-a^2)} \leq \dot \f_3 \leq \sqrt{M/(1-a^2)}$. Решения этого уравнения имеют вид:
\begin{enumerate}
\item в $C_1$:
\begin{align*}
\f_3 &= \sqrt{\frac{1-a^2(1-k^2)}{a^2(1-a^2)}}\left( \Pi\left( \frac{a^2k^2}{a^2-1}; \am(a\theta,k),k \right)\right. 
- \left. \Pi\left( \frac{a^2k^2}{a^2-1}; \am(a\theta_0,k),k \right) \right);
\end{align*}
\item в $C_2$:
\begin{align*}
\f_3 
&= \sqrt{\frac{k^2+a^2(1-k^2)}{a^2(1-a^2)}}\left( \Pi\left( \frac{a^2}{a^2-1}; \am\left(\dfrac{a\theta}{k},k\right),k \right)\right. 
- \left.\Pi\left( \frac{a^2}{a^2-1}; \am\left(\dfrac{a\theta_0}{k},k\right),k\right) \right);
\end{align*}
\item в $C_3$:
$$
\f_3 =\sqrt{1-a^2}t  + \left( \arctan\left( \frac{a}{\sqrt{1-a^2}}\th a\theta \right) - \arctan\left( \frac{a}{\sqrt{1-a^2}}\th a\theta_0 \right)  \right);
$$
\item в $C_4$:
$$
\f_3 = t;
$$
\item в $C_5$:
$$
\f_3 = \sqrt{1-a^2}t.
$$
\end{enumerate}
Здесь $\am(\varphi,k)$ -- амплитуда Якоби, а $\Pi(n;\varphi,k)$ -- эллиптический интеграл третьего рода. Заметим, что из последних двух выражений видно, что геодезические, которые соответствуют областям $C_4$ и $C_5$, являются вращениями вокруг горизонтальных базисных векторов $e_1 = (1, 0, 0), e_2 = (0, 1, 0) \in \R^3$.

\subsubsection{Периодические геодезические}

\begin{proposition}
Для любого $a \in (0, 1)$  в соответствующей субримановой задаче на группе $\SO(3)$  существует бесконечное количество геодезических.

В случае $\lam \in C_1$ ($\lam \in C_2$)  периодическая геодезическая может иметь только период $T = \frac{4 K(k)}{a}$ ($T = \frac{4 k K(k)}{a}$), и такие траектории существуют тогда и только тогда, когда  для некоторых $n ,m \in \N$  выполнено равенство $\f_3(m T) = 2 \pi n$. Это равенство выполняется вдоль некоторой геодезической тогда и только тогда, когда 
\be{so3gen_nm}
\frac nm > \frac 1a
\ee
 в случае $C_1$ и 
\be{so3gen_nm1}
\frac nm > 1
\ee
 в случае $C_2$.  Различным несократимым дробям $\frac nm \in \Q_+$    соответствуют различные периодические геодезические.
\end{proposition}

\begin{proposition}
Любая периодическая геодезическая для $\lam \in C_1$ ($\lam \in C_2$)  однозначно определяется несократимой дробью $\frac nm \in \Q_+$, удовлетворяющей условию \eq{so3gen_nm} (соотв. \eq{so3gen_nm1}).

При $\lam \in C_3$  геодезические непериодичны.

При $\lam \in C_4 \cup C_5$  геодезические периодичны.
\end{proposition}

Так как $\pi_1(\SO(3)) = \Z_2$, то существуют только два гомотопических класса замкнутых путей на $\SO(3)$. Следующее утверждение показывает, какие из периодических геодезических стягиваемы (нуль-гомотопны).

\begin{proposition}
Рассмотрим периодическую геодезическую $g_t \in \SO(3)$, которая является проекцией экстремали $\lam_t \in C_1$ ($\lam_t \in C_2$)  и которая задана своей несократимой дробью $\frac nm \in \Q_+$, удовлетворяющей \eq{so3gen_nm} 
(соотв. \eq{so3gen_nm1}).  В этом случае геодезическая $g_t$  стягиваема тогда и только тогда, когда $n$  четно.

Все геодезические, соответствующие $\lam \in C_4 \cup C_5$, нестягиваемы.

\end{proposition}

\subsubsection{Условия оптимальности}
Рассмотрим трехмерную единичную сферу в алгебре кватернионов
$$
S^3 = \{ q = q^0 + i q^1 + j q^2 + k q^3 \in \H \mid (q^0)^2 + (q^1)^2 + (q^2)^2 + (q^3)^2 = 1 \}.
$$
Сфера $S^3$ односвязна и образует двулистное накрытие группы $\SO(3)$. Геодезическая $g_t \in \SO(3)$  имеет лифт $q_t \in S^3$, $q_0 = 1$,  вида
$$
q_t = \exp\left( - \frac{\f_1(0)}{2} k \right) \exp\left( - \frac{\f_2(0)}{2} i \right) \exp\left( \frac{\f_3(t)}{2} k \right) \exp\left(  \frac{\f_2(t)}{2} i \right) \exp\left( \frac{\f_1(t)}{2} k \right),    
$$  
где углы Эйлера $\f_i(t)$  совпадают с аналогичными углами в \eq{so3gen_f2}--\eq{so3gen_f3}.

\begin{theorem}
Пусть $g_t \in \SO(3)$, $t \in [0, t_1]$, есть геодезическая, а $q_t \in S^3$, $q_0 = 1$,  есть ее лифт на $S^3$. Пусть $\t_t$  есть соответствующая выпрямленная координата маятника \eq{so3gen_pend},  и $\tau = a(\t_0 + \frac t2)$. 

Тогда кривая $g_t$  неоптимальна, если для некоторого $t \in (0, t_1)$  выполняется одно из следующих условий:
\begin{itemize}
\item[$(1)$]
$q_t^0 = 0$,
\item[$(2)$]
$q_t^1 = 0$ и $\sn \tau \neq 0$, если $\lam_0 \in C_1 \cup C_2$, или  $\tau \neq 0$, если $\lam_0 \in C_3$,
\item[$(3)$]
$q_t^2 = 0$ и $\cn \tau \neq 0$, если $\lam_0 \in C_1$,
\item[$(3)$]
$q_t^3 = 0$ и $\cn \tau \neq 0$, если $\lam_0 \in C_2$.
\end{itemize}
\end{theorem}

\subsubsection{Библиографические комментарии}
Результаты  этого раздела  получены в работе~\cite{BS16}.

\subsection[Задача о качении шара по плоскости без прокручивания и проскальзывания]{Задача о качении шара по плоскости \\без прокручивания и проскальзывания}\label{subsec:roll}
\subsubsection{История задачи}\label{subsubsec:roll_history}
В 1983 году Дж.~Хаммерсли \cite{hammersley} рассмотрел следующую  {\em оксфордскую задачу о шаре}. Шар единичного радиуса лежит на бесконечной горизонтальной плоскости. {Состояние } шара определяется его пространственной ориентацией и положением на плоскости. Требуется перевести шар из заданного начального состояния в заданное конечное состояние с помощью последовательности качений. Каждое качение выполняется вдоль некоторой прямой на плоскости: длина  и направление качений выбираются нами, но качение должно выполняться \ddef{без прокручиваний и проскальзываний}, то есть ось вращения должна быть горизонтальной  и скорость шара в точке касания с плоскостью должна  быть нулевой. Какое наименьшее число качений $N$ необходимо для достижения любого конечного состояния? С использованием кватернионов Хаммерсли показал, что $N \in \{3, 4  \}$. Далее, он поставил две континуальные версии задачи о шаре:
\begin{itemize}
\item[$(a)$]
найти кривую $\G$ на плоскости минимальной длины $T$, переводящую шар в заданное конечное состояние;
\item[$(b)$]
перевести шар просто в некоторую заданную ориентацию, не заботясь о ее положении на плоскости.
\end{itemize}
Для задачи $(b)$ Хаммерсли указывает, что оптимальная кривая $\G$ есть отрезок или дуга окружности, и $0 \leq T \leq \pi\sqrt 3$, где  верхняя граница достигается, только если требуемая переориентация сферы есть ее поворот на $\pi$ вокруг  вертикальной оси.

В заключительном разделе статьи \cite{hammersley} <<Варианты для двадцать первого века>>  Хаммерсли ставит ряд вариаций и обобщений указанных задач о шаре, остающихся открытыми до сих пор.

В 1986 году А.~Артурс и Дж.~Уолш \cite{arthur_walsh} исследовали задачу $(a)$. С использованием кватернионов и принципа максимума Понтрягина они доказали, что точка контакта  шара и плоскости $(x, y)$
удовлетворяет уравнениям:
\begin{align*}
&\dx = \sin \p, \qquad \dy = - \cos \p,\\
&\ddot \p = \lam \cos(\p + \eps), \qquad \lam, \eps \equiv \const.
\end{align*}
Артурс и Уолш указали, что эти дифференциальные уравнения интегрируются в эллиптических интегралах первого и третьего рода, и оставили задачу оптимального управления для численного исследования.

Независимо от этих работ, в 1993 году Р.~Брокетт и Л.~Даи \cite{brock_dai} поставили <<задачу о пластинах и шаре>> (The Plate-Ball Problem). Они рассмотрели шар, катящийся без прокручивания и проскальзывания между двумя плоскими горизонтальными пластинами, расстояние между которыми равно диаметру шара. Брокетт и Даи записали управляемую систему для шара в форме \eq{roll_dxy}--\eq{roll_J} и показали, что нильпотентная аппроксимация  этой системы эквивалентна управляемой  системе \eq{cartanp21} на группе Картана (см. раздел  \ref{subsec:cartan}).

В том же 1993 году В.~Джурджевич \cite{jurd_ball} подробно исследовал \ddef{задачу об оптимальном качении шара по плоскости без прокручиваний и проскальзываний}, опираясь на постановку Брокетта и Даи \cite{brock_dai}, и независимо от работ \cite{hammersley, arthur_walsh}.  Джурджевич рассмотрел эту задачу как левоинвариантную  задачу оптимального управления на группе Ли $G = \R^2\times\SO(3)$:
\begin{align}
&\dx = u_1, \qquad \dy = u_2, \label{roll_dxy}\\
&\dot R = R
\begin{pmatrix}
0 & 0 & - u_1\\
0 & 0 & -u_2\\
u_1 & u_2 & 0
\end{pmatrix}, \label{roll_dR}\\
&g = (x, y, R) \in G, \quad u = (u_1, u_2) \in \R^2, \label{roll_gu} \\
&g(0) = \Id = (0, 0, E_{11} + E_{22 } + E_{33}), \qquad g(t_1) = g_1\label{roll_g(0)},\\
&J = \frac 12 \int_0^{t_1}(u_1^2 + u_2^2)\,dt \to \min.\label{roll_J}
\end{align}
Далее, он применил  к этой задаче принцип максимума Понтрягина в инвариантной формулировке для групп Ли (см. раздел \ref{subsubsec:PMP}), и получил следующие результаты. Оптимальные анормальные управления  постоянны и порождают качение шара по прямой; эти управления нестрого анормальны. Нормальные экстремали суть траектории гамильтоновой системы с гамильтонианом $H = \frac 12(h_1 - H_2)^2 + \frac 12(h_2 + H_1)^2$, где гамильтонианы $h_1$ и $h_2$ соответствуют векторным полям $\frac{\partial}{\partial x}$ и $\frac{\partial}{\partial y}$, а гамильтонианы $H_1$, $H_2$, $H_3$ соответствуют 
левоинвариантным полям на $\SO(3)$, задающим вращение трехмерного пространства  с генераторами
$$
A_1 = E_{32} - E_{23}, \qquad A_2 = E_{13} - E_{31}, \qquad A_3 = E_{21} - E_{12}. 
$$
Вертикальная подсистема этой гамильтоновой системы есть
\begin{align*}
&\dh_1 = \dh_2 = 0,\\
&\dot H_1 = (h_1 - H_2) H_3, \quad \dot H_2 = (h_2 + H_1)H_3,\\
&\dot H_3 = - h_1 H_1 - h_2 H_2.
\end{align*}
Эта подсистема имеет интегралы $h_1$, $h_2$, $H$ и $M = H_1^2 + H_2^2 + H_3^2$, поэтому интегрируема. Более того, эта подсистема сведена к уравнению маятника. Для интегрирования уравнения  для ориентации шара $R(t) \in\SO(3)$ вводятся углы Эйлера $\f_1$, $\f_2$, $\f_3$, для этих углов получены дифференциальные уравнения, которые качественно исследованы и частично проинтегрированы. Показано, что траектория точки контакта шара и плоскости $(x(t), y(t))$ есть эйлерова эластика, см.~раздел \ref{subsec:elastica}. Получена связь между типом  пересечения цилиндра $\{H = \const\}$ и сферы $\{M = \const  \}$, типом эластик и качественным поведением углов Эйлера $\f_1$, $\f_2$, $\f_3$.

Дальнейшее изложение в этом разделе опирается на \cite{s2r2_sym, s2r2}.

\subsubsection{Постановка задачи}\label{subsubsec:roll_state}
\paragraph{Механическая постановка}

Рассматривается механическая система, состоящая из двух горизонтальных плоскостей и сферы, касающейся этих плоскостей. Нижняя плоскость неподвижна, а сфера катится без прокручивания и проскальзывания благодаря горизонтальному движению верхней плоскости.
 Состояние такой системы описывается точкой контакта сферы с нижней плоскостью и ориентацией сферы в трехмерном пространстве. Требуется перекатить сферу из заданного начального состояния в заданное терминальное состояние так, чтобы кривая, пробегаемая точкой контакта на плоскости, имела минимальную длину. Управлением является  скорость верхней плоскости, или, что эквивалентно, скорость центра сферы.
 
Рассматривается кинематика данной системы, поэтому наличие верхней плоскости можно игнорировать и  изучать качение сферы по (нижней) плоскости без прокручивания и проскальзывания. Отсутствие проскальзывания означает, что мгновенная скорость точки контакта сферы и плоскости равна нулю, а отсутствие прокручивания означает, что вектор угловой скорости сферы горизонтален. Качение одной поверхности по другой  без прокручивания и проскальзывания моделирует работу руки робота-манипулятора, и задачи о таком движении вызывают большой интерес в механике, робототехнике и теории управления (см., например, работы~\cite{laumond, li_canny, bicchi_prat_sast,  marigo_bicchi, notes}).

\paragraph{Математическая постановка}
Пусть $e_1$, $e_2$, $e_3$  --- неподвижный правый репер в пространстве $\R^3$, такой, что векторы $e_1$, $e_2$  лежат  в  плоскости $\R^2 \cong (\R^2, 0) \subset \R^3$, по которой катится сфера $S^2$ единичного радиуса, а вектор $e_3$  направлен в полупространство, содержащее эту сферу. Репер  $e_1$, $e_2$, $e_3$ закреплен в точке $O \in (\R^2, 0)$. Пусть $f_1$, $f_2$, $f_3$ ---  подвижный правый репер, закрепленный в катящейся сфере $S^2$. Обозначим координаты точки в $\R^3$ в базисе $e_1$, $e_2$, $e_3$ как $(x,y,z)$, а координаты этой точки в базисе $f_1$, $f_2$, $f_3$, перенесенном в точку $O$, как $(X,Y,Z)$. Таким образом,
$$
x e_1 + y e_2 + z e_3 = X f_1 + Y f_2 + Z f_3.
$$
Пусть  матрица $R \in \SO(3)$  переводит координаты точки в неподвижном репере $e_1$, $e_2$, $e_3$ в ее координаты в подвижном репере $f_1$, $f_2$, $f_3$, т.е.
$$
\vect{X \\ Y \\ Z } = R \vect{x \\ y \\ z }. 
$$
 Состояние системы <<сфера $S^2$ и плоскость $\R^2$>>  задается координатами $(x,y)$ точки контакта $S^2$ и   $\R^2$, и матрицей вращения $R$. В качестве управлений будем использовать вектор $(u_1, u_2)$  скорости центра сферы.
Задача об оптимальном качении сферы по плоскости формализуется как следующая задача оптимального управления: 
\begin{align}
&\dot x = u_1, \label{sysx}\\
&\dot y = u_2, \label{sysy}\\
&\dot R = R(u_2 A_1 - u_1 A_2), \label{sysR}\\
&Q = (x,y,R) \in G = \R^2 \times \SO(3), \label{sysQ}\\
&u = (u_1, u_2) \in \R^2, \label{sysu}\\
&Q(0) = Q_0 = (0, 0, \Id), \qquad Q(t_1) = Q_1, \label{sysQ0}\\
&l = \int_0^{t_1} \sqrt{u_1^2 + u_2^2} \, dt \to \min. \label{sysl}
\end{align}
Здесь и далее мы используем базисные матрицы в алгебре Ли $\so(3)$:
\be{A123}
A_1 = 
\left(\begin{array}{ccc}
0 & 0 & 0 \\
0 & 0 & -1 \\
0 & 1 & 0 
\end{array}\right),
\quad
A_2 = 
\left(\begin{array}{ccc}
0 & 0 & 1 \\
0 & 0 & 0 \\
-1 & 0 & 0 
\end{array}\right),
\quad
A_3 = 
\left(\begin{array}{ccc}
0 & -1 & 0 \\
1 & 0 & 0 \\
0 & 0 & 0 
\end{array}\right).
\ee

\paragraph{Левоинвариантная субриманова задача}

Задача~\eq{sysx}--\eq{sysl} есть левоинвариантная субриманова задача на группе Ли $G = \R^2 \times \SO(3)$. 
Введем следующий левоинвариантный репер на этой группе Ли:
$$
e_1 = \pder{}{x}, \quad e_2 = \pder{}{y}, \quad 
V_i(R) = R A _i, \quad i = 1, 2, 3.
$$
В терминах левоинвариантных полей
$$
X_1 = e _1 - V_2, \qquad X_2 = e_2 + V_1,
$$
управляемая система~\eq{sysx}--\eq{sysu} принимает вид
\be{dotQ}
\dot Q =  u_1 X_1 (Q) + u_2 X_2(Q), \qquad Q \in G = \R^2 \times \SO(3), \quad (u_1, u_2) \in \R^2.
\ee
Функционал~\eq{sysl} есть функционал субримановой длины для левоинвариантной субримановой структуры, заданной полями $X_1$, $X_2$ как ортонормированным базисом:
\begin{align}
&l = \int_0^{t_1} \langle \dot Q, \dot Q\rangle^{1/2} \, dt \to \min, \label{lsubriem}\\
&\langle X_i, X_j\rangle = \delta_{ij}, \qquad i, \ j = 1, 2. \nonumber
\end{align} 

\paragraph{Существование оптимальных управлений }
Матричные коммутаторы $[A_i, A_j] = A_i A_j - A_j A_i$ вычисляются следующим образом:
$$
[A_1, A_2] = A_3, \qquad
[A_2, A_3] = A_1, \qquad
[A_3, A_1] = A_2.
$$
Таблица умножения в алгебре Ли $\gg = \R^2 \oplus \so(3) = \spann(e_1, e_2, V_1, V_2, V_3)$  группы Ли $G$  имеет вид:
\begin{align*}
&\ad e_i = 0, \qquad
[V_1, V_2] = V_3, \qquad
[V_2, V_3] = V_1, \qquad
[V_3, V_1] = V_2.
\end{align*}

В силу равенств 
$$
[X_1, X_2] = V_3, \qquad
[X_1, V_3] = - V_1, \qquad
[X_2, V_3] = - V_2, 
$$
векторные поля $X_1$, $X_2$ в правой части системы~\eq{dotQ}  порождают алгебру Ли $\gg$. По теореме Рашевского-Чжоу, система~\eq{dotQ} вполне управляема. Из теоремы Филиппова следует существование оптимальных управлений в задаче~\eq{sysx}--\eq{sysl} для любых $Q_0, Q_1 \in G$ в классе существенно ограниченных измеримых управлений.

\subsubsection{Экстремали}\label{subsubsec:roll_extrem}
Введем линейные на слоях $T^*G$ гамильтонианы:
\begin{align*}
&h_i(\lam) = \langle \lam, e_i \rangle, \quad
i = 1, 2,\\
&H_i(\lam) = \langle \lam, V_i \rangle, \quad
i = 1, 2, 3.
\end{align*}
\paragraph{Анормальные траектории}
Анормальные траектории постоянной скорости имеют вид
\begin{align*}
&x_t = u_1 t, \qquad y_t = u_2 t, \\
&R_t = \exp(t(u_2 A_1 - u_1 A_2)).
\end{align*}
Они нестрого анормальны и оптимальны. В анормальном случае сфера равномерно катится по прямой.

\paragraph{Нормальная гамильтонова система}
 В нормальном случае гамильтонова система с гамильтонианом 
$$
H = \frac 12 ((h_1 - H_2)^2 + (h_2 + H_1)^2)
$$
записывается в координатах так
\begin{align}
&\dot h_1 = \dot h_2 = 0, \label{ham31}\\
&\dot H_1 = (h_1 - H_2)H_3, \label{ham32}\\
&\dot H_2 = (h_2 + H_1) H_3, \label{ham33}\\
&\dot H_3 = - h_1 H_1 - h_2 H_2, \label{ham34}\\
&\dot Q = (h_1 - H_2) X_1 + (h_2 + H_1) X_2. \label{ham35}
\end{align}
Как всегда в субримановых задачах, можно ограничиться геодезическими единичной скорости, т.е. экстремальными траекториями, вдоль которых
$H  \equiv \frac 12$. При таком ограничении удобно перейти в сопряженном пространстве от координат $(h_1, h_2, H_1, H_2, H_3)$  к новым координатам $(r, \a, \theta, c)$: 
\begin{align}
&h_1 = r \cos \a, \qquad h_2 = r \sin \a, \label{h1h2}\\
&h_1 - H_2 = \cos(\t + \a), \qquad h_2 + H_1 = \sin (\t + \a), \label{h1-H2}\\
&c = H_3. \nonumber
\end{align}
После этого гамильтонова система для нормальных экстремалей~\eq{ham31}--\eq{ham35} принимает следующую форму:
\begin{align}
&\dot \t = c, \label{roll_pend1}\\
&\dot c = -r \sin \t , \label{roll_pend2}\\
&\dot \a = \dot r = 0, \label{roll_pend3} \\
&\dot x = \cos(\t + \a), \label{roll_x}\\
&\dot y = \sin(\t + \a), \label{roll_dy}\\
&\dot R =  R \Om, \qquad \Om = \sin(\t + \a) A_1 - \cos(\t + \a) A_2. \label{roll_dRR}
\end{align}
Семейство нормальных экстремалей $\lam_t$ параметризуется цилиндром $C$, состоящим из начальных точек $\lam = \restr{\lam_t}{t = 0}$: 
\begin{align*}
C &= \{ \lam \in\gg^* \mid H(\lam) = 1/2 \} \\
& \cong \{(h_1, h_2, H_1, H_2, H_3) \in \R^5 \mid (h_1-H_2)^2 + (h_2 + H_1)^2 = 1 \}\\
& \cong \{(\t, c, \a, r) \mid \t \in S^1, \ c \in \R, \ \a \in S^1, \ r \geq 0 \}.
\end{align*}
Экспоненциальное отображение определяется как
\begin{align*}
&\Exp(\lam,t) = \pi  \circ e^{t \vH}(\lam) = Q_t,\\
&\map{\Exp}{N}{M}, \\
&N = C \times \R_+ = \{ (\lam,t) \mid \lam \in C, \ t > 0 \}.
\end{align*}

В случае $r = 0$  эластика $(x_t,y_t)$ есть прямая (при $H_3 = c = 0$)  или окружность (при $H_3 = c \neq 0$), будем называть такие эластики \ddef{вырожденными}.

В случае $r \neq 0$  эластика $(x_t, y_t)$ принадлежит одному из четырех классов в зависимости от полной энергии $E = c^2/2 - r \cos \t$ маятника~\eq{roll_pend1}, \eq{roll_pend2}, см. раздел \ref{subsec:elastica}:
\begin{enumerate}
\item
инфлексионный при   $E \in (-r, r)$, \\
\item
неинфлексионный при   $E \in (r, + \infty)$, \\
\item
критический при   $E = r$, $c \neq 0$, \\
\item
прямая при    $E = -r$ и при $E = r$,  $c = 0$.
\end{enumerate}
Эластики классов 1--3 будем называть \ddef{невырожденными}.

\paragraph{Симплектическое слоение}
На коалгебре Ли $\gg^*$ имеются функции Казимира $h_1$, $h_2$, $M = H_1^2 +H_2^2 +H_3^2$. Симплектическое слоение состоит из:
\begin{itemize}
\item
сфер $\{h_1, h_2 = \const, \quad M = \const>0  \}$,
\item
точек $\{h_1, h_2 = \const, \quad H_1 = H_2 = H_3 = 0  \}$.
\end{itemize}

Нормальная гамильтонова система имеет интегралы $h_1$, $h_2$, $M$, $E = \frac 12(M + h_1^2 + h_2^2) - H$ и интегрируема  в эллиптических функциях и интегралах.

Различные типы геодезических, проецирующихся в эйлеровы эластики $(x_t, y_t)$, соответствуют разным типам пересечения поверхности  уровня гамильтониана $\{H = \const\}$ с  симплектическими листами.

\paragraph{Выпрямляющие координаты }
Цилиндр $C = \{\lam \in  \gg^* \mid H(\lam) = \frac 12   \}$ стратифицируется согласно разным типам движения маятника \eq{roll_pend1}, \eq{roll_pend2}:
\begin{align*}
&C = \sqcup_{i=1}^7 C_i, \\
&C_1 = \{\lam \in C \mid E \in (-r, r), \ r > 0\}, \\
&C_2 = \{\lam \in C \mid E \in (r, + \infty), \ r > 0\}, \\
&C_3 = \{\lam \in C \mid E = r > 0, \ c \neq 0\}, \\
&C_4 = \{\lam \in C \mid E = -r, \ r > 0\}, \\
&C_5 = \{\lam \in C \mid E = r > 0, \ c = 0\}, \\
&C_6 = \{\lam \in C \mid  r = 0, \ c \neq 0\}, \\
&C_7 = \{\lam \in C \mid  r = 0, \ c = 0\}.
\end{align*}

В области $\cup_{i=1}^3 C_i$ введем координаты $(\f, k, \a, r)$, выпрямляющие уравнения маятника \eq{roll_pend1}, \eq{roll_pend2}.

Если $\lam = (\t, c, \a, r)\in C_1$, то
$$\sin (\t/2) = k \sn(\sr \f, k), \quad \cos (\t/2) = \dn(\sr \f, k), \quad {c}/{2} =  k \sqrt{r} \cn(\sr \f, k),$$ 
при этом
$k = \sqrt{(E + r)/(2r)}  \in (0, 1)$, $\sr \f \pmod {4 K} \in [0, 4 K]$.

Если
$\lam = (\t, c, \a, r)\in C_2$, то
$$\sin (\t/2) = \pm \sn (\sr \f/k, k), \quad \cos (\t/2) = \cn (\sr \f/k, 
k), \quad c/2 =  \pm \sr/k  \dn (\sr \f/k, k), $$
где $ \pm = \sgn c$, при этом 
$k =  \sqrt{2 r/(E + r)}  \in (0, 1)$,  $\sr \f \pmod {2 k K}  \in [0, 2kK]$. 

Если $\lam \in C_3$, то
$$\sin (\t/2) = \pm \th(\sr \f), \quad \cos (\t/2) = 1/\ch(\sr \f), \quad c/2 =  \pm \sr/\ch(\sr \f),$$ 
где $\pm = \sgn c$, при этом
$k = 1$, $\f \in (- \infty, + \infty)$.

В новых координатах уравнения маятника \eq{roll_pend1}, \eq{roll_pend2} принимают форму:
\be{}\nonumber
\dot \f = 1, \quad \dot k = 0, \quad  \dot \a = 0, \quad  \dot r = 0,
\ee
откуда $\f_t = \f + t$; $k, \a, r = \const$.

\paragraph{Интегрирование вертикальной подсистемы гамильтоновой системы ПМП}
Если
$\lam  \in C_1$, то
$$\sin (\t_t/2) = k \sn(\sr \f_t, k), \quad \cos (\t_t/2) = \dn(\sr \f_t, k), \quad {c_t}/{2} =  k \sqrt{r} \cn(\sr \f_t, k).$$
Если
$\lam  \in C_2$, то
$$\sin (\t_t/2) = \pm \sn (\sr \f_t/k, k), \quad \cos (\t_t/2) = \cn (\sr 
\f_t/k, k), \quad c_t/2 =  \pm \sr/k  \dn (\sr \f_t/k, k),$$ 
где $\pm = \sgn c$.

Если
$\lam \in C_3$, то
$$\sin (\t_t/2) = \pm \th(\sr \f_t), \quad \cos (\t_t/2) = 1/\ch(\sr 
\f_t), \quad c_t/2 =  \pm \sr/\ch(\sr \f_t),$$ 
где $\pm = \sgn c$.

Для случаев $\lam \in \cup_{i=4}^7 C_i$  система~\eq{roll_pend1}--\eq{roll_pend3} 
интегрируется непосредственно:
$\t_t \equiv 0$, $c_t \equiv 0$ при $\lam \in C_4$; 
$\t_t \equiv \pi$, $c_t \equiv 0$ при $\lam \in C_5$;
$\t_t = c t + \t$, $c_t \equiv c \neq 0$ при $\lam \in C_6$;
$\t_t \equiv \t$, $c_t \equiv 0$ при $\lam \in C_7$.

\paragraph{Интегрирование уравнений для $x$, $y$}
Для интегрирования уравнений	\eq{roll_x}, \eq{roll_dy} с начальным условием $x_0 = y_0 = 0$ воспользуемся симметрией задачи --- поворотом
$$
x = \bar x\cos \a - \bar y \sin \a, \quad y = \bar x \sin \a + \bar y \cos \a.
$$
В новых переменных получаем задачу Коши
\be{dotbarxt}
\dot{\bar x}_t = \cos \t_t, \quad  \dot{\bar y}_t = \sin \t_t, \qquad \bar x_0 
= \bar y_0 = 0,
\ee
решения которой параметризуются следующим образом.

Если $\lam  \in C_1$, то
\begin{align*}
&\bx_t = (2(\E(\sr \f_t) - \E(\sr\f)) - \sr t)/\sr,  \\
& \by_t = 2k(\cn (\sr \f) - \cn (\sr \f_t))/\sr.
\end{align*}

Если 
$\lam  \in C_2$, то
\begin{align*}
&\bx_t = 2 \left(\E(\sr \f_t/k) - \E(\sr \f/k) - (2 - k^2) \sr t 
/(2k)\right) /(k \sr),
\\
&\by_t = \pm 2 (\dn (\sr \f/k) - \dn (\sr \f_t/k))/(k \sr), \qquad \pm 
= \sgn c.
\end{align*}

Если 
$\lam  \in C_3$, то
\begin{align*}
&\bx_t = (2(\th (\sr \f_t) - \th(\sr \f)) - \sr t)/\sr,    \\
&\by_t =  \pm 2 ( 1/\ch(\sr  \f) - 1/\ch(\sr \f_t))/\sr, \qquad 
\pm = \sgn c.
\end{align*}

При $\lam \in \cup_{i=4}^7 C_i$  уравнения~\eq{dotbarxt} интегрируются 
непосредственно:
$\bx_t = t$, $\by_t = 0$  при $\lam \in C_4$;
$\bx_t = -t$, $\by_t = 0$ при $\lam \in C_5$; 
$\bx_t = (\sin(c t + \t) - \sin \t)/c$, $\by_t = (\cos \t - \cos(c t + 
\t))/c$ при $\lam \in C_6$; 
$\bx_t = t \cos \t$, $\by_t = t \sin \t$  при $\lam \in C_7$. 

\paragraph{Интегрирование уравнений для $R$}
Пусть $M = H_1^2 + H_2^2 +H_3^2 > 0$. Тогда
\be{R(t)}
R(t) =  e^{(\a - \f_3^0)  A_3} e^{-\f_2^0  A_2} e^{\f_1(t)  A_3}  
e^{\f_2(t)  A_2}  e^{(\f_3(t)- \a)  A_3},
\ee
где углы $\f_i$  определяются из соотношений~\eq{cosphi2}--\eq{el_int_III}  при $r \neq 1$  и~\eq{cosphi2r1}--\eq{phi_t}   
 при $r = 1$,  а угол $\f_1$  удовлетворяет начальному условию  $ \f_1^0 = 0$.

Входящие в разложение~\eq{R(t)} экспоненты матриц, содержащие  $\f_2$, 
$\f_3$, выражаются через $\cos \f_2$, $\sin \f_2$,  $\cos \f_3$, $\sin 
\f_3$, которые с помощью соотношений~\eq{cosphi2}, \eq{cosphi3}, 
\eq{cosphi2r1}, \eq{cosphi3r1} выражены через переменные $c$, 
$\cos(\t/2)$, $\sin(\t/2)$, которые, в свою очередь представлены выше как функции эллиптических координат или 
непосредственно. При $r = 1$  имеем $\f_1(t) = \sqrt M t/2$. 
Интегрирование уравнения \eq{dotphi1} при $r \neq 1$ вынесено в 
следующий пункт.

В случае $M = 0$  имеем $r = 1$, $c = 0$, $\t = 0$, откуда $u_1 = \cos 
\a$, $u_2 = \sin \a$. Поэтому $\Om = u_2 A_1 - u_1 A_2 \equiv \const$  
и $R(t) = e^{t \Om}$. 

\paragraph{Интегрирование уравнений для $\f_1$}
Вдоль нормальных геодезических углы $\f_i$  удовлетворяют при $r \neq 1$ равенствам:
\begin{align}
&\cos \f_2 = c/\sM, \qquad \sin \f_2 = \pm \sqrt{M - c^2}/\sM, 
\label{cosphi2} \\ 
&\cos \f_3 = \mp \sin \t/\sqrt{M - c^2}, \qquad \sin \f_3 = \pm (r - 
\cos \t)/\sqrt{M - c^2}, \label{cosphi3} \\
&\dot\f_1 = \sM(1 - r \cos \t)/(M - c^2), \label{dotphi1}
\end{align}
а при $r = 1$ равенствам:
\begin{align}
&\cos \f_2 = c/\sM, \qquad \sin \f_2 = \pm 2 \sin(\t/2)/\sM, 
\label{cosphi2r1} \\ 
&\cos \f_3 = \mp \cos(\t/2), \qquad \sin \f_3 = \pm \sin(\t/2), 
\label{cosphi3r1} \\
&\dot\f_1 = \sM/2.  \label{dotphi1r1}
\end{align}

Введем в рассмотрение \ddef{эллиптический интеграл III рода} в следующей 
форме:
\be{el_int_III}
\Pi(n, u, k) = \int_0^u \frac{dt}{(1 - n \sin^2 t)\sqrt{1 - k^2 \sin^2 
t}} = \int_0^{F(u,k)} \frac{dv}{1 - n \sn^2 v}.
\ee

Пусть $r \neq 1$.
Если $\lam_1 \in C_1$, то
\begin{align}
&\f_1(t) = \frac{\sM}{2} t + \frac{\sM(1 + r)}{2 \sr(1-r)} (\Pi(l, 
\am(\sr(\f + t)), k) - \Pi(l, \am(\sr\f ), k)),\label{phi_t}
\end{align}
где $l = -\frac{4 k^2 r}{(1-r)^2}$. 

Если $\lam_1 \in C_2$, то
\begin{align*}
&\f_1(t) = \frac{\sM}{2} t + \frac{\sM k(1 + r)}{2 \sr(1-r)} (\Pi(l, 
\am(\sr(\f + t)/k), k) - \Pi(l, \am(\sr\f/k), k)),
\end{align*}
где $l = -\frac{4 r}{(1-r)^2}$.

Если $\lam_1 \in C_3$, то
\begin{align*}
&\f_1(t) = \frac{\sM}{2} t + \frac{\sM k(1 - r^2)}{8 r^{3/2}} 
(I(\sr(\f + t), a) - I( \sr\f, a)), \\
&I(v,a) = \int_0^v \frac{dt}{a^2 + \th^2t} = 
\frac{a t - \arctg a + \arctg (e^t(a^2 \ch t + \sh t)/a)}{a + a^3}, 
\end{align*}
где $ a = (1-r)/(2 \sr)$.

Если $\lam_1 \in C_6$, то $\f_1(t) = \sqrt{1 
+ c^2}\, t$.

Если $\lam \in C_4 \cup C_5 \cup C_7$, то
$\t_t \equiv \const = \t$, 
$\Om = \sin(\a + \t) A_1 - \cos(\a + \t) A_2 \equiv \const$, $R(t) = 
e^{t \Om}$.

\paragraph{Управляемая система в терминах кватернионов}
Для описания ориентации катящейся сферы удобно, наряду с 
матрицей вращения $R$, использовать кватернионы.

Пусть $\H = \{q = q_0 + i q_1 + j q_2 + k q_3 \mid q_0,\dots, q_3 \in \R 
\}$ есть алгебра кватернионов, $S^3 = \{q \in \H \mid |q|^2 = q_0^2 + 
q_1^2 + q_2^2 + q_3^2 = 1\}$ --- единичная сфера, $I = \{q \in 
\H \mid \mathrm{Re}\, q = q_0 = 0\}$ --- подпространство чисто мнимых 
кватернионов. Любой кватернион $q \in S^3$ задает вращение евклидова 
пространства $I$:
\[q \in S^3 \quad\Rightarrow \quad R_q(a) = q a q^{-1}, \quad a \in I, \quad R_q \in \SO(3)\cong \SO(I).\]

Соответствие между кватернионом $q = q_0 + i q_1 + j q_2 + k q_3 \in S^3$ и матрицей $R \in \SO(3)$ имеет вид:
\begin{equation}
R = \left(\begin{array}{ccc}
q_0^2 + q_1^2 - q_2^2 - q_3^2 & 2\ q_1\ q_2 - 2\ q_0\ q_3& 2\ q_0\ q_2 
+ 2\ q_1\ q_3 \\
2\ q_1\ q_2 + 2\ q_0\ q_3 &q_0^2 - q_1^2 + q_2^2 - q_3^2 & -2\ q_0\ 
q_1 + 2\ q_2\ q_3 \\
-2\ q_0\ q_2 + 2\ q_1\ q_3 & 2\ q_0\ q_1 + 2\ q_2\ q_3 & q_0^2 - q_1^2 
- q_2^2 + q_3^2 
\end{array}\right).
\end{equation}

Управляемая система \eq{roll_dR} в терминах кватернионов принимает форму
\begin{eqnarray}
\begin{cases}
\dot{q}_0 = \frac{1}{2}(q_2 u_1 - q_1 u_2),\\
\dot{q}_1 = \frac{1}{2}(q_3 u_1 + q_0 u_2),\\
\dot{q}_2 = \frac{1}{2}(-q_0 u_1 + q_3 u_2),\\
\dot{q}_3 = \frac{1}{2}(-q_1 u_1 - q_2 u_2),
\label{roll_dq}
\end{cases}
q = q_0 + i q_1 + j q_2 + k q_3 \in S^3, \quad   (u_1, u_2) \in \R^2.
\end{eqnarray}

Управляемая система на $\R^2\times\SO(3)$ \eq{roll_dxy}, \eq{roll_dR} с начальным условием $g(0) = \Id$ имеет лифт на $\R^2\times S^3$ вида \eq{roll_dxy}, \eq{roll_dq} с начальным условием $(x, y)(0) = (0, 0)$, $q(0) = 1$.

\subsubsection{Симметрии}\label{subsubsec:roll_max}
\paragraph{Симметрии семейства экстремальных траекторий}
Вращения эластик $(x_s, y_s)$  вокруг начала координат в плоскости $(x,y)$  порождают однопараметрическую группу симметрий траекторий гамильтоновой системы~\eq{roll_pend1}--\eq{roll_dRR}:
$$
\{ \F^{\b} \mid \b \in S^1 \},
$$
где вращение $\F^{\b}$  определяется следующим образом:
\begin{align}
&\map{\F^{\b}}{\{\lam_s \mid s \in [0,t]\}}{\{\lam_s^{\b} \mid s \in [0,t]\}}, \label{roll_rot1}\\
&\lam_s = (\t_s, c_s, \a, r, Q_s), \qquad Q_s = (x_s, y_s, R_s),  \label{roll_rot2}\\
&\lam_s^{\b} = (\t_s^{\b}, c_s^{\b}, \a^{\b}, r, Q_s^{\b}), \qquad Q_s^{\b} = (x_s^{\b}, y_s^{\b}, R_s^{\b}), \label{roll_rot3} \\
&\t_s^{\b} = \t_s, \qquad c_s^{\b} = c_s, \qquad \a^{\b} = \a + \b, \label{roll_rot4} \\
&\vect{x_s^{\b} \\ y_s^{\b}}
= \left(\begin{array}{cc}
\cos \b & - \sin \b \\
\sin \b &  \cos \b
\end{array}\right)
\vect{x_s \\ y_s},  \label{roll_rot5}\\
&R_s^{\b} = e^{\b A_3} R_s e^{-\b A_3}, \qquad \Om_s^{\b} = e^{\b A_3} \Om_s e^{-\b A_3}. \label{roll_rot6}
\end{align}

\begin{proposition}
\label{propos:rotlams}
Если $\{\lam_s \mid s \in [0, t]\}$  есть траектория системы~\eq{roll_pend1}--\eq{roll_dRR}, то для любого $\b \in S^1$ кривая $\{\lam_s^{\b} \mid s \in [0, t]\}$ есть также траектория этой системы.
\end{proposition}
Отражения траекторий $(\t_s, c_s)$ маятника~\eq{roll_pend1}, \eq{roll_pend2}   в осях координат $\t$, $c$ и в начале координат  продолжаются до  дискретных симметрий $\eps^1$, $\eps^2$, $\eps^3$  семейства траекторий гамильтоновой системы~\eq{roll_pend1}--\eq{roll_dRR}:
\begin{align*}
&\map{\eps^i}{\{\lam_s \mid s \in [0,t]\}}{\{\lam_s^{i} \mid s \in [0,t]\}}, \qquad i = 1, \ 2, \ 3, \\
&\lam_s = (\t_s, c_s, \a, r, Q_s), \qquad Q_s = (x_s, y_s, R_s), \\
&\lam_s^{i} = (\t_s^{i}, c_s^{i}, \a^{i}, r, Q_s^{i}), \qquad Q_s^{i} = (x_s^{i}, y_s^{i}, R_s^{i}).
\end{align*}

Отражению траекторий $(\t_s, c_s)$ маятника~\eq{roll_pend1}, \eq{roll_pend2}   в оси координат $\t$  соответствует  дискретная симметрия $\eps^1$ семейства экстремальных траекторий:
\begin{align*}
&\t_s^{1} = \t_{t-s}, \qquad c_s^{1} = -c_{t-s}, \qquad \a^{1} = \a + \pi, \\
&x_s^1 = x_{t-s} - x_t, \qquad y_s^1 = y_{t-s} - y_t, \\
&R_s^{1} = (R_t)^{-1} R_{t-s}, \qquad \Om_s^{1} = -\Om_{t-s}. 
\end{align*}

Отражение траекторий маятника $(\t_s, c_s)$ в оси координат $c$  порождает  симметрию $\eps^2$  экстремальных траекторий:
\begin{align*}
&\t_s^{2} = -\t_{t-s}, \qquad c_s^{2} = c_{t-s}, \qquad \a^{2} = \pi - \a, \\
&x_s^2 = x_{t-s} - x_t, \qquad y_s^2 = y_t - y_{t-s}, \\
&R_s^{2} = I_2 (R_t)^{-1} R_{t-s} I_2, \qquad \Om_s^{2} = - I_2 \Om_{t-s} I_2,  \\
&I_2 = I_2^{-1} = e^{\pi A_2} = 
\left(\begin{array}{ccc}
- 1 & 0 & 0 \\
0 & 1 & 0 \\
0 & 0 & -1 
\end{array}\right). 
\end{align*}

Отражение траекторий маятника $(\t_s, c_s)$ в начале координат  $(\t,c) = (0,0)$  продолжается до  симметрии $\eps^3$  экстремальных траекторий:
\begin{align*}
&\t_s^{3} = -\t_{s}, \qquad c_s^{3} = -c_{s}, \qquad \a^{3} = - \a, \\
&x_s^3 = x_{s}, \qquad y_s^3 = - y_{s}, \\
&R_s^{3} = I_2 R_{s} I_2, \qquad \Om_s^{3} = I_2 \Om_{s} I_2. 
\end{align*}
\begin{proposition}
Если $\{\lam_s \mid s \in [0, t]\}$  есть траектория системы~\eq{roll_pend1}--\eq{roll_dRR}, то  кривые $\{\lam_s^{i} \mid s \in [0, t]\}$, $i = 1, 2, 3$, суть также траектории этой системы.
\end{proposition}

\paragraph{Симметрии экспоненциального отображения}
Действие вращений $\Fb$ и отражений $\eps^i$ в прообразе и образе экспоненциального отображения определяется так, чтобы они коммутировали с действием экспоненциального отображения.

Вращения $\mapto{\Fb}{\lam}{\lam^{\b}}$~\eq{roll_rot1}--\eq{roll_rot6} являются симметриями гамильтоновой системы, поэтому их действие в $T^*G$  естественно распадается в прямую сумму действий в 
$N = \gg^* \times \R_+$ (на $(\lam,t)$, где $\lam$ --- начало экстремали) и в $G$ (на $Q_t$ ---  конец соответствующей экстремальной траектории):
\begin{align*}
&\map{\Fb}{N}{N}, \qquad (\lam, t) \mapsto (\lam^{\b}, t), \\
&\lam = (\t, c, \a, r), \qquad \lam^{\b} = (\t, c, \a^{\b}, r), \\
&\a^{\b} = \a + \b,\\
\intertext{и}
&\map{\Fb}{G}{G}, \qquad Q \mapsto Q^{\b}, \\
&Q = (x, y, R), \qquad Q^{\b} = (x^{\b}, y^{\b}, R^{\b}), \\
&\vect{x^{\b} \\ y^{\b}} = 
\left(\begin{array}{cc}
\cos \b & - \sin \b\\
\sin \b & \cos \b 
\end{array}\right)
\vect{x \\ y}, 
\qquad R^{\b} = e^{\b A_3} R e^{-\b A_3}.
\end{align*}

Действие отражений $\eps^i$  в $N$ определяется ограничением их действия на вертикальные составляющие экстремальных траекторий в начальный момент времени $s = 0$:
\begin{align*}
&\map{\eps^i}{N}{N}, \qquad (\lam, t) \mapsto (\lam^{i}, t), \qquad i = 1, 2, 3,\\
&\lam = (\t, c, \a, r), \qquad \lam^{i} = (\t^i, c^i, \a^{i}, r), 
\end{align*}
где $\lam = \restr{\lam_s}{s = 0}$, $\lam^i = \restr{\lam_s^i}{s = 0}$.   Явные выражения для действия $\eps^i$  в $N$:
\begin{align*}
&\mapto{\eps^1}{(\t, c, \a, r, t)}{(\t^1, c^1, \a^1, r, t) = (\t_t, -c_t, \a+ \pi, r, t)}, \\ 
&\mapto{\eps^2}{(\t, c, \a, r, t)}{(\t^2, c^2, \a^2, r, t) = (-\t_t, c_t, \pi - \a, r, t)}, \\ 
&\mapto{\eps^3}{(\t, c, \a, r, t)}{(\t^3, c^3, \a^3, r, t) = (-\t, -c, -\a, r, t)}. 
\end{align*}

Действие отражений в $G$  определяется их действием на экстремальные траектории в конечный момент времени $s = t$:
\begin{align*}
&\map{\eps^i}{G}{G}, \qquad Q \mapsto Q^{i}, \qquad i = 1, 2, 3,\\
&Q = (x, y, R), \qquad Q^i = (x^i, y^i, R^i), 
\end{align*}
где $Q = \restr{Q_s}{s = t}$, $Q^i = \restr{Q_s^i}{s = t}$. Явные формулы:
\begin{align*}
&\mapto{\eps^1}{(x, y, R)}{(x^1, y^1, R^1) = (-x, -y, (R)^{-1})}, \\ 
&\mapto{\eps^2}{(x, y, R)}{(x^2, y^2, R^2) = (-x, y, I_2 (R)^{-1}I_2}), \\ 
&\mapto{\eps^3}{(x, y, R)}{(x^3, y^3, R^3) = (x, -y, I_2 R I_2 )}. 
\end{align*}

Итак, определено действие вращений и отражений в прообразе и образе экспоненциального отображения:
\begin{align}
&\map{\Fb, \ \eps^i}{N}{N}, \qquad 
(\lam, t) \mapsto (\lam^{\b}, t), \ (\lam^i,t), \label{PhibetaepsiN}\\
&\map{\Fb, \ \eps^i}{G}{G}, \qquad 
Q \mapsto Q^{\b}, \ Q^i. \label{PhibetaepsiM}
\end{align}
Существенно, что образ $Q^i = \eps^i(Q)$ зависит лишь от прообраза $Q$, но не от момента времени $t$.

\begin{proposition}
\label{propos:symExp}
Отображения $\Fb$, $\eps^i$  являются симметриями экспоненциального отображения.
\end{proposition}

Рассмотрим группу симметрий экспоненциального отображения, порожденную вращениями и отражениями:
$$
\Sym = \langle \Fb, \ \eps^1, \ \eps^2, \ \eps^3 \rangle.
$$
Таблица умножения в этой группе имеет следующий вид:
$$
\begin{array}{|c|c|c|c|c|}
\hline
\cdot \, \circ \, \cdot  & \eps^1 & \eps^2 & \eps^3 & \Fb\\
  \hline
\eps^1 & \Id & \eps^3 & \eps^2 & \Fb \circ \eps^1 \\
  \hline
\eps^2  & \eps^3 & \Id & \eps^1 & \F^{-\b} \circ \eps^2\\
  \hline
\eps^3  & \eps^2 & \eps^1 & \Id & \F^{-\b} \circ \eps^3\\
  \hline
  \F^{\g}  & \eps^1 \circ \F^{\g} & \eps^2 \circ \F^{-\g} & \eps^3 \circ \F^{-\g}  & \F^{\b + \g}\\
  \hline
  \end{array}
$$
Отсюда получаем явное описание группы симметрий экспоненциального отображения: 
$$
\Sym = \{ \Fb, \ \Fb \circ \eps^i \mid \b \in S^1, \ i = 1, 2, 3\} \cong \SO(2)\times(\Z_2\times\Z_2).
$$

Определим \ddef{множество Максвелла, соответствующее группе $\langle \eps^i, \ \Fb\rangle$, $i = 1, 2, 3$}:
\begin{multline*}
\MAX^i = \{(\lam,t) \in N \mid \exists \ \b \in S^1 \ : \   (\tilde \lam, t) = \eps^i \circ \Fb(\lam,t), \  
\Exp(\lam, s) \not\equiv \Exp(\tilde\lam, s), \ \Exp(\lam, t) =  \Exp(\tilde \lam, t)\}.
\end{multline*}

\subsubsection{Условия оптимальности}
\begin{theorem}
\label{th:Max1}
Пусть $t > 0$ и  $Q_s = (x_s, y_s, R_s) = \Exp(\lam,s)$   есть такая экстремальная траектория, что:
\begin{enumerate}
\item[$(1)$]
$q_3(t) = 0$,
\item[$(2)$]
эластика $\{(x_s, y_s) \mid s \in [0,t]\}$  невырождена и не центрирована в точке перегиба.
\end{enumerate}
Тогда $(\lam,t) \in \MAX^1$,  поэтому для любого $t_1 > t$  траектория $Q_s$, $s \in [0, t_1]$, неоптимальна.
\end{theorem}
\begin{theorem}
\label{th:Max2}
Пусть $t > 0$ и  $Q_s = (x_s, y_s, R_s) = \Exp(\lam,s)$   есть такая экстремальная траектория, что:
\begin{enumerate}
\item[$(1)$]
$(x q_1 + y q_2)(t) = 0$,
\item[$(2)$]
эластика $\{(x_s, y_s) \mid s \in [0,t]\}$  невырождена и не центрирована в вершине.
\end{enumerate}
Тогда $(\lam,t) \in \MAX^2$,  поэтому для любого $t_1 > t$  траектория $Q_s$, $s \in [0, t_1]$, неоптимальна.
\end{theorem}
\begin{theorem}
\label{th:Max3}
Пусть $t > 0$ и  $Q_s = (x_s, y_s, R_s) = \Exp(\lam,s)$   есть такая экстремальная траектория, что:
\begin{enumerate}
\item[$(1)$]
$(x q_1 + y q_2)(t) =  q_3(t) = 0$  или $(y q_1 - x q_2)(t) =  q_0(t) = 0$.
\item[$(2)$]
эластика $\{(x_s, y_s) \mid s \in [0,t]\}$  невырождена.
\end{enumerate}
Тогда $(\lam,t) \in \MAX^3$,  поэтому для любого $t_1 > t$  траектория $Q_s$, $s \in [0, t_1]$, неоптимальна.
\end{theorem}
\begin{remark}
Учитывая то, что для любого кватерниона $q = q_0 + iq_1 +jq_2 + kq_3 \in S^3$, соответствующее движение $R_q: \R^3 \to \R^3$ есть вращение вокруг вектора $(q_1, q_2, q_3) \in\R^3$, можно дать следующую наглядную интерпретацию условию (1) теорем \ref{th:Max1}--\ref{th:Max3}:
\begin{enumerate}
\item
Условие (1) теоремы \ref{th:Max1} означает,  что вращение сферы $R_t$ есть поворот вокруг некоторой горизонтальной оси;
\item
Условие (1) теоремы \ref{th:Max2} означает,  что вращение  $R_t$ есть поворот вокруг некоторой  оси, ортогональной вектору перемещения точки контакта сферы и плоскости $(x_t, y_t, 0)$;
\item
Условие (1) теоремы \ref{th:Max3} означает,  что вращение  $R_t$ есть поворот вокруг горизонтальной  оси, ортогональной вектору $(x_t, y_t, 0)$, или что $R_t$ есть поворот на угол $\pi$ вокруг некоторой оси, лежащей в вертикальной плоскости, которая содержит вектор $(x_t, y_t, 0)$.
\end{enumerate}
\end{remark}

\subsubsection{Библиографические комментарии}
Раздел \ref{subsubsec:roll_history} опирается на работы \cite{hammersley, arthur_walsh, brock_dai, jurd_ball}, раздел \ref{subsubsec:roll_state} --- на \cite{s2r2_sym}, раздел \ref{subsubsec:roll_extrem} --- на \cite{s2r2_sym, s2r2},  раздел \ref{subsubsec:roll_max} --- на \cite{s2r2_sym}.

\subsection{Субриманова задача на группе Энгеля}\label{subsec:engel}
\subsubsection{Постановка задачи}\label{subsubsec:engel_state}
\paragraph{Геометрическая постановка}
Пусть на евклидовой плоскости заданы точки $a_0, \, a_1 \in \R^2$, соединенные кривой $\g_0 \subset \R^2$. Пусть также заданы число $S \in \R$ и прямая $L \subset \R^2$. Требуется соединить точки $a_0$, $a_1$
кратчайшей кривой $\g \subset \R^2$ так, чтобы кривые $\g_0$ и $\g$ ограничивали на плоскости область алгебраической площади $S$ с центром масс, принадлежащим прямой $L$.
Таким образом, это некоторое обобщение  (усложнение) задачи Дидоны, см.~раздел~\ref{subsec:heis}.

\paragraph{Задача оптимального управления}
Поставленную геометрическую задачу можно переформулировать как задачу оптимального управления
\begin{align}
&\dg = u_1 X_1(g) + u_2 X_2(g), \quad g = (x, y, z, v) \in \R^4, \label{engelp21}\\
&g(0) = g_0, \quad g(t_1) = g_1, \label{engelp22}\\
&l=\int_0^{t_1}\sqrt{u_1^2 + u_2^2}\,dt \to \min, \label{engelp23}\\
&X_1 = \dfrac{\partial}{\partial x} - \dfrac{y}{2}\dfrac{\partial}{\partial z} , \quad 
X_2 = \dfrac{\partial}{\partial y} + \dfrac{x}{2}\dfrac{\partial}{\partial z} + \dfrac{x^2 +y^2}{2}\dfrac{\partial}{\partial v}. \label{engelX12}
\end{align}
Эта задача --- субриманова  для субримановой структуры на $\R^4$, заданной векторными полями $X_1$, $X_2$ как ортонормированным репером.

\paragraph{Алгебра Энгеля и группа Энгеля}
\ddef{Алгеброй Энгеля} называется  алгебра Ли  $\gg$, в которой существует базис $(X_1, \dots, X_4)$, в котором ненулевые коммутаторы  суть
$$
[X_1, X_2] = X_3, \quad [X_1, X_3] = X_4,
$$
см.~Рис.~\ref{fig:engel}.

\begin{figure}[htb]
\setlength{\unitlength}{1cm}

\begin{center}
\begin{picture}(4, 4)
\put(1.1, 2.9){ \vector(1, -1){0.8}}
\put(1, 2.9){ \vector(0, -1){1.8}}
\put(2.9, 2.9){ \vector(-1, -1){0.8}}
\put(1.9, 1.9){ \vector(-1, -1){0.8}}

\put(1, 1){ \circle*{0.15}}
\put(1, 3){ \circle*{0.15}}
\put(3, 3){ \circle*{0.15}}
\put(2, 2){ \circle*{0.15}}

\put(1, 0.5) {$X_4$}
\put(1, 3.3) {$X_1$}
\put(3, 3.3) {$X_2$}
\put(2, 2.4) {$X_3$}

\end{picture}

\caption{ Алгебра Энгеля\label{fig:engel}}

\end{center}
\end{figure}
Алгебра Энгеля есть нильпотентная алгебра Ли с  градуировкой
$\gg = \gg^{(1)} \oplus\gg^{(2)} \oplus \gg^{(3)}$, $ \gg^{(1)} = \spann(X_1, X_2)$, $ \gg^{(2)} = \R X_3$, $ \gg^{(3)}=\R X_4$, $[\gg^{(1)}, \gg^{(i)}] =  \gg^{(i+1)} $, 
$\gg^{(4)} = \{0\}$,
поэтому она является алгеброй Карно. Соответствующая связная односвязная группа Ли $G$ называется \ddef{группой Энгеля}.

Группа Энгеля имеет линейное представление:
\begin{align*}
\left\{
\begin{pmatrix}
1 & b & c & d \\
0 & 1 & a & a^2/2\\
0 & 0 & 1 & a \\
0 & 0 & 0 & 1
\end{pmatrix}
\mid a, b, c, d \in \R
\right\}.
\end{align*}
На пространстве $\R^4_{x, y, z, v}$ можно ввести закон умножения
\begin{align*}
\begin{pmatrix}
x_1\\
y_1\\
z_1\\
v_1
\end{pmatrix} \cdot
\begin{pmatrix}
x_2\\
y_2\\
z_2\\
v_2
\end{pmatrix} = 
\begin{pmatrix}
x_1 + x_2\\
y_1 + y_2\\
z_1 + z_2 + (x_1 y_2 - x_2 y_1)/2\\
v_1 + v_2 + y_1y_2(y_1 + y_2)/2 + x_1z_2 + x_1y_2(x_1 + x_2)/2
\end{pmatrix},
\end{align*}
превращающий это пространство в группу Энгеля: $G \cong \R^4_{x,y,z,v}$, а поля \eq{engelX12} в левоинвариантные поля на этой группе. Таким образом, задача \eq{engelp21}--\eq{engelp23} есть левоинвариантная субриманова задача на группе Энгеля. Поэтому можно считать, что начальная точка в \eq{engelp22} есть единица группы Энгеля: $g_0 = \Id = (0, 0, 0, 0)$.

Все вполне неголономные левоинвариантные субримановы задачи ранга 2 на группе Энгеля переводятся друг в друга изоморфизмом этой группы \cite{symmetry}.

\paragraph{Особенности задачи}
Субриманова задача на группе Энгеля есть простейшая левоинвариантная субриманова задача со следующими свойствами:
\begin{itemize}
\item
она имеет глубину $3$,  ее вектор роста равен $(2, 3, 4)$,
\item
она имеет нетривиальные анормальные кратчайшие,
\item
ее геодезические параметризуются неэлементарными функциями (эллиптическими функциями Якоби),
\item
ее субриманова сфера несубаналитична.
\end{itemize}

Эта задача доставляет нильпотентную аппроксимацию любой \ddef{субримановой задачи энгелева типа} (то есть с вектором роста $(2, 3, 4)$, см. раздел \ref{subsec:class4}), в частности, для мобильного робота с прицепом.

\subsubsection{Симметрии распределения и субримановой структуры}\label{subsubsec:engel_inf_sym}
\begin{theorem}
Алгебра Ли инфинитезимальных симметрий распределения $\spann(X_1, X_2)$ на группе Энгеля  параметризуется гладкими функциями на этой группе, постоянными вдоль поля $X_2$.
\end{theorem}

\begin{theorem}
Алгебра Ли инфинитезимальных симметрий нильпотентной  субримановой структуры на группе Энгеля 
изоморфна алгебре Энгеля и состоит из правоинвариантных векторных полей на этой группе.
\end{theorem}

\subsubsection{Геодезические}\label{subsubsec:engel_geod}
Существование оптимальных управлений в задаче \eq{engelp21}--\eq{engelp23} следует из теорем Рашевского-Чжоу и Филиппова.

\paragraph{Принцип максимума Понтрягина}
Перейдем от задачи минимизации длины \eq{engelp23} к эквивалентной задаче минимизации энергии 
\begin{align}\label{engelp24}
&J = \frac{1}{2}\int_0^{t_1} (u_1^2 + u_2^2) dt \to \min.
\end{align}

Введем линейные на слоях $T^*G$ гамильтонианы 
$h_i(\lam) = \langle\lam, X_i\rangle$, $i = 1, \dots, 4$. Тогда  принцип максимума Понтрягина для задачи \eq{engelp21}, \eq{engelp22}, \eq{engelp24} принимает форму:
\begin{align*}
&\dh_1 = -u_2 h_3,\\
&\dh_2 = u_1 h_3,\\
&\dh_3 = u_1 h_4,\\
&\dh_4 = 0,\\
&\dg = u_1 X_1 + u_2 X_2,\\
&u_1 h_1 + u_2 h_2 + \dfrac{\nu}{2}(u_1^2 + u_2^2) \to \max_{(u_1, u_2) \in \R^2},\\
& \nu \leq 0,\\
& (h_1, \dots, h_4, \nu) \neq 0.
\end{align*}
\paragraph{Анормальные экстремали}
Анормальные экстремали постоянной скорости могут быть параметризованы как
\begin{align}
&h_1= h_2 = h_3 = 0, \quad h_4 \equiv \const \neq 0, \nonumber\\
&u_1\equiv 0, \quad u_2\equiv \pm 1, \nonumber\\
&x = z \equiv 0,  \quad y = \pm \, t,  \quad v = \pm \,\dfrac{t^3}{6}. \label{engel_abn}
\end{align}
Анормальные траектории \eq{engel_abn} суть однопараметрические подгруппы $g(t) = e^{\pm \, tX_2}$. Они проецируются на плоскость $(x, y)$ в прямые, потому являются субримановыми кратчайшими. Анормальное множество есть одномерное гладкое многообразие, диффеоморфное прямой:
$$
\Abn = \{g \in G \mid x = z = v - y^3/6 = 0 \}.
$$
\paragraph{Нормальные экстремали}
Нормальные экстремали являются траекториями нормальной  гамильтоновой системы
\begin{align}\label{engel_Ham}
\dot{\lam} = \Vec{H}(\lam), \quad \lam \in T^*G, 
\end{align}
с гамильтонианом $H = \dfrac{1}{2}(h_1^2 + h_2^2).$ Введем на поверхности уровня $\{H =1/2\}$ координаты $(\t, c, \a)$: 
\be{}\nonumber
h_1 = -\sin\t, \quad  h_2 = \cos\t, \quad   h_3 = c, \quad h_4 = \a,
\ee
тогда гамильтонова система \eq{engel_Ham} примет форму
\begin{align}
&\dot{\t} = c, \quad \dc = - \a \sin\t, \quad \dot{\a} = 0, \label{engel_vert}\\
&\dg =  -\sin \t\, X_1 + \cos\t\,X_2.\label{engel_hor}
\end{align}
Вертикальная подсистема \eq{engel_vert} есть уравнение маятника в поле силы тяжести с ускорением  свободного падения $g = \a l$, где $l$ --- длина маятника. Таким образом, при $\a> 0$ ($\a < 0$) сила тяжести направлена  вниз (вверх) относительно оси, от которой отсчитывается угол $\t$, а при $\a = 0$ маятник движется в невесомости.

Проекции нормальных экстремалей на плоскость $(x, y)$ суть эйлеровы эластики, см.~раздел~\ref{subsec:elastica}.

Анормальные кратчайшие  удовлетворяют нормальной гамильтоновой системе \eq{engel_vert}, \eq{engel_hor} при $\t = \pi + 2\pi n$, $c = 0$, поэтому они нестрого анормальны.
\paragraph{Симплектическое слоение и функции Казимира}
На коалгебре Ли $\gg^*$ существуют 2 независимые функции Казимира:
\begin{align*}
h_4,  \quad E = \dfrac{h_3^2}{2}  - h_2 h_4,
\end{align*}
где
$E$ есть полная энергия маятника \eq{engel_vert}.

Симплектическое слоение на $\gg^*$ состоит из:
\begin{itemize}
\item
параболических цилиндров
$$
\{E = \const, \quad h_4 = \const \neq 0, \quad h_3^2 + h_4^2 \neq 0  \},
$$
\item
пар плоскостей
$$
\{E = \const, \quad h_4 = 0, \quad h_3 \neq 0  \},
$$
\item
точек
$$
\{h_i = \const, \quad i = 1, \dots, 4, \quad h_3^2 + h_4^2 =0  \}.
$$
\end{itemize}
Симплектические листы $2$-мерны и $0$-мерны, потому вертикальная подсистема \eq{engel_vert} интегрируема по Лиувиллю. Фазовый портрет гамильтоновой системы \eq{engel_Ham} на цилиндре 
$C\cap \{h_4 = \const\}$, где $C = \gg^* \cap \left\{H = \frac{1}{2} \right\}$,
получается пересечением этого цилиндра с поверхностью уровня энергии $E$.
\paragraph{Параметризация нормальных геодезических}
Семейство нормальных экстремалей на поверхности уровня $\{H = \frac{1}{2} \}$ параметризуется начальными точками, принадлежащими цилиндру $C$.

Рассмотрим стратификацию цилиндра $C$ на подмногообразия, соответствующие разным типам траекторий маятника \eq{engel_vert}:
\begin{align*}
&C = \bigsqcup_{i=1}^7 C_i,  \\
&C_1 = \{\lambda \in C \mid \alpha \neq 0, E\in(- |\alpha|, 
|\alpha|)\}, \\
&C_2 = \{\lambda \in C \mid \alpha \neq 0, E\in(|\alpha|,+\infty)\}, 
\\
&C_3 = \{\lambda \in C \mid \alpha \neq 0, E=|\alpha|, c \neq 0 \}, \\
&C_4 = \{\lambda \in C \mid \alpha \neq 0, E=-|\alpha|\}, \\
&C_5 = \{\lambda \in C \mid \alpha \neq 0, E=|\alpha|, c = 0\}, \\
&C_{6} = \{\lambda \in C \mid \alpha = 0, \ c \neq 0\}, \\
&C_7 = \{\lambda \in C \mid \alpha = c = 0\}. 
\end{align*}
Далее, множества $C_i, \, i=1, \dots, 5,$ разбиваются на подмножества 
в зависимости от знака переменной $\alpha$:
\begin{align*}
&C_i^+ = C_i \cap \{\alpha>0\}, \qquad C_i^- = C_i \cap 
\{\alpha<0\}, \qquad i\in\{1, \dots, 5\}.
\end{align*}
 
Более того, подмножества $C_6$, $C_2^{\pm}, C_3^{\pm}$ разбиваются на 
связные компоненты в зависимости от знака переменной $c$:
\begin{align*}
&C_{6+} = C_6 \cap \{c>0\}, \quad C_{6-} = C_6 \cap \{c<0\}, \\
&C_{i+}^{\pm} = C_i^{\pm} \cap \{c>0\}, \quad C_{i-}^{\pm} = C_i^{\pm} 
\cap \{c<0\}, \quad i\in\{2,3\}. 
\end{align*}

Для нормализации нормальных  геодезических введем на стратах $C_1$, $C_2$, $C_3$ эллиптические координаты $(\f, k, \a)$, в которых уравнение маятника \eq{engel_vert} выпрямляется.

\medskip\noindent
В области $C_1^+$
\begin{align*}
&k=\sqrt{\frac{E+\alpha}{2 \alpha}} = \sqrt{\frac{c^2}{4 
\alpha}+\sin^2 \frac{\theta}{2}}\in (0,1),\\
&\sin\frac{\theta}{2} = k \sn (\sqrt{\alpha} \varphi), 
&&\cos\frac{\theta}{2} = \dn (\sqrt{\alpha} \varphi), \\
&\frac {c}{2} = k \sqrt{\alpha} \cn (\sqrt{\alpha} \varphi), 
&&\varphi 
\in [0,4 K].
\end{align*}

\medskip\noindent
В области $C_2^+$
\begin{align*}
&k=\sqrt{\frac{2 \alpha}{E+\alpha}} = \frac{1}{\sqrt{\frac{c^2} {4 
\alpha}+\sin^2 \frac{\theta}{2} }}\in (0,1),\\
&\sin\frac{\theta}{2} = \sgn c \, \sn \frac{\sqrt{\alpha} \varphi}{k}, 
&&\cos\frac{\theta}{2} = \cn \frac{\sqrt{\alpha} \varphi}{k}, \\
&\frac{c}{2} = \sgn c \, \frac {\sqrt{\alpha}}{k}  \dn \frac 
{\sqrt{\alpha} \varphi}{k}, &&\varphi \in [0,2 k K],\\
&\psi = \frac {\varphi}{k}. 
\end{align*}

\medskip\noindent
На множестве $C_3^+$
\begin{align*}
&k=1,\\
&\sin\frac{\theta}{2} = \sgn c \,\th (\sqrt{\alpha} \varphi), 
&&\cos\frac{\theta}{2} = \frac {1}{\ch (\sqrt{\alpha} \varphi)}, \\
&\frac{c}{2} = \sgn c \, \frac {\sqrt{\alpha}}{\ch (\sqrt{\alpha} 
\varphi)}, &&\varphi \in (-\infty, +\infty).
\end{align*}

На множествах $C_1^-, C_2^-, C_3^-$ определим новые координаты 
следующим образом:
\begin{align}
&\varphi (\theta, c, \alpha) = \varphi (\theta - \pi, c, -\alpha), 
\label{phi_a}\\
&k(\theta, c, \alpha) = k (\theta-\pi, c, -\alpha). \label{k_a}
\end{align}

Вертикальная подсистема \eq{engel_vert} принимает в новых 
координатах следующую форму:
\begin{align*}
&\dot{\varphi}=1, \qquad \dot{k}=0, \qquad  \dot{\alpha}=0,
\end{align*}
поэтому ее решения имеют вид
\begin{align}
&\varphi_t = \varphi + t, \qquad  k=\const, \qquad 
\alpha=\const. \label{phit}
\end{align}
Задача инвариантна относительно левых сдвигов на группе Энгеля, а также дилатаций
\begin{align}
&\d_s : \, (t, x, y, z, v) \mapsto (e^s t, e^s x, e^s y, e^{2s} z, e^{3s} v), \label{engel_dil1}\\
&(\t, c, \a) \mapsto (\t, e^{-s} c, e^{-2s}\a), \quad (\f, k, \a) \mapsto (e^s \f, k, e^{-2s} \a),\label{engel_dil2}
\end{align}
и отражений
\begin{align*}
&(t, x, y, z, v) \mapsto (t, -x, -y, z, -v),\\
&(\t, c, \a) \mapsto (\t - \pi, c, -\a), \quad (\f, k, \a) \mapsto (\f, k, -\a).
\end{align*}
Дилатации на группе Энгеля задают поток векторного поля
$$
Y = x \dfrac{\partial}{\partial x} + y \dfrac{\partial}{\partial y} + 2z \dfrac{\partial}{\partial z} + 3v \dfrac{\partial}{\partial v}.
$$

При $\lam = (\f, k, \a) \in \cup_{i =1}^3 C_i$, $\a = 1$, геодезические параметризуются следующим образом.

Если $\lam \in C_1$, то 
\begin{align}
&x_t = 2 k (\cn \varphi_t - \cn \varphi), \nonumber \\
&y_t = 2 \big(\E(\varphi_t) - \E(\varphi)\big)-t, \nonumber\\
&z_t = 2 k \big(\sn \varphi_t \dn \varphi_t - \sn \varphi \dn \varphi 
- \frac{y_t}{2} (\cn \varphi_t + \cn \varphi)\big), \nonumber\\
&v_t= \frac{y_t^3}{6} + 2 k^2 \cn^2 \varphi y_t - 4 k^2 \cn \varphi 
(\sn \varphi_t \dn \varphi_t - \sn \varphi \dn \varphi) +  \nonumber \\
&\quad+ 2 k^2 \bigg(\frac{2}{3} \cn \varphi_t \dn \varphi_t \sn 
\varphi_t- \frac{2}{3} \cn \varphi \dn \varphi \sn \varphi + \frac{1-
k^2}{3 k^2} t  +\frac{2 k^2 -1}{3 k^2}\Big(\E(\varphi_t)-\E(\varphi)\Big)\bigg).\label{ExpC1}
\end{align}

\medskip\noindent
Если $\lam \in C_2$, то 
\begin{align}
&x_t = \frac {2 \sgn c} {k} \Big(\dn \psi_t - \dn \psi \Big), \nonumber\\
&y_t = \frac {k^2-2}{k^2} t + \frac {2}{k} \Big(\E(\psi_t) - \E(\psi)\Big), \nonumber\\
&z_t = - \frac{x_t y_t}{2} - \frac{2 \sgn c \dn \psi}{k} y_t + 2 \sgn c \, (\cn \psi_t \sn \psi_t - \cn \psi \sn \psi), \nonumber\\
&v_t= \frac {4}{k} \bigg( \frac {1}{3}\cn \psi_t \dn \psi_t \sn \psi_t 
- \frac{1}{3} \cn \psi \dn \psi \sn \psi -\frac{1-k^2}{3 k^3} t - \frac {k^2-2}{6 k^2} \Big(\E(\psi_t)-\E(\psi)\Big)\bigg) +\nonumber\\
&\quad + \frac {y_t^3}{6} + \frac {2 y_t}{k^2} \dn^2 \psi - \frac{4}{k} \dn \psi \big(\cn \psi_t \sn \psi_t - \cn \psi \sn 
\psi \big), \nonumber \\
&\psi=\frac{\varphi}{k}, \quad \psi_t = \psi + \frac{t}{k}. \label{ExpC2}
\end{align}

\medskip\noindent
Если $\lam \in C_3$, то 
\begin{align}
&x_t = 2 \sgn c \left(\frac {1}{\ch \varphi_t} - \frac{1}{\ch \varphi}\right), \nonumber\\
&y_t = 2 (\tangh \varphi_t - \tangh \varphi) - t , \nonumber\\
&z_t = - \frac{x_t y_t}{2} - \frac{2 \sgn c}{\ch \varphi}\, y_t + 2 \sgn 
c\left(\frac{\tangh \varphi_t}{\ch \varphi_t} - \frac{\tangh \varphi}{\ch 
\varphi}\right), \nonumber \\
&v_t= \frac{2}{3}\left(\tangh \varphi_t - \tangh \varphi + 2 \frac{\tangh 
\varphi_t}{\ch^2 \varphi_t}-2\frac{\tangh \varphi}{\ch^2 \varphi}\right) 
+ \frac{y_t^3}{6} + \frac {2 y_t}{\ch^2 \varphi}  -  \frac{4}{\ch \varphi} \left(\frac{\tangh \varphi_t}{\ch 
\varphi_t} - \frac{\tangh \varphi}{\ch \varphi}\right). \label{ExpC3}
\end{align}
Параметризация геодезических для произвольных $\lam = (\f, k, \a) \in \cup_{i =1}^3 C_i$  получается из случая $\a = 1$ с помощью дилатаций и отражения:
\begin{itemize}
\item
если $\a > 0$, то
\begin{align*}
(x_t, y_t, z_t, v_t)(\f, k, \a) = \left(\dfrac{x_{t'}}{\a^{1/2}}, \dfrac{y_{t'}}{\a^{1/2}}, \dfrac{z_{t'}}{\a}, \dfrac{v_{t'}}{\a^{3/2}}\right)(\sqrt \a \, \f, k, 1), \quad t' = t\sqrt \a,
\end{align*}
\item
если $\a < 0$, то
\begin{align*}
(x_t, y_t, z_t, v_t)(\f, k, \a) = (-x_t, -y_t, z_t, -v_t)(\f, k, -\a).
\end{align*}
\end{itemize}
В оставшихся случаях $\lam \in \cup_{i=4}^7C_i$ геодезические параметризуются элементарными функциями.

\medskip\noindent
Если $\lam \in C_4$, то
\begin{align*}
&x_t = 0, \qquad y_t = t \sgn \alpha , \qquad z_t = 0, \qquad v_t = 
\frac{t^3}{6} \sgn \alpha. 
\end{align*}

\medskip\noindent
Если $\lam \in C_5$, то 
\begin{align*}
&x_t = 0, \qquad y_t = - t \sgn \alpha, \qquad z_t = 0, \qquad v_t = - 
\frac{t^3}{6} \sgn \alpha. 
\end{align*}

\medskip\noindent
Если $\lam \in C_6$, то 
\begin{align*}
&x_t = \frac{\cos (c t + \theta) - \cos \theta}{c}, &&y_t = \frac{\sin 
(c t + \theta) - \sin \theta}{c}, \nonumber\\
&z_t = \frac{c t - \sin c t}{2 c^2}, 
&&v_t = \frac{3 \cos \theta - 2 c t \sin \theta - 4 \cos (c t + \theta) + \cos (2 c t + \theta)}{4 c^3}.  
\end{align*}

\medskip\noindent
Если $\lam \in C_7$, то 
\begin{align*}
&x_t = - t \sin \theta, \qquad y_t = t \cos \theta, \qquad z_t = 0, \qquad v_t = \frac{t^3}{6} \cos \theta. 
\end{align*}

Проекции геодезических на плоскость $(x, y)$ суть эйлеровы эластики (см. раздел \ref{subsec:elastica}): инфлексионные при $\lam \in C_1$,  неинфлексионные при $\lam \in C_2$, критические при $\lam \in C_3$, прямые при $\lam \in C_4 \cup C_5 \cup C_7$, и окружности при $\lam \in C_6$.

Семейство всех геодезических параметризуется экспоненциальным отображением
\begin{eqnarray*}
&&\Exp \colon N = C \times \R_+ \to M,\\
&&\Exp (\lambda, t) = g_t = (x_t, y_t, z_t, v_t).
\end{eqnarray*}

\subsubsection{Симметрии экспоненциального отображения и время Максвелла}\label{subsubsec:engel_max}
Дилатации \eq{engel_dil1}, \eq{engel_dil2} образуют однопараметрическую группу симметрий экспоненциального отображения. Имеется также дискретная группа симметрий, образованная отражениями:
\begin{align*}
&\Sym = \{\Id, \eps^1, \dots, \eps^7 \} \cong \Z_2\times\Z_2\times\Z_2.
\end{align*}

Обозначим через $\vH_v = c\dfrac{\partial}{\partial \t} - \a \sin \t \dfrac{\partial}{\partial c} \in \VVec(C)$ вертикальную часть нормального гамильтонова поля~$\vH$.
Следующие отображения $\eps^i: C \to C$ сохраняют поле направлений векторного поля $\vH_v$:
\begin{align*}
&\varepsilon^1 : (\theta,c, \alpha) ~ \mapsto ~ (\theta, -c, \alpha), &&\varepsilon^2 : (\theta,c, \alpha) ~ \mapsto ~ (-\theta, c, \alpha), \\
&\varepsilon^3 : (\theta,c, \alpha) ~ \mapsto ~ (-\theta, -c, \alpha),  &&\varepsilon^4 : (\theta,c, \alpha) ~ \mapsto ~ (\theta+\pi, c, -\alpha), \\
&\varepsilon^5 : (\theta,c, \alpha) ~ \mapsto ~ (\theta+\pi, -c, -\alpha), &&\varepsilon^6 : (\theta,c, \alpha) ~ \mapsto ~ (-\theta+\pi, c, -\alpha), \\
&\varepsilon^7 : (\theta,c, \alpha) ~ \mapsto ~ (-\theta+\pi, -c, -\alpha). 
\end{align*}
А именно: $\eps^i_*\vH_v = \vH_v$ при $i = 3, 4, 7$, и $\eps^i_*\vH_v = -\vH_v$ при $i =1, 2, 5, 6$. Действие отражений $\eps^i: C \to C$ продолжается до симметрий экспоненциального отображения следующим образом.

Действие $\eps^i: N \to N$, $N = C \times \R_{+}$, определяется как 
\[
{\varepsilon}^i (\lambda, t) = 
 \begin{cases}
 \big({\varepsilon}^i (\lambda), t\big), & \textrm{если }{\varepsilon}^i_*  \vec{H}_v = \vec{H}_v, \\
 \big({\varepsilon}^i \circ e^{t \vec{H}_v} (\lambda), t\big), & \textrm{если } {\varepsilon}^i_* \vec{H}_v = - \vec{H}_v.
\end{cases}
\]

Действие $\eps^i: G \to G$ определяется как 
\begin{align*}
&\varepsilon^i(q) = \varepsilon^i(x,y,z,v) = g^i = (x^i, y^i, z^i, v^i), \\
&(x^1, y^1, z^1, v^1) = (x, y, -z, v - x z), \\
&(x^2, y^2, z^2, v^2) = (-x, y, z, v - x z), \\
&(x^3, y^3, z^3, v^3) = (-x, y, -z, v), \\
&(x^4, y^4, z^4, v^4) = (-x, -y, z, -v), \\
&(x^5, y^5, z^5, v^5) = (-x, -y, -z, -v + x z), \\
&(x^6, y^6, z^6, v^6) = (x, -y, z, -v + x z), \\
&(x^7, y^7, z^7, v^7) = (x, -y, -z, -v). 
\end{align*}

\begin{proposition}
Группа $\Sym$  есть подгруппа группы симметрий экспоненциального отображения.
\end{proposition}

\begin{theorem}
Первое время Максвелла, соответствующее  группе симметрий $\Sym$, для почти всех геодезических выражается следующим образом:
\begin{align}
&\lambda \in C_1 & \then &\tmax = \min \big(2 p_z^1(k), 4 K(k)\big)/\sigma, \label{tmaxC1}\\
&\lambda \in C_2 & \then &\tmax = 2 k K(k) /\sigma, \label{tmaxC2}\\
&\lambda \in C_6 & \then &\tmax = 2 \pi / |c|, \label{tmaxC6}\\
&\lambda \in C_3 \cup C_4 \cup C_5 \cup C_7 & \then &\tmax=+\infty, \label{tmaxC4}
\end{align}
где 
$\sigma = \sqrt{|\alpha|}$, 
$p^1_z(k)\in \big(K(k), 3K(k)\big)$ есть первый положительный корень функции $f_z(p,k)=\dn p \,\sn p+ (p-2\E(p))\cn p$.
\end{theorem}
\begin{remark}
Для тех геодезических, для которых первое время Максвелла не равно $\tmax$, оно больше этого значения, а $\tmax$ есть первое сопряженное время.
\end{remark}
\begin{theorem}
Функция $\tmax : C \to (0, + \infty]$ имеет следующие свойства инвариантности:
\begin{itemize}
\item[$(1)$]
$\tmax(\lam)$ зависит только от значений $E$ и $|\a|$,
\item[$(2)$]
$\tmax(\lam)$  есть первый интеграл поля $\vH_v$,
\item[$(3)$]
$\tmax(\lam)$  инвариантно относительно отражений: если $(\lam, t) \in C\times \R_+$, $(\lam^i, t) = \eps^i(\lam, t)$, то $\tmax(\lam^i) = \tmax(\lam)$,
\item[$(4)$]
$\tmax(\lam)$ однородна относительно дилатаций: если $\lam \in C$, $\lam_s = \d_s(\lam) \in C$, то $\tmax(\lam_s) = e^s\tmax(\lam)$, $s \in \R$.
\end{itemize}
\end{theorem}

\subsubsection{Нижняя оценка сопряженного времени}\label{subsubsec:engel_conj}
\begin{theorem}
Для любого  $\lam \in C$ 
$$
\tconj(\lam) \geq \tmax(\lam).
$$
\end{theorem}

\subsubsection{Диффеоморфная структура экспоненциального отображения}\label{subsubsec:engel_cut}
Рассмотрим подмножество в пространстве состояний, не содержащее неподвижных точек симметрий $\eps^1$, $\eps^2$:
\begin{align*}
&\widetilde G = \{g \in G \mid \eps^1(g) \neq g \neq \eps^2(g) \}= \{g \in G \mid xz \neq 0  \},
\end{align*}
и его связные компоненты:
\begin{align*}
& G_1 = \{ g \in G ~ | ~ x < 0, z > 0 \}, \\ 
& G_2 = \{ g \in G ~ | ~ x < 0, z < 0 \},  \\ 
& G_3 = \{ g \in G ~ | ~ x > 0, z < 0 \},  \\ 
& G_4 = \{ g \in G ~ | ~ x > 0, z > 0 \}.  
\end{align*}
Также рассмотрим открытое плотное подмножество в пространстве всех потенциально оптимальных  геодезических:
\begin{align*}
 &\widetilde{N} = \{ (\lambda, t) \in N ~ | ~ t < \tmax (\lambda), \ c_{t/2} \sin \theta_{t/2} \ne 0 \}, 
\end{align*}
и его связные компоненты:
\begin{align*}
& D_1 = \{ (\lambda, t) \in N ~ | ~ t \in \big(0, \tmax(\lambda)\big), \ \sin\theta_{t/2} > 0, \ c_{t/2} > 0 \}, \\
& D_2 = \{ (\lambda, t) \in N ~ | ~ t \in \big(0, \tmax(\lambda)\big), \ \sin\theta_{t/2} > 0, \ c_{t/2} < 0 \},  \\
& D_3 = \{ (\lambda, t) \in N ~ | ~ t \in \big(0, \tmax(\lambda)\big), \ \sin\theta_{t/2} < 0, \ c_{t/2} < 0 \},  \\
& D_4 = \{ (\lambda, t) \in N ~ | ~ t \in \big(0, \tmax(\lambda)\big), \ \sin\theta_{t/2} < 0, \ c_{t/2} > 0 \}.  
\end{align*}
\begin{theorem}
Следующие отображения являются диффеоморфизмами:
\begin{align*}
&\Exp : D_i \to G_i, \quad i=1, \dots, 4, \\
&\Exp : \widetilde N \to \widetilde G.
\end{align*}
\end{theorem}

\subsubsection{Время разреза}
\begin{theorem}
Для любого $\lam \in C$
$$
\tcut(\lam)= \tmax(\lam).
$$
\end{theorem}

\subsubsection{Множество разреза и его стратификация}\label{subsubsec:engel_synth}
\begin{theorem}
Множество разреза $\Cut$ содержится в объединении координатных подпространств 
 $\{x = 0\}$ и $\{z = 0\}$. Оно инвариантно относительно дилатаций  и дискретных симметрий: 
\begin{align*}
e^{t Y} (\Cut) &= \Cut, \qquad t \in \R,\\
\varepsilon^i(\Cut) &= \Cut, \qquad i = 1, \dots, 7.
\end{align*}
\end{theorem}
\begin{theorem}
\label{th:cut_strat}
Множество разреза имеет стратификацию
\begin{align*}
\Cut =& (\I_{x+}  \sqcup \I_{x-}) \sqcup (\NN_{x+} \sqcup \NN_{x-}) \sqcup (\I_{z+}  \sqcup \I_{z-}) \sqcup\\ 
&\sqcup (\CI_{x+}^{+} \sqcup \CI_{x+}^{-} \sqcup \CI_{x-}^{+} \sqcup \CI_{x-}^{-}) 
 \sqcup (\CN_{x+}^{+} \sqcup \CN_{x+}^{-} \sqcup \CN_{x-}^{+} \sqcup \CN_{x-}^{-})\sqcup\\
&\sqcup (\CI_{z+}^{+} \sqcup \CI_{z+}^{-} \sqcup \CI_{z-}^{+} \sqcup \CI_{z-}^{-}) \sqcup\\
									 &\sqcup (\EE_{+} \sqcup \EE_{-}).
\end{align*}	
Пересечения множества разреза с координатными подпространствами имеют стратификации
\begin{align*}
\Cut \cap \{z=0\} =& (\I_{z+}  \sqcup \I_{z-}) \sqcup (\CI_{z+}^{+} \sqcup \CI_{z+}^{-} \sqcup \CI_{z-}^{+} \sqcup \CI_{z-}^{-}) \\
									 &\sqcup (\I_{x+}^0 \sqcup \I_{x-}^0)\sqcup (\EE_{+} \sqcup \EE_{-}), \\
\Cut \cap \{x=0\} =& (\I_{x+}  \sqcup \I_{x-})  \sqcup (\CI_{x+}^{+} \sqcup \CI_{x+}^{-} \sqcup \CI_{x-}^{+} \sqcup \CI_{x-}^{-}) \\
									&\sqcup (\NN_{x+} \sqcup \NN_{x-})\sqcup (\CN_{x+}^{+} \sqcup \CN_{x+}^{-} \sqcup \CN_{x-}^{+} \sqcup \CN_{x-}^{-}) \\
									&\sqcup (\I_{z+}^0 \sqcup \I_{z-}^0)\sqcup  (\EE_{+} \sqcup \EE_{-}), \\
\Cut \cap \{x=z=0\} =& (\I_{z+}^0 \sqcup \I_{z-}^0) \sqcup (\I_{x+}^0 \sqcup \I_{x-}^0) \sqcup (\EE_{+} \sqcup \EE_{-}).
\end{align*}
При этом $\I_{x\pm}^0 \subset \I_{x\pm}$, $\I_{z\pm}^0 \subset \I_{z\pm}$, а также
\begin{align*}
\I_{z+} &= \Big\{g \in G \mid z=0, y > Y_0^1 |x|, \ w < G_1(x,y)\Big\} \simeq \R^3, \\
\I_{x+} &= \Big\{g \in G \mid x=0, y > 0, \ w > G_2(z,y)\Big\} \simeq \R^3, \\
\NN_{x\pm} &= \Big\{g \in G \mid x=0, \sgn z =\pm 1, \ - G_3(z,-y) < w < G_3(z,y)\Big\} \simeq \R^3, \\
\CI_{z+}^{\pm} &= \Big\{g \in G \mid z=0, y > Y_0^1 |x|, \ w = G_1(x,y), \sgn x = \pm 1\Big\} \simeq \R^2, \\
\CI_{x+}^{\pm} &= \Big\{g \in G \mid z=0, y > 0, \ w = G_2(x,y), \sgn z = \pm 1\Big\} \simeq \R^2, \\
\CN_{x\pm}^{+} &= \Big\{g \in G \mid x=0, \sgn z =\pm 1, \ w = G_3(z,y)\Big\} \simeq \R^2, \\
\I_{z\pm}^0 &= \Big\{g \in G \mid x=z=0, \ y w < 0, \sgn y = \pm 1 \Big\} \simeq \R^2,\\
\I_{x\pm}^0 &= \Big\{g \in G \mid x=z=0, \ y w > 0, \sgn y = \pm 1 \Big\} \simeq \R^2,\\
\EE_{\pm} &= \Big\{g \in G \mid x=y=z=0, \ \sgn w = \pm 1 \Big\} \simeq \R^1, \\
\I_{z-} &= \varepsilon^4(\I_{z+}), \quad \I_{x-} = \varepsilon^4(\I_{z+}),  \\ 
\CI_{z-}^{\pm}&=\varepsilon^4(\CI_{z+}^{\pm}), \quad  \CI_{x-}^{\pm}=\varepsilon^4(\CI_{x+}^{\pm}), \quad  \CN_{x\pm}^{-} = \varepsilon^4(\CN_{x\pm}^{+}),
\end{align*}
где $Y_0^1 < 0$, а $G_i$, $i=1,2,3,$ --- некоторые гладкие функции, удовлетворяющие свойствам: 
\begin{align*}
G_1(0,y) = 0, \ &G_1(-x,y) = G_1(x,y), \  G_1(\rho x, \rho y) = \rho^3 G_1(x, y), \rho>0, \\
G_2(0,y) = 0, \ &G_2(-z,y) = G_2(z,y), \  G_2(\rho^2 z, \rho y) = \rho^3 G_2(z, y), \rho>0, \\
&G_3(-z,y) = G_3(z,y), \  G_3(\rho^2 z, \rho y) = \rho^3 G_3(z, y), \rho>0.
\end{align*}
\end{theorem}
Трехмерные страты $\I_{x\pm}$, $\I_{z\pm}$ (соотв. $\NN_{x\pm}$) состоят из точек, для которых проекции кратчайших на плоскость $(x,y)$ суть инфлексионные, т.е. имеющие точки перегиба (соотв. неинфлексионные, т.е. не имеющие точек перегиба) эластики, см.~раздел \ref{subsec:elastica}. Для одномерных стратов $\EE_{\pm}$ соответствующие эластики замкнуты (имеют форму восьмерки, <<figure-of-eight elastica>>). 

На Рис.~\ref{fig:dec1}, \ref{fig:dec2} изображены стратификации множества разреза и его пересечения с координатными подпространствами.   На Рис.~\ref{fig:dec1} показана топология примыкания стратов множества разреза в факторе по дилатациям $Y$. 
На Рис.~\ref{fig:dec2} представлено пересечение $\Cut \cap \{x=z=0\}$.

\figout{
\twofiglabelsize
{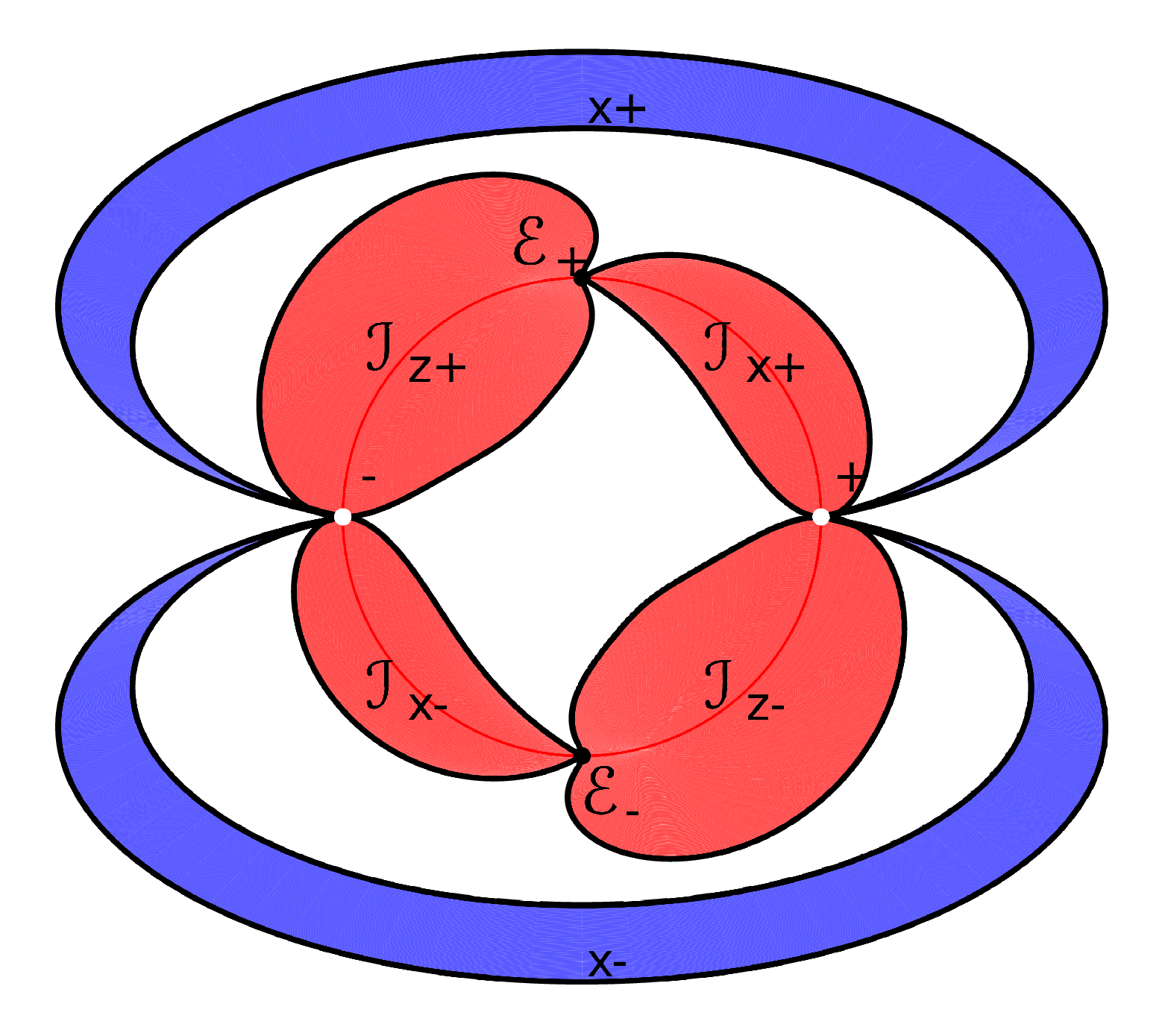}{Стратификация множества разреза: глобальная структура}{fig:dec1}{0.4}
{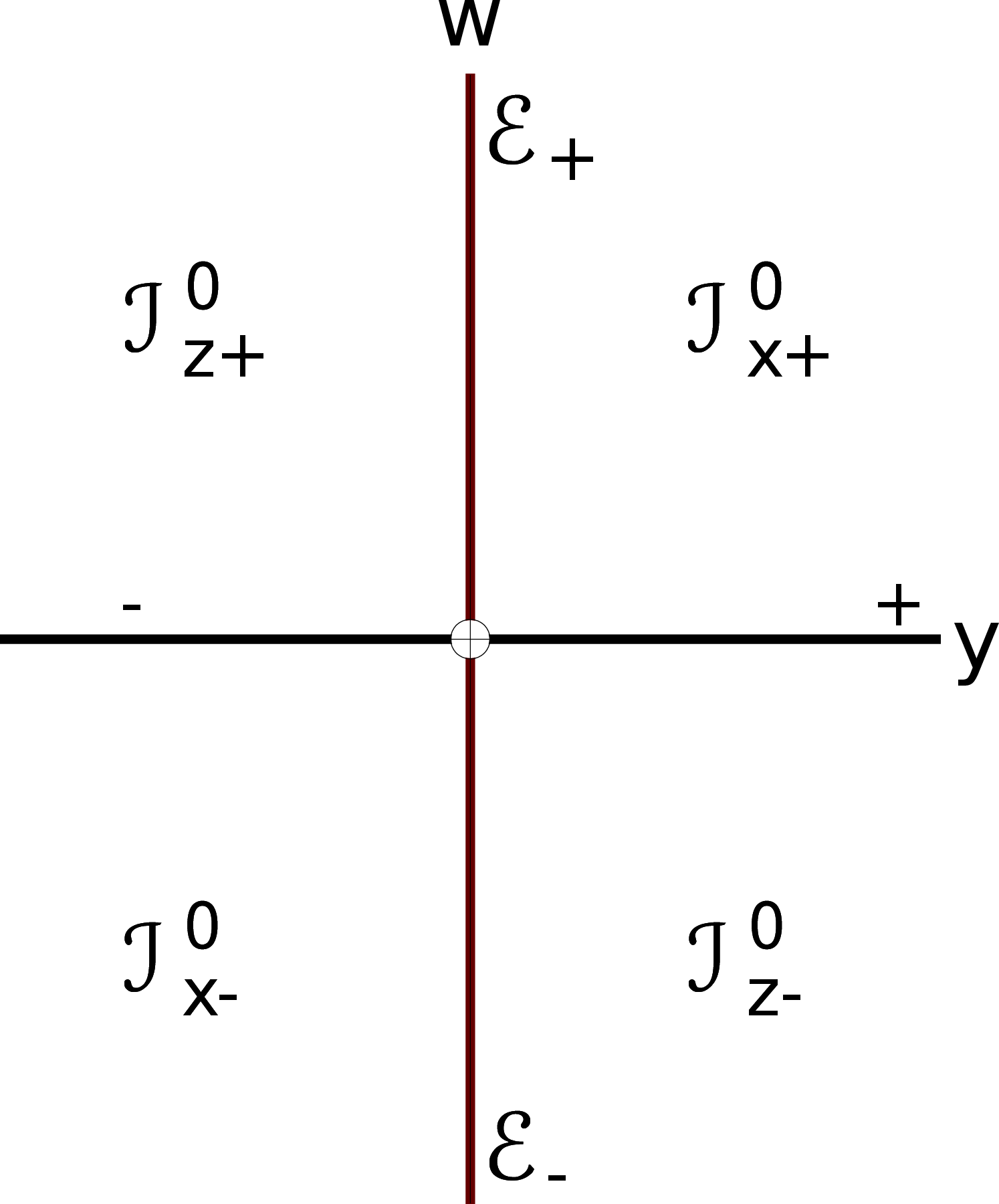}{Пересечение множества разреза с подпространством $\{x = z = 0\}$}{fig:dec2}{0.4}
}

На Рис.~\ref{fig:dec3} изображено множество $\Cut \cap \{z=0\}$ после факторизации по дилатациям $Y$; фактор $\{z=0\} / e^{\R Y}$ представлен топологической сферой $\{g \in G \mid x^6 + y^6 + w^2 = 1\}$.
Аналогично  на Рис.~\ref{fig:dec4} изображен фактор $(\Cut \cap \{x=0\}) / e^{\R Y}$ на топологической сфере $\{g \in G \mid y^6 + |z|^3 + w^2 = 1\}$. 

\figout{
\twofiglabelsize
{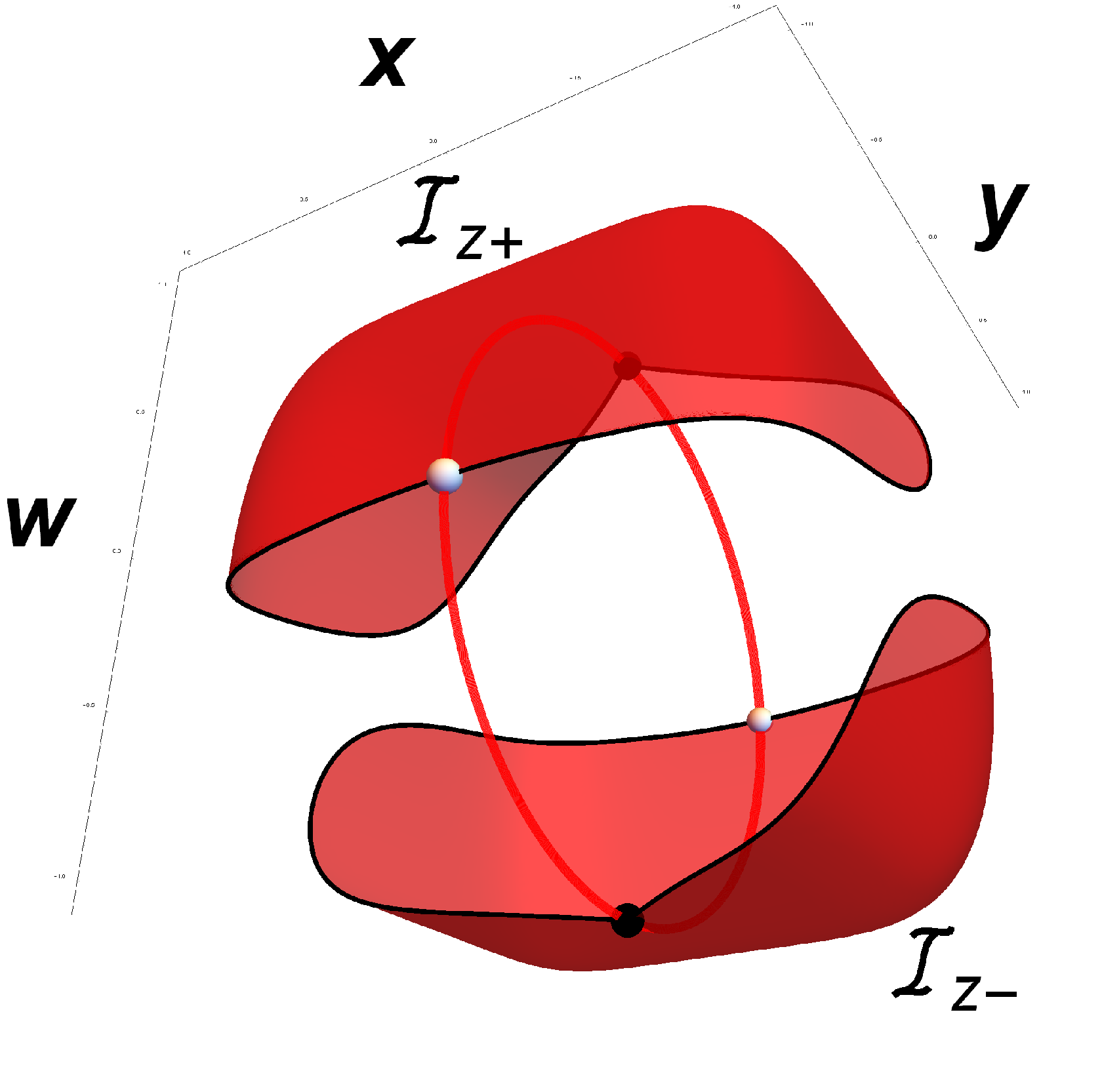}{Пересечение $\Cut \cap \{z=0\}$}{fig:dec3}{0.4}
{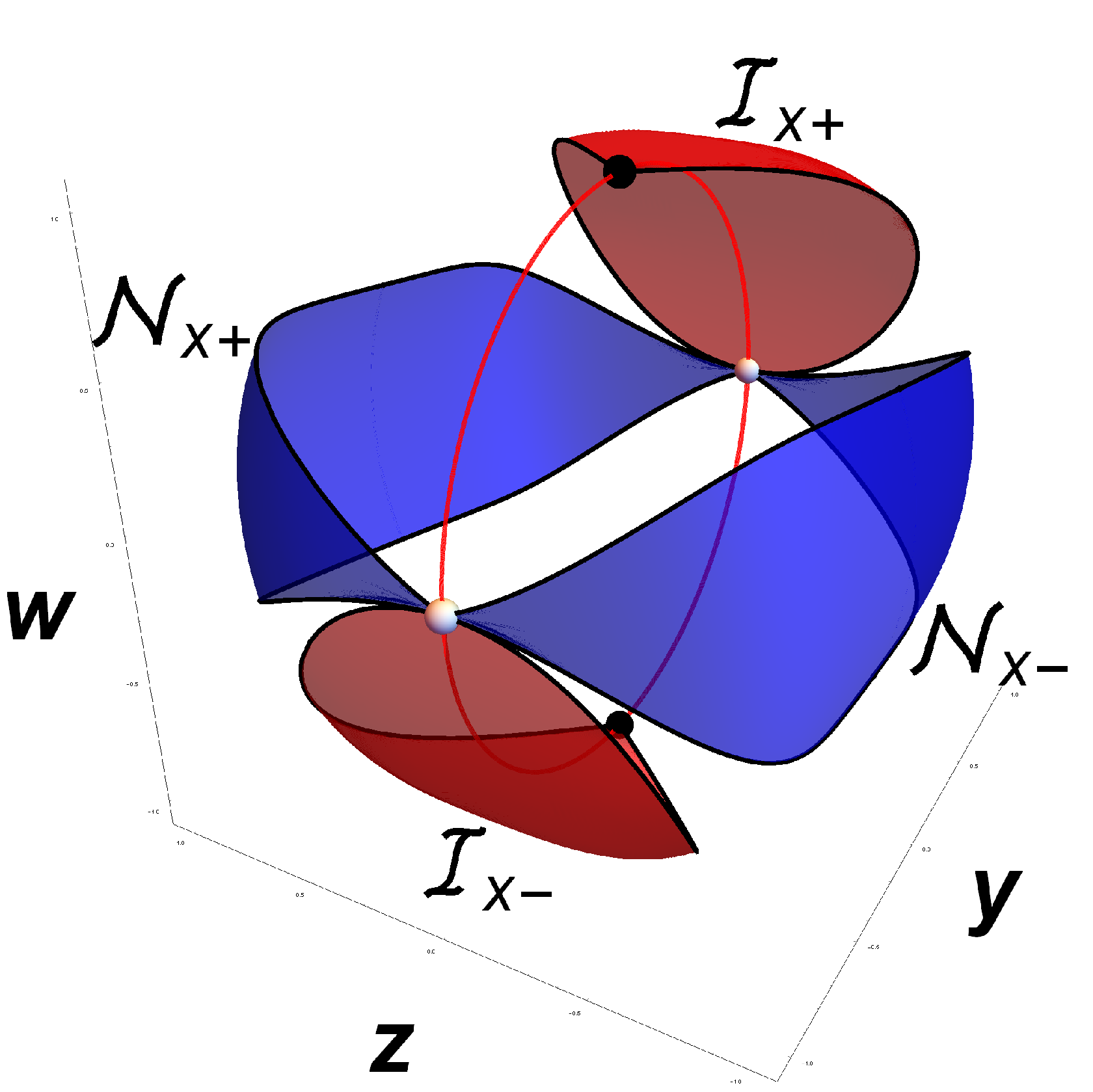}{Пересечение $\Cut \cap \{x=0\}$}{fig:dec4}{0.4}
}

Очевидно, что в каждую точку  $g_1 \in G \setminus \Cut$ приходит ровно одна субриманова кратчайшая. Ниже аналогичное свойство описано для точек $g_1 \in \Cut$. 

\begin{theorem}
\begin{itemize}
\item[$(1)$]
В каждую точку трехмерных стратов множества разреза приходят ровно две кратчайшие (эти страты состоят из точек Максвелла, не являющихся сопряженными). 
\item[$(2)$]
В каждую точку двумерных стратов приходит единственная кратчайшая (эти страты состоят из сопряженных точек, не являющихся точками Максвелла). 
\item[$(3)$]
В каждую точку одномерных стратов приходит однопараметрическое  семейство кратчайших (эти страты состоят из точек Максвелла, являющихся одновременно сопряженными точками).
\end{itemize}
\end{theorem}

Множество разреза незамкнуто т.к. оно содержит точки, сколь угодно близкие к начальной точке $q_0$, но не саму эту точку (это общий факт субримановой геометрии). Замыкание множества разреза в субримановой задаче на группе Энгеля допускает следующее простое описание.
\begin{theorem}
$\cl (\Cut) = \Cut  \sqcup \A_+ \sqcup \A_- \sqcup \{g_0\} $.
\end{theorem}

Примыкание анормальных траекторий $\A_{\pm}$  к стратам множества разреза изображено на Рис.~\ref{fig:dec1} слева.

\begin{theorem}
Имеют место стратификации
\begin{align*}
\Cut \cap \Conj &= \bigsqcup_{i \in \{+,-\}, j\in\{+,-\}} \Big(\mathcal{CI}_{zi}^j \sqcup \mathcal{CI}_{xi}^j \sqcup \mathcal{CN}_{xi}^j\Big) \sqcup \mathcal{E}_+\sqcup \mathcal{E}_-, \\
\Cut \cap \Max &= \bigsqcup_{i \in \{+,-\}} \bigg(\mathcal{I}_{zi} \sqcup \mathcal{I}_{xi} \sqcup \mathcal{N}_x^{i} \sqcup \mathcal{E}_i\bigg).
\end{align*}
\end{theorem}

\subsubsection{Сфера}\label{subsubsec:engel_sphere}
Субримановы сферы переходят друг в друга при левых сдвигах
$$
L_g(S_R(g_0)) = S_R(g g_0)
$$
и дилатациях
$$
\d_s(S_R(\Id))= S_{R'}(\Id), \quad R' = e^s R,
$$
поэтому достаточно исследовать единичную сферу $S = S_1(\Id)$.

Единичная сфера инвариантна относительно отражений:
$$
\eps^i(S)  = S, \quad i = 1, \dots,7.
$$
Рассмотрим сечение единичной сферы двумерным инвариантным многообразием основных симметрий $\eps^1$, $\eps^2$:
$$
\tS = \{ g \in S \mid \eps^1(g) = \eps^2(g) = g\} = S \cap \{ x = z = 0\},
$$
см. Рис. \ref{fig:engel_sphere}.

\figout{
\onefiglabelsize{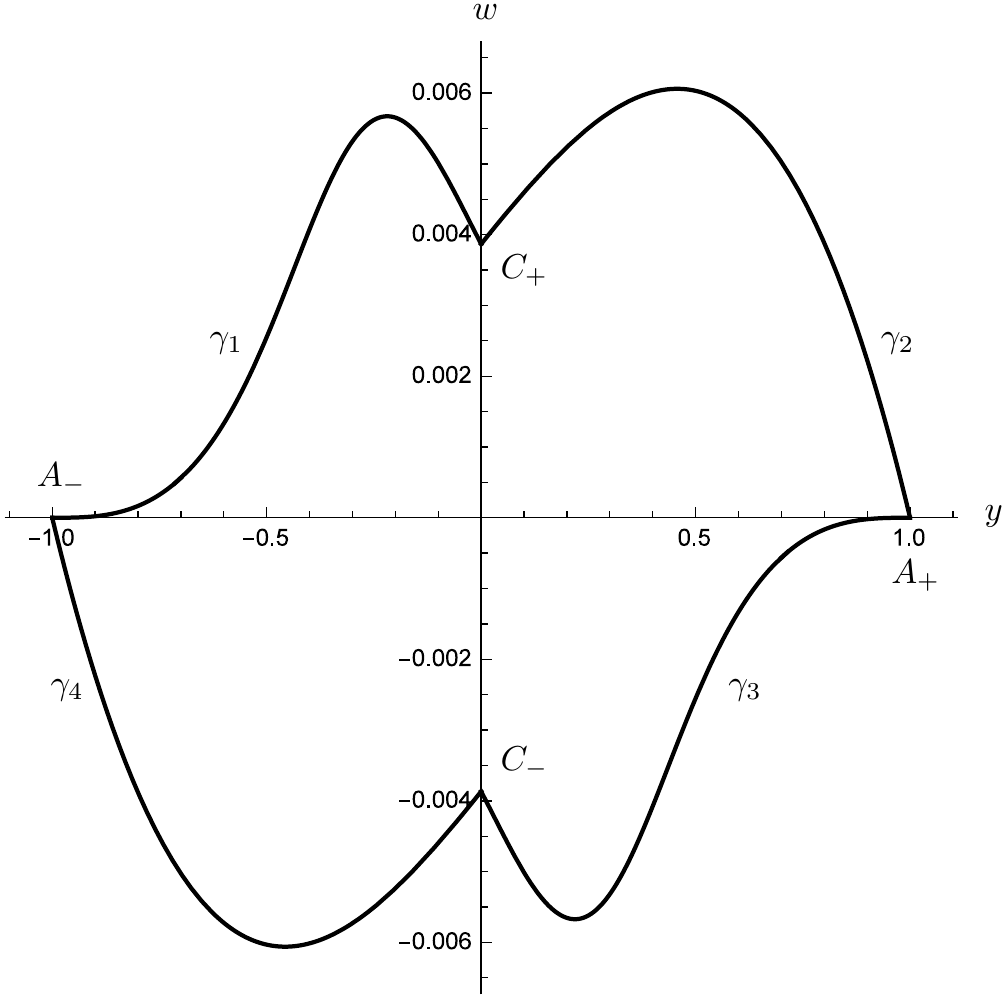}{Сечение сферы $\tS = S \cap \{x=z=0\}$}{fig:engel_sphere}{0.5}
}

Сечение $\tS$ центрально-симметрично в силу отражения $\eps^4$:
\begin{align*}
&\eps^4(\gamma_i)  = \gamma_{i+2}, \qquad i = 1, 2,\\
&\eps^4(A_+) = A_-, \qquad \eps^4(C_+) = C_-.
\end{align*}
Различные точки сечения $\tS$ можно охарактеризовать следующим образом:
\begin{itemize}
\item
$A_{\pm}$ --- точки на анормальных кратчайших,
\item
$C_{\pm}$ --- сопряженные точки, точки Максвелла, точки разреза,
\item
$g \in \gamma_i$ --- точки Максвелла, точки разреза.
\end{itemize}
Точки сечения $\tS$ имеют следующую кратность $\mu$  (количество кратчайших, приходящих из $\Id$  в эту точку):
\begin{itemize}
\item
$\mu(A_{\pm}) = 1$,
\item
$\mu(C_{\pm}) = \mathfrak{c}$ (континуум $\cong S^1$),
\item
$g \in \gamma_i \then \mu(g) = 2$.
\end{itemize}
\begin{theorem}
Сечение $\tS$  имеет следующую регулярность в различных своих точках:
\begin{itemize}
\item[$(1)$]
кривые $\gamma_i$  аналитичны и регулярны,
\item[$(2)$]
$A_{\pm}$,
$C_{\pm}$ --- особые точки, в них $\tS$  негладкая, но липшицева,
\item[$(3)$]
 $\overline{\gamma}_2 = \gamma_2 \cup \{ C_+, A_+\}$  гладкая класса $C^{\infty}$, 
\item[$(4)$]
 $\gamma_1 \cup \{ C_+\}$  гладкая класса $C^{\infty}$, 
\item[$(5)$]
 $\gamma_1 \cup \{ A_-\}$  гладкая класса $C^{1}$. 
\end{itemize}
\end{theorem}
\begin{theorem}
\begin{itemize}
\item[$(1)$]
Множество $\tS \setminus \{ A_+, A_-\}$  полуаналитично, потому субаналитично.
\item[$(2)$]
В окрестности точки $A_-$  кривая $\gamma_1$  есть график неаналитической функции 
$$
w = \frac 16 Y^3 - 4 Y^3 \exp(-2/Y) (1+ o(1)), \qquad Y = (y+1)/2 \to 0.
$$
\item[$(3)$]
Поэтому множество $\tS$ неполуаналитично, следовательно, несубаналитично.
\item[$(4)$]
Следовательно, сфера $S$  несубаналитична.
\end{itemize}
\end{theorem}

\begin{remark}
Утверждение о несубаналитичности сферы Энгеля $S$ следует также из проекции сферы Энгеля на (несубаналитическую) сферу Мартине (см.~раздел~\ref{subsec:martinet}).
\end{remark}
\begin{theorem}
В окрестности точки $A_-$  кривая $\gamma_1$   есть график   функции из $\exp$-$\log$  категории:
$$
w = F\left(Y, \frac{e^{-1/Y}}{Y}\right), \qquad Y = (y+1)/2 \to 0,
$$
где $F(\xi, \eta)$  есть аналитическая функция в окрестности точки $(\xi, \eta ) = (0, 0)$.

Поэтому множество $\tS$  принадлежит $\exp$-$\log$ категории.
\end{theorem}
\begin{theorem}
Разбиение $$\tS = \cup_{i=1}^4 \gamma_i \cup \{A_+, A_-, C_+, C_-\}$$  есть стратификация 
 Уитни.
\end{theorem}

\subsubsection{Явные выражения для субриманова расстояния}
Для некоторых точек группы Энгеля известно их субриманово расстояние  до единичного элемента:
\begin{itemize}
\item
Анормальная кратчайшая $g(t) = e^{\pm t X_2}$, $x = z= w = 0$, $y = \pm t$:
$$
d(\Id, g(t)) = t.
$$
\item
Центральный элемент группы $g(t) = e^{\pm t X_4}$, $x = y = z = 0$, $w =  \pm t$:
\begin{align*}
&d(\Id, g(t)) = C \sqrt[3]{t}, \\
&C = \sqrt[3]{48 K^2(k_0)} \approx 6,\!37, \qquad K(k_0) - 2 E(k_0) = 0, \quad k_0 \approx 0,\!91.
\end{align*}
\end{itemize}

\subsubsection{Метрические прямые}\label{subsubsec:engel_mline}
\begin{theorem}
Натурально параметризованными метрическими прямыми на группе Энгеля являются следующие геодезические (и только они):
\begin{itemize}
\item[$(1)$]
однопараметрические подгруппы, касающиеся распределения:
\begin{align}
&e^{(u_1X_1 + u_2X_2) t} = \Exp(\lam, t), \quad t \in \R, \label{engel_mline1}\\
&u_1 = - \sin \t, \quad u_2 = \cos \t, \quad \lam = (\t, c=0, \a) \in C_4\cup C_5,\nonumber
\end{align}
\item[$(2)$]
критические геодезические:
\begin{align}
&\Exp(\lam, t), \quad \lam \in C_3, \quad t\in \R.\label{subsubsec:engel_mline2}
\end{align}
\end{itemize}
\end{theorem}
\begin{remark}
Геодезические \eq{engel_mline1} проецируются на плоскость $(x, y)$ в евклидовы прямые, из них анормальными являются только кривые
$$
e^{X_2 t} = \Exp(\lam, t), \quad \lam=(\t = 0, c=0, \a)\in C_4 \cup C_5.
$$
Геодезические \eq{engel_mline1} проецируются на плоскость $(x, y)$ в критические эйлеровы эластики (см.~Рис.~\ref{fig:el7}), так называемые солитоны Эйлера.

\subsubsection{Библиографические комментарии}
Разделы \ref{subsubsec:engel_state}, \ref{subsubsec:engel_geod}, \ref{subsubsec:engel_max} опираются на \cite{engel};
раздел
\ref{subsubsec:engel_inf_sym} --- на \cite{symmetry};
раздел \ref{subsubsec:engel_conj} --- на \cite{engel_conj};
разделы \ref{subsubsec:engel_cut}, \ref{subsubsec:engel_mline} --- на \cite{engel_cut};
раздел \ref{subsubsec:engel_synth} --- на \cite{engel_synth};
раздел \ref{subsubsec:engel_sphere} --- на \cite{engel_sphere}.

Параметризация субримановых геодезических на группе Энгеля впервые получена в работе \cite{vesrh_gran91}.

\end{remark}

\subsection{Субриманова задача на группе Картана}\label{subsec:cartan}
\subsubsection{Постановка задачи}\label{subsubsec:cartan_state}
\paragraph{Геометрическая постановка}
Рассмотрим следующее обобщение (усложнение) задач на группе Гейзенберга (раздел \ref{subsec:heis}) и группе 
Энгеля (раздел \ref{subsec:engel}) --- \ddef{обобщенную задачу Дидоны}. Пусть на евклидовой плоскости заданы точки $a_0, a_1 \in \R^2$, 
соединенные кривой $\gamma_0 \subset \R^2$. Пусть также заданы число $S \in \R$ и точка $c\in \R^2$. Требуется соединить точки $a_0$, $a_1$ кратчайшей кривой $\gamma \subset \R^2$
так, чтобы кривые $\gamma_0$ и $\gamma$ ограничивали на плоскости область алгебраической площади $S$, с центром масс $c$.

\paragraph{Задача оптимального управления}
Эту геометрическую задачу можно переформулировать как задачу оптимального управления
\begin{align}
&\dg = u_1 X_1(g) + u_2 X_2(g), \quad g = (x, y, z, v, w) \in \R^5, \label{cartanp21}\\
&g(0) = g_0, \quad g(t_1) = g_1, \label{cartanp22}\\
&l=\int_0^{t_1}\sqrt{u_1^2 + u_2^2}\,dt \to \min, \label{cartanp23}\\
&X_1 = \dfrac{\partial}{\partial x} - \dfrac{y}{2}\dfrac{\partial}{\partial z} - \dfrac{x^2 +y^2}{2}\dfrac{\partial}{\partial w}, \quad 
X_2 = \dfrac{\partial}{\partial y} + \dfrac{x}{2}\dfrac{\partial}{\partial z} + \dfrac{x^2 +y^2}{2}\dfrac{\partial}{\partial v}. \label{cartanX12}
\end{align}
Это субриманова задача для субримановой структуры на $\R^5$, заданной векторными полями $X_1$, $X_2$ как ортонормированным репером.

\paragraph{Алгебра Картана и группа Картана}
\ddef{Алгеброй Картана} называется пятимерная свободная нильпотентная алгебра Ли  $\gg$ с двумя образующими, глубины 3. Существует базис $\gg = \spann(X_1, \dots, X_5)$, в котором ненулевые скобки Ли суть
$$
[X_1, X_2] = X_3, \quad [X_1, X_3] = X_4, \quad [X_2, X_3] = X_5,
$$
см.~Рис.~\ref{fig:cartan}.

\begin{figure}[htb]
\setlength{\unitlength}{1cm}

\begin{center}
\begin{picture}(4, 4)
\put(1.1, 2.9){ \vector(1, -1){0.8}}
\put(1, 2.9){ \vector(0, -1){1.8}}
\put(2.9, 2.9){ \vector(-1, -1){0.8}}
\put(3, 2.9){ \vector(0, -1){1.8}}
\put(1.9, 1.9){ \vector(-1, -1){0.8}}
\put(2.1, 1.9){ \vector(1, -1){0.8}}

\put(1, 1){ \circle*{0.15}}
\put(1, 3){ \circle*{0.15}}
\put(3, 1){ \circle*{0.15}}
\put(3, 3){ \circle*{0.15}}
\put(2, 2){ \circle*{0.15}}

\put(1, 0.5) {$X_4$}
\put(3, 0.5) {$X_5$}
\put(1, 3.3) {$X_1$}
\put(3, 3.3) {$X_2$}
\put(2, 2.4) {$X_3$}

\end{picture}

\caption{ Алгебра Картана \label{fig:cartan}}

\end{center}
\end{figure}
Алгебра Картана имеет градуировку
$\gg = \gg^{(1)} \oplus\gg^{(2)} \oplus \gg^{(3)}$, $ \gg^{(1)} = \spann(X_1, X_2)$, $ \gg^{(2)} = \R X_3$, $ \gg^{(3)} = \spann(X_4, X_5)$, $[\gg^{(1)}, \gg^{(i)}] =  \gg^{(i+1)} $, 
$\gg^{(4)} = \gg^{(5)} = \{0\}$,
поэтому она является алгеброй Карно. Соответствующая связная односвязная группа Ли $G$ называется \ddef{группой Картана}.

На пространстве $\R^5_{x,y,z,v,w}$ можно ввести закон умножения
$$
\vect{x_1 \\ y_1 \\ z_1 \\ v_1 \\ w_1}
\cdot
\vect{x_2 \\ y_2 \\ z_2 \\ v_2 \\ w_2}
=
\vect{x_1 + x_2 \\ y_1 + y_2 \\ z_1 + z_2  + \frac 12 (x_1 y_2 - y_1 x_2)\\
v_1 +  v_2 + \frac 12 (x_1^2 + y_1^2 + x_1 x_2 + y_1 y_2) y_2 + x_1 z_2\\
w_1 + w_2 - \frac 12 (x_1^2 + y_1^2 + x_1x_2 + y_1 y_2)x_2 + y_1z_2},
$$
превращающий это пространство в группу Картана: $G \cong \R^5_{x,y,z,v,w}$, а поля (\ref{cartanX12}) в левоинвариантные поля на этой группе. 
Поэтому задача \eq{cartanp21}--\eq{cartanp23} есть   \ddef{левоинвариантная субриманова задача на группе Картана}. Следовательно, можно считать, что $g_0 = \Id = (0,\dots, 0)$.

Помимо модели \eq{cartanX12}, известны и другие модели субримановой  задачи на группе Картана \cite{brock_dai, monroy, symmetry}.

Левоинвариантная субриманова задача с вектором роста $(2,3,5)$ на группе Картана единственна, с точностью до изоморфизма этой группы \cite{symmetry}. 

\paragraph{Особенности задачи}
Субриманова задача на группе Картана есть простейшая левоинвариантная задача со следующими свойствами:
\begin{itemize}
\item
она имеет анормальные кратчайшие, касающиеся каждого вектора распределения,
\item
это  следующая по сложности после задачи Дидоны задача на свободной группе Карно максимального роста (ее вектор роста равен $(2, 3, 5)$).
\end{itemize}

Эта задача --- единственная свободная нильпотентная субриманова задача глубины 3 с интегрируемым по Лиувиллю нормальным гамильтоновым полем принципа максимума Понтрягина (неинтегрируемыми по Лиувиллю являются свободные нильпотентные задачи глубины 3, ранга более 2 \cite{borisov}, 
а также глубины более 3, ранга не менее 2 \cite{LS2018}.

Распределение $\Delta = \spann(X_1, X_2)$ имеет 14-мерную алгебру инфинитезимальных симметрий --- особую алгебру $\gg_2$, этот факт восходит к знаменитой  пятимерной работе Эли Картана \cite{cartan1910},  см.~также далее п.~\ref{subsubsec:cartan_inf_sym}.

Наконец, субриманова задача на группе Картана доставляет нильпотентную аппроксимацию любой задачи с вектором роста $(2, 3, 5)$, в частности:
\begin{itemize}
\item
задачи о качении двух твердых тел друг по другу без прокручивания и проскальзывания \cite{li_canny, cdc99, marigo_bicchi},
\item
машины с двумя прицепами \cite{laumond},
\item
задачи о движении электрического заряда в плоскости под действием магнитного поля \cite{monroy}.
\end{itemize}
Любой из этих причин достаточно для детального исследования субримановой задачи на группе Картана.

\subsubsection{Симметрии распределения и субримановой структуры}\label{subsubsec:cartan_inf_sym}
\begin{theorem}
Алгебра Ли инфинитезимальных симметрий распределения $\D$ на группе Картана есть $14$-мерная алгебра $\gg_2$ --- некомпактная вещественная форма комплексной особой алгебры Ли $\gg_2^{\C}$.
\end{theorem}

\begin{theorem}
Алгебра Ли инфинитезимальных симметрий нильпотентной  субримановой структуры на группе Картана есть $6$-мерная алгебра Ли,  в которой можно выбрать базис $X_0, Y_1, \dots, Y_5$ с ненулевыми скобками
\begin{align*}
&[X_0, Y_1] = - Y_2, \quad [X_0, Y_2] = Y_1,\\
&[X_0, Y_4] = -Y_5, \quad [X_0, Y_5] = Y_4,\\
&[Y_1, Y_2] = Y_3,\\
&[Y_1, Y_3] = Y_4, \quad [Y_2, Y_3] = Y_5. 
\end{align*}
Векторные поля $Y_1, \dots, Y_5$ --- правоинвариантные поля на группе $G$, а поле $X_0$ обращается в нуль в единице этой группы. Коммутаторы  симметрий с базисными полями субримановой структуры имеют вид:
\begin{align*}
&[Y_i, X_j] = 0, \quad i, j = 1, \dots, 5,\\
&[X_0, X_1]= -X_2, \quad [X_0, X_2] = X_1, \quad [X_0, X_3] = 0,\\
&[X_0, X_4] = -X_5, \quad [X_0, X_5] = X_4.
\end{align*}
В модели \eq{cartanX12}
$$
X_0 = -y\dfrac{\partial}{\partial x} + x \dfrac{\partial}{\partial y} - w \dfrac{\partial}{\partial v} + v \dfrac{\partial}{\partial w}.
$$
\end{theorem}
Представление  алгебры Ли симметрий распределения и субримановой структуры векторными полями в $\R^5$ приведено в работе \cite{symmetry}.
\subsubsection{Геодезические}\label{subsubsec:cartan_geod}
Существование оптимальных управлений в задаче \eq{cartanp21}--\eq{cartanp23} следует из теорем Рашевского-Чжоу и Филиппова.

\paragraph{Принцип максимума Понтрягина}
Переходя от минимизации длины \eq{cartanp23} к минимизации энергии $J = \frac{1}{2}\int_0^{t_1} (u_1^2 + u_2^2) dt$ и используя линейные на слоях $T^*G$ гамильтонианы 
$h_i(\lam) = \langle\lam, X_i\rangle$, $i = 1, \dots, 5$, получаем условия принципа максимума Понтрягина:
\begin{align*}
&\dh_1 = -u_2 h_3,\\
&\dh_2 = u_1 h_3,\\
&\dh_3 = u_1 h_4 + u_2 h_5,\\
&\dh_4 = 0,\\
&\dh_5 = 0,\\
&\dg = u_1 X_1 + u_2 X_2,\\
&u_1 h_1 + u_2 h_2 + \dfrac{\nu}{2}(u_1^2 + u_2^2) \to \max_{(u_1, u_2) \in \R^2},\\
& \nu \leq 0,\\
& (h_1, \dots, h_5, \nu) \neq 0.
\end{align*}
\paragraph{Анормальные экстремали}
Анормальные экстремали постоянной скорости могут быть параметризованы как
\begin{align}
&h_1= h_2 = h_3 = 0, \quad (h_4, h_5) \equiv \const \neq 0, \nonumber\\
&(u_1, u_2) \equiv \const, \nonumber\\
&x = u_1 t, \label{cartan_abn1}\\
&y = u_2 t, \label{cartan_abn2}\\
&z = 0, \label{cartan_abn3}\\
&v = (u_1^2 + u_2^2)u_1 t^3/6, \label{cartan_abn4}\\
&w = -(u_1^2 + u_2^2)u_2 t^3/6.\label{cartan_abn5}
\end{align}
Анормальные траектории \eq{cartan_abn1}--\eq{cartan_abn5} суть однопараметрические подгруппы 
 $g_t = e^{t(u_1X_1 + u_2X_2) }$, касающиеся распределения $\D$. Они проецируются на плоскость $(x, y)$ в прямые, потому являются кратчайшими.

Анормальное множество есть двумерное гладкое многообразие, диффеоморфное $\R^2$:
$$
\Abn = \{g \in G \mid z = v - (x^2 + y^2)x/6 = w + (x^2 + y^2)y/6  =0 \}.
$$

\paragraph{Нормальные экстремали}
Нормальные экстремали удовлетворяют гамильтоновой системе 
\begin{align}\label{cartan_Ham}
\dot{\lam} = \Vec{H}(\lam), \quad \lam \in T^*G, 
\end{align}
с гамильтонианом $H = \dfrac{1}{2}(h_1^2 + h_2^2).$ Введем на поверхности уровня $\{H =1/2\}$ координаты $(\t, c, \a, \b) \in S^1\times\R\times\R_+\times S^1$: 
\be{}   
h_1 = \cos \t, \quad  h_2 = \sin \t, \quad   h_3 = c, \quad h_4 = \a \sin \b, \quad h_5 = -\a \cos \b,
\ee
тогда гамильтонова система \eq{cartan_Ham} примет форму
\begin{align}
&\dot{\t} = c, \quad \dc = - \a \sin(\t - \b), \quad \dot{\a} = \dot{\b} = 0, \label{cartan_vert}\\
&\dg = \cos \t\, X_1 + \sin \t\,X_2.\label{cartan_hor}
\end{align}
Вертикальная подсистема \eq{cartan_vert} есть уравнение маятника.

Проекции нормальных геодезических на плоскость $(x, y)$ суть эйлеровы эластики, см.~раздел~\ref{subsec:elastica}.

Анормальные кратчайшие \eq{cartan_abn1}--\eq{cartan_abn5} удовлетворяют нормальной гамильтоновой системе \eq{cartan_vert}, \eq{cartan_hor} при $\t = \b$, $c = 0$, поэтому они нестрого анормальны.

\paragraph{Симплектическое слоение и функции Казимира}
На коалгебре Ли $\gg^*$ существуют 3 независимые функции Казимира:
\begin{align*}
h_4, \quad h_5, \quad E = \dfrac{h_3^2}{2} + h_1 h_5 - h_2 h_4.
\end{align*}
Симплектическое слоение на $\gg^*$ состоит из:
\begin{itemize}
\item
2-мерных параболических цилиндров 
$$
\{h_4 = \const, \quad h_5 = \const, \quad E = \const, \quad h_4^2 + h_5^2 \neq 0 \},
$$
\item
2-мерных аффинных плоскостей
$$
\{h_4 = h_5 = 0, \quad h_3 = \const  \neq 0\},
$$
\item
точек
$$
\{h_1 = \const, \quad h_2 = \const, \quad h_3 = h_4 = h_5 = 0  \}.
$$
\end{itemize}
Размерность симплектических листов не больше 2, поэтому вертикальная подсистема \eq{cartan_vert} интегрируема по Лиувиллю.

\paragraph{Параметризация нормальных геодезических}
Семейство нормальных экстремалей на поверхности уровня $\{H = \frac{1}{2} \}$ параметризуется начальными точками, принадлежащими цилиндру
$$
C = \gg^* \cap \left\{H = \frac{1}{2} \right\}.
$$
Этот цилиндр стратифицируется в зависимости от разных типов траекторий маятника \eq{cartan_vert}:
\begin{align*}
&C = \bigsqcup_{i=1}^7 C_i,  \\
&C_1 = \{ \lam \in C \mid \a > 0, \ E \in (-\a, \a) \}, \\
&C_2 = \{ \lam \in C \mid \a > 0, \ E \in (\a,  + \infty) \}, \\
&C_3 = \{ \lam \in C \mid \a > 0, \ E =  \a , \ \t - \b \neq \pi  \}, \\
&C_4 = \{ \lam \in C \mid \a > 0, \ E  = -\a  \}, \\
&C_5 = \{ \lam \in C \mid \a > 0, \ E =  \a, \ \t - \b = \pi \}, \\
&C_6 = \{ \lam \in C \mid \a =  0, \ c \neq 0 \}, \\
&C_7 = \{ \lam \in C \mid \a = c = 0 \}.
\end{align*}
Для параметризации нормальных геодезических введем на стратах $C_1$, $C_2$, $C_3$ эллиптические координаты $(\f, k, \a, \b)$, в которых  уравнение маятника \eq{cartan_vert} выпрямляется:\\
если {$\lam \in C_1$}, то
\begin{align*}
&k = \sqrt{\frac{E + \a}{2\a}} = \sqrt{\sin^2 \frac{\t - \b}{2} +
\frac{c^2}{4\a}} \in (0, 1), \\
&\f \in [0, 4K], \\
&\begin{cases}
\sin \ds \frac{\t - \b}{2} = k \sn(\sqrt{\a} \f), \\
\ds \frac{c}{2} =  k \sqrt{\a} \cn(\sqrt{\a} \f),
\end{cases}
\end{align*}
если { $\lam \in C_2$}, то
\begin{align*}
&k =  \sqrt{\frac{2 \a}{E + \a}} = \frac{1}{\sqrt{\sin^2 \frac{\t - \b}{2} +
\frac{c^2}{4\a}}}
         \in (0, 1), \\
&\f \in [0, 2kK], \\
&\begin{cases}
\sin \ds \frac{\t - \b}{2} = \pm \sn \frac{\sqrt{\a}\f}{k}, \\
\ds \frac{c}{2} =  \pm \frac{\sqrt{\a}}{k} \dn \frac{\sqrt{\a}\f}{k},
\end{cases} \\
&\pm = \sgn c, \\
&\p = \frac{\f}{k},
\end{align*}
если {$\lam \in C_3$}, то
\begin{align*}
&k = 1, \\
&\f \in (- \infty, + \infty), \\
&\begin{cases}
\sin \ds \frac{\t - \b}{2} = \pm \tanh(\sqrt{\a} \f), \\
\ds \frac{c}{2} =  \pm \frac{\sqrt{\a}}{\cosh(\sqrt{\a} \f)},
\end{cases} \\
&\pm = \sgn c.
\end{align*}
Тогда
$$
\dot{\f} = 1, \quad \dot{k} = \dot{\a} = \dot{\b} = 0.
$$
Задача инвариантна относительно левых сдвигов на группе Картана, дилатаций
\begin{align*}
e^{sY}: \, (t, x, y, z, v, w)&\mapsto (e^s t, e^s x, e^s y, e^{2s} z, e^{3s} v, e^{3s} w ), \nonumber \\ 
(\t, c, \a, \b)&\mapsto (\t, e^{-s} c, e^{-2s} \a, \b), \nonumber \\ 
(\f, k, \a, \b)&\mapsto (e^s \f, k, e^{-2s} \a, \b), \nonumber\\
Y = x \dfrac{\partial}{\partial x} + y \dfrac{\partial}{\partial y} &+ 2z \dfrac{\partial}{\partial z} + 3 v \dfrac{\partial}{\partial v} + 3w \dfrac{\partial}{\partial w},
\end{align*}
и вращений
\begin{multline}
e^{rX_0}: \, (x, y, z, v, w)\mapsto (x \cos r - y\sin r, x \sin r + y\cos r, 
z,\, v\cos r - w\sin r, v\sin r + w\cos r).
\end{multline}
С помощью вращений и дилатаций любой ковектор $\lam = (\f,k,\a,\b) \in\cup_{i=1}^3C_i$ переводится в фундаментальное множество 
$\{\a =1, \, \b = 0 \}$.
При $\a =1$, $\b =0$, $\lam \in\cup_{i=1}^3C_i $, геодезические $g_t = (x_t, y_t,z_t, v_t, w_t)$ параметризуются следующим образом.

Если $\lam \in C_1$, то
\begin{align*}
&x_t = 2(\E(\f_t) - \E(\f)) - (\f_t - \f), \\
&y_t = 2k(\cn \f - \cn \f_t), \\
&z_t = 2k(\sn \f_t \dn \f_t - \sn \f \dn \f) - k(\cn \f + \cn \f_t) x_t, \\
&v_t = 2k \sn \f_t \dn \f_t x_t - k \cn \f_t x_t^2 - (1 - 2k^2 + 2k^2 \cn \f \cn
\f_t) y_t, \\
&w_t = - \frac 16 \left(
x_t^3 + 2(2k^2 - 1 + 6k^2 \cn^2 \f) x_t + 2 (\f_t - \f)     \right. \\
&\qquad\qquad + 8k^2(\sn \f_t \cn \f_t \dn \f_t - \sn \f \cn \f \dn \f) \\
&\qquad\qquad \left. - 24 k^2 \cn \f (\sn \f_t \dn \f_t - \sn \f \dn \f)
\right),
\end{align*}
где $\f_t = \f + t$.

Если $\lam \in C_2$, то
\begin{align*}
&x_t = \frac 2k \left(\E(\p_t) - \E(\p) - \frac{2 - k^2}{2}(\p_t - \p)\right),
\\
&y_t = \pm \frac 2k (\dn \p - \dn \p_t), \\
&z_t = \pm \left( 2(\sn \p_t \cn \p_t - \sn \p \cn \p) - \frac 1k (\dn \p + \dn
\p_t) x_t\right), \\
&v_t = \pm \left(2 \sn \p_t \cn \p_t x_t - \frac 1k \dn \p_t x_t^2 \right)
              + \frac{1}{k^2}(2 - k^2 - 2 \dn \p \dn \p_t) y_t, \\
&w_t = - \frac 16 \left(
            x_t^3 + \frac{2}{k^2} (2 - k^2  + 6  \dn^2 \p) x_t + 2 k (\p_t - \p)
\right. \\
&\qquad\qquad    +  \frac{8}{k} (\sn \p_t \cn \p_t \dn \p_t - \sn \p \cn \p \dn
\p)  \\
&\qquad\qquad \left. -
         \frac{24}{k} \dn \p (\sn \p_t \cn \p_t - \sn \p \cn \p) \right), \\
&\pm = \sgn c,
\end{align*}
где $\p_t = \p + \frac tk$.

Если $\lam \in C_3$, то
\begin{align*}
&x_t = 2(\tanh \f_t - \tanh \f) - (\f_t - \f),    \\
&y_t =  \pm 2 \left( \frac{1}{\cosh \f} - \frac{1}{\cosh \f_t} \right), \\
&z_t =  \pm \left( 2 \left( \frac{\sinh \f_t}{\cosh^2 \f_t} -  \frac{\sinh
\f}{\cosh^2 \f} \right)
         - \left( \frac{1}{\cosh \f} + \frac{1}{\cosh \f_t} \right) x_t \right),
\\
&v_t = \pm \left( \frac{2}{\sinh \f_t} x_t - \frac{1}{\cosh \f_t} x_t^2\right) +
    \left( 1 - \frac{2}{\cosh \f \cosh \f_t}\right) y_t, \\
&w_t = - \frac 16 \left( x_t^3 + 6 \frac{2 + \cosh^2 \f}{\cosh^2 \f} x_t +
6(\f_t - \f)   \right. \\
&\qquad\qquad \left.    -  \frac{24}{\cosh \f}
        \left( \frac{\sinh \f_t}{\cosh^2 \f_t} -  \frac{\sinh \f}{\cosh^2 \f}
\right)
      - 8(\tanh^3 \f_t - \tanh^3 \f) \right), \\
&\pm = \sgn c,
\end{align*}
где $\f_t = \f + t$.

Параметризация геодезических при произвольных $\lam = (\f,k,\a,\b) \in\cup_{i=1}^3C_i$ получается из случая $\a = 1$, $\b = 0$ с помощью вращений и дилатаций:
\begin{align*}
&g_t(\f,k,\a,\b) = e^{-r X_0}\circ e^{-s Y} (g_{t'}(\f', k, \a' = 1, \b' = 0)),\\
&t' = t\sqrt{\a}, \quad \f' = \f \sqrt{\a}, \quad r = -\b, \quad s = \frac 12 \ln {\a}.
\end{align*}
В оставшихся случаях $\lam = (\t, c, \a, \b) \in\cup_{i=4}^7C_i$ геодезические параметризуются элементарными функциями.

Если $\lam = (\t, c, \a, \b) \in C_4 \cup C_5 \cup C_7$ и $\b = 0$, то 
$$(x_t, y_t, z_t, v_t, w_t) = (t, 0, 0, 0, -t^3/6).$$
В общем случае $\lam \in C_4 \cup C_5 \cup C_7$
$$
g_t(\t, c, \a, \b) = e^{-r X_0}(g_t(\t', c, \a, \b' =0)), \quad \t' = \t - \b, \quad r = -\b.
$$
Если $\lam = (\t = 0, c, \a =0) \in C_6$, то
\begin{align*}
&x_t = \frac{\sin \tau}{c}, \\
&y_t = \frac{1 - \cos \tau}{c}, \\
&z_t = \frac{\tau - \sin \tau}{2c^2}, \\
&v_t = \frac{\cos 2\tau - 4 \cos \tau + 3}{4c^3}, \\
&w_t = \frac{\sin 2\tau - 4 \sin \tau + 2\tau}{4c^3}, \\
&\tau = ct.
\end{align*}
В общем случае $\lam \in C_6$
$$
g_t(\t, c, \a =0, t) = e^{\t X_0}(g_t(\t' =0, c, \a=0, t)).
$$

Семейство всех геодезических параметризуется экспоненциальным  отображением
$$
\Exp: \, (\lam, t) \mapsto g_t = \pi \circ e^{t \vec H}(\lam), \quad C \times \R_+ \to G.
$$

\subsubsection{Симметрии и страты Максвелла}\label{subsubsec:cartan_max}
\paragraph{Непрерывные симметрии}
Дилатации и вращения образуют двухпараметрическую группу непрерывных симметрий экспоненциального отображения.

Введем линейные на слоях $T^*G$ гамильтонианы
$$
h_0(\lam) = \lan \lam, X_0(g) \ran, \qquad   h_Y(\lam) = \lan \lam, Y(g) \ran,
\qquad \lam \in T^*G,
$$
и соответствующие гамильтоновы векторные поля:
$$
\vh_0, \ \vh_Y \in \VVec(T^* G).
$$
Тогда
\begin{align*}
&[\vh_0, \vH] = 0,       &&\vh_0 H = 0, \\
&[\vh_Y, \vH] = - 2 \vH, &&\vh_Y H = - 2 H.
\end{align*}
Обозначим также вертикальное эйлерово поле на $T^*G$:
$e = \sum_{i = 1}^5 h_i \pder{}{h_i}$.
Так как гамильтониан
$H$ квадратичен на слоях, гамильтоново поле
$\vH$ линейно на слоях, поэтому
$$
[e, \vH] = \vH, \qquad e H = 2 H.
$$
Следовательно, векторное поле
$Z = \vh_Y + e$
удовлетворяет равенствам
\be{}\nonumber
[Z, \vH] = - \vH, \qquad Z H = 0.
\end{equation}
Более того,
\be{}\nonumber
[\vh_0, Z] = 0.
\end{equation}
\begin{proposition}
\label{prop:flow_commut}
Для любых $t,  s,  r \in \R$, $\lam \in T^*G$
$$
e^{rZ} \circ e^{s \vh_0} \circ e^{t \vH}(\lam) =
e^{t' \vH} \circ e^{r Z} \circ e^{s \vh_0}(\lam),
\qquad
\text{где  } t' = t e^r.
$$
\end{proposition}

\paragraph{Дискретные симметрии} 
Вертикальная подсистема \eq{cartan_vert} факторизуется по вращениям $X_0$ и дилатациям $Y$ в стандартное уравнение маятника
\be{}\nonumber
\dot \t = c, \quad \dot c = - \sin \t, \qquad (\t, c) \in S^1 \times \R.
\ee
Поле направлений этого уравнения имеет очевидные дискретные симметрии --- отражения в координатных осях и в начале координат
\begin{align*}
&\mapto{\eps^1}{(\t, c)}{(\t, - c)},\\
&\mapto{\eps^2}{(\t, c)}{(-\t, c)},\\
&\mapto{\eps^3}{(\t, c)}{(-\t, - c)}.
\end{align*}
Эти отражения порождают группу диэдра
$$
D_2 = \{\Id, \eps^1, \eps^2, \eps^3\} = \Z_2 \times \Z_2.
$$
Действие отражений естественно продолжается на эйлеровы эластики $(x_t, y_t)$, так что по модулю  вращений в плоскости $(x, y)$:
\begin{itemize}
\item
$\eps^1$ есть отражение эластики в центре ее хорды, 
\item
$\eps^2$ есть отражение эластики в серединном перпендикуляре ее хорды, 
\item
$\eps^3$ есть отражение эластики в ее хорде.  
\end{itemize}
Действие отражений также естественно продолжается в прообраз экспоненциального отображения:
$$
\map{\eps^i}{C \times\R_+}{C \times \R_+}, \qquad i = 1, 2, 3, 
$$
и в его образ:
$$
\map{\eps^i}{G}{G}, \qquad i = 1, 2, 3, 
$$
так что
\be{}\nonumber
\eps^i \circ \Exp(\lam, t) = \Exp \circ \, \eps^i(\lam,t), \qquad (\lam, t) \in C \times \R_+, \quad i = 1, 2, 3.
\ee
В явном виде:
\begin{align*}
&\eps^1 \ : (\t, c, \a, \b, t)  \mapsto (\t^1, c^1, \a, \b^1, t) = 
(\t_t, -c_t, \a, \b, t), \\ 
&\eps^2 \ : (\t, c, \a, \b, t)  \mapsto (\t^2, c^2, \a, \b^2, t) = 
(-\t_t, c_t, \a, -\b, t), \\ 
&\eps^3 \ : (\t, c, \a, \b, t)  \mapsto (\t^3, c^3, \a, \b^3, t) = 
(-\t, -c, \a, -\b, t), 
\end{align*}
\begin{align*}
&\eps^1: \, (x, y, z, v, w) \mapsto (x, y, -z, v - xz, w - yz),\\
&\eps^2: \, (x, y, z, v, w) \mapsto (x, -y, z, -v + xz, w - yz),\\
&\eps^3: \, (x, y, z, v, w) \mapsto (x, -y, -z, -v, w).
\end{align*}
Группа $\Sym$ симметрий  экспоненциального отображения состоит из вращений,  отражений и их композиций:
\begin{align*}
&e^{s \vh_0}, \quad \map{e^{s\vh_0} \circ \eps^i}{C \times \R_+}{C \times \R_+}, \\ 
&e^{s X_0}, \quad \map{e^{sX_0} \circ \eps^i}{G}{G}. 
\end{align*}
\begin{theorem}
Пусть $\lam \in C$. Первое время Максвелла, соответствующее группе $\Sym$ симметрий экспоненциального отображения, для почти всех геодезических $\Exp(\lam, t)$ выражается следующим образом:
\begin{align*}
&\lam \in C_1 \then \tmax(\lam) = 
\min \left(\frac{2}{\sqrt{\a}}p_1^z(k), 
\frac{2}{\sqrt{\a}}p_1^V(k)\right),  \\
&\lam \in C_2 \then \tmax(\lam) =  
\frac{2k}{\sqrt{\a}}p_1^V(k), \nonumber \\
&\lam \in C_6 \then \tmax(\lam) =  \frac{4}{|c|}p_1^V(0),  \nonumber\\
&\lam \in C_i, \ i = 3, 4, 5, 7 \then \tmax(\lam) =  + \infty. \nonumber
\end{align*}
Здесь $p = p_1^z(k) \in (K, 3 K)$ есть первый положительный корень функции
\be{}\nonumber
f_z(p,k) = \sn p \dn p - (2 \E(p) -  p) \cn p,
\ee
$p = p_1^V(k)$ есть первый положительный корень функции 
\begin{align*}
f_V(p) &= 
\dfrac 43 \, \ss \, \dd \, (-p - 2(1 - 2 k^2 + 6 k^2 \cc^2)(2 \E(p) - p) + (2 \E(p) - p)^3 \nonumber\\
&\qquad + 8 k^2  \cc \, \ss \, \dd) +   4  \cc \, (1 - 2 k^2  \ss^2)(2 \E(p) - p)^2, \quad  \\
p_1^V(k) &\in [2K, 4 K) \qquad \text{при } \lam \in C_1, \nonumber
\intertext{и}
f_V(p) &= 
\dfrac 43 \{ 3\,\dd\,(2\E(p) - (2 - k^2)p)^2   + \cc\,
[8\E^3(p) - 4\E(p)(4 + k^2) \nonumber\\
&\qquad -
 12\E^2(p)(2-k^2)p 
 + 6\E(p)(2-k^2)^2 p^2 \nonumber\\
&\qquad + 
 p(16 - 4k^2-3k^4-(2-k^2)^3p^2)]\,\ss -  
 \nonumber \\
 &\qquad -
2\,\dd\,(-4k^2+3(2\E(p)-(2-k^2)p)^2)\,\ss^2 + 
\nonumber\\
& \qquad 
+ 12 k^2 \cc(2 \E(p) - (2 - k^2)p)\ss^3 
- 8 k^2\,\ss^4\dd\},
\quad   \\
p_1^V(k) &\in (K, 2 K) \qquad \text{при } \lam \in C_2, \nonumber
\end{align*}
а $p=p_1^V(0) \in (\pi/2, \pi)$ есть первый положительный корень функции
$$
f_V^0(p) = [(32 p^2 - 1) \cos 2 p - 8 p \sin 2 p + \cos 6 p]/512.
$$
\end{theorem}
\begin{remark}
Для тех геодезических, для которых первое время Максвелла, соответствующее группе симметрий $\Sym$, не равно ${\tmax}$,  оно больше этого значения, а ${\tmax}$ есть первое сопряженное время.
\end{remark}
\begin{theorem}
Функция $\tmax: \, C \to (0, \, + \infty] $ имеет
следующие свойства инвариантности:
\begin{itemize}
\item[$(1)$]
$\tmax(\lam)$ зависит только от значений $E$ и $|\a|$,
\item[$(2)$]
$\tmax(\lam)$ есть первый интеграл поля $\vH_v$,
\item[$(3)$]
$\tmax(\lam)$ инвариантна относительно отражений: если $(\lam, t) \in C \times \R_+$, $(\lam^i, t) = \eps^i(\lam, t)$, то $\tmax(\lam^i)$ = $\tmax(\lam)$,
\item[$(4)$]
$\tmax(\lam)$ однородна относительно дилатаций: если $\lam \in C$, $\lam_s = \d_s(\lam) \in C$, то $\tmax(\lam_s) = e^s\, \tmax(\lam)$, $s \in \R$.
\end{itemize}
\end{theorem}

\subsubsection{Нижняя оценка сопряженного времени}\label{subsubsec:cartan_conj}
\begin{theorem}
Для любого $\lam \in C$
$$
\tconj(\lam) \geq \tmax(\lam).
$$
\end{theorem}
\subsubsection{Время разреза и кратчайшие}\label{subsubsec:cartan_cut}
\begin{theorem}
Для любого $\lam \in C$
$$
\tcut(\lam) = \tmax(\lam).
$$
\end{theorem}

\begin{theorem}
Пусть $g_1=(x_1,y_1,z_1,v_1,w_1) \in G$.
Если $z_1\neq 0$ и $x_1 v_1 + y_1 w_1 - {(x_1^2 + y_1^2)
z_1}/{2}\neq 0$, то существует единственная кратчайшая, соединяющая
$g_0=\Id$ c $g_1$.
\end{theorem}

\subsubsection{Метрические прямые}\label{subsubsec:cartan_mline}
\begin{theorem}
Натурально параметризованными метрическими прямыми на группе Картана являются следующие геодезические (и только они):
\begin{itemize}
\item[$(1)$]
однопараметрические подгруппы, касающиеся распределения:
\begin{align}
&e^{t(u_1 X_1 + u_2 X_2) } = \Exp(\lam, t), \label{cartan_mline1}\\
&u_1 = \cos \t, \quad u_2 = \sin \t, \quad \lam=(\t, c=0, \a, \b)\in C_4\cup C_5\cup C_7,\nonumber
\end{align}
\item[$(2)$]
критические геодезические
\begin{align}
&\Exp(\lam, t), \quad \lam \in C_3.                   \label{cartan_mline2}
\end{align}
\end{itemize}
\end{theorem}
\begin{remark}
Геодезические \eq{cartan_mline1}   проецируются на плоскость $(x, y)$ в евклидовы прямые, а геодезические \eq{cartan_mline2} --- в критические эйлеровы эластики (см.~Рис.~\ref{fig:el7}), так называемые солитоны Эйлера.
\end{remark}

\subsubsection{Библиографические комментарии}
Субриманова задача на группе Картана впервые рассматривалась в работе Р.~Брокетта и Л.~Даи \cite{brock_dai}, где показана интегрируемость геодезических в эллиптических  функциях.

Разделы \ref{subsubsec:cartan_state}, \ref{subsubsec:cartan_geod} опираются на \cite{dido_exp}; раздел \ref{subsubsec:cartan_inf_sym} --- на \cite{symmetry}; раздел \ref{subsubsec:cartan_max} --- 
на \cite{max1, max2, max3}; раздел \ref{subsubsec:cartan_conj} --- на \cite{cartan_conj}; разделы \ref{subsubsec:cartan_cut}, \ref{subsubsec:cartan_mline}  --- на \cite{cartan_cut}.

\section{Вместо заключения: некоторые неохваченные вопросы}\label{sec:uncovered}

Некоторые вопросы, близкие к рассмотренным выше, остались неохваченными из-за большого объема обзора. Перечислим их здесь:
\begin{enumerate}
\item
левоинвариантные субфинслеровы задачи \cite{buseman,b1, b2, ber94_2, berzub20, berzub20_2, convex,BBLDS,boscain3level,S18,ADS19, AS19_1, AS19_2,AS19_3,(36)},
\item
левоинвариантные сублоренцевы задачи \cite{groch06,grong_vas,kor_mar,CHSY, ASHY},
\item
левоинвариантные субримановы задачи с неинтегрируемым геодезическим потоком \cite{LS2017,LS2018,SS17,gole_karidi,borisov},
\item
приложения левоинвариантных задач к нильпотентной аппроксимации и конструктивному решению двухточечной задачи управления \cite{agrach_sarych,masht12,agrachev_marigo,bellaiche,Bellaiche11,Bellaiche10,BS,Jean2013,Fer35,gromov,b:hermes,b:lafsus,suss2,laumond,l81,stefani,TilSas119,vendit_laumond_oriolo,b:venditelli,walsh},
\item
приложения левоинвариантных задач к обработке изображений и моделям зрения \cite{MDSBB,benyosef,BDMG15,BDMS17,Boscain2014,Gauthier,citti_sarti,DuitsJMIV2014,DuitsIJCV2007,JDCS16,FMCS,DGA19,Franken2009,MAS13,MDSBB17,petitot,petitot_book},
\item
приложения левоинвариантных задач к робототехнике \cite{Ar-Eu-mob-rob,ard_gub,ArdentovRCD,EU-robot,laum98}.
\end{enumerate}

\newpage
\listoffigures
\addcontentsline{toc}{section}{Список иллюстраций}


\begin{thebibliography}{99}\addcontentsline{toc}{section}{Список литературы}

\subsection*{Книги и обзоры по теории управления, вариационному исчислению и \\субримановой геометрии}



\bibitem{agr_UMN}
А.А.Аграчев, Некоторые вопросы субримановой геометрии, {\em Успехи математических наук}, 71:6 (432), 2016, 3--36.

\bibitem{notes}
А.А.~Аграчев, Ю. Л. Сачков,  
{\em Геометрическая теория управления},
Физматлит, 2005.

\bibitem{ATF}
В.М. Алексеев, В.М. Тихомиров, С.В. Фомин,  
{\em Оптимальное управление},
Физматлит, 2005.

\bibitem{griffith}
Ф. Гриффитс, {\em  Внешние дифференциальные системы и вариационное исчисление},
М.: Мир, 1986. 




\bibitem{versh_gersh}
Вершик А.М., Гершкович В.Я.
{\em Неголономные динамические системы. Геометрия распределений и вариационные
задачи}.
Итоги науки и техники: Современные проблемы математики, Фундаментальные
направления, т. 16, ВИНИТИ, Москва, 1987, 5--85.

\bibitem{versh_gersh2}
Вершик А.М., Гершкович В.Я.
{\em Неголономные задачи и геометрия распределений}.
Добавление к книге \cite{griffith}, с. 318--349.

\bibitem{zelikin}
М.И. Зеликин, {\em Оптимальное управление и вариационное исчисление}, М.: URSS, 2014.

\bibitem{PGBM}
Понтрягин Л.\,С., Болтянский В.\,Г., Гамкрелидзе Р.\,В., Мищенко Е.\,Ф.
{\em Математическая теория оптимальных процессов},
М.: Наука, 1961.

\bibitem{JMS07}	
Сачков Ю.Л. Теория управления на группах Ли, {\em Современная математика. Фундаментальные направления}, 2007, т. 27, 5--59.

\bibitem{mono}
Сачков Ю.Л. {\em Управляемость и симметрии инвариантных систем на группах Ли и однородных пространствах}. ---М.: Физматлит, 2007.

\bibitem{intro}
Сачков Ю.Л. {\em Введение в геометрическую теорию управления}, М.: URSS, 2021, 160 C.

\bibitem{cime}
Agrachev A.A.,  
Geometry of optimal control problems and Hamiltonian systems. In: Nonlinear and Optimal Control Theory,  Lecture Notes in Mathematics. CIME, 1932, Springer Verlag,  2008, 1-59.

\bibitem{agrachev_open}
A. A. Agrachev. {\em Some open problems}, pp. 1–13. Springer International
Publishing, Cham, 2014.

\bibitem{ABB_book}
A. Agrachev, D. Barilari, U. Boscain, {\em A Comprehensive Introduction to
sub-Riemannian Geometry from Hamiltonian viewpoint}, Cambridge Studies in Advanced Mathematics, Cambridge Univ. Press, 2019

\bibitem{agr_gam}
A. A. Agrachev and R. V. Gamkrelidze, Symplectic geometry for optimal control, Nonlinear controllability
and optimal control, Monogr. Textbooks Pure Appl. Math., vol. 133, Dekker, New York, 1990, pp. 263--
277.

\bibitem{SRG} A.~Bellaiche, J.~Risler, Eds., Sub-Riemannian geometry.
Birkh\"aser, Progress in Math., 1996, v.144

\bibitem{bloch} 
A. Bloch,
{\em Nonholonomic Mechanics and Control},
Interdisciplinary Applied Mathematics, Volume 24, Springer, 2003.

\item U.~Boscain, B.~Piccoli, Optimal synthesis for control
systems on 2-D manifolds. Springer SMAI, v.43, 2004.

\bibitem{sussmann}
R.~W.~Brockett, R.~S.~Millman, H.~J.~Sussmann, Eds.,
Differential geometric control theory. Birkh\"auser Boston, 1983

\bibitem{capogna}
Luca Capogna, Donatella Danielli, Scott D. Pauls, and Jeremy T. Tyson.
{\em An introduction to the Heisenberg group and the sub-Riemannian
isoperimetric problem}, volume 259 of Progress in Mathematics.
Birkh\"auser Verlag, Basel, 2007.

\bibitem{euler}
Euler L., \textit{Methodus inveniendi lineas curvas maximi minimive proprietate gaudentes, sive Solutio problematis isoperimitrici latissimo sensu accepti}, Lausanne, Geneva, 1744.


\bibitem{isidori}
A.~Isidori,
{\em Nonlinear control systems: an introduction},
Springer-Verlag,  1985.


\item B.~Jakubczyk, W.~Respondek, Eds.,
Geometry of feedback and optimal control.
Marcel Dekker, 1998

\item
Fr\'ed\'eric Jean. {\em Control of nonholonomic systems: from sub-Riemannian
geometry to motion planning}. SpringerBriefs in Mathematics. Springer,
Cham, 2014.

\bibitem{jurd_book}
V.~Jurdjevic,
{\em Geometric Control Theory},
Cambridge University Press, 1997.

\bibitem{jurd16}
Velimir Jurdjevic. Optimal control and geometry: integrable systems.
Cambridge University Press, Cambridge, 2016.


\bibitem {montgomery_book}
R. Montgomery, {\em A tour of subriemannnian geometries, their geodesics and applications}, Amer. Math. Soc., 2002


\bibitem{nijmeijer_vanderschaft}
H.~Nijmeijer, A. ~van der Schaft,
{\em Nonlinear dynamical control systems},
Springer-Verlag, 1990.

\bibitem{rifford}
L. Rifford. {\em Sub-Riemannian geometry and Optimal Transport}.
SpringerBriefs in Mathematics, 2014.

\bibitem{sontag}
E.\,D.~Sontag,
{\em Mathematical control theory : deterministic finite dimensional systems},
Springer-Verlag, 1990.


\bibitem{sus_nonlin} 
H.~J.~Sussmann, Ed.,
Nonlinear controllability and optimal control.
Marcel Dekker, 1990

\bibitem{zelikin_borisov}
M.\,I.~Zelikin, V.\,F.~Borisov,
{\em Theory of chattering control with applications to astronautics, robotics,
economics, and engineering},
Basel: Birkhauser, 1994.



\subsection*{Другие книги и  статьи}
\bibitem{tikhomirov_stories}
В.М. Тихомиров, {\em Рассказы о максимумах и минимумах}, М.: МЦНМО, 2006

\bibitem{whit_watson}
Э.Т. Уиттекер, Дж.Н. Ватсон,
{\em 
Курс современного анализа},
М.: УРСС, 2002.

\bibitem{warner}
Ф. Уорнер,
{\em 
Основы теории гладких многообразий и групп Ли},
М.: Мир, 1987.

\bibitem{hadamard}
Krantz~S.~G., Parks~H.~R., \textit{The Implicit Function Theorem: History, Theory, and Applications}, Birkauser, 2001.

\bibitem{lawden}
D. F. Lawden,  {\em  Elliptic Functions and Applications}, Springer, 1989

\bibitem{akhiezer}
Н.И. Ахиезер,
{\em 
Элементы теории эллиптических функций},
М.: Наука, 1970.

\bibitem{vilenkin}
Н.Я. Виленкин, {\em Специальные функции и теория представлений групп}, М.: Мир, 1965. 

\bibitem{lan_liv}
Ландау Л. Д., Лифшиц Е. М. Механика. — 5-е изд., стереотип. — М.: ФИЗМАТЛИТ, 2012. — 224 с.

\bibitem{filippov}
А.Ф. Филиппов, О некоторых вопросах теории оптимального регулирования, {\em Вестник Московского университета, Сер. мат., мех., астрон., физ., хим.}, 1959, No. 2, 25--32. 


\subsection*{Классификации левоинвариантных субримановых задач}
\bibitem{versh_gersh88}
А. М. Вершик, В. Я. Гершкович, Геометрия неголономной сферы трехмерных групп Ли, Геометрия и теория особенностей в нелинейных уравнениях. Новое в глобальном анализе, 61—75. Воронеж, 1987. Перевод на англ. язык:  Lect. Notes in Math. 1334, 309-331 (1988).

\bibitem{FG}
Elisha Falbel and Claudio Gorodski. Sub-Riemannian homogeneous
spaces in dimensions 3 and 4. Geom. Dedicata, 62(3): 227–252, 1996.

\bibitem{agr_bar}
A.Agrachev, D. Barilari, 
Sub-Riemaian structures on 3D Lie groups. {\em J. Dynamical and Control Systems}, 2012, v.18, 21--44

\bibitem{almeida1}
D. Almeida, Sub-riemannian symmetric spaces of Engel type, {\em Mat. contemp.}, 17 (1998), pp. 45–58.


\bibitem{almeida2}
D. Almeida, Sub-riemannian homogeneous spaces of Engel type, {\em J. Dynamic. and Control Syst.}, 20 (2014), pp. 149–
166.

\bibitem{BM}
I. Beschastnyi, A. Medvedev, Left-invariant Sub-Riemannian Engel structures: abnormal geodesics and integrability, 
{\em SIAM Journal of Control and Optimization}, 2018, 56 (5), 3524--3537.


\subsection*{Субриманова задача на группе Гейзенберга}
\bibitem{brock}
Brockett R.
Control theory and singular Riemannian geometry//
In: New Directions in Applied Mathematics, (P.~Hilton and G.~Young eds.),
Springer-Verlag, New York, 11--27, 1982.

\subsection*{Машины Маркова-Дубинса и Ридса-Шеппа}
\bibitem{markov}
А.А. Марков, Несколько примеров решения особого рода задач о наибольших и наименьших величинах // Сообщ. Харьков. матем. общ. Вторая сер., 1:2 (1889), 250–276.

\bibitem{dubins}
{\it Dubins L.E.}
On curves of minimal length with a constraint on average curvature, and with prescribed initial and terminal positions and tangents // American Journal of Mathematics. 1957. V. 79. No. 3. P. 497--516.

\bibitem{sussmann_tang}
Sussmann, H.J. and Tang, G., Shortest paths for the Reeds-Shepp car: a~worked out example of the use of geometric techniques in nonlinear optimal control, \textit{Report SYCON-91-10}, 1991.

\bibitem{reeds}
Reeds, J.A. and Shepp, L.A., Optimal paths for a car that goes both forwards and backwards, \textit{Pacific J. Math.}, 1990. vol.~145, no.~2, pp.~367--393.

\bibitem{soueres_laumond}
P. So\`eres and J.-P. Laumond, Shortest paths synthesis for a car-like robot, {\em IEEE Transactions on Automatic Control}, vol. 41, No. 5, May 1996, 672--688.

\bibitem{soueres}
P. So\`eres, Commande optimale et robots mobiles non holonomes, Ph.D. dissertation, Universit\'e Paul Sabatier, no. 1554, 1993.

\bibitem{boissonat}
J.D. Boissonat, A. Cerezo, J. Leblong, Shortest paths of bounded curvature in the plane, in {\em Proc. IEEE Int. Conf. Robotics Automat.}, Nice, France, 1992.

\bibitem{pecsvaradi}
Th. Pecsvaradi, Optimal horizontal guidance law for aircraft in the terminal area,
{\em IEEE Transactions on Automatic Control},  vol. AC-17, No. 6, December 1972, 763--772.

\subsection*{Субриманова задача на двухступенных свободных нильпотентных группах Ли}
\bibitem{brock80}
R.W. Brockett. “Control theory and singular Riemannian geometry”,
New Directions in Applied Mathematics, Occasion of the Case Centennial
Celebration (April 25/26, 1980), eds. P.J. Hilton, G.S. Young, 1982,
pp. 11–27.

\bibitem{gaveau}
B. Gaveau. “Principe de moindre action, propagation de la chaleur et
estimees sous elliptiques sur certains groupes nilpotents”, Acta Mathematica,
139 (1977), pp. 95–153.

\bibitem{liusus}
W.S. Liu, H. Sussmann, 
Shortest paths for sub-Riemannian metrics on rank-2 distributions, 
 Memoirs of the American Mathematical Society, No. 564, Vol. 118, November 1995

\bibitem{monroy_meneses}
F. Monroy-Pérez, A. Anzaldo-Meneses. “The step-2 nilpotent $(n, n(n + 1)/2)$
sub-Riemannian geometry”, Journal of Dynamical and Control Systems, 12:2
(2006), pp. 185–216.

\bibitem{montanari_morbidelli}
A. Montanari, D. Morbidelli. On the sub-Riemannian cut locus in a model of free
two-step Carnot group, Calc. Var. 56, 36 (2017), 1 - 26.

\bibitem{myasn36}
O. Myasnichenko. “Nilpotent $(3, 6)$ sub-Riemannian problem”, Journal of
Dynamical and Control Systems, 8:4 (2002), pp. 573–597

\bibitem{myasn06}
Oleg Myasnichenko. Nilpotent $(n, n(n + 1)/2)$ sub-Riemannian
problem. J. Dyn. Control Syst., 12(1): 87–95, 2006.

\bibitem{rizzi_serres}
L. Rizzi, U. Serres. On the cut locus of free, step two Carnot groups, Proc. Amer. Math. Soc. 145 (2017), 5341-5357

\bibitem{ledonne_sard}
E. Le Donne, R. Montgomery, A. Ottazzi, P. Pansu, D.Vittone,
Sard property for the endpoint map on some Carnot groups,
Annales de l'Institut Henri Poincar\'e C, Analyse non lin\'eaire
Volume 33, Issue 6, November–December 2016, Pages 1639-1666

\subsection*{Двухступенные  субримановы задачи коранга 1 и 2}
\bibitem{corank1}
A. Agrachev, D. Barilari, U. Boscain, On the Hausdorff volume in sub-Riemannian geometry, 
{\em Calculus of variations and partial differential equations}, 
43 (2012), 3--4, 355--388.

\bibitem{mon_anz_heis}
F. Monroy-P\'erez, A. Anzaldo-Meneses, Optimal control on the Heisenberg group, {\em Journal of Dynamical and Control Systems},  vol. 5, No. 4, 1999, 473--499.

\bibitem{corank2}
D. Barilari, U. Boscain, J.-P. Gauthier,
On 2-step, corank 2 nilpotent sub-Riemannian metrics,
SIAM Journal on Control and Optimization 50 (2012), 1:559--582

\subsection*{Субримановы $\ddd\oplus\sss$ структуры} 
\bibitem{agr95}
A.A. Agrachev, Methods of control theory in nonholonomic geometry. In {\em Proceedings of the International Congress of Mathematicians, Vols. 1, 2}, pp. 1473--1483, Basel, Birkh\"auser, 1995.

\bibitem{brock99}
R.W. Brockett, Explicitly solvable control problems with nonholonomic constraints. In {\em Decision and control, 1999, Proc. 38th IEEE Conf. on}, vol. 1, pp. 13--16, IEEE, 1999.


\bibitem{brock73}
R.W. Brockett, Lie theory and control systems defined on spheres. {\em SIAM J. Appl. Math.}, 25: 213--225, 1973

\bibitem{jurd99}
V. Jurdjevic, Optimal control, geometry, and mechanics. In {\em Mathematical control theory}, pp. 227--267, Springer, New York, 1999.

\bibitem{jurd01}
V. Jurdjevic, Hamiltonian point of view of non-Euclidean geometry and elliptic functions. {\em Systems Control Lett.}, 43 (1): 25--41, 2001.

\bibitem{BCG02a}
U. Boscain, T. Chambrion, J.-P. Gauthier, On the $K \ + P$   problem for a three-level quantum system: optimality implies resonance. {\em J. Dyn. Control Syst.}, 8 (4): 547--572, 2002.

\subsection*{Симметричные субримановы задачи на группах $\SO(3)$, $\SU(2)$, и $\SL(2)$}
\bibitem{gersh84}
В.Я. Гершкович, Вариационная задача с неголономной связью на $\SO(3)$. --- В сб.: Новое в глобальном анализе. --- Воронеж: Изд-во ВГУ, 1984, с. 149--152.

\bibitem{versh_gersh86}
А. М. Вершик, В. Я. Гершкович, “Геодезический поток на $\SL(2,\R)$ с неголономными ограничениями”, Записки научных семинаров ЛОМИ / Акад. наук СССР, Мат. ин-т им. В. А. Стеклова, Ленингр. отд-ние. - Л. : Наука, Т. 155 : Дифференциальная геометрия, группы Ли и механика : 7 : сб. работ / под ред. Л. Д. Фаддеева. – 1986, 7–17

\bibitem{ber01}
В. Н. Берестовский, И. А. Зубарева, “Формы сфер специальных неголономных левоинвариантных внутренних метрик на некоторых группах Ли”, Сибирский математический журнал, 42:4 (2001), 731–748

\bibitem{ber14}
В. Н. Берестовский, “Универсальные методы поиска нормальных геодезических на группах Ли с левоинвариантной субримановой метрикой”, Сиб. матем. журн., 55:5 (2014), 959–970

\bibitem{ber151}
	В. Н. Берестовский, “(Локально) кратчайшие специальной субримановой метрики на группе Ли $\SO 0(2,1)$”, Алгебра и анализ, 27:1 (2015), 3–22

\bibitem{BZ15a}
В. Н. Берестовский, И. А. Зубарева, “Геодезические и кратчайшие специальной субримановой метрики на группе Ли $\SO(3)$”, Сиб. матем. журн., 56:4 (2015), 762–774

\bibitem{BZ15b} 
В. Н. Берестовский, И. А. Грибанова, “Субриманово расстояние в группах Ли $\SU(2)$ и $\SO(3)$”, Матем. тр., 18:2 (2015), 3–21

\bibitem{BZ16}
В. Н. Берестовский, И. А. Зубарева, “Геодезические и кратчайшие специальной субримановой метрики на группе Ли $\SL(2)$”, Сиб. матем. журн., 57:3 (2016), 527–542

\bibitem{BZ16_1}
В. Н. Берестовский, И. А. Зубарева, “Локально изометричные накрытия группы Ли $\SO 0(2,1)$ со специальной субримановой метрикой”, Матем. сб., 207:9 (2016), 35–56  

\bibitem{BZ16_2}
В. Н. Берестовский, И. А. Зубарева, “Субриманово расстояние в группе Ли $\SO 0(2,1)$”, Алгебра и анализ, 28:4 (2016), 62–79

\bibitem{BZ17}
В. Н. Берестовский, И. А. Зубарева, “Субриманово расстояние в группе Ли $\SL(2)$”, Сиб. матем. журн., 58:1 (2017), 22–35 

\bibitem{boscain_rossi}
U. Boscain, F. Rossi. “Invariant Carnot–Caratheodory metrics on 
$S^3$, $\SO(3)$, $\SL(2)$ and lens spaces”, SIAM J. Control Optim., 47 (2008), pp. 1851–1878

\bibitem{chang_markina_vasiliev}
D.Ch. Chang, I. Markina, A. Vasiliev,
Sub-riemannian geodesics on the 3-sphere,
{\em Compl. anal. oper. theory}, 3 (2009), 361--377.


\bibitem{BS16}
И. Ю. Бесчастный, Ю. Л. Сачков, “Геодезические в субримановой задаче на группе $\SO(3)$”, Матем. сб., 207:7 (2016), 29–56 

\subsection*{Римановы задачи на группах $\SO(3)$, $\SU(2)$, $\SL(2)$ и  $\PSL(2)$}

\bibitem{bates-fasso}
{\it L.~Bates, F.~Fass\`{o}.} The Conjugate Locus for the Euler Top. I. The Axisymmetric Case. //
International Mathematical Forum. 2007. 2, 43. 2109--2139.

\bibitem{sakai}
{\it T.~Sakai.} Cut loci of Berger's sphere. // Hokkaido Mathematical Journal. 1981. 10. 143--155.


\bibitem{so3_sl2_dan}
	А.В.Подобряев, Ю.Л.Сачков, Левоинвариантные симметричные римановы задачи на группах собственных движений плоскости Лобачевского и сферы,  {\em Доклады Академии Наук}, 2017, 473, No. 6, с. 640--642
	
\bibitem{so3}
A.Podobryaev, Yu. L. Sachkov,
Cut locus of a left invariant Riemannian metric on $\SO(3)$ in the axisymmetric case,
{\em Journal of Geometry and Physics},    110 (2016) 436--453.

\bibitem{diamsu2}
А. В. Подобряев, Диаметр сферы Берже, Матем. заметки, 103:5 (2018), 779–784



\bibitem{sl2}
A.Podobryaev, Yu. L. Sachkov,
Symmetric Riemannian Problem on the Group of Proper Isometries of Hyperbolic Plane, J Dyn Control Syst (2018) 24: 391 

\subsection*{Субриманова задача в плоском случае Мартине}

\bibitem{martinet}
A.Agrachev, B. Bonnard, M. Chyba, I. Kupka,
Sub-Riemannian sphere in Martinet flat case. J. ESAIM: Control, Optimisation and Calculus of Variations, 1997, v.2,
377--448 

\subsection*{Субриманова задача на группе $\SE(2)$}

\bibitem{ber94}
В. Н. Берестовский, “Геодезические левоинвариантной неголономной римановой метрики на группе движений евклидовой плоскости”, Сиб. матем. журн., 35:6 (1994), 1223–1229

\bibitem{max_sre}
I. Moiseev, Yu. L. Sachkov,
Maxwell strata in sub-Riemannian problem  on the group of motions of a plane,
{\em ESAIM: COCV}, 16 (2010), 380--399.

\bibitem{cut_sre1}
Yu. L. Sachkov,
Conjugate and cut time in the sub-Riemannian problem on the group of motions of a 
plane, {\em ESAIM: COCV}, 16 (2010), 1018--1039.

\bibitem{cut_sre2}
Yu. L. Sachkov,
Cut locus and optimal synthesis in the sub-Riemannian problem  on the group of motions of a 
plane, {\em ESAIM: COCV}, 17 (2011), 293--321.

\bibitem{bicycle}
A. Ardentov,  G. Bor, E. Le Donne,  R. Montgomery, Yu. Sachkov,
 Bicycle paths, elasticae and sub-Riemannian geometry, 
{\em Nonlinearity},  принята к публикации. 

\bibitem{se2_GO}
Ю.Л. Сачков, Однородные субримановы геодезические на группе движений плоскости, {\em направлена для  публикации}.


\subsection*{Субриманова задача на группе $\SH(2)$}

\bibitem{sh2_1}
Y.A.Butt, A.I. Bhatti, Yu. L. Sachkov,
Parametrization of Extremal Trajectories in Sub-Riemannian Problem on Group of
Motions of Pseudo Euclidean Plane,
{\em Journal of Dynamical and Control Systems},  
Vol. 20 (2014), No. 3 (July), 341--364.

\bibitem{sh2_2}
Y.A.Butt, A.I. Bhatti, Yu. L. Sachkov,
Maxwell strata and conjugate points in sub-Riemannian problem on the group $\SH(2)$,
{\em Journal of Dynamical and Control Systems}, 22 (2016),  747--770

\bibitem{sh2_3}
Y.A.Butt, A.I. Bhatti, Yu. L. Sachkov,
Cut Locus and Optimal Synthesis in Sub-Riemannian Problem on the Lie
Group $\SH(2)$,
{\em Journal of Dynamical and Control Systems},  23 (2017), 155--195




\subsection*{Задача Эйлера об эластиках}

\bibitem{el_AiT}
А.А.Ардентов, Ю.Л.Сачков. Решение задачи Эйлера об эластиках.  {\em Автоматика и телемеханика}, 2009, No.  4, 78--88.



\bibitem{love}
{\sl Ляв~А.} Математическая теория упругости. Москва-Ленинград:
ОНТИ, 1935.

\bibitem{el_dan}
Сачков Ю.Л. Оптимальность эйлеровых эластик, {\em Доклады Академии Наук}, том 417, No. 1, ноябрь 2007, С. 23--25.

\bibitem{euler_rus}
{\sl Эйлер~Л.} Метод нахождения кривых линий, обладающих
свойствами максимума или минимума, или решение изопериметрической
задачи, взятой в самом широком смысле, Леонарда Эйлера,
королевского профессора и члена Императорской Петербургской
Академии Наук. Приложение I. <<Об упругих кривых>>.
Москва-Ленинград: ГТТИ, 1934. С. 447-572.

\bibitem{el_stable}
	Сачков Ю.Л., Левяков С.В. Устойчивость инфлексионных эластик, центрированных в вершинах или точках перегиба, {\em Труды МИАН}, 2010,  T. 271, 187--203.
	
\bibitem{antman}
S.S.Antman,
The influence of elasticity on analysis: Modern developments,
{\em Bulletin American Math. Society}, 1983, v. 9, No. 3, 267--291.

\bibitem{DBernoulli}
D.Bernoulli, 26th letter to L. Euler (October, 1742), In: Fuss, {\em Correspondance math\'ematique et physique}, t.2, St. Petersburg, 1843.

\bibitem{JBernoulli}
J.Bernoulli,
V\'eritable hypoth\`ese de la r\'esistance des solides, avec la demonstration de la corbure des corps qui font ressort,
In: {\em Collected works}, t.2, Geneva, 1744.
	

\bibitem{birkhoff64}
{\sl Birkhoff~G., de Boor~C.R.} Piecewise polynomial interpolation
and approximation // Approximation of Functions. Proc. Sympos.
General Motors Res. Lab., 1964. Elsevier. Amsterdam: 1965. P.
164-190.

\bibitem{born}
Max Born. Stabilitat der elastischen Linie in Ebene und Raum.
Preisschrift Und Dissertation, 1: 5–101, 1906.	
	
\bibitem{jurd_noneuclid}
V.Jurdjevic,
Non-Euclidean elastica,
{\em Am. J. Math.}, v. 117 (1995), 93--125.

\bibitem{levien}
{\it Levien R.} The elastica: a mathematical history // Technical Report No. UCB\slash EECS-2008-103,	2008, P. 1--25.



\bibitem{manning96}
R.S.Manning, J.H.Maddocks, J.D.Kahn,
A continuum rod model of sequence-dependent DNA structure,
{\em J. Chem. Phys.} v. 105 (1996), 5626--5646.

\bibitem{manning98}
R.S.Manning, K.A.Rogers, J.H.Maddocks,
Isoperimetric conjugate points with application to the stability of DNA minicircles,
{\em Proc. R. Soc. Lond. A}, v. 454 (1998), 3047--3074.

\bibitem{el_SE2}
A. Mashtakov, A. Ardentov, Yu. Sachkov, Relation between Euler’s elasticae and sub-Riemannian geodesics on SE(2), {\em Regular and Chaotic Dynamics}, December 2016, Volume 21, Issue 7, pp 832--839.	

\bibitem{mumford}
D.Mumford,
Elastica and computer vision, In:
{\em Algebraic geometry and its applications},
C.L.Bajaj, Ed., Springer-Verlag, New-York, 1994, 491--506.

\bibitem{saalchutz}
L.Saalsch\"utz,
{\em Der belastete Stab}, Leipzig, 1880.
	
\bibitem{el_max}
Yu. L. Sachkov,
Maxwell strata in Euler's elastic problem,
Journal of Dynamical and Control Systems,
Vol. 14 (2008), No. 2  (April), pp. 169--234. 


\bibitem{el_conj}
Yu. L. Sachkov,
Conjugate points in Euler's elastic problem,
{\em Journal of Dynamical and Control Systems}, vol. 14 (2008), No. 3 (July).

\bibitem{el_closed}
Yu. Sachkov, Closed Euler Elasticae, {\em Proceedings of the Steklov Institute of Mathematics}, V. 278 (2012), 218--232.

\bibitem{el_cut}
Yu. L. Sachkov, E.Sachkova,
Exponential mapping in Euler's elastic problem,
{\em Journal of Dynamical and Control Systems}, Vol. 20 (2014),    No. 4, 443--464.

\bibitem{timoshenko}
S.Timoshenko,
{\em History of Strength of Materials},
McGraw-Hill, New-York, 1953.

\bibitem{truesdell}
C.Truesdell,
The Influence of Elasticity on Analysis: The Classic Heritage,
{\em Bulletin American Math. Society}, 1983, v. 9, No. 3, 293--310.

\bibitem{Greenhill}
A. G. Greenhill, 
{\em The applications of elliptic functions} (New York, Macmillan, 1892).

\bibitem{Frisch}
R. Frisch Fay, 
{\em Flexible Bars}, Butterworths, London (1962).

\bibitem{Bisshopp}
K.E. Bisshopp, D.C. Drucker, Large deflection of cantilever beams. 
{\em Quart. Appl. Math.} 3 (1945) 272-275.

\bibitem{Lardner}
T. J. Lardner, A note on the elastica with large loads, 
{\em Int. J. Solids Struct.}, Vol. 21. No. 1. pp. 21-26, I985.

\bibitem{Wang}
C. Y. Wang, Post-buckling of a clamped-simply supported elastica, 
{\em Int. J. Non-Linear Mechanics}, Vol. 32. No. 6, pp. I 115-l 122, 1997.

\bibitem{Panayotounakos}
D. E. Panayotounakos, P. S. Theocaris, Analytic solutions for nonlinear differential equations describing the elastica of straight bars : Theory, 
{\em Journal of the Franklin Institute}, Vol. 325, No. 5, pp. 621433, 1988.

\bibitem{Seide}
P. Seide. Large deflections of a simply supported beam subjected to moment at one end, 
{\em Trans. ASME, J. Appl. Mech.}, Vol. 51, 519-525, 1984.

\bibitem{Naschie}
M. S. El Naschie, Thermal initial post buckling of the extensional elastica, 
{\em Int. J. Mechanical Science}, 1976, Vol. 18, pp. 321-324.

\bibitem{Stampouloglou}
I. H. Stampouloglou, E. E. Theotokoglou, P. N. Andriotaki , Asymptotic solutions to the non-linear cantilever elastica, 
{\em Int. J. Non-Linear Mechanics}, 40 (2005) 1252 - 1262.

\bibitem{Glassmaker}
N. J. Glassmaker, C. Y. Hui, Elastica solution for a nanotube formed by self-adhesion of a folded thin film, 
{\em J. Appl. Physics}, V. 96, N. 6, 3429-3434.

\bibitem{Tang}
T. Tang, N. J. Glassmaker, On the inextensible elastica model for the collapse of nanotubes, {\em Mathematics and Mechanics of Solids}, 2009, doi:10.1177/1081286509105923.

\bibitem{Mikata}
Y. Mikata, Complete solution of elastica for a clamped-hinged beam, and its applications to a carbon nanotube, {\em Acta Mechanica} 190, 133-150 (2007).

\bibitem{Heijden}
Van der Heijden, G.H.M., Neukirch, S., Goss, V.G.A., Thompson, J.M.T.: Instability and self-contact phenomena in the writhing of clamped rods. 
{\em Int. J. Mech. Sci.} 45, 161-196 (2003).



\subsection*{Задача о качении сферы по плоскости}

\bibitem{s2r2}
Маштаков А.П., Сачков Ю.Л.,  Экстремальные траектории и асимптотика времени Максвелла в задаче об оптимальном качении сферы по
плоскости, {\em Математический сборник} (2011),  Т.  202, No. 9, C. 97--120.

\bibitem{s2r2_sym}
	Сачков Ю.Л. Симметрии и страты Максвелла и симметрии в задаче об оптимальном качении сферы по плоскости, {\em Мат. Сборник}, 2010,  T. 201, N 7, С. 99--120.

\bibitem{arthur_walsh} A.M. Arthurs, G.R.Walsh, On the Hammersley’s minimum problem for a rolling
sphere, {\em Math. Proc. Cambridge Phil. Soc.}, 99 (1986), 529–534.

\bibitem{bicchi_prat_sast}
Bicchi A., Prattichizzo D., Sastry S., Planning motions of rolling surfaces// IEEE Conf. on Decision and Control, 1995

\bibitem{hammersley}
J.M. Hammersley, Oxford commemoration ball, Probability, Statistics and Analysis
(London Math. Soc. lecture notes), 79 (1983), 112–142.



\bibitem{jurd_ball} 
V. Jurdjevic, The geometry of the plate-ball problem, {\em Arch. Rat. Mech. Anal.}, 124
(1993), 305–328.

\bibitem{roll_besch}
И. Ю. Бесчастный,
Об оптимальном качении сферы с прокручиванием, без проскальзывания, {\em Матем. сб.}, 205:2 (2014),  3–38




\subsection*{Субриманова задача на группе Энгеля}
\bibitem{vesrh_gran91}
Вершик А.М., Граничина О.А.,  Редукция неголономных вариационных задач к изопериметрическим и связности в главных расслоениях, {\em Мат. заметки} (1991), Т.  202, No. 11, C. 37--44.

\bibitem{engel}
Ардентов А.А., Сачков Ю.Л.,  Экстремальные траектории в нильпотентной субримановой задаче на группе Энгеля, {\em Мат. сборник} (2011), Т.  49, No. 5, C. 31--54.

\bibitem{engel_dan}	А.А.Ардентов, Ю.Л.Сачков, Множество разреза в субримановой задаче на группе Энгеля, {\em Доклады Академии Наук}, Выпуск 6, 2018, Том 478,  623--626.

\bibitem{engel_conj}
Ardentov, A.\,A. and Sachkov, Yu.\,L., Conjugate points in nilpotent sub-Riemannian problem on the Engel group, \textit{Journal of Mathematical Sciences}, 2013, vol.\,195, no.\,3, pp.\,369--390.

\bibitem{engel_cut}
Ardentov, A.\,A. and Sachkov, Yu.\,L., Cut time in sub-Riemannian problem on Engel group, \textit{ESAIM: COCV}, 2015, vol.\,21, no.\,4, pp.\,958--988.

\bibitem{engel_synth}
A.A. Ardentov, Yu. L. Sachkov,
Maxwell Strata and Cut Locus in the Sub-Riemannian Problem on the Engel Group, Regular and Chaotic Dynamics, December 2017, Vol. 22, Issue 8, pp 909--936.

\bibitem{engel_sphere}
Yu. L. Sachkov, On the structure of Engel sub-Riemannian sphere, {\em в работе}.

\subsection*{Субриманова задача на группе Картана}

\bibitem{cartan1910}
E.~Cartan,
L\`es systemes de Pfaff a cinque variables et l\`es equations aux
derivees partielles du second ordre,
{\em Ann. Sci. \`Ecole Normale} {\bf 27} (1910), 3: 109--192.


\bibitem{dido_exp}
Сачков Ю.Л., Экспоненциальное отображение в обобщенной задаче Дидоны, {\em Математический сборник}, 2003, 194 (9),  63--90.

\bibitem{max1}
Сачков Ю.Л.  Дискретные симметрии в обобщенной задаче Дидоны, {\em Мат. Сборник}, 2006, T. 197, N 2, S. 95--116.


\bibitem{max2}
Сачков Ю.Л. Множество Максвелла в обобщенной задаче Дидоны, {\em Мат. Сборник}, 2006, Т. 197,  № 4, С. 123--150.

\bibitem{max3}
Сачков Ю.Л. Полное описание стратов Максвелла в обобщенной задаче Дидоны, {\em Мат. Сборник}, 2006,  T. 197, N 6, С. 111--160.

\bibitem{brock_dai}
Brockett R., Dai L.
Non-holonomic kinematics and the role of el\-lip\-tic functions in
constructive controllability//
In: Nonholonomic Motion Planning, Z.~Li and J.~Canny, Eds.,
Kluwer, Boston, 1993, 1--21.

\bibitem {symmetry}
Yu. L. Sachkov, Symmetries of flat rank two distributions and sub-Riemannian structures, {\em Transactions of the AMS}, 2004, 2, 457-494

\bibitem{cartan_conj}
Yu. L. Sachkov, Conjugate time in the  sub-Riemannian problem on the  Cartan group, {\em Journal of Dynamical and Control Systems}, accepted.

\bibitem{cartan_cut}
A. Ardentov, E. Hakavuori,  Cut time in the  sub-Riemannian problem on the  Cartan group, {\em in preparation}.






\subsection*{$\exp$-$\log$  категория}
\bibitem{dries}
L.V.D. Dries, A. Macintyre, D. Marker, The elementary theory of restricted analytic fields with exponentiation, {\em Annals of Mathematics}, 140, 1994, 183--205.
\bibitem{lion}
J.M. Lion, J.P. Rolin, Th\'eor\`emes de pr\'eparation pur les fonctions logarithmo-exponentielles, {\em Annales de l'Institut Fourier}, 47, 1997, 859--884.

\subsection*{Субримановы задачи с неинтегрируемым геодезическим потоком}

\bibitem{LS2017}		Л. В. Локуциевский, Ю. Л. Сачков, Неинтегрируемость по Лиувиллю субримановых задач на свободных группах Карно глубины 4,  {\em Доклады Академии Наук}, 2017, т. 95, no. 3, С. 211–213.

\bibitem{LS2018}	
Л. В. Локуциевский, Ю. Л. Сачков, Об интегрируемости по Лиувиллю субримановых задач на группах Карно глубины 4 и больше, {\em Матем. сб.}, 209:5 (2018), 74–119 


\bibitem{SS17}
Ю. Л. Сачков, Е.Ф. Сачкова, Вырожденные анормальные траектории в субримановой задаче с вектором роста (2,3,5,8), Дифференциальные уравнения, 2017, 3, 362--374.

\bibitem{gole_karidi}
Gol\'e, C., Karidi, R. A note on Carnot geodesics in nilpotent Lie groups. {\em Journal of Dynamical and Control Systems} 1, 535–549 (1995).

\bibitem{borisov}
Bizyaev I. A.,  Borisov A. V.,  Kilin A. A.,  Mamaev I. S., Integrability and Nonintegrability of Sub-Riemannian Geodesic Flows on Carnot Groups, Regular and Chaotic Dynamics, 2016, vol. 21, no. 6, pp. 759-774



\subsection*{Левоинвариантные субфинслеровы задачи}
\bibitem{ADS19}
Ардентов А.А., Э. Ле Донне, Сачков Ю.Л., Субфинслеровы геодезические на группе Картана, Regular and Chaotic Dynamics, 2019,  vol. 24, no. 1, pp. 36-60

\bibitem{AS19_1}
Ардентов А.А., Э. Ле Донне, Сачков Ю.Л., Экстремальные траектории в субфинслеровой задаче на группе Картана, Труды МИАН, т. 304 (2019), 49–67.

\bibitem{AS19_2}
Ардентов А.А., Сачков Ю.Л., Субфинслерова задача на группе Картана, Доклады Академии Наук, 2019, 484, No. 2, 138-141

\bibitem{AS19_3}
Ардентов А.А., Сачков Ю.Л., Субфинслеровы структуры на группе Энгеля, Доклады Академии Наук, 2019, 485, No. 4, 395—398

\bibitem
{b1}
Берестовский В.Н., \emph{Однородные многообразия с внутренней метрикой. {II}}, Сибирский математический журнал, \textbf{30} (1989), no.~2, 14--28,
  225.

\bibitem
{b2}
Берестовский В.Н., \emph{О структуре однородных локально компактных пространств с внутренней метрикой}, Сибирский математический журнал, \textbf{30} (1989), no.~1, 23--34.

\bibitem{ber94_2}
В. Н. Берестовский, “Геодезические неголономных левоинвариантных внутренних метрик на группе Гейзенберга и изопериметриксы плоскости Минковского”, Сиб. матем. журн., 35:1 (1994), 3–11
 
\bibitem{berzub20}
В. Н. Берестовский, И. А. Зубарева, “Экстремали левоинвариантной субфинслеровой метрики на группе Энгеля”, Сиб. матем. журн., 61:4 (2020), 735–751

\bibitem{berzub20_2}
	В. Н. Берестовский, И. А. Зубарева, “ПМП, (ко)присоединённое представление и нормальные геодезические левоинвариантных (суб)финслеровых метрик на группах Ли”, Чебышевский сб., 21:2 (2020), 43–64 

\bibitem{convex}
Локуциевский Л.В. 
Выпуклая тригонометрия с приложениями к субфинслеровой геометрии // Матем. сб. 2019. Т. 210. № 8. С. 120--148.

\bibitem
{BBLDS}
Davide Barilari, Ugo Boscain, Enrico {L}{e Donne}, and Mario Sigalotti,
  \emph{Sub-{F}insler structures from the time-optimal control viewpoint for
  some nilpotent distributions}, J. Dyn. Control Syst. \textbf{23} (2017),
  no.~3, 547--575.
	
\bibitem
{boscain3level}
Ugo Boscain, Thomas Chambrion, and Gr{\'e}goire Charlot, \emph{Nonisotropic
  3-level quantum systems: complete solutions for minimum time and minimum
  energy}, Discrete Contin. Dyn. Syst. Ser. B \textbf{5} (2005), no.~4,
  957--990 (electronic).
	
\bibitem{buseman}
 H. Busemann,
The Isoperimetric Problem in the Minkowski Plane,
{\em American Journal of Mathematics}, Vol. 69, No. 4 (Oct., 1947), pp. 863--871

\bibitem{S18}
	Yu. L. Sachkov, Optimal bang-bang trajectories in sub-Finsler problem on the Cartan group, Russian Journal of Nonlinear Dynamics, 2018, vol.4, n.4, 583-593. 

\bibitem{(36)}
	Yu. L. Sachkov, Periodic controls in step 2 strictly convex sub-Finsler problems,
Regular and Chaotic Dynamics, 2020, Vol. 25, No. 1, pp. 33--39.

\subsection*{Левоинвариантные сублоренцевы задачи}
\bibitem{groch06}
M. Grochowski, Reachable Sets for the Heisenberg sub-Lorentzian Metric on $\R^3$. An Estimate
for the Distance Function, {\em Journal of Dynamical and Control Systems}, Vol. 12, No. 2, 2006.

\bibitem{grong_vas}
Grong E., Vasil’ev A. Sub-Riemannian and sub-Lorentzian geometry on $\SU(1, 1)$ and on its
universal cover, {\em J. Geom. Mech.} 3 (2011), no. 2, 225-260.

\bibitem{kor_mar}
A. Korolko, I. Markina, Nonholonomic Lorentzian Geometry on Some $H$-type Groups, 
{\em Journal of Geometric Analysis}, Vol. 19, No. 4, 2009.

\bibitem{CHSY}
Qihui Cai, Tiren Huang, Yu. Sachkov, Xiaoping Yang, Geodesics in the Engel group with a sub-Lorentzian metric, {\em Journal of Dynamical and Control Systems}, Vol. 22 (2016), 465--484.

\bibitem{ASHY}
А. А. Ардентов, Ю. Л. Сачков, Т. Хуанг, К. Янг, “Экстремальные траектории в сублоренцевой задаче на группе Энгеля”, Матем. сб., 209:11 (2018), 3–31

\subsection*{Разные статьи по левоивариантным задачам} 



\bibitem{ber18}
В. Н. Берестовский, “Геодезические и кривизны специальных субримановых метрик на группах Ли”, Сиб. матем. журн., 59:1 (2018), 41–55 


\bibitem{MS15}
	А.П. Маштаков, Ю.Л. Сачков, Суперинтегрируемость левоинвариантных субримановых структур на унимодулярных трехмерных группах Ли, Дифференциальные уравнения, 2015, 11, 1482—1488.

\subsection*{Нильпотентная аппроксимация и конструктивное решение \\задачи управления}
\bibitem{agrach_sarych}
Аграчев А.А., Сарычев А.А.
Фильтрация алгебры Ли векторных полей и нильпотентная аппроксимация управляемых
систем//
Докл. АН СССР.
1987.
Т. 295, No. 4, С. 777–781

\bibitem{masht12}
А. П. Маштаков, Алгоритмическое и программное обеспечение решения
конструктивной задачи управления неголономными пятимерными
системами, Программные системы: теория и приложения, 3:1 (2012),
с. 3–29



\bibitem{agrachev_marigo}
A. Agrachev and A. Marigo. Nonholonomic tangent spaces: intrinsic
construction and rigid dimensions. Electron. Res. Announc. Amer. Math.
Soc., 9: 111–120, 2003.

\bibitem{bellaiche}
Bellaiche A.
The tangent space in sub-Riemannian geometry//
In: Sub-Riemannian geometry, A.~Bellaiche and J.-J.~Risler,
Eds., Birkh\"auser, Basel, Swizerland, 1996, 1--78.

\bibitem{Bellaiche11}  
{\it Bellaiche~A., Laumond~J. P., Chyba~M.}
Canonical nilpotent approximation of control systems: application to nonholonomic motion planning. 32nd IEEE CDC, 1993.

\bibitem{Bellaiche10} 
{\it Bellaiche~A., Laumond~J.P., Risler~J.J.}
Nilpotent infinetisimal approximations to a control Lie algebra. IFAC NCSDS, Bordeaux, 1992. P.~174--181.

\bibitem{BS}
Rosa Maria Bianchini and Gianna Stefani. Graded approximations
and controllability along a trajectory. SIAM J. Control Optim., 28(4):
903–924, 1990.

\bibitem{Jean2013}
{\it Chitour~Y., Jean~F., Long~R.}
A Global Steering Method for Nonholonomic Systems, Journal of Differential Equations, 2013, V. 254, P. 1903--1956.

\bibitem{Fer35} {\it Fernandes~C., Gurvits~L., Li~Z.X.} 
A variational approach to optimal nonholonomic motion planning. IEEE ICRA, Sacramento, 1991. P. 680--685.

\bibitem{gromov}
Mikhael Gromov. Carnot–Carath\'eodory spaces seen from within. In
Sub-Riemannian geometry, volume 144 of Progress in Mathematics,
pp. 79–323. Birkh\"auser, Basel, 1996.

\bibitem{b:hermes}
{\it Hermes H.}
Nilpotent and high-order approximations of vector fields systems //  SIAM. 1991. V. 33. P. 238--264.

\bibitem{b:lafsus}
{\it Lafferriere~G., Sussmann~H.J.}
A differential geometric approach to motion planning. 
Nonholonomic Motion Planing. 1992. Editors: Zexiang Li, J.F. Canny.

\bibitem{suss2}
Laferriere, G. and Sussmann, H.J., A differential geometric approach to motion planning, \textit{Nonholonomic motion planning}, Kluwer, 1993, pp.~235--270.

\bibitem{laumond}
Laumond J.P.
Nonholonomic motion planning for mobile robots// Preprint No.~98211.
 Toulouse, France: LAAS-CNRS, 1998.

\bibitem{l81} Murray R. M., {\em Nilpotent bases for a class on nonintegrable distributions with applications to trajectory generation for nonholonomic systems}// Math. Control Signal Syst., university of California, Berkeley, 1990.

\bibitem{stefani}
{\it Stefani G.}
Polynomial approximations to control systems and local controllability // Proc. 24th. I.E.E.E. Conference on Decision and
Control, Ft. Lauderdale, Fla., 1985, P. 33--38.

\bibitem{TilSas119} Tilbury D., Murray R., Sastry S., {\em Trajectory generation for the n-trailer problem using Goursat normal form}// IEEE TAC, 1995, V. 40, No. 5, P. 802--819.

\bibitem{vendit_laumond_oriolo}
Vendittelli M., Laumond J.P., Oriolo G.
Steering nonholonomic systems via nilpotent approximations: The general two-
trailer
system//
1999 IEEE International Confer. on Robotics and Automation, May 10--15,
1999,
Detroit, MI.

\bibitem{b:venditelli}
{\it Venditelli~M., Oriolo~G., Jea~F., Laumond~J.P.}
Nonhomogeneous nilpotent approximations for nonholonomic systems with singularities // Transactions on Automatic Control. 2004. P. 261--266.

\bibitem{walsh}
Walsh, G.C., Montgomery, R., and Sastry, S.S., Optimal path planning on matrix Lie group, in \textit{Proc. of the 33rd IEEE Conference on Decision and Control}, 1994, vol.~2, pp.~1258--1263. 


\subsection*{Приложения субримановой геометрии к моделям зрения и \\обработке изображений}
\bibitem{MDSBB}
	А. П. Маштаков, P. Дайтс, Ю. Л. Сачков, E. Беккерс, И. Ю. Бесчастный, Субримановы геодезические на группе SO(3) в задаче поиска кровеносных сосудов на сферических изображениях сетчатки,  Доклады Академии Наук, 2017, 473, № 5, с. 521–524

\bibitem{benyosef}
G.~Ben-Yosef and O.~Ben-Shahar, {\em A tangent bundle theory for visual
curve completion}. PAMI 34(7), pp.1263–-1280, 2012.

\bibitem{BDMG15}
E .J. Bekkers, R. Duits, A. Mashtakov, G. R. Sanguinetti. “A PDE approach
to data-driven sub-Riemannian geodesics in SE(2)”, SIAM Journal on
Imaging Sciences, 8:4 (2015), pp. 2740–2770

\bibitem{BDMS17}
E. J. Bekkers, R. Duits, A. Mashtakov, Y. Sachkov. “Vessel tracking via
sub-Riemannian geodesics on the projective line bundle”, Geometric Science
of Information, GSI 2017, LNCS, vol. 10589, eds. F. Nielsen, F. Barbaresco,
Springer, Cham, 2017, pp. 773–781.

\bibitem{Boscain2014}
U.~Boscain, R.~Duits, F.~Rossi and Y.~L.~Sachkov,
Curve Cuspless Reconstruction via sub-Riemannian Geometry,
ESAIM: Control, Optimisation and Calculus of Variations, 20 (2014), 748--770,
doi:10.1051/cocv/2013082

\bibitem{Gauthier}
U. Boscain, J.-P. Gauthier, R. Chertovskih, A. Remizov
Hypoelliptic diffusion and human vision: a semidiscrete new twist //
SIAM J. Imaging Sci., 7 (2), 2014, pp. 669--695.

\bibitem{citti_sarti}
G. Citti, A. Sarti. A cortical based model of perceptual completion in the
roto-translation space, Journal of Mathematical Imaging and Vision, 24:3
(2006), pp. 307–326.

\bibitem{DuitsJMIV2014}
R.~Duits, U.~Boscain, F.~Rossi and Y.~L.~Sachkov,
Association Fields via Cuspless Sub-Riemannian Geodesics in SE(2) //
JMIV, 49 (2), 2014, pp. 384--417, doi: 10.1007/s10851-013-0475-y

\bibitem{DuitsIJCV2007}
R.~Duits, M.~Felsberg, G.~Granlund, B.H.~Romeny, {\em Image Analysis and Reconstruction using a Wavelet Transform Constructed from a Reducible Representation of the Euclidean Motion Group.} IJCV, vol.72 (1), pp. 79--102, 2007.

\bibitem{JDCS16}
R.~Duits, A.~Ghosh,  T.\,C.\,J.~Dela Haije, A.~Mashtakov,
On Sub-Riemannian Geodesics in SE(3) Whose Spatial Projections do not Have Cusps // Journal of Dynamical and Control Systems, 2016, 22(4), pp:771--805, doi:10.1007/s10883-016-9329-4

\bibitem{FMCS}
B. Franceschiello, A. Mashtakov, G. Citti, A. Sarti. “Modelling of the
Poggendorff illusion via sub-Riemannian geodesics in the roto-translation
group”, New Trends in Image Analysis and Processing, ICIAP 2017, LNCS,
vol. 10590, eds. S. Battiato, G. M. Farinella, M. Leo, G. Gallo, Springer,
2017, pp. 37–47.

\bibitem{DGA19}
B. Franceschiello, A. Mashtakov, G. Citti, A. Sarti, Geometrical optical illusion via sub-Riemannian geodesics in the roto-translation group // Differential Geometry and its Applications, Volume 65, 2019, pp. 55--77, doi:10.1016/j.difgeo.2019.03.007

\bibitem{Franken2009}
E.~Franken and R.~Duits,
{\em Crossing-Preserving Coherence-Enhancing Diffusion on Invertible Orientation Scores.}
IJCV, 85(3), pp. 253-278, 2009.

\bibitem{MAS13}
A. P. Mashtakov, A. A. Ardentov, Yu. L. Sachkov. “Parallel algorithm
and software for image inpainting via sub-Riemannian minimizers on the
group of rototranslations”, Numerical Mathematics: Theory, Methods and
Applications, 6:1 (2013), pp. 95–115.

\bibitem{MDSBB17}
A. Mashtakov, R. Duits, Y. Sachkov, E. J. Bekkers, I. Beschastnyi. “Tracking
of lines in spherical images via sub-Riemannian geodesics in SO(3)”, JMIV,
58:2 (2017), pp. 239–264

\bibitem{petitot}
J.Petitot, The neurogeometry
	of pinwheels as a sub-Riemannian contact structure, 
	{\em J. Physiology - Paris} 97 (2003), 265--309.

\bibitem{petitot_book}
J.Petitot, 
{\em Neurogeometrie de la vision --- Modeles mathematiques et physiques des architectures fonctionnelles}, 2008, Editions de l'Ecole Polytechnique. 



\subsection*{Приложения субримановой геометрии к робототехнике}
\bibitem{Ar-Eu-mob-rob} А. А. Ардентов, А. В. Смирнов, Управление мобильным роботом вдоль эластик Эйлера, Программные системы: теория и приложения, 2017, том 8, выпуск 4, 163–178

\bibitem{ard_gub}
А. А. Ардентов, И. С. Губанов. Моделирование парковки автомобиля
с прицепом вдоль путей Маркова-Дубинса и Ридса-Шеппа. {\em Программные
системы: теория и приложения}, 2019, 10:4(43), с. 97–110.



\bibitem{ArdentovRCD}
{\it Ardentov A.A.} 
Controlling of a mobile robot with a trailer and its nilpotent approximation // Regular and Chaotic Dynamics. 2016. V. 21. No. 7--8, P. 775--791.

\bibitem{EU-robot}
Ardentov A.A., Karavaev Y.L., Yefremov K.S.
Euler Elasticas for Optimal Control of the Motion of Mobile Wheeled Robots: the Problem of Experimental Realization // RCD.
2019. V. 24. No. 3. P. 312--328.

\bibitem{laum98}
Laumond, J.-P., \textit{Nonholonomic motion planning for mobile robots}, Tutorial notes, 1998.

\bibitem{monroy}
Anzaldo-Menezes A., Monroy-P\'erez F.
Charges in magnetic fields and sub-Riemannian geodesics//
Contemporary trends in nonlinear geometric control theory and its
applications. Singapore:  World Scientific, 2002. P.~183--202.



\bibitem{li_canny}
Li Z., Canny J.
Motion of two rigid bodies with rolling constraint//
IEEE Trans. Robot. Automat. (1). 1990. V.~6. P.~62--72.


\bibitem{cdc99}
Agrachev A.\,A., Sachkov Yu.L.
An intrinsic approach to the control of rolling bodies//
Proceedings of the 38-th
IEEE Conference on Decision and Control, Phoenix, Arizona, USA, December
7--10, 1999. V.~1. P.~431--435.

\bibitem{marigo_bicchi}
Marigo A., Bicchi A.
Rolling bodies with regular surface: the holonomic case//
 Differential geometry and control: Summer Research Institute
on Differential Geometry and Control,  Univ.
Colorado, Boulder, June 29--July 19, 1997, / ed.~G.~Ferreyra et al. Providence, RI: Amer. Math.~Soc., 1999. P.~241--256 (Proc.~Sympos.~Pure Math. V.~64).



\end{thebibliography}
\end{document}